\documentclass[12pt]{amsart}
\usepackage[a4paper]{anysize}
\usepackage{amsthm}
\usepackage{amssymb}
\usepackage{times}
\usepackage{mathrsfs}

\input xy
\xyoption{v2}
\xyoption{all}

    \newtheorem{theorem}{Theorem}[subsection]
    \newtheorem{proposition}[theorem]{Proposition}
    \newtheorem{lemma}[theorem]{Lemma}
    \newtheorem{corollary}[theorem]{Corollary}
    
    \theoremstyle{definition}
    \newtheorem{example}[theorem]{Example}
    \newtheorem{definition}[theorem]{Definition}
    \newtheorem{remark}[theorem]{Remark}
    \theoremstyle{remark}
    \newtheorem{claim}[theorem]{Claim}

    \def\set{\setcounter{equation}
             {\value{theorem}}\addtocounter{theorem}{1}}
    \numberwithin{equation}{subsection}
    \def\sset{\setcounter{subsubsection}
             {\value{theorem}}\addtocounter{theorem}{1}}

\def\A{{\mathbb A}}

\def\cB{{{\mathscr B}}}
\def\C{{\mathbb C}}

\def\cC{{{\mathscr C}}}  
\def\sC{{\text{\sf C}}}
\def\cD{{{\mathscr D}}}

\def\D{{\mathbb D}}

\def\sD{{\text{\sf D}}}
\def\sE{{\text{\sf E}}}
\def\E{{\mathbb E}}
\def\cE{{{\mathscr E}}}
\def\F{{\mathbb F}}
\def\cF{{{\mathscr F}}}
\def\fF{{\mathfrak F}}
\def\G{{\mathbb G}}
\def\cG{{{\mathscr G}}}

\def\cK{{{\mathscr K}}}
\def\cL{{{\mathscr L}}}

\def\cM{{{\mathscr M}}}

\def\N{{\mathbb N}}
\def\cO{{{\mathscr O}}}

\def\fp{{\mathfrak p}}
\def\P{{\mathbb P}}

\def\fq{{\mathfrak q}}
\def\Q{{\mathbb Q}}
\def\R{{\mathbb R}}
\def\SS{{\mathbb S}}

\def\Z{{\mathbb Z}}

\def\cK{{{\mathscr K}}}

\newcommand{\colim}[1]
{\underset{#1}{\mathrm{colim}\,}}

\newcommand{\liminv}[1]
{\underset{#1}{\lim}}

\newcommand{\derotimes}
{\overset{\mathbf{L}}{\otimes}}

\renewcommand{\labelenumi}{(\roman{enumi})}

\newenvironment{pfclaim}[1][$\Diamond$]
{\def\claimQED{{#1}}\noindent {\em Proof of the claim. }}
{\leavevmode\unskip\penalty9999 \hbox{}\nobreak\hfill
    \quad\hbox{\claimQED}{\smallskip}}

\def\Tr{\mathrm{Tr}}
\def\tr{\mathrm{tr}}
\def\gr{\mathrm{gr}}

\def\Ann{\mathrm{Ann}}

\def\Gal{\mathrm{Gal}}
\def\Spec{\mathrm{Spec}}
\def\Hom{\mathrm{Hom}}

\def\Tor{\mathrm{Tor}}
\def\End{\mathrm{End}}

\def\Ind{\mathrm{Ind}}
\def\Spf{\mathrm{Spf}}
\def\Coker{\mathrm{Coker}}
\def\Ker{\mathrm{Ker}}
\def\Img{\mathrm{Im}}
\def\Aut{\mathrm{Aut}}
\def\Alg{{{\text{-}}\mathbf{Alg}}}

\def\Loc{\text{-}\mathbf{Loc}}
\def\Mod{\text{-}\mathbf{Mod}}

\def\Set{\mathbf{Set}}
\def\chara{\mathrm{char}}

\def\rk{\mathrm{rk}}
\def\one{\mathbf{1}}
\def\bar#1{\overline{#1}}
\def\hat{\widehat}  
\def\fm{\mathfrak m}
\def\fp{\mathfrak p}
\def\eps{\varepsilon}

\def\et{\mathrm{\acute{e}t}}

\def\bmu{\boldsymbol\mu}

\def\sm{\mathrm{sm}}
\def\sing{\mathrm{sing}}

\renewcommand\emptyset{\varnothing}

\def\sA{{\mathsf A}}
\def\sB{{\mathsf B}}
\def\sC{{\mathsf C}}
\def\sT{{\mathsf T}}
\def\sU{{\mathsf U}}
\def\sV{{\mathsf V}}
\def\sX{{\mathsf X}}
\def\sY{{\mathsf Y}}
\def\sZ{{\mathsf Z}}
\def\fP{{\mathfrak P}}
\def\fV{{\mathfrak V}}
\def\fU{{\mathfrak U}}
\def\fX{{\mathfrak X}}
\def\fY{{\mathfrak Y}}
\def\fW{{\mathfrak W}}

\def\fd{\mathfrak d}

\def\fn{\mathfrak n}
\def\sp{\mathrm{sp}}
\def\Spa{\mathrm{Spa}}

\def\an{\mathrm{an}}

\def\Fil{\mathrm{Fil}}

\def\trg{\triangleright}
\def\ad{\mathrm{ad}}
\renewcommand\theta{\vartheta}
\def\Nm{\mathrm{Nm}}
\def\sff{{\mathsf f}}
\def\sg{{\mathsf g}}
\def\sa{{\mathsf a}}
\def\ssw{\mathsf{sw}}
\def\su{{\mathsf u}}
\def\sreg{\mathsf{reg}}
\def\spsi{{\mathsf\Psi}}
\def\isom{\stackrel{\sim}{\to}}
\def\fpet{{\mathrm{fp}\et}}
\def\fget{{\mathrm{fg}\et}}

\def\alphaenu{\renewcommand{\labelenumi}{(\alph{enumi})}}
\def\romanenu{\renewcommand{\labelenumi}{(\roman{enumi})}}
\def\addenu{\addtocounter{enumi}{1}}

\font\prelim=pagko

\makeindex

\begin{document}
\title{Local monodromy in non-archimedean analytic geometry}

\author{Lorenzo Ramero}

\maketitle

\centerline{\prelim fifth release}
\medskip

%\centerline\today

\tableofcontents

\vskip 0.5cm

    Lorenzo Ramero 

    Universit\'e de Bordeaux I 

    Institut de Math\'ematiques 

    351, cours de la Liberation 

    F-33405 Talence Cedex

{\em e-mail address:} {\ttfamily ramero@math.u-bordeaux.fr}

{\em web page:} {\ttfamily http://www.math.u-bordeaux.fr/$\sim$ramero}

\newpage

\setcounter{section}{-1}

\hfill{\small\it ``Der Weg ist das Ziel''}

\smallskip

\hfill{\small\rm Bierdeckel im Caf\'e Pendel}

\hfill{\small\rm Bonn, ca. Mai 2004}

\bigskip

\section{Introduction (Snapshots from an Expedition)}

\subsection{Where we came from}
Let $F:=k((T))$ be the field of Laurent power series with
coefficients in a field $k$; if $k$ has positive characteristic,
denote by $G_F$ the Galois group of a separable closure of $F$.
As it is well known, every $\ell$-adic representation admits
a natural {\em break decomposition\/} as a finite direct sum
of $\Q_\ell[G_F]$-submodules. These decompositions
are nicely compatible with the standard operations defined
on the category $\Q_\ell[G_F]\Mod$ of $\Q_\ell[G_F]$-modules,
such as tensor products and $\Hom$ functors.

On the other hand, if $k$ has characteristic zero, one may consider
the category $\sD.\sE.(F/k)$ of finite dimensional $F$-vector spaces
endowed with a $T$-adically continuous $k$-linear connection.
Then the theory of Levelt-Turritin says that every object
of $\sD.\sE.(F/k)$ admits a natural decomposition, satisfying wholly
analogous compatibilities (see \cite{Kat2}).

The parallelisms revealed by the study of $\Q_\ell[G_F]\Mod$ and
$\sD.\sE.(F/k)$ are aboundant and striking, to the point that one can
write down a heuristic dictionary to translate definitions and
theorems back and forth between them (see \cite[Appendix]{Kat2}).
More recently, Yves Andr\'e has introduced notions of
{\em slope filtration} and {\em Hasse-Arf filtration} for a general
tannakian category
(see \cite[D\'ef.1.1.1 and D\'ef.2.2.1]{And}), which extract the
basic features that are common to these two categories
(and indeed, to others as well). He has shown that, quite
generally, the existence of a Hasse-Arf filtration imposes
very binding restrictions on the structure of a tannakian category.

\subsection{Where we were heading}
The aim of this work is to exhibit another specimen of the
same sort as those considered by Andr\'e. Namely, let
$(K,|\cdot|)$ be an algebraically closed valued field
of mixed characteristic $(0,p)$, complete for its rank one
valuation $|\cdot|:K\to\Gamma_K\cup\{0\}$ (we may view the
value group $\Gamma_K$ as a subgroup of $\R_{>0}$). Let also
$\Lambda$ be a local ring which is a filtered union
of finite rings on which $p$ is invertible.
Our objects of study are the  locally constant sheaves of free
$\Lambda$-modules of finite rank on the \'etale site of the
punctured disc :
$$
\D(r)^*:=\{x\in K~|~0<|x|\leq r\}\qquad
\text{(for any $r\in\Gamma_K$)}.
$$
These modules form a category $\Lambda\Loc(r)$, which is
tannakian if $\Lambda$ is a field. However, we are really
interested in describing the {\em local monodromy\/} of these
sheaves, {\em i.e.} their behaviour in an arbitrarily
small neighborhood of the missing center of the disc, hence we
do not distinguish two local systems $\cF$ and $\cF'$ on the \'etale
sites $\D(r)^*_\et$, respectively $\D(r')^*_\et$, if they become
isomorphic after restriction to some smaller disc $\D(r'')^*$
(with $0<r''\leq r,r'$). Hence, we are really concerned with
the $2$-colimit category :
$$
\Lambda\Loc(0^+):=\colim{r\to 0^+}\Lambda\Loc(r).
$$

\subsection{What we hoped to find there}
It is instructive to consider first the case where the monodromy
is finite, {\em i.e.} the local system $\cF$ under consideration
becomes constant on some finite Galois \'etale covering
$X\to\D(r)^*$, say with Galois group $G$. This case is already
highly non-trivial : $\cF$ is the same as a $\Lambda[G]$-module
of finite type, and if we insisted on a complete description
of the {\em global\/} monodromy of $\cF$, we would have to classify
all such representations, so essentially all possible finite
coverings of $\D(r)^*$ -- a task that is probably beyond the reach
of current techniques. On the other hand, the {\em germs\/} of
finite coverings of $\D(r)^*$ are completely classified by the
so-called $p$-adic Riemann existence theorem, proved by O.Gabber
around 1982 (unpublished) and by W.L\"utkebohmert about ten years
later (\cite{Lut}). Explicitly, after restriction to a smaller
disc, every such finite covering becomes a disjoint union of 
cyclic coverings of Kummer type; therefore the local monodromy
of our $\cF$ will be a representation of a finite cyclic group.

Though no such description is known for general
\'etale coverings of $\D(1)^*$ ({\em i.e.} for those of infinite degree),
this provides some evidence for the thesis that, by shifting the focus
to germs of coverings, one should expect to reach a new, but substantially
tamer mathematical territory -- one in which some general geographical
features are discernible and can be used as worthwhile reference points.

\sset\subsection{Why we were not disappointed}
This expectation is largely borne out by our main theorem \ref{th_decompose},
which can be stated as follows. Suppose that $H^1(\D(r)_\et^*,\cF)$
is a $\Lambda$-module of finite type, in which case we say that
$\cF$ has {\em bounded ramification}; this condition is independent
of the value $r\in\Gamma_K$. Then there exists a connected
open subset $U\subset\D(r)^*$ such that $U\cap\D(\eps)^*\neq\emptyset$
for every $\eps>0$, and such that the restriction of $F$ to $U_\et$
admits a {\em break decomposition} as a direct sum of locally constant
subsheaves :
\set\begin{equation}\label{eq_decompose}
\cF_{|U}\simeq\bigoplus_{\gamma\in\Gamma_0}\cF(\gamma)
\end{equation}
indexed by the ordered group $\Gamma_0$, which is the product
of ordered groups $\Q\times\R$, endowed with the lexicographic
ordering (of course $F(\gamma)\neq 0$ for only finitely many
$\gamma\in\Gamma$). This decomposition is compatible in the usual
way with tensor products and $\Hom$ functors; moreover, one
may define Swan conductors for $\cF$, and there is also an
adequate analogue of the Hasse-Arf theorem. Since the Swan
conductor determines the Euler-Poincar\'e characteristic
of $R\Gamma(\D(r)_\et^*,\cF)$, it follows easily that, in
case $\Lambda$ is a field, the subcategory
$\Lambda\Loc(0^+)_\mathrm{b.r.}$ of $\Lambda\Loc(0^+)$
consisting of $\Lambda$-modules with bounded ramification,
is tannakian. Therefore, the break decomposition in
$\Lambda\Loc(0^+)_\mathrm{b.r.}$ would be precisely what is
needed to define a filtration of Hasse-Arf type, in the
sense of \cite{And}, if it were not for the following two
short-comings. First, the filtration is not indexed by the
real numbers, but by the more complicated group $\Gamma_0$.
This is however a minor divergence, which can be cured,
for instance, by generalizing a little the definition of slope
filtration in a tannakian category. More seriously, the
break decomposition of an object of $\Lambda\Loc(0^+)_\mathrm{b.r.}$
is defined (a priori) only in a strictly larger tannakian category;
this is because the open subset $U$ may not contain any punctured
disc $\D(\eps)^*$ (though it intersects all of them).

I expect that actually every local system with bounded
ramification admits a break decomposition already over
some small punctured disc, and I hope to address this
question in a future release of this work. Once this
result is available, it will be possible to apply the
tannakian machinery of \cite{And} to study the structure of
such local systems.

\subsection{Planning for the journey}
The proof of theorem \ref{th_decompose} is divided into two separate
steps. The first step consists in describing the monodromy of $\cF$
{\em around the border\/} of a disc $\D(r)^*$. This is one of the
points where the standard topological intuition may be misleading:
a non-archimedean punctured disc is far from being ``homotopically
equivalent'' to an annulus, for any reasonable notion of homotopy
equivalence. Indeed, it is easy to construct examples of (finite
or infinite) connected Galois coverings of $\D(r)^*$ that are
completely split at the border, that is over every annulus of
the form:
$$
\D(a,r):=\{x\in K~|~a\leq|x|\leq r\}
$$
with $a\in\Gamma_K$ sufficiently close to $r$. (And conversely,
an \'etale covering of $\D(r)^*$ may be completely split
near the center, and connected on every annulus $\D(a,r)$ with
$a$ sufficiently close to $r$.) But the crucial difference is
that the monodromy around the border is {\em always finite}.
As a consequence, the study of this monodromy is an essentially
algebraic affair that can be carried out by a suitable extension of
classical ramification theory for henselian discrete valuation rings
with perfect residue field.

This extension has been developed by R.Huber in \cite{Hu3}: the
main tool is a certain rank two valuation $\eta(r)$, defined
on the ring of analytic functions $\cO(\D(a,r))$ on $\D(a,r)$
(for any $a<r$), and continuous for the $p$-adic topology.
In a precise sense (best expressed in the language of adic spaces),
this valuation is localized at the border of the disc. Moreover,
the value group $\Gamma_{\eta(r)}$ of $\eta(r)$ is naturally isomorphic
to $\Gamma_K\times\Z$, ordered lexicographically. 
Now, let $f:X\to\D(a,r)_\et$ be a finite, connected, \'etale and
Galois covering of Galois group $G$, such that $f^*\cF$ is a constant
sheaf. Then $\eta(r)$ admits finitely many extensions $x_1,\dots,x_n$ to
$\cO(X)$, and the action of $G$ permutes transitively these extensions.
Fix one of these valuations $x:=x_i$; the stabilizer $St_x\subset G$
can be naturally identified with the Galois group of a finite 
separable extension 
$$
\kappa(\eta(r))^{\wedge h}\subset\kappa(x)^{\wedge h}
$$
of henselian valued fields, obtained by suitably henselizing the
completions (for the valuation topologies) $\kappa(\eta(r))^\wedge$,
$\kappa(x)^\wedge$ of the fields of fractions of $\cO(\D(a,r))$
and respectively $\cO(X)$. Let $\pi_1(r)$ be the Galois group of
a separable closure of $\kappa(\eta(r))^{\wedge h}$; we may regard
$\cF$ as a $\Lambda[G]$-module, hence as a $\Lambda[St_x]$-module,
by restriction, and then as a $\Lambda[\pi_1(r)]$-module, by inflation.

The group $\pi_1(r)$ admits a natural (upper numbering) higher
ramification filtration, wholly analogous to the standard one for
discrete valued fields, except that it is indexed by the ordered
group $\Gamma_{\eta(r)}\otimes_\Z\Q$. Therefore, when $\Lambda$
is a field, the tannakian category $\Lambda[\pi_1(r)]\Mod$
is yet another example of a category with a slope filtration,
except that the filtration is indexed by $\Gamma_{\eta(r)}\otimes_\Z\Q$,
rather than by $\R$. Moreover, this slope filtration  is even
of {\em Hasse-Arf type}, provided we redefine appropriately the
Newton polygon of a representation.

\subsection{Division of labour}
R.Huber's ramification theory yields, for every radius $r\in\Gamma_K$,
a $\Lambda[\pi_1(r)]$-equivariant break decomposition of the stalk
$\cF_{\eta(r)}$. The second step of the proof of theorem
\ref{th_decompose} consists in describing how this decomposition
evolves as $r$ changes. This step presents in turn two aspects:
on the one hand, we have to examine the {\em continuity\/} properties
of the breaks, {\em i.e.} the way the decomposition varies
in a neighborhood of a given radius $r$; on the other hand, we
have to make an {\em asymptotic\/} study, to determine the
behaviour of the decomposition for $r\to 0^+$. The upshot
is that, for large values of $r$, the breaks of $\cF_{\eta(r)}$
vary in a continuous, but apparently patternless manner; but,
as we approach the missing center of the disc, eventually a
coherent order emerges: the decompositions fall into step, and
they give rise to the asymptotic decomposition \eqref{eq_decompose}.

\subsection{Surveyor's gear}
Both the continuity and the aymptotic study ultimately
rely on the remarkable properties of certain conductor
functions attached to our local system $\cF$. To define
these conductors we may assume, without loss of generality,
that $\Gamma_K=\R_{>0}$. Suppose that $f:X\to\D(a,b)$ is
a finite Galois \'etale covering, with group $G$, such that
$f^*\cF$ is constant. For every $r\in[a,b]$, we have the
$\pi_1(r)$-equivariant break decomposition
$$
\cF_{\eta(r)}\simeq M_1(\gamma_1(r))\oplus\cdots\oplus M_n(\gamma_n(r))
$$
where $n$ depends also on $r$, and the breaks $\gamma_i(r)$
live in $\R_{>0}\times\Q$. For any $\gamma\in\R_{>0}\times\Q$,
let us denote by $\gamma^\flat$ and $\gamma^\natural$ the
projections of $\gamma$ on $\R_{>0}$ and respectively $\Q$.
Set also $m_i:=\rk_\Lambda M_i(\gamma_i(r))$ for every $i=1,\dots,n$.
Then we may consider the  conductor functions :
$$
\delta_\cF:[\log 1/b,\log 1/a]\to\R_{\geq 0}
\quad\text{and}\quad
\ssw^\natural(\cF,\cdot):[a,b]\to\Z
$$
defined by letting :
$$
\delta_\cF(-\log r):=
-\sum_{i=1}^n\log\gamma_i(r)^\flat\cdot m_i
\quad\text{and}\quad
\ssw^\natural(\cF,r^+):=
\sum_{i=1}^n\gamma_i(r)^\natural\cdot m_i.
$$
We show that $\delta_\cF$ is a piecewise linear, continuous
and convex function, and moreover the right derivative of $\delta_\cF$
is computed by $\ssw^\natural(\cF,\cdot)$ (see proposition
\ref{prop_invoke}). The function $\ssw^\natural(\cF,\cdot)$
can be characterized in terms of the Swan conductor of the
covering $X$. Namely, for every $r\in[a,b]$, choose a valuation
$x$ of $\cO(X)$ extending the valuation $\eta(r)$; then the higher
ramification filtration of $St_x$ determines, in the usual way, a
$\Z$-valued Swan character $\ssw_x$ of $St_x$ (see \ref{sec_conductors}),
which is the character of an element of $K^0(\Z_\ell[St_x])$.
We induce to get a virtual character of $G$ :
$$
\ssw_G^\natural(r^+):=\Ind^G_{St_x}\ssw^\natural_x
$$
which is independent of the choice of $x$. Let now
$\rho\in K_0(\Lambda[G])$ be the $\Lambda[G]$-module corresponding
to $\cF$; since $\ssw_G^\natural(r^+)$ lies in $K^0(\Z_\ell[G])$,
we may apply the natural pairing
$$
\langle\cdot,\cdot\rangle_G:K^0(\Z_\ell[G])\times K_0(\Lambda[G])\to\Z
$$
and we obtain the identity :
$$
\ssw^\natural(\cF,r^+)=\langle\ssw_G^\natural(r^+),\rho\rangle_G.
$$
On the other hand, for every $r\in[a,b]$ one can consider
the preimage $X(r)\subset X$ of the annulus $\D(r,r)\subset\D(a,b)$.
The ring of analytic functions ${\mathscr O}^+_X(r)$ on $X(r)$ whose
sup-norm is $\leq 1$ is a finite free module over the analogous ring
${\mathscr O}^+_\D(r)$ of bounded functions on $\D(r,r)$. Hence the
discriminant $\fd^\flat_f(r)$ of this ring extension is well-defined,
and it is an invertible function on $\D(r,r)$, since $f$ is \'etale.
Its sup-norm $|\fd^\flat_f(r)|$ is a real number in $]0,1]$,
and $\delta_f(-\log r):=-\log|\fd^\flat_f(r)|\in\R_{\geq 0}$.
Let now $\cG:=f_*\Lambda_X$ be the direct image of the constant
sheaf $\Lambda_X$ on $X_\et$; then $\cG$ corresponds to the regular
representation of $G$, and we have the identity :
$$
\delta_\cG=\delta_f.
$$
Hence, the right (logarithmic) derivative of the discriminant
is the Swan conductor of the regular representation of $G$ : this
is our analogue of Hasse's F\"uhrerdiskriminantenproduktformel.

\subsection{Detour to visit a relative}
The proof of the convexity of $\delta_\cF$ is accomplished
by a rather technical argument, involving semi-stable
reduction and a vanishing cycle calculation. As a corollary,
one derives a proof of the convexity of the discriminant
function $\delta_f$. However, the convexity of $\delta_f$
can also be shown by a completely elementary argument
that uses little more than some valuation theory, the
first rudiments of the theory of adic spaces, and some
simple tools from $p$-adic analysis borrowed from
\cite{Bo-Gun} and the first chapter of \cite{Fr-VdP}.
We present this argument in section \ref{sec_convexity},
since it is of independent interest : indeed, as explained
 in section \ref{sec_Riemann}, with its aid one may quickly
derive a proof of the $p$-adic Riemann existence theorem.
This proof is not only much more elementary than L\"utkebomert's;
it is also significantly simpler than Gabber's original
argument\footnote{Of course, the reader will have to take my word
for it, since Gabber never published his proof.}. All in all,
I believe it is a convincing demonstration of the new possibilities
opened up by the theory of adic spaces.

\subsection{Dulcis in fundo}\label{sec_dulcis} The ideas which
enable to tackle successfully the asymptotic study of the break
are developed in section \ref{sec_bound-ram}, and find their roots
in harmonic analysis techniques, such as Fourier transform and the
allied method of stationary phase; this should come as no surprise to
any reader familiar with the works of Katz ({\em e.g.} \cite{Kat3})
or Laumon (\cite{Lau2}). Closer to home, these ideas represent an
extension (perhaps, a vindication) of my previous work \cite{Ra},
where I introduced a Fourier transform for sheaves of
$\Lambda$-modules on the \'etale site of the analytification
$(\A^1_K)^\ad$ of the affine line. For more details, we refer
the reader to remark \ref{rem_station}. This intrusion of concepts
and viewpoints originating from such a seemingly far removed
area of mathematics, is for me one of the most appealing
aspects of this project. It was already one of the main themes
in \cite{Ra}, and I believe that it runs deeper than a mere
formal analogy : for instance, from this perspective, formula
\eqref{eq_decompose} is none else than a spectral decomposition
of the local system $\cF$. Whereas \cite{Ra} dealt only with a
suggestive, but {\em ad hoc} class of local systems, we have
now a good grasp of all those local systems $\cF$ whose ramification
is bounded. This boundedness condition can also be restated in
terms of Swan conductors, hence it is, on the one hand a purely
local condition that serves to circumscribe {\em the\/} good
class of local systems that should be stable under the usual
Yoga of cohomological operations. On the other hand, it is a
finiteness condition on the cohomology of $\cF$, hence --
from the harmonic viewpoint -- it is essentially like confining
our attention to the class of ``integrable sheaves''; {\em i.e.}
we are really doing harmonic analysis in the space $L_1$ : a most
natural restriction.

We cannot resist ending on a more speculative note.
As it has been seen, in many situations local monodromy
is described via the higher ramification filtration on
a Galois group, defined by an appropriate valuation.
This cannot be literally true for the local monodromy
theory of the punctured disc, since the trivializing
covering $X\to\D(r)^*$ of a local system may have infinite
degree, in which case the field of fractions of $\cO(X)$ has
infinite transcendence degree over the field of fractions
$\cO(\D(r)^*)$. Nevertheless, one may ask whether there
exists a valuation ``localized at the origin'', which governs,
in some mysterious way, the local monodromy theory  of $\D(r)^*$.
It turns out that there exists one natural candidate, which
is well-defined on the ring $A:=K[T,T^{-1}]$ (of regular functions
on the ``algebraic punctured disc'') : namely, the
rank two specialization $w$ of the degree valuation
$v:A\to\Z\cup\{\infty\}$ such that $v(T^n):=n$ for every $n\in\Z$.
The value group of $w$ is the lexicographically ordered group
$\Delta:=\Z\times\Gamma_K$, and one has the rule :
$w(a\cdot T^n):=(n,|a|)$, for every $n\in\Z$ and every
$a\in K$. Notice that $\Delta\otimes_\Z\Q$ is isomorphic --
as an ordered group -- to the group $\Gamma_0$ that indexes
our break decomposition. However, this valuation $w$ is
not $p$-adically continuous, hence it does not lie in the
adic spectrum $\Spa\,A$, but only in the larger valuation
spectrum $\mathrm{Spv}\,A$.

\medskip

{\bf Acknowledgements}: I am hugely indebted with Ofer
Gabber for explaining me the proof of his unpublished
theorem and for many discussions, through which I have
learned many of the techniques that play an essential
role in this work. Especially, I have learned from him an idea,
due originally to Deligne (\cite{Lau}), that allows to calculate
the rank of vanishing cycles; the same idea is recycled in the
proof of proposition \ref{prop_vanish}. I thank Isabelle Vidal
for sending me a copy of her thesis \cite{Vid}, where she uses
de Jong's method of alterations to deduce consequences for the
\'etale cohomology of schemes; her argument is applied here to
the study of vanishing cycles (theorem \ref{th_rationality}).
I also wholeheartedly acknowledge the support of the IHES and
the Max-Planck Institute in Bonn, where large parts of this
research have been carried out.
This paper was stimulated and inpired by Roland Huber's work
\cite{Hu3}.

\section{Algebraic preliminaries}

\subsection{Power-multiplicative seminorms}
Real-valued valuations on fields and topological algebras have
been standard tools in $p$-adic analytic geometry since the
earliest infancy of the subject; by contrast, the role played
by higher rank valuations in several fundamental questions has
been recognized only in recent times.

In non-archimedean analysis one encounters more generally certain
ultrametric real-valued norms (or seminorms) that are not multiplicative,
but only power-multiplicative. We shall see that higher rank
power-multiplicative seminorms appear just as naturally, and
are just as useful.

\sset\subsubsection{}\label{subsec_gamma-plus}
In this section we let $(\Gamma,<)$ be a totally ordered abelian
group, whose neutral element is denoted by $1$ and whose composition
law we write multiplicatively. As customary, we shall extend the
ordering and the composition law of $\Gamma$ to the set $\Gamma\cup\{0\}$,
in such a way that $\gamma>0$ for every $\gamma\in\Gamma$, and
$\gamma\cdot 0=0\cdot\gamma=0$ for every $\gamma\in\Gamma\cup\{0\}$.
Notice also that the ordering of $\Gamma$ extends uniquely to
$\Gamma_\Q:=\Gamma\otimes_\Z\Q$.
Finally we let $\Gamma^+:=\{\gamma\in\Gamma~|~\gamma\leq 1\}$

\begin{definition}\label{def_Gamma-norm}
Let $A$ be a ring. A {\em power-multiplicative $\Gamma$-valued seminorm}
on $A$ is a map
$$
|\cdot|:A\to\Gamma\cup\{0\}
$$
satisfying the following conditions:
\begin{enumerate}
\renewcommand{\labelenumi}{(\alph{enumi})}
\item
$|0|=0$ and $|1|=1$.
\item
$|x-y|\leq\max(|x|,|y|)$ for all $x,y\in A$.
\item
$|x\cdot y|\leq|x|\cdot|y|$ for all $x,y\in A$.
\item
$|x^n|=|x|^n$ for every $x\in A$ and every $n\in\N$.
\renewcommand{\labelenumi}{(\roman{enumi})}
\end{enumerate}
One says also that $(A,|\cdot|)$ is a {\em $\Gamma$-seminormed ring}.
If $|x|\neq 0$ whenever $x\neq 0$, one says that $|\cdot|$ is a
{\em $\Gamma$-valued norm\/} and correspondingly one talks of
$\Gamma$-normed rings. If in (c) we have equality for every
$x,y\in A$, then we say that $(A,|\cdot|)$ is a {\em $\Gamma$-valued
ring\/} and that $|\cdot|$ is a {\em $\Gamma$-valued valuation}.
A morphism $\phi:(A,|\cdot|)\to(A',|\cdot|')$ of $\Gamma$-seminormed
rings is a ring homomorphism $\phi:A\to A'$ such that
$|\phi(a)|'=|a|$ for every $a\in A$. Notice that the subset
\set\begin{equation}\label{eq_notation+}
A^+:=\{a\in A~|~|a|\leq 1\}\subset A
\end{equation}
is a subring of $A$. The {\em support} of $|\cdot|$ is the ideal
$$
\mathrm{supp}(|\cdot|):=\{a\in A~|~|a|=0\}\subset A.
$$
If $|\cdot|$ is a valuation, $\mathrm{supp}(|\cdot|)$ is a prime
ideal.
\end{definition}

\begin{lemma}\label{lem_strict.ineq}
Let $(A,|\cdot|)$ be a seminormed ring and $x,y\in A$ any two
elements such that $|x|\neq|y|$. Then $|x+y|=\max(|x|,|y|)$.
\end{lemma}
\begin{proof} Let us say that $|x|<|y|$. By (b) of definition
\ref{def_Gamma-norm} we have:
$$
|y|\leq\max(|x+y|,|x|)\leq\max(|y|,|x|)=|y|.
$$
Hence $|y|=|x+y|$, which is the claim.
\end{proof}

\sset\subsubsection{}
Let $(A,|\cdot|)$ be a semi-normed ring. For every monic polynomial
$p(T):=T^m+a_1T^{m-1}+\cdots+a_m\in A[T]$ we set
$$
\sigma(p):=\max(|a_i|^{1/i}~|~i=1,\dots,m)\in\Gamma_\Q\cup\{0\}
$$
and call $\sigma(p)$ the {\em spectral value\/} of $p(T)$.

\begin{lemma}\label{lem_specter} Let $p,q\in A[T]$ be monic
polynomials. Then $\sigma(pq)\leq\max(\sigma(p),\sigma(q))$.
If $|\cdot|$ is a valuation, the above inequality is, in fact,
an equality.
\end{lemma}
\begin{proof} {\em Mutatis mutandis}, this is the same as the
proof of \cite[\S1.5.4, Prop.1]{Bo-Gun}.
\end{proof}

\sset\subsubsection{}\label{subsec_minimal}
Let $A$ be a normal domain, $K$ the field of fractions of $A$,
and $A\to B$ an injective integral ring homomorphism such that
$B$ is torsion-free as an $A$-module. For any
$b\in B\otimes_AK$ the set of polynomials $P(T)\in K[T]$ such
that $P(b)=0$ is an ideal, whose monic generator $\mu_b(T)$
is the {\em minimal polynomial} of $b$.

\begin{lemma}\label{lem_integrality} Keep the assumptions of
\eqref{subsec_minimal} and suppose that $b\in B$. Then
$\mu_b(T)\in A[T]$.
\end{lemma}
\begin{proof} By \cite[Th.10.4]{Mat} we have $A=\bigcap_vA_v$
where $v$ ranges over all the valuations of $K$ whose valuation
ring $A_v$ contains $A$. We can therefore suppose that $A$ is
the valuation ring of one such $v:A\to\Gamma_v\cup\{0\}$. For
any polynomial $p(T):=\sum_{i=0}^na_iT^i\in K[T]$
let
$$
|p|_v:=\max(v(a_i)~|~i=0,\dots,n)\in\Gamma_v\cup\{0\}.
$$
By assumption, there is a monic polynomial $P(T)\in A[T]$ such
that $P(b)=0$, hence $\mu_b$ divides $P$ in $K[T]$. The assertion
is therefore a consequence of the following:
\begin{claim} $|p\cdot q|_v=|p|_v\cdot |q|_v$ for all $p,q\in K[T]$.
\end{claim}
\begin{pfclaim}[] We leave it as an exercise for the reader: one
can easily adapt the direct argument used in the classical proof
of the Gauss lemma (cp. \cite[p.44]{Bo-Gun}).
\end{pfclaim}
\end{proof}

\sset\subsubsection{}\label{subsec_spectralnorm}
In the situation of \eqref{subsec_minimal}, let $b\in B$ and
say that
$$
\mu_b(T)=T^n+a_1T^{n-1}+a_2T^{n-2}+\cdots+a_n
$$
for some $n\in\N$ and elements $a_1,\dots,a_n\in A$, in view of lemma
\ref{lem_integrality}. Suppose now that $|\cdot|:A\to\Gamma\cup\{0\}$
is a power multiplicative $\Gamma$-valued seminorm on $A$.
Then the {\em spectral seminorm\/} of $b$ is defined as
$$
|b|_\sp:=\sigma(\mu_b(T))\in\Gamma_\Q\cup\{0\}.
$$
The name of $|\cdot|_\sp$ is justified by the following:

\begin{proposition}\label{prop_spectral-norm}
With the notation of \eqref{subsec_spectralnorm},
the pair $(B,|\cdot|_\sp)$ is a $\Gamma_\Q$-seminormed ring.
\end{proposition}
\begin{proof} We consider the $A$-algebra $A':=A[\Gamma_\Q]$.
Hence $A'$ is generated as an $A$-algebra by elements $[\gamma]$,
for all $\gamma\in\Gamma_\Q$, subject to the relations
$$
[\gamma]\cdot[\delta]=[\gamma\cdot\delta]\qquad
\text{for all $\gamma,\delta\in\Gamma_\Q$}.
$$
Clearly, every element of $A'$ admits a unique expression
of the form $\sum_{\gamma\in\Gamma_\Q}a_\gamma\cdot[\gamma]$
where $a_\gamma=0$ for all but finitely many $\gamma\in\Gamma_\Q$.
Set also $B':=B\otimes_AA'$.
\begin{claim}\label{cl_A'.is.good}
$A'$ is a normal domain, flat over $A$, and $B'$ is a torsion-free
$A'$-module.
\end{claim}
\begin{pfclaim} We can write
$\Gamma_\Q=\bigcup_{i\in I}\Gamma_i$, the filtered union
of all the finitely generated subgroups $\Gamma_i\subset\Gamma_\Q$.
Then $A'=\bigcup_{i\in I}A[\Gamma_i]$, and it suffices to show the
claim for the subalgebras $A'_i:=A[\Gamma_i]$. Since $\Gamma$
is torsion-free, we have (non-canonical) isomorphisms
$\Gamma_i\simeq\Z^{N_i}$, whence $A'_i\simeq
A[T^{\pm 1}_1,\dots,T^{\pm 1}_N]\simeq
A\otimes_\Z\Z[T^{\pm 1}_1,\dots,T^{\pm 1}_N]$,
from which flatness is clear. Normality follows as well, in
view of \cite[Ch.IV, Prop.6.14.1]{EGAIV-2}. Likewise, $B'$
is the filtered union of the $A'_i$-algebra
$B[T^{\pm 1}_1,\dots,T^{\pm 1}_N]$, and the latter are
clearly torsion-free over $A'_i$.
\end{pfclaim}

We define a map $|\cdot|:A'\to\Gamma_\Q$ by the rule:
$$
\sum_{\gamma\in\Gamma_\Q}a_\gamma\cdot[\gamma]\mapsto
\max(|a_\gamma|\cdot\gamma~|~\gamma\in\Gamma_\Q).
$$
\begin{claim}\label{cl_replace.by.A'}
$(A',|\cdot|)$ is a $\Gamma_\Q$-seminormed ring.
\end{claim}
\begin{pfclaim} All conditions of definition \ref{def_Gamma-norm}
are clearly fulfilled, except possibly for (d). However, for given
$x=\sum_\gamma a_\gamma\cdot[\gamma]$, let $\delta\in\Gamma$
be the minimal element such that $|a_\delta|\cdot\delta=|x|$.
Say that $x^n=\sum_\gamma b_\gamma\cdot[\gamma]$; we have
$$
b_{\delta^n}\cdot[\delta^n]=
\sum_{\gamma_1{\cdot}\dots{\cdot}\gamma_n=\delta^n}
a_{\gamma_1}\cdot[\gamma_1]{\cdot}\dots{\cdot} a_{\gamma_n}\cdot[\gamma_n].
$$
If now $\gamma_1{\cdot}\dots{\cdot}\gamma_n=\delta^n$ and the $\gamma_i$
are not all equal to $\delta$, then necessarily $\gamma_i<\delta$
for some $i\leq n$. By the choice of $\delta$ it follows that
$|a_{\gamma_i}|\cdot\gamma_i<|a_\delta|\cdot\delta$, hence
$|a_{\gamma_1}\cdot[\gamma_1]{\cdot}\dots{\cdot}
a_{\gamma_n}\cdot[\gamma_n]|<|x|^n$.
In view of lemma \ref{lem_strict.ineq} we deduce that
$|b_{\delta^n}|\cdot\delta^n=|x|^n$, so (d) holds as well.
\end{pfclaim}

From claim \ref{cl_A'.is.good} it follows that the induced map
$A'\to B'$ is injective, and by claim \ref{cl_replace.by.A'}
we can replace $A$, $B$ and $\Gamma$ by respectively $A'$, $B'$
and $\Gamma_\Q$, which allows us to assume that
\set\begin{equation}\label{eq_monoids}
\parbox{14cm}{$|A|=|B|_\sp=\Gamma\cup\{0\}=\Gamma_\Q\cup\{0\}$
and there is a group homomorphism $[\cdot]:\Gamma\to A^\times$ which
is a left inverse for $|\cdot|:A\to\Gamma\cup\{0\}$.
}\end{equation}
Let $B^+:=\{x\in B~|~|x|_\sp\leq 1\}$.

\begin{claim}\label{cl_integral.plus}
$A^+$ is normal and $B^+$ is the integral closure of $A^+$
in $B$.
\end{claim}
\begin{pfclaim} To show that $A^+$ is normal, it suffices to
prove that $A^+$ is integrally closed in $A$. Hence, suppose that
$x\in A$ satisfies an equation of the form $x^n+a_1x^{n-1}+\cdots+a_n=0$,
where $|a_i|\leq 1$ for $i=1,\dots,n$. It follows that
$|x|^n\leq\max(|x|^{n-i}~|~i=1,\dots,n)$, which is possible only when
$|x|\leq 1$, as required. Next, let $b\in B^+$; by definition this
means that $\mu_b(T)\in A^+[T]$, so $x$ is integral over $A^+$.
Conversely, we apply lemma \ref{lem_integrality} with $A^+$ in place
of $A$, to see that $|b|_\sp\leq 1$ for every $b\in B$ which is
integral over $A^+$.
\end{pfclaim}

Finally, we verify conditions (a)-(d) of definition \ref{def_Gamma-norm}.
(a) is obvious. Let $x,y\in B$ and say that $\delta:=|x|\leq\gamma:=|y|$;
thanks to \eqref{eq_monoids} we have 
$$
|x\cdot[\gamma^{-1}]|\leq|x|\cdot\gamma^{-1}\leq 1
$$
and likewise $|y\cdot[\gamma^{-1}]|\leq 1$. By claim \ref{cl_integral.plus}
it follows that $x\cdot[\gamma^{-1}]$ and $y\cdot[\gamma^{-1}]$
are integral over $A^+$, hence the same holds for their sum and
again claim \ref{cl_integral.plus} implies that
$|(x+y)\cdot[\gamma^{-1}]|\leq 1$. Consequently
$|x+y|\leq|(x+y)\cdot[\gamma^{-1}]|\cdot\gamma\leq|y|$, which is (b).
For (c) one considers the product 
$x\cdot[\delta^{-1}]\cdot y\cdot[\gamma^{-1}]$ which is integral over
$A^+$ by an analogous argument; then 
$|x\cdot y|\leq|x\cdot[\delta^{-1}]\cdot
y\cdot[\gamma^{-1}]|\cdot\delta\cdot\gamma\leq|x|\cdot|y|$, which is (c).
Finally, suppose that $|x^n|=\eps<\delta^n$; by \eqref{eq_monoids}
the value group $\Gamma$ is divisible, hence we can consider the
element $z:=x\cdot[\eps^{-1/n}]$ and in fact $|z^n|\leq 1$, hence
$z^n$ is integral over $A^+$, so the same holds for $z$, and again
$|z|\leq 1$, therefore $|x|\leq\eps^{1/n}<\delta$, a contradiction
that shows (d). 
\end{proof}

\begin{remark} (i) Proposition \ref{prop_spectral-norm} generalizes
\cite[\S3.2.2, Th.2]{Bo-Gun}, which deals with the special case
of real-valued norms. The proof of {\em loc.cit.} does not extend
to the present case, since it is based on a smoothing technique
that makes sense only for real-valued seminorms.

(ii) In the situations encountered in later sections,
it is probably not too hard to verify directly that the spectral norm is
power-multiplicative (the same can already be said for most applications
of \cite[\S3.2.2, Th.2]{Bo-Gun}). However, it seems desirable to
have a general statement such as proposition \ref{prop_spectral-norm}.
\end{remark}

\begin{lemma}\label{lem_nobigdeal}
In the situation of \eqref{subsec_spectralnorm} :
\begin{enumerate}
\item
Let $(A,|\cdot|)\to(A',|\cdot|')$ be a flat morphism of\/
$\Gamma$-seminormed normal domains, suppose that $B':=A'\otimes_AB$
is torsion-free over $A'$, and endow it with the spectral seminorm
$|\cdot|'_\sp$ relative to the induced injective ring homomorphism
$A'\to B'$. Then:
\begin{enumerate}
\item
The natural map $(B,|\cdot|_\sp)\to(B',|\cdot|'_\sp)$ is a morphism
of\/ $\Gamma_\Q$-seminormed rings. 
\item
If\/ $|\cdot|:A\to\Gamma\cup\{0\}$ is a valuation, we
have $|ab|_\sp=|a|\cdot|b|_\sp$ for every $a\in A$ and $b\in B$.
\end{enumerate}
\item
Suppose that $B=B_1\times\cdots\times B_r$, where each $B_i$ is
an $A$-algebra fulfilling the conditions of \eqref{subsec_minimal},
and for $i=1,\dots,r$, denote by $|\cdot|_{\sp,i}$ the spectral
norm of $B_i$. Then, for every $b:=(b_1,\dots,b_r)\in B$ we have:
$$
|b|_\sp\leq\max(|b_i|_{\sp,i}~|~i=1,\dots,r)
$$
and if the norm of $A$ is a valuation, the inequality is actually
an equality.
\item
Suppose that $(A,|\cdot|)\to(B,|\cdot|_B)$ is an extension of
valuation rings, such that $B$ is integral over $A$. Then the
spectral norm $|\cdot|_\sp$ is a valuation equivalent to $|\cdot|_B$.
\end{enumerate}
\end{lemma}
\begin{proof} (i.a): For given $b\in B$, let $C\subset B$ be a
finite $A$-subalgebra with $b\in C$; then $C':=A'\otimes_AC\subset B'$.
Let $K$ and $K'$ be the fields of fraction of $A$ and respectively
$A'$; the element $b$ induces a $K$-linear (resp. $K'$-linear)
endomorphisms on the finite dimensional $K$-vector space
(resp. $K'$-vector space) $C\otimes_AK$ (resp. $C'\otimes_{A'}K'$),
and the spectral seminorms of $b$ in $B$ and $B'$ are defined
in terms of the minimal polynomials of these endomorphisms.
Hence the assertion boils down to the invariance of the minimal
polynomial under base field extensions.

(i.b): If $\mu_b(T)=T^n+a_1T^{n-1}+a_2T^{n-2}+\cdots+a_n$, then 
$\mu_{ab}(T)=
T^n+a\cdot a_1T^{n-1}+a^2\cdot a_2T^{n-2}+\cdots+a^n\cdot a_n$,
from which the assertion follows easily.

(ii): The minimal polynomial of $b$ is the least common multiple
of the minimal polynomials of $b_1,\dots,b_r$; hence the claim
follows from lemma \ref{lem_specter}.

(iii): Let $b\in B$ such that $|b|_B=1$, and 
$\mu_b(T):=T^n+\sum^n_{i=1}a_iT^{n-i}$ the minimal polynomial
of $b$ over $\mathrm{Frac}(A)$. Since $b$ is integral over $A$, we
have $a_i\in A$ for every $i=1,\dots,n$, hence $|b|_\sp\leq 1$;
on the other hand, since :
$$
1=|b|^n_B=
|\sum^n_{i=1}a_ib^{n-i}|_B\leq\max(|a_i|~|~i=1,\dots,n)
$$
we have as well : $1\leq|b|_\sp$, so that $|b|_\sp=1$. Finally,
for a general element $b\in B\setminus\{0\}$, we may find
$s\in\N\setminus\{0\}$ and $a\in A$ such that $|b^s\cdot a|_B=1$,
hence $|b|_B^s\cdot|a|=1=|b|^s_\sp\cdot|a|$, in view of (i.b).
The assertion follows.
\end{proof}

\subsection{Normed modules}
Throughout this section $(A,|\cdot|)$ denotes a $\Gamma$-normed
ring (for some ordered abelian group $\Gamma$). Following
(\cite[\S2.1.1, Def.1]{Bo-Gun}), we shall say that a {\em faithfully
$\Gamma$-seminormed $A$-module\/} is a pair $(V,|\cdot|_V)$ consisting of an
$A$-module $V$ and a map $|\cdot|_V:V\to\Gamma\cup\{0\}$
such that:
\begin{enumerate}
\item
$|x-y|_V\leq\max(|x|_V,|y|_V)$ for every $x,y\in V$.
\item
$|ax|=|a|\cdot|x|_V$ for every $a\in A$ and $x\in V$.
\end{enumerate}
If moreover $|\cdot|_V$ satisfies also the axiom:

(iii)\;\ $|x|_V=0$ if and only if $x=0$ \\
then we say that $(V,|\cdot|_V)$ is a faithfully $\Gamma$-normed
$A$-module.
In the following we will suppose that all the $\Gamma$-seminormed
$A$-modules under consideration are faithfully seminormed, so we
shall refer to them simply as ``seminormed $A$-modules''
(or ``normed $A$-modules'' if (iii) holds).

\sset\subsubsection{}\label{subsec_A-cartesian}
Let $V:=(V,|\cdot|_V)$ be a free seminormed $A$-module of finite rank.
Following \cite[\S2.4.1, Def.1]{Bo-Gun}, we say that
$V$ is {\em $A$-cartesian\/} if there exists a basis $\{v_1,\dots,v_n\}$
of $V$ such that
$$
\left|\sum_{i=1}^na_iv_i\right|_V=\max_{1\leq i\leq n}|a_i|\cdot|v_i|_V
$$
for all $a_1,\dots,a_n\in A$. A basis with this property is called
{\em $A$-orthogonal} (or just orthogonal).

\sset\subsubsection{}\label{subsec_immediate}
Suppose that $(A,|\cdot|)$ is a valuation ring. Recall that
an {\em immediate extension\/} of $A$ is a flat morphism of valuation
rings $(A,|\cdot|)\to(A',|\cdot|')$ inducing an isomorphism of
value groups $\Gamma\stackrel{\sim}{\to}\Gamma'$ and residue fields
$A^\sim\stackrel{\sim}{\to} A^{\prime\sim}$. For instance, the
henselization $(A^h,|\cdot|^h)$ of $(A,|\cdot|)$ is an immediate
extension; also the completion $(A^\wedge,|\cdot|^\wedge)$ of
$(A,|\cdot|)$ relative to its valuation topology, is an immediate
extension.

\begin{lemma}\label{lem_immediate}
Let $(A,|\cdot|)$ be a valuation ring, $(A',|\cdot|')$ an
immediate extension of $A$, and $B$ a flat $A$-algebra; set
$B':=A^\wedge\otimes_AB$. We endow $B$ (resp. $B'$) with the
spectral seminorm $|\cdot|_\sp$ (resp. $|\cdot|'_\sp$) relative
to the valuation of $A$ (resp. of $A'$). Suppose furthermore
that both $B$ and $B'$ are reduced. Then:
\begin{enumerate}
\item
If $(B,|\cdot|_\sp)$ is $A$-cartesian, then $(B',|\cdot|'_\sp)$
is $A'$-cartesian. 
More precisely, a subset $\{b_1,\dots,b_d\}$ is an orthogonal
basis of $B$ if and only if $\{1\otimes b_1,\dots,1\otimes b_d\}$
is an orthogonal basis of $B'$.
\item
Conversely, suppose that $(B',|\cdot|'_\sp)$ is $A'$-cartesian,
and assume that $(A',|\cdot|')\subset(A^\wedge,|\cdot|^\wedge)$,
the completion of $A$ for the valuation topology.
Then $(B,|\cdot|_\sp)$ is $A$-cartesian.
\end{enumerate}
\end{lemma}
\begin{proof} Notice that $B$ is free of finite rank over $A$
if and only if $B^\wedge$ is free of finite rank over $A^\wedge$
(\cite[Rem.3.2.26(ii)]{Ga-Ra}), hence we may assume from start
that $B$ is free of finite rank.

(i): Suppose that $\{b_1,\dots,b_d\}$ is an orthogonal
basis of $B$; for every $a'_1,\dots,a'_d\in A^\wedge$ such that
$x:=\sum_{i=1}^d a_i'\otimes b_i\neq 0$, we have $|x|'_\sp\neq 0$,
since by assumption $B^\wedge$ is reduced. Since $A'$ is an immediate
extension of $A$, we can find $a_1,\dots,a_d\in A$ such that
\set\begin{equation}\label{eq_slight-approx}
\text{either}\quad a_i=0\qquad\text{or}\qquad
|a_i-a'_i|<|a'_i|\qquad\text{ for every $i\leq d$}
\end{equation}
and especially, $|a_i|=|a'_i|$ for every $i\leq d$. 
Set $y:=\sum_{i=1}^da_ib_i$; we deduce :
$|x-1\otimes y|'_\sp<\max(|a_i|\cdot|b_i|_\sp)=|y|_\sp=|1\otimes y|'_\sp$,
by lemma \ref{lem_nobigdeal}(i.a), whence :
$$
|x|'_\sp=|y|_\sp=\max(|a'_i|\cdot|1\otimes b_i|'_\sp~|~i=1,\dots,d)
$$
{\em i.e.} $\{1\otimes b_1,\dots,1\otimes b_d\}$ is an orthogonal basis.
Conversely, if $\{1\otimes b_1,\dots,1\otimes b_d\}$ is orthogonal,
obviously $\{b_1,\dots,b_d\}$ is orthogonal in $B$.

(ii): In view of (i), we can assume that $B^\wedge$ is
$A^\wedge$-cartesian, and it remains to show that $B$ is $A$-cartesian.
Choose an orthogonal basis $e'_1,\dots,e'_d$ for $B^\wedge$.
We shall use the following analogue of \eqref{eq_slight-approx} :

\begin{claim} We can find $e_1,\dots,e_n\in B$ such that
$|e_i-e'_i|'_\sp<|e'_i|'_\sp$ for $i=1,\dots,d$.
\end{claim}
\begin{pfclaim} Write $e'_i=\sum_{j=1}^d a'_j\otimes b_j$ for
some $a'_1,\dots,a'_d\in A^\wedge$ and $b_1,\dots,b_d\in B$;
choose approximations $a_1,\dots,a_d\in A^\wedge$ of these
elements and set $e_i:=\sum^d_{j=1}a_j b_j$. By lemma
\ref{lem_nobigdeal}(i.b) we have :
$|e_i-e'_i|_\sp\leq\max(|a_j-a'_j|\cdot|b_j|_\sp~|~j\leq d)$,
which can be made arbitrarily small.
\end{pfclaim}

It follows that
$|e_i|_\sp=|e'_i|'_\sp$ and $(e_i~|~i=1,\dots,d)$ is a basis of $B$
(by Nakayama's lemma); furthermore, for every $a_1,\dots,a_d\in A$
we have:
$$
\begin{array}{r@{\:}l}
|\sum\limits_{i=1}^d a_ie_i-\sum\limits_{i=1}^d a_ie'_i|'_\sp=
|\sum\limits_{i=1}^d a_i(e_i-e'_i)|'_\sp & \leq
\max(|a_i|\cdot|e_i-e'_i|'_\sp~|~i=1,\dots,d) \\
& <\max(|a_i|\cdot|e'_i|'_\sp~|~i=1,\dots,d)=
   |\sum\limits_{i=1}^d a_ie'_i|'_\sp.
\end{array}
$$
Hence:
$$
\begin{array}{c}
|\sum\limits_{i=1}^d a_ie_i|_\sp=|\sum\limits_{i=1}^d a_ie'_i|'_\sp=
\max(|a_i|\cdot|e'_i|'_\sp~|~1\leq i\leq d)=
\max(|a_i|\cdot|e_i|_\sp~|~1\leq i\leq d).
\end{array}
$$
In other words, the basis $(e_i~|~i=1,\dots,d)$ is orthogonal.
\end{proof}

\begin{remark} I do not know whether lemma \ref{lem_immediate}(ii)
holds for an arbitrary immediate extension $(A,|\cdot|)\to(A',|\cdot|')$.
If the rank of the valuation ring $A$ is greater than one,
this seriously limits the usefulness of lemma \ref{lem_immediate},
since for instance, for such valuations, the henselization $A^h$
of $A$ is not necessarily contained in $A^\wedge$.
\end{remark}

\sset\subsubsection{}\label{subsec_artificial}
Let $(K,|\cdot|)$ be a valued field with value group $\Gamma$
and residue field $K^\sim$, $L$ a finite extension of $K$ and 
$$
|\cdot|_i:L\to\Gamma_i\cup\{0\}\qquad i=1,\dots,k
$$
the finitely many extensions of $|\cdot|$. For every $i\leq k$,
let $L^\sim_i$ be the residue field of the valuation ring $L^+_i$
of $(L,|\cdot|_i)$, and set $f_i:=[L^\sim_i:K^\sim]$,
$e_i:=(\Gamma_i:\Gamma)$. Furthermore, for every pair of
integers $i,j\leq k$ let $\Gamma_{ij}$ be the value group
of the valuation ring $L^+_{ij}:=L^+_i\cdot L_j^+$; the embedding
$L_i\subset L_{ij}$ induces a natural surjective order-preserving
group homomorphisms $\Gamma_i\to\Gamma_{ij}$, whose kernel we
denote $\Delta_{ij}$. Then we have natural isomorphisms of
ordered groups $\Gamma_i/\Delta_{ij}\simeq\Gamma_j/\Delta_{ji}$,
for every such pair $(i,j)$. For every $i\leq k$ set also
$\Delta_i:=\bigcap_{j\neq i}\Delta_{ij}$.

\begin{proposition}\label{prop_mess}
In the situation of \eqref{subsec_artificial}, endow the
$K$-algebra $L$ with its spectral norm $|\cdot|_\sp$ (relative
to the norm $|\cdot|$ on $K$), and suppose moreover that :
\begin{enumerate}
\alphaenu
\item
The extension $K\subset L$ is {\em defectless}, {\em i.e.}
$\sum_{i=1}^d e_i\cdot f_i=[L:K]$.
\item
$\Gamma+\Delta_i=\Gamma_i$\ \ for every $i\leq k$.
\item
For every $i\leq k$, the quotient
$\Delta_i/(\Gamma\cap\Delta_i)$ consists of equivalence
classes $\bar\alpha_{i1},\dots,\bar\alpha_{ie_i}$ of elements
of $\Delta_i$ :
$$
\alpha_{i1}:=1>\alpha_{i2}>\cdots>\alpha_{ie_i}
$$
such that $\alpha_{ij}>\gamma$ for every $i\leq k$, every
$j\leq e_i$, and every $\gamma\in\Gamma^+\setminus\{1\}$.
\item
$\Delta_i\neq\{1\}$ for every $i\leq k$. 
\end{enumerate}
Then $(L,|\cdot|_\sp)$ is a cartesian $(K,|\cdot|)$-module.
\end{proposition}
\begin{proof} Assumption (c) and \cite[Th.5]{Rib} imply
that for every $i\leq k$ and every $j\leq e_i$ we may find
$x_{ij}\in L$ such that :
\set\begin{equation}\label{eq_approx-1}
|x_{ij}|_i=\alpha_{ij}
\qquad\text{and}\qquad
|x_{ij}|_l=1\quad\text{for every $l\neq i$}.
\end{equation}
Next, for every $i\leq k$, let $c_{i1},\dots,c_{if_i}\in L_i^+$ whose
equivalence classes form a basis of the $K^\sim$-vector space
$L^\sim_i$. In view of assumption (d) we may find, for every
$i\leq k$, a non-zero element $b_i\in K$ such that :
$$
|b_i|\in\Gamma^+\cap\Delta_i\setminus\{1\}.
$$
Then, according to \cite[Lemme 9]{Rib} we may find, for every
$i\leq k$ and every $j\leq f_i$, elements $y_{ij}\in L$ such that :
\set\begin{equation}\label{eq_approx-2}
|y_{ij}-c_{ij}|_i\leq|b_i|\qquad\text{and}\qquad
|y_{ij}|_l\leq|b_i|\quad\text{for every $l\neq i$}.
\end{equation}
The proposition now follows from assumptions (a) and (b),
and the following :

\begin{claim} The family :
$$
\Sigma:=(x_{ij}y_{il}~|~i=1,\dots,k;~j=1,\dots,e_i;~l=1,\dots,f_i)
$$
is orthogonal relative to the spectral norm $|\cdot|_\sp$.
\end{claim}
\begin{pfclaim}[] Let $a\in L$ be an element which can be written
in the form :
$$
a=\sum_{i=1}^k a_i\qquad\text{where :}\qquad
a_i=\sum_{j=1}^{e_i}\sum_{l=1}^{f_i}a_{ijl}x_{ij}y_{il}
\quad\text{for every $i=1,\dots,k$}
$$
with $a_{ijl}\in K$. We have to show that :
$$
|a|_\sp=\max(|a_{ijl}|\cdot|x_{ij}y_{il}|_\sp
~|~i=1,\dots,k;~j=1,\dots,e_i;~l=1,\dots,f_i).
$$
Let $(K^h,L^h)$ be the henselization of
$(K,|\cdot|)$, and set $L^{h+}:=L^+\otimes_{K^+}K^{h+}$,
where $L^+:=\{x\in L~|~|x|_\sp\leq 1\}$ is the integral
closure of $K^+$ in $L$. In view of lemma \ref{lem_nobigdeal}(i.a),
it suffices to verify that $\Sigma$ is an orthogonal system
of elements of $(L^h,|\cdot|^h_\sp)$, where $L^h=L^{h+}\otimes_{K^+}K$
and $|\cdot|_\sp$ is the spectral norm of $L^h$ relative to
$|\cdot|^h$. However, $L^h=L^h_i\times\cdots\times L^h_k$,
where $(L^h_i,|\cdot|^h_i)$ denotes the henselization of
$(L_i,|\cdot|_i)$, for every $i\leq k$, hence lemma
\ref{lem_nobigdeal}(ii) yields the identity :
\set\begin{equation}\label{eq_yields-the-ident}
|a|_\sp=\max(|a|_i~|~i=1,\dots,k).
\end{equation}
For every $i\leq k$ let us set :
$$
\gamma_i:=\max(|a_{ijl}|~|~j=1,\dots,e_i;~l=1,\dots,f_i).
$$
In view of lemma \ref{lem_nobigdeal}(i.b) we may assume that :
$$
\max(\gamma_i~|~i=1,\dots,k)=1.
$$
Fix $i\leq k$. By inspecting \eqref{eq_approx-1} and
\eqref{eq_approx-2} we see that in that case :
$$
|a_t|_i\leq|b_i|\qquad\text{for every $t\neq i$}.
$$
whence :
\set\begin{equation}\label{eq_go-for-it}
|a|_i\leq\max(|a_i|_i,|b_i|)
\qquad\text{and equality holds if $|a_i|_i>|b_i|$.}
\end{equation} 

$\bullet$\ \ 
Now, suppose first that $\gamma_i=1$, and denote by $j_0$ the
smallest integer $j\leq e_i$ such that $|a_{ijl}|=1$ for some
$l\leq f_i$. With this notation, we may decompose :
$$
a_i=a'_i+x_{ij_0}\cdot\sum_{l=1}^{f_i}a_{ij_0l}y_{il}
$$
where $a'_i$ is the sum of the terms $a_{ijl}x_{ij}y_{il}$
with $j\neq j_0$. By \eqref{eq_approx-2}, the residue classes
in $L^\sim_i$ of the elements $y_{i1},\dots,y_{if_i}$ are all
distinct, hence \eqref{eq_approx-1} yields :
$$
|x_{ij_0}\cdot\sum_{l=1}^{f_i}a_{ij_0l}y_{il}|_i=
\alpha_{ij_0}
$$
and, on the other hand :
$$
|a'_i|_i\leq\alpha_{i,j_0+1}
$$
since $\alpha_{i,j_0+1}$ is greater than any element in
$\Gamma^+\setminus\{1\}$. Summing up we obtain :
\set\begin{equation}\label{eq_summing-uppo}
|a_i|_i=\alpha_{ij_0}=
\max(|a_{ijl}|\cdot|x_{ij}y_{il}|_\sp~|~j=1,\dots,e_i;~l=1,\dots,f_i)
\qquad\text{if $\gamma_i=1$}
\end{equation}
and combining with \eqref{eq_go-for-it} :
\set\begin{equation}\label{eq_eval=1}
|a|_i=\alpha_{ij_0}\qquad\text{whenever $\gamma_i=1$}.
\end{equation}
$\bullet$\ \  Suppose next that $\gamma_i<1$. Then
$|a_i|_i\leq\gamma_i$, and in view of \eqref{eq_go-for-it}
we deduce :
\set\begin{equation}\label{eq_eval<1}
|a|_i\leq\max(|\gamma_i|,|b_i|)\qquad\text{whenever $\gamma_i<1$}.
\end{equation}
Finally, from \eqref{eq_yields-the-ident}, \eqref{eq_eval=1}
and \eqref{eq_eval<1} we conclude that, in order to evaluate
$|a|_\sp$ we may neglect all the terms $a_{ijl}$ such that $\gamma_i<1$,
and then the sought identity follows from \eqref{eq_summing-uppo}.
\end{pfclaim}
\end{proof}

\begin{remark} Of the four conditions of proposition \ref{prop_mess},
the first one is very natural, and in fact characterizes cartesian
extensions of a valued fields, in the rank one case (see
\cite[\S3.6.2, Prop.5]{Bo-Gun}). On the other hand, conditions
(b), (c) and (d) appear (to me) as artifical, and certainly
leave room for improvements.
\end{remark}

\sset\subsubsection{}
Suppose that $V:=(V,|\cdot|_V)\neq 0$ and $W:=(W,|\cdot|_W)$ are
two normed $A$-modules; let $\psi:V\to W$ be an $A$-linear homomorphism.
We say that $\psi$ is {\em bounded\/} if there exists $\gamma\in\Gamma$
such that 
$$
|\psi(x)|_W/|x|_V\leq\gamma\qquad \text{for every $x\in V\setminus\{0\}$}.
$$
We denote by $\cL(V,W)$ the $A$-module of all bounded $A$-linear
homomorphisms $V\to W$.
If $\Gamma\subset\R$, then one can define the norm of $\psi$ as the
supremum of $|\psi(x)|_W/|x|_V$ for $x$ ranging over all
$x\in V\setminus\{0\}$ (\cite[\S2.1.6]{Bo-Gun}).
For more general groups $\Gamma$, this quantity is not necessarily
defined. Hence, for $\psi\in\cL(V,W)$ we shall set
$$
|\psi|_\cL:=\sup_{x\in V\setminus\{0\}}\frac{|\psi(x)|_W}{|x|_V}
$$
{\em whenever this is well defined as an element of $\Gamma\cup\{0\}$}.
Lemma \ref{lem_descartes} shows that, if $V$ and $W$ are $A$-cartesian,
the norm $|\psi|_\cL$ is well defined for every $A$-linear map $\psi$.

\begin{lemma}\label{lem_descartes}
Let $V:=(V,|\cdot|_V)$ and $W:=(W,|\cdot|_W)$ be two free 
$A$-cartesian normed $A$-modules of finite rank. Then:
\begin{enumerate}
\item
$\cL(V,W)=\Hom_k(V,W)$ and the pair $(\cL(V,W),|\cdot|_\cL)$
is an $A$-cartesian normed $A$-module.
\item
If $(v_i~|~i=1,\dots,n)$ and $(w_j~|~j=1,\dots,m)$ are orthogonal basis
of\/ $V$, resp. $W$, then the basis
$(v_i^*\otimes w_j~|~i=1,\dots,n; j=1,\dots,m)$ of\/ $\cL(V,W)$ is
orthogonal.
\end{enumerate}
\end{lemma}
\begin{proof} (i) will follow from the more precise assertion (ii).
The basis in (ii) is characterized by the identities
$$
v_i^*\otimes w_j(v_k)=\delta_{ik}\cdot w_j\qquad
\text{for every $i=1,\dots,n$ and $j=1,\dots,m$.}
$$
\begin{claim}\label{cl_v_i.w_k}
$|v_i^*\otimes w_j|_\cL=|w_j|_W/|v_i|_V$.
\end{claim}
\begin{pfclaim} By definition we have
$$
|v_i^*\otimes w_j|_\cL=\sup_{\underline b\in A^n\setminus\{0\}}
\frac{|b_i|\cdot|w_j|_W}{\displaystyle\max_{1\leq k\leq n}|b_k|\cdot|v_k|_V}.
$$
For given $\underline b:=(b_1,\dots,b_n)$, in order for the
expression on the right-hand side to be non-zero, it is
necessary that $b_i\neq 0$; in that case the denominator
of the right-hand side cannot be made lower than $|b_i|\cdot|v_i|_V$,
so the claim follows.
\end{pfclaim}

Taking into account claim \ref{cl_v_i.w_k}, the lemma boils down
to the following
\begin{claim} For every $n\times m$ matrix $(\alpha_{ij})$ with
coefficients in $A$ we have:
$$
\sup_{\underline b\in A^n\setminus\{0\}}
\frac{{\displaystyle\max_{1\leq j\leq m}}|\sum_i\alpha_{ij}b_i|\cdot|w_j|_W}
{\displaystyle\max_{1\leq k\leq n}|b_k|\cdot|v_k|_V}=
\max_{ij}|\alpha_{ij}|\cdot\frac{|w_j|_W}{|v_i|_V}.
$$
\end{claim}
\begin{pfclaim}[] The inequality $\geq$ can be shown by choosing,
for every $r\leq n$, the vector
$\underline b{}_r:=(b_{1r},\dots,b_{nr})$ such that $b_{ir}=0$ for
$i\neq r$ and $b_{rr}=1$. For the inequality $\leq$ one remarks
that
$$
\frac{{\displaystyle\max_{1\leq j\leq m}}|\sum_i\alpha_{ij}b_i|\cdot|w_j|_W}
{\displaystyle\max_{1\leq k\leq n}|b_k|\cdot|v_k|_V}\leq
\frac{{\displaystyle\max_{ij}}\,|\alpha_{ij}b_i|\cdot|w_j|_W}
{\displaystyle\max_{1\leq k\leq n}|b_k|\cdot|v_k|_V}\leq
\max_{ij}\frac{|\alpha_{ij}b_i|\cdot|w_j|_W}{|b_i|\cdot|v_i|_V}
$$
from which the claim follows easily.
\end{pfclaim}
\end{proof}

\sset\subsubsection{}\label{subsec_tensornorm}
Let $V$ and $W$ be as in lemma \ref{lem_descartes}. By
lemma \ref{lem_descartes}(i) there is a natural $A$-linear
isomorphism
$$
\cL(V\otimes_AW,A)\simeq\cL(V,\cL(W,A))
$$
whence a natural structure of $A$-cartesian normed $A$-module
on $\cL(V\otimes_kW,A)$. After dualizing (and applying
again lemma \ref{lem_descartes}) we deduce that $V\otimes_AW$
carries a natural structure of $A$-cartesian normed $A$-module.
Furthermore, let $(v_i~|~i=1,\dots,n)$ and $(w_j~|~j=1,\dots,m)$
be orthogonal bases for $V$ and respectively $W$; using repeatedly
lemma \ref{lem_descartes}(ii) one sees easily that
$(v_i\otimes w_j~|~i=1,\dots,n; j=1,\dots,m)$ is an orthogonal
basis for $V\otimes_AW$ and moreover
$$
|v_i\otimes w_j|=|v_i|_V\cdot|w_j|_W\qquad
\text{for every $i=1,\dots,n$ and $j=1,\dots,m$.}
$$

\begin{remark} At least when $\Gamma=\R$, and $A$ is a field,
it should be possible to use the characterization of
\cite[\S2.1.7, Cor.3]{Bo-Gun} to see that the above normed $A$-module
structure on $V\otimes_AW$ agrees with the one defined on the complete
tensor product $V\hat\otimes_AW$ as in \cite[\S2.1.7]{Bo-Gun}.
\end{remark}

\sset\subsubsection{}\label{subsec_ext-powers}
As a special case of \eqref{subsec_tensornorm}, we deduce a
natural norm on every tensor power $V^{\otimes k}$ of $V$.
All these $A$-modules are $A$-cartesian. For every $k\in\N$
we have a natural imbedding of $A$-modules:
$\Lambda^k_AV\hookrightarrow V^{\otimes k}$
induced by the antisymmetrizer operator
(\cite[Ch.III, \S7.4, Remarque]{Bourbaki})
$$
V^{\otimes k}\to V^{\otimes k}\quad :\quad
v_1\otimes\cdots\otimes v_k\mapsto
\sum_{\sigma\in S_k}v_{\sigma(1)}\otimes\cdots\otimes v_{\sigma(k)}.
$$
Hence the norm of $V^{\otimes k}$ restricts to a natural
norm on $\Lambda^k_AV$. Let $\{v_1,\dots,v_n\}$ be an orthogonal
basis for $V$. For every subset $I\subset\{1,\dots,n\}$ of
cardinality $|I|=k$ we set $v_I:=v_{i_1}\wedge\cdots\wedge v_{i_k}$,
where $i_1<\cdots<i_k$ are the elements of $I$. One checks easily
that
$$
|v_I|=|v_{i_1}|_V{\cdot}\dots{\cdot}|v_{i_k}|_V
$$
and the basis $(v_I~|~I\subset\{1,\dots,n\}, |I|=k)$ is orthogonal.

\sset\subsubsection{}\label{subsec_A^+-submodules}
In the situation of \eqref{subsec_ext-powers}, consider a
free $A^+$-submodule $V^+\subset V$ such that the natural map
$A\otimes_{A^+}V^+\to V$ is an isomorphism. The highest exterior
power $\Lambda^n_{A^+}V^+$ is a rank one free $A^+$-submodule
of the $A$-cartesian module $\Lambda^n_AV$. Pick any generator
$e$ of $\Lambda^n_{A^+}V^+$; one sees easily that the value
$$
|V^+|:=|e|
$$
is independent of the choice of $e$. Especially, if $A$ is an
integral domain and $I\subset A^+$ is any principal ideal, then
$|I|$ is well defined.

\begin{lemma}\label{lem_Fitting}
Suppose that $A$ is an integral domain and let
$V^+_1\subset V^+_2\subset V$ be two $A^+$-submodules of the
free cartesian $A$-module $V$ of finite rank, fulfilling the
conditions of \eqref{subsec_A^+-submodules}. Then we have:
$$
|V^+_1|=|F_0(V^+_2/V^+_1)|\cdot|V^+_2|.
$$
\end{lemma}
\begin{proof} Here $F_0$ denotes the Fitting ideal (see
\cite[Ch.XIX]{Lan} for generalities on Fitting ideals).
Let $n$ be the rank of $V$; more or less by definition we have 
$F_0(V^+_2/V^+_1)=F_0(\Lambda^n_{A^+}V^+_2/\Lambda^n_{A^+}V^+_1)$,
from which the assertion follows easily.
\end{proof}

\subsection{Henselian algebras and complete algebras}\label{sec_prelim}
Let $(K,|\cdot|)$ be a complete valued field of rank one, $\fm$
the maximal ideal of the valuation ring $K^+$ of $K$,
$K^\sim:=K^+/\fm$ the residue field and $\Gamma_K$
the value group. Let also $\pi\in\fm$ be a fixed non-zero element.

\sset\subsubsection{}\label{subsec_normal-compl}
For any $K^+$-algebra $R$, let us denote by $R\Alg_{\fpet/K}$
(resp. $R\Alg_{\fget/K}$) the category of $R$-algebras $B$ that
are finitely presented (resp. finitely generated) as $R$-modules
and such that $B_K:=B\otimes_{K^+}K$ is \'etale over $R_K:=R\otimes_{K^+}K$.
Furthermore, for any object $B$ of $R\Alg_{\fget/K}$, let $B^\nu$
be the integral closure in $B_K$ of the image of $B$.

\begin{proposition}\label{prop_normal}
With the notation of \eqref{subsec_normal-compl}, suppose that
$R$ is $K^+$-flat and henselian along its ideal $\pi R$, and
denote by $R^\wedge$ the $\pi$-adic completion of $R$. Then :
\begin{enumerate}
\item
The base change functor $B\mapsto B^\wedge:=R^\wedge\otimes_RB$
induces equivalences :
$$
R\Alg_{\fpet/K}\isom R^\wedge\Alg_{\fpet/K}
\qquad
R\Alg_{\fget/K}\isom R^\wedge\Alg_{\fget/K}.
$$
\item
$\Ann_{B^\wedge}(\pi)=\Ann_B(\pi)$\ \  for every object
$B$ of $R\Alg_{\fget/K}$.
\item
Suppose furthermore that $R_K$ and $R^\wedge_K$ are normal domains.
Then the natural map :
$$
B^\nu\otimes_RR^\wedge\to(B\otimes_RR^\wedge)^\nu
$$
is an isomorphism for every object $B$ of $R\Alg_{\fget/K}$.
\end{enumerate}
\end{proposition}
\begin{proof} (i): The assertion for $R\Alg_{\fpet/K}$ follows from
\cite[Ch.III, \S3, Th.5 and \S4, Rem.2, p.587]{Elk}.
Next, let $B^\wedge$ be any object of $R^\wedge\Alg_{\fget/K}$;
we may find a filtered system $(B^\wedge_\lambda~|~\lambda\in\Lambda)$
of finitely presented $R^\wedge$-algebras, with surjective transition
maps $\phi^\wedge_{\lambda\mu}:B^\wedge_\lambda\to B^\wedge_\mu$,
whose colimit is $B$.
Since $B^\wedge$ is integral over $R^\wedge$, we may also
arrange that $B^\wedge_\lambda$ is integral over $R^\wedge$ for every
$\lambda\in\Lambda$, and then every $B^\wedge_\lambda$ is of finite
presentation as $R^\wedge$-module. Furthermore, since $B^\wedge_K$
is a finitely presented $R^\wedge_K$-algebra, we may assume
that
\set\begin{equation}\label{eq_no-variation}
B^\wedge_\lambda\otimes_{K^+}K\simeq B^\wedge
\qquad\text{for every $\lambda\in\Lambda$}.
\end{equation}
Therefore every $B^\wedge_\lambda$ is an object of
$R^\wedge\Alg_{\fpet/K}$. In this case, under the foregoing
equivalence, this family comes from a filtered system
$(B_\lambda~|~\lambda\in\Lambda)$ of objects of $R\Alg_{\fpet/K}$.
Let $B$ be the colimit of the latter family; then
$B\otimes_RR^\wedge\simeq B^\wedge$. Moreover :

\begin{claim}\label{cl_ess-surj}
(i)\ \ $B$ is a finitely generated $R$-module.
\begin{enumerate}
\addenu
\item
The induced map : $B_\lambda\otimes_{K^+}K\to B_K$ is bijective
for every $\lambda\in\Lambda$.
\end{enumerate}
\end{claim}
\begin{pfclaim} For (i), it suffices to show that the transition
maps $\phi_{\lambda\mu}:B_\lambda\to B_\mu$ are still surjective.
Indeed, let $C_{\lambda\mu}:=\Coker\,\phi_{\lambda\mu}$; then
$C_{\lambda\mu}\otimes_RR^\wedge=0$, hence
$C_{\lambda\mu}/\pi C_{\lambda\mu}=0$, and therefore $C_{\lambda\mu}=0$,
by Nakayama's lemma.

(ii): Let $C_\lambda:=B/B_\lambda$; we have to show that
$\pi^nC_\lambda=0$ for large enough $n\in\N$. However,
$C_\lambda/\pi^n C_\lambda\simeq C_\lambda\otimes_R(R/\pi^nR)\simeq
C_\lambda\otimes_R(R^\wedge/\pi^nR^\wedge)\simeq
C_\lambda\otimes_RR^\wedge$ for every sufficiently large $n\in\N$,
by \eqref{eq_no-variation}.
It follows that $C_\lambda/\pi^nC_\lambda=C_\lambda/\pi^mC_\lambda$
for every $m>n$, {\em i.e.} $\pi^nC_\lambda=\pi^mC_\lambda$, whence
$\pi^mC_\lambda=0$ by (i) and Nakayama's lemma.
\end{pfclaim}

Claim \ref{cl_ess-surj} implies that $B$ is an object of $R\Alg_{\fget/K}$,
hence the base change functor is essentially surjective on the subcategory
$R^\wedge\Alg_{\fget/K}$. Next, let $C$ be any other object of
$R\Alg_{\fget/K}$, and $\alpha^\wedge:B^\wedge\to C^\wedge$ a
ring homomorphism. Choose a filtered system $(C_\mu~|~\mu\in\Lambda')$
of objects of $R\Alg_{\fpet/K}$ whose colimit is $C$; for every
$\lambda\in\Lambda$, let $\psi_\lambda:B^\wedge_\lambda\to B^\wedge$
be the natural map; we may find $\mu\in\Lambda'$ such that the composition
$\alpha^\wedge\circ\phi_\lambda$ factors though a morphism
$\alpha^\wedge_{\mu\lambda}:B^\wedge_\lambda\to C^\wedge_\mu$.
Let $\alpha_{\mu\lambda}:B_\lambda\to C_\mu$ be the corresponding
morphism in $R\Alg_{\fpet/K}$ and $\alpha_\lambda:B_\lambda\to C$
its composition with the natural map $C_\mu\to C$. One checks easily
that $\alpha_\lambda$ does not depend on the choice of
$\alpha_{\mu\lambda}$, and moreover, the collection
$(\alpha_\lambda~|~\lambda\in\Lambda)$ is compatible with
the transition morphisms
$\phi_{\lambda\lambda'}:B_\lambda\to B_{\lambda'}$ of the filtered
system $(B_\lambda~|~\lambda\in\Lambda)$, therefore it gives rise
to a map $\alpha:B\to C$, and by construction it is clear that
$\alpha\otimes_R\one_{R^\wedge}=\alpha^\wedge$. This shows that
the base change is a full functor
$R\Alg_{\fget/K}\to R^\wedge\Alg_{\fget/K}$. A similar argument,
by reduction to finitely presented algebras, yields also the
faithfulness of the functor, thus concluding the proof of
assertion (i).

(ii): To start out, we claim that the natural map
$B\to B^\wedge$ is injective for every object $B$ of
$R\Alg_{\fget/K}$. Indeed, let $I\subset B$ be the kernel
of this map; then both $B^\wedge$ and $(B/I)^\wedge$
are the coproduct of $B$ and $R^\wedge$ in the category of
$R$-algebras, hence the natural map $B^\wedge\to(B/I)^\wedge$
is an isomorphism, so the same holds for the map $B\to B/I$,
in view of (i). Now, let $T:=\Ann_B(\pi)$; clearly
$T\otimes_RR^\wedge=T$, whence a right exact sequence :
$T\to B^\wedge\stackrel{\pi}{\to} B^\wedge\to 0$. However,
$T$ injects into $B^\wedge$, since it injects into $B$ and
$B$ injects into $B^\wedge$; hence the foregoing sequence
is short exact, which is assertion (ii).

(iii): Under the standing assumptions, $B^\nu$ is the
colimit of the filtered family $(B_\lambda~|~\lambda\in\Lambda)$
consisting of all the objects of $R\Alg_{\fget/K}$ such that
$B_\lambda\otimes_{K^+}K=B$ and such that $\pi$ is
regular in $B_\lambda$. In view of (ii), the family
$(B^\wedge_\lambda~|~\lambda\in\Lambda)$ consists of
all the objects of $R^\wedge\Alg_{\fget/K}$ such that
$B^\wedge_\lambda\otimes_{K^+}K=B$ and such that $\pi$ is
regular in $B^\wedge_\lambda$, hence its colimit is $(B^\wedge)^\nu$.
\end{proof}

\sset\subsubsection{}\label{subsec_future}
Let $A$ a $K^+$-algebra of finite presentation, and $A^h$
(resp. $A^\wedge$) the henselization of $A$ along its ideal
$\pi A$ (resp. the $\pi$-adic completion of $A$).
Recall that the $\pi$-adic completion of $A^h$ is naturally
isomorphic to $A^\wedge$.
(indeed, $A^h/\pi^n A^h$ is the henselization of $A/\pi^n A$
along its ideal $\pi A/\pi^n A$, for every $n\in\N$
\cite[Ch.XI, \S2, Prop.2]{Ray}; therefore
the natural map $A/\pi^n A\to A^h/\pi^n A^h$ is an isomorphism
for every $n\in\N$, which implies the claim).

\begin{lemma}\label{lem_future} In the situation of
\eqref{subsec_future}, suppose that $A$ is flat over
$K^+$. Then $A^\wedge$ is faithfully flat over $A^h$.
\end{lemma}
\begin{proof} To begin with, we claim that $A^\wedge$ is
flat over $K^+$. Indeed, suppose that $\pi x=0$ for some
$x\in A^\wedge$; choose a sequence $(x_k~|~k\in\N)$ of elements
of $A$ converging to $x$ in the $\pi$-adic topology of $A^\wedge$.
Then the sequence $(\pi x_k~|~k\in\N)$ converges to $0$, and
since $A$ has no $\pi$-torsion, it follows easily that the
sequence $(x_k~|~k\in\N)$ also converges to $0$, so $x=0$.

\begin{claim} In order to show the lemma, it suffices to
prove that $A^\wedge$ is flat over $A$.
\end{claim}
\begin{pfclaim} First, if $A^\wedge$ is flat over $A$, then
$A^\wedge$ is flat over $A^h$ as well. To conclude,
by standard reductions, it suffices to show that a finitely
generated $A^h$-module $M$ vanishes if and only if
$M^\wedge:=A^\wedge\otimes_{A^h}M$ vanishes. But if $M^\wedge=0$,
it follows that $M/\pi M=0$, and then we invoke Nakayama's lemma
to see that $M=0$.
\end{pfclaim}

Hence, let us that show $A^\wedge$ is flat over $A$; to this
aim, since $A/\pi A\simeq A^\wedge/\pi A^\wedge$,
\cite[Lemma 5.2.1]{Ga-Ra} says that it suffices to show that
$A^\wedge_K:=A^\wedge\otimes_{K^+}K$ is flat over
$A_K:=A\otimes_{K^+}K$. Say that $A=K^+[T_1,\dots,T_n]/I$
for some finitely generated ideal $I$.

\begin{claim}\label{cl_induced-top}
$\pi^nI=I\cap\pi^nK^+[T_1,\dots,T_n]$ for every $n\in\N$.
\end{claim}
\begin{pfclaim} By assumption $\Tor^{K^+}_1(A,K^+/\pi^nK^+)=0$;
hence
$$
I/\pi^nI=\Ker(K^+[T_1,\dots,T_n]/\pi^nK^+[T_1,\dots,T_n]\to
A/\pi^nA).
$$
The claim follows easily.
\end{pfclaim}

By \cite[Th.8.4]{Mat} we have
$A^\wedge\simeq K^+\langle T_1,\dots,T_n\rangle/I^\wedge$, where
$I^\wedge$ is the completion of $I$ for the subspace topology
as a submodule of $K^+[T_1,\dots,T_n]$; by claim \ref{cl_induced-top}
the subspace topology is nothing else than the $\pi$-adic topology.
It then follows from \cite[Prop.7.1.1]{Ga-Ra}(iv) that :
$$
I^\wedge=IK^+\langle T_1,\dots,T_n\rangle
$$
hence $A^\wedge\simeq
K^+\langle T_1,\dots,T_n\rangle\otimes_{K^+[T_1,\dots,T_n]}A$
and $A^\wedge_K\simeq 
K\langle T_1,\dots,T_n\rangle\otimes_{K[T_1,\dots,T_n]}A_K$,
hence we are reduced to the case where $A=K^+[T_1,\dots,T_n]$.
Let $\fn\subset A^\wedge_K$ be any maximal ideal, and set
$\fq:=\fn\cap A_K$; it suffices to show that $A^\wedge_{K,\fn}$
is flat over $A_{K,\fq}$. However, it is well known that
$E:=A^\wedge_K/\fn$ is a finite extension of $K$, hence
the same holds for $A_K/\fq\subset E$. Choose any maximal
ideal $\fn_E\subset A_E^\wedge:=E\otimes_KA^\wedge_K\simeq
E\langle T_1,\dots,T_n\rangle$ lying over $\fn$ and let
$\fq_E$ be the preimage of $\fn_E$ in $A_E:=E[T_1,\dots,T_n]$.
Since the extension $A^\wedge_{K,\fn}\to A^\wedge_{E,\fn_E}$
is faithfully flat, we are reduced to showing that
$A^\wedge_{E,\fn_E}$ is flat over $A_{E,\fq_E}$.
Hence we can replace $K$ by $E$ and assume from start that
$A/\fn\simeq K$, in which case
$\fn=(T_1-a_1,\dots,T_n-a_n)$ for some $a_1,\dots,a_n\in K^+$.
Clearly the $\fn$-adic completions of $A_K$ and $A^\wedge_K$
are both isomorphic to $K[[T_1-a_1,\dots,T_n-a_n]]$, and
by \cite[Th.8.8]{Mat} this latter ring is faithfully flat
over both $A_{K,\fq}$ and $A^\wedge_{K,\fn}$. The claim
follows easily.
\end{proof}

\section{Study of the discriminant}\label{chap_study-disc}

\subsection{Discriminant}\label{sec_discrim}
Let $R\to S$ be a ring homomorphism such that $S$ is a free $R$-module
of finite rank. Every element $a\in S$ defines an $R$-linear endomorphism
$$
\mu_a:S\to S\qquad b\mapsto ab
$$
whose trace and determinant we denote respectively by $\tr_{S/R}(a)$ and
$\Nm_{S/R}(a)$. There follows a well-defined $R$-bilinear {\em trace form}
$$
\Tr_{S/R}:S\otimes_RS\to R\qquad a\otimes b\mapsto\tr_{S/R}(ab).
$$
It is well known (see {\em e.g.} \cite[Th.4.1.14]{Ga-Ra}) that $\Tr_{S/R}$
is a perfect pairing if and only if $S$ is \'etale over $R$. Pick a basis
$e_1,\dots,e_d$ of $S$; one defines the {\em discriminant\/} of $S$ over
$R$ as the element
$$
\fd_{S/R}:=\det(\Tr_{S/R}(e_i\otimes e_j)~|~1\leq i,j\leq d)\in R.
$$
One verifies easily that $\fd_{S/R}$ is well defined ({\em i.e.} independent
of the choice of the basis) up to the square of an invertible element of
$R$.

\sset\subsubsection{}\label{subsec_ring-of-rings}
Let $R$ be a (not necessarily commutative) unitary ring; for any integer
$m>0$ we let $M_m(R)$ be the unitary ring of all $m\times m$ matrices with
entries in $R$. For every $a\in R$ and every pair of integers $i,j\leq m$
we denote by $E_{ij}(a)\in M_m(R)$ the elementary matrix whose $(i,j)$-entry
equals $a$, and whose other entries vanish; moreover, sometimes we may denote
by $1_m$ the unit of $M_m(R)$. If $n>0$ is any other integer,
we let
$$
\alpha_n:M_n(M_m(R))\stackrel{\sim}{\to}M_{nm}(R)
$$
be the unique ring isomorphism such that
$$
E_{ij}(E_{kl}(a))\mapsto E_{i(m-1)+k,j(m-1)+l}(a) \qquad
\text{for all $a\in R$,\quad $1\leq i,j\leq n$,\quad $1\leq k,l\leq m$.}
$$
Suppose now that $t:=(t_{ij})\in M_n(M_m(R))$ is a matrix whose entries
$t_{ij}$ commute pairwise; let $T\subset M_m(R)$ be the commutative ring
generated by all the $t_{ij}$; we can then view $t$ as an element of $M_n(T)$
so that its determinant is well-defined as an element of $T$.
To avoid ambiguities, we shall write $\det_n(t_{ij}~|~1\leq i,j\leq n)$
for this determinant.

\begin{lemma}\label{lem_blocks} With the notation of
\eqref{subsec_ring-of-rings}, suppose that the ring $R$ is commutative
and let $t:=(t_{ij}~|~1\leq i,j\leq n)\in M_n(M_m(R))$ be an element
such that all the matrices $t_{ij}\in M_m(R)$ commute pairwise.
We have the identity:
\set\begin{equation}\label{eq_double-det}
\det(\mathrm{det}_n(t_{ij}~|~1\leq i,j\leq n))=\det(\alpha_n(t)).
\end{equation}
\end{lemma}
\begin{proof} We proceed by induction on $n$, the case $n=1$ being
trivial. Hence, assume $n>1$; suppose first that $t_{11}\in M_m(R)$
is an invertible matrix. It follows that the matrix
$$
s:=E_{11}(t_{11})+\sum_{k=2}^nE_{kk}(1_m)
$$
is invertible in $M_n(M_m(R))$, and obviously its entries commute pairwise and
with the entries of $t$; furthermore the sought identity is easily verified for
$s$. Since both sides of \eqref{eq_double-det} are multiplicative in $t$,
we are therefore reduced to verifying the identity for $s^{-1}\cdot t$, hence
we can assume that $t_{11}=1_m$. Next, let
$$
e:=(1_m)_n-\sum_{k=2}^nE_{1k}(t_{1k}).
$$
Clearly $e$ is invertible in $M_n(M_m(R))$, and again its entries commute both
pairwise and with the entries of $t$; thus it suffices to show the sought
identity for $e$ and for $e^{-1}\cdot t$. The identity is obvious for $e$,
therefore we can replace $t$ by $e^{-1}\cdot t$ and assume that
$t_{1j}=\delta_{1j}\cdot 1_m$ for every $j\leq n$. In this case,
$\det_n(t_{ij}~|~1\leq i,j\leq n)$ equals the determinant of the
$(n-1)\times(n-1)$-minor $t'$ obtained by omitting the first row
and the first column of $t$, and $\det(\alpha_n(t))=\det(\alpha_{n-1}(t'))$.
By inductive assumption, the sought identity is already know for such a minor,
so the proof is complete in case $t_{11}$ is invertible. For a general $t$,
we notice that $\det(t_{11}+\lambda 1_m)$ is a non-zero-divisor in the free
polynomial $R$-algebra $R[\lambda]$ and consider the localization
$S:=R[\lambda,\det(t_{11}+\lambda 1_m)^{-1}]$. The assumptions of the lemma
are verified by the matrix $t'':=t+E_{11}(\lambda 1_m)\in M_n(M_m(S))$ and
moreover $t''_{11}$ is invertible in $M_n(S)$, so \eqref{eq_double-det} holds
for $t''$, and actually yield an identity in the subring $R[\lambda]$ of $S$.
After specializing the latter identity in $\lambda=0$, we deduce that
\eqref{eq_double-det} holds for $t$ as well.
\end{proof}

\begin{proposition}\label{prop_tower-discr}
Let $A\to B\to C$ be maps of commutative rings and suppose that $B$
(resp. $C$) is a free $A$-module (resp. $B$-module) of finite rank.
Let $r:=\rk_BC$. Then we have:
$$
\fd_{C/A}=(\fd_{B/A})^r\cdot\Nm_{B/A}(\fd_{C/B}).
$$
\end{proposition}
\begin{proof} Let $d:=\rk_AB$; pick bases $e_1,\dots,e_d\in B$ of the
$A$-module $B$ and $f_1,\dots,f_r\in C$ of the $B$-module $C$; clearly
the system $(e_if_j~|~i\leq d, j\leq r)$ is a basis for the free
$A$-module $C$ of rank $dr$. We let $T\in M_r(M_d(A))$ be the element
whose $(j,j')$-entry is the matrix $T_{jj'}$ such that
$$
(T_{jj'})_{ii'}:=\tr_{C/A}(e_if_je_{i'}f_{j'})
\qquad\text{for every $1\leq i,i'\leq d$.}
$$
In the notation of \eqref{subsec_ring-of-rings}, we have
$\fd_{C/A}=\det(\alpha_r(T))$. Moreover, let $M\in M_r(B)$
(resp. $N\in M_d(A)$) be the matrix such that 
$M_{jj'}:=\tr_{C/B}(f_jf_{j'})$ (resp. such that
$N_{ii'}:=\tr_{B/A}(e_ie_{i'})$); by the transitivity of
the trace, we can write
\set\begin{equation}\label{eq_taking.into}
(T_{jj'})_{ii'}=\tr_{B/A}(e_ie_{i'}\cdot M_{jj'}).
\end{equation}
Let $\mu:B\to M_d(A)$ be the unique ring homomorphism such that
$$
be_i=\sum_{k=1}^d\mu(b)_{ki}e_k\qquad\text{for all $b\in B$}.
$$
Especially: $e_{i'}M_{jj'}=\sum_{k=1}^d\mu(M_{jj'})_{ki'}e_k$ and consequenly:
\set\begin{equation}\label{eq_account}
\tr_{B/A}(e_ie_{i'}M_{jj'})=\sum_{k=1}^d\mu(M_{jj'})_{ki'}\cdot\tr_{B/A}(e_ke_i)
                           =\sum_{k=1}^dN_{ik}\cdot\mu(M_{jj'})_{ki'}.
\end{equation}
Finally, let $\Delta(N),\mu(M)\in M_r(M_d(A))$ be the matrices such that
$$
\Delta(N)_{jj'}:=N\cdot\delta_{jj'}\qquad
\mu(M)_{jj'}:=\mu(M_{jj'})\qquad\text{for all $1\leq j,j'\leq r$.}
$$
Taking into account \eqref{eq_taking.into} and \eqref{eq_account} we see that
$$
T=\Delta(N)\cdot\mu(M)
$$
whence, an application of lemma \ref{lem_blocks} delivers the sought identity.
\end{proof}

\sset\subsubsection{}\label{subsec_different}
In the situation of \eqref{sec_discrim} we let
$$
\tau_{S/R}:S\to S^*:=\Hom_R(S,R)
$$
be the map such that $\tau_{S/R}(b)(b')=\Tr_{S/R}(b\otimes b')$ for
every $b,b'\in S$. Notice that $S^*$ is an $S$-module with the natural
scalar multiplication defined by the rule: $(b\cdot\phi)(b'):=\phi(bb')$
for every $b,b'\in S$ and $\phi\in S^*$. With respect to this $S$-module
structure, $\tau$ is $S$-linear; thus, we can define the {\em different ideal}
$$
\cD_{S/R}:=\Ann_S(\Coker\,\tau_{S/R})\subset S.
$$
In this generality, not much can be said about the ideal $\cD_{S/R}$.
However, suppose furthermore that there is an isomorphism of $S$-modules
$\omega:S^*\stackrel{\sim}{\to}S$; it follows easily that $\cD_{S/R}$
is the principal ideal generated by $\delta:=\omega\circ\tau(1)$. Denote
by $\Nm_{S/R}(\cD_{S/R})\subset R$ the principal ideal generated
by $\Nm_{S/R}(\delta)$.

\begin{lemma}\label{lem_different}
Under the assumptions of \eqref{subsec_different} we have the identity:
$$
\Nm_{S/R}(\cD_{S/R})=\fd_{S/R}.
$$
\end{lemma}
\begin{proof} Indeed, let $r:=\rk_RS$; directly from the
definition we deduce that
$$
\fd_{S/R}=\Ann_R(\Coker\,\Lambda^r_R\tau_{S/R})=
\Ann_R(\Coker\,\Lambda^r_R(\omega\circ\tau_{S/R})).
$$
which implies straightforwardly the assertion.
\end{proof}

\begin{example}\label{ex_different}
Suppose that $R$ is a henselian valuation ring and $S$ is the
integral closure of $R$ in a finite extension of the field of
fractions of $R$; moreover suppose that $S$ is a finitely
presented $R$-module. Then actually $S$ is a free $R$-module
of finite rank. Notice that $S$ is a valuation ring and an
$S$-module is $S$-torsion-free if and only if it is
$R$-torsion-free; in particular we see that $S^*$ is a finitely
presented $S$-torsion-free $S$-module, hence it is free over
$S$ (see {\em e.g.} \cite[lemma 6.1.14]{Ga-Ra}) and clearly
$\rk_SS^*=1$, so lemma \ref{lem_different} applies to the extension
$R\subset S$. Moreover :
\end{example}

\begin{lemma}\label{lem_bound-discr}
Keep the assumptions of example {\em\ref{ex_different}}
and suppose moreover that the valuation $|\cdot|_R$ of $R$ has
rank one; let $d:=\rk_RS$. Then:
\begin{enumerate}
\item
In case the valuation $|\cdot|_R$ is discrete, 
$|\fd_{S/R}|\geq|\pi|_R^{d(1-1/e)}\cdot|d|^d$, where $\pi\in R$
is a uniformizer of $R$ and $e$ is the ramification index of
$S$ over $R$.
\item
In case the valuation $|\cdot|_R$ is not discrete,
$|\fd_{S/R}|\geq|d|^d$.
\end{enumerate}
\end{lemma}
\begin{proof} Obviously $\tr_{S/R}(d^{-1})=1$ in either case. If
the valuation of $R$ is discrete, we can write
$\tr_{S/R}(\pi^{-1}\cdot d^{-1})=\pi^{-1}$, hence
$\pi^{-1}\cdot d^{-1}\notin\cD_{S/R}^{-1}$ and consequently
$\cD_{S/R}^{-1}\subset\pi_S\cdot\pi^{-1}\cdot d^{-1}S$
(inclusion of fractional ideals, where $\pi_S$ is a uniformizer for $S$).
The bound then follows easily from lemma \ref{lem_different}.
In case the valuation is non-discrete, the same argument yields
the weaker estimate: $\cD_{S/R}^{-1}\subset x\cdot d^{-1}S$,
for every $x$ with $|x|_R<1$, whence the sought inequality,
again by lemma \ref{lem_different}.
\end{proof}

\begin{lemma}\label{lem_Fitt-discr}
Let $R$ be a valuation ring and $S_2\subset S_1$ two
$R$-algebras that are both free of the same finite rank
as $R$-modules. Then we have:
$$
\fd_{S_2/R}=F_0(S_1/S_2)^2\cdot \fd_{S_1/R}.
$$
\end{lemma}
\begin{proof} (Here $F_0$ denotes the Fitting ideal of the
torsion $R$-module $S_1/S_2$.) This is a special case of
\cite[Lemma 7.5.4]{Ga-Ra}.
\end{proof}

The following example will play a key role in later sections.

\begin{example}\label{ex_discriminant}
Let $K$, $K^+$ and $\pi$ be as in \eqref{sec_prelim}.
As usual, one lets $K^+\langle T_1,\dots,T_n\rangle$ be the
$\pi$-adic completion of $K^+[T_1,\dots,T_n]$.
We consider the (continuous) ring homomorphism
$$
\psi:K^+\langle\xi\rangle\to K^+\langle S,T\rangle/(ST-\pi^2)
\quad :\quad \xi\mapsto S+T.
$$
Notice that $K^+\langle S,T\rangle/(ST-\pi^2)$ is generated
over $K^+\langle\xi\rangle$ by the class of $S$, which
satisfies the integral equation
$$
S^2-S\xi+\pi^2=0.
$$
The matrix of the trace form for the morphism $\psi$, relative to
the basis $(1,S)$ is:
$$
\left(
\begin{array}{cc}
2 & \xi \\
\xi & \xi^2-2\pi^2
\end{array}
\right).
$$
Finally, the discriminant of $\psi$ is $\fd_\psi:=\xi^2-4\pi^2$.
\end{example}

\subsection{Finite ramified coverings of annuli}\label{sec_annuli}
We keep the notation and assumptions of \eqref{sec_prelim},
and we suppose additionally that $K$ is algebraically closed.

\sset\subsubsection{}
We shall use rather freely the language of adic
spaces of \cite{Hu1} and \cite{Hu2}.
For a quick review of the main definitions, we refer also to
\cite[\S7.2.15-27]{Ga-Ra}. Recall that such an adic space is
a datum of the form $(X,\cO_X,\cO^+_X)$, where $(X,\cO_X)$ is
a locally ringed space and $\cO^+_X\subset\cO_X$ is a subsheaf
of rings satisfying certain natural conditions (see {\em loc.cit.});
moreover such an adic space is obtained by gluing {\em affinoid}
open subspaces that are {\em adic spectra\/} $\Spa\,A$
attached to certain pairs $A:=(A^\trg,A^+)$ consisting of
a ring and an integrally closed subring $A^+\subset A^\trg$.

For every $x\in X$ we shall denote :
\begin{itemize}
\item
$\kappa(x)$ the residue field of $\cO_{X,x}$, which is
a valued field whose valuation we denote $|\cdot|_x$.
\item
$(\kappa(x)^\wedge,|\cdot|^\wedge_x)$ (resp. $(\kappa(x)^h,|\cdot|^h_x)$)
the completion for the valuation topology (resp. the henselization)
of $(\kappa(x),|\cdot|_x)$.
\item
$(\kappa(x)^{\wedge h},|\cdot|_x^{\wedge h})$ the henselization
of $(\kappa(x)^\wedge,|\cdot|^\wedge_x)$.
\end{itemize}
In agreement with \eqref{eq_notation+}, we shall write $\kappa(x)^+$
for the valuation ring of the valuation $|\cdot|_x$, and likewise we
define $\kappa(x)^{\wedge+}$ and so on. For future reference we point out :

\begin{lemma}\label{lem_completing-kappa}
Let $X$ be an analytic adic space, $x\in X$ any point and $y\in X$
a specialization of $x$. Then the induced map
$\kappa(y)^\wedge\to\kappa(x)^\wedge$ is an isomorphism of complete
topological fields.
\end{lemma}
\begin{proof} It follows directly from \cite[Lemma 1.1.10(iii)]{Hu2}.
\end{proof}

\begin{remark}\label{rem_comp-hensel} Let $f:X\to Y$ be
a finite map of analytic adic spaces, $x\in X$ any point and
$y:=f(x)$. Notice that the extension of valued fields
$(\kappa(y),|\cdot|_y)\subset(\kappa(x),|\cdot|_x)$
is usually {\em not\/} algebraic, whereas the induced map on
completions :
$$
(\kappa(y)^\wedge,|\cdot|^\wedge_y)\to (\kappa(x)^\wedge,|\cdot|^\wedge_x)
$$
is always a finite algebraic extension (\cite[Lemma 1.5.2]{Hu2}),
but the latter does not necessarily induce a finite map between
the corresponding valuation rings (see \eqref{subsec_not-hensel}).
That is why it is useful to take henselizations : the induced ring
homomorphism $\kappa(y)^{\wedge h+}\to\kappa(x)^{\wedge h+}$ is finite.
\end{remark}

\sset\subsubsection{}\label{subsec_conflict}
Following R.Huber (\cite{Hu1}, \cite{Hu2}), we call a
{\em $K$-affinoid algebra\/} a pair $A:=(A^\trg,A^+)$,
where $A^\trg$ is  a $K$-algebra of topologically finite type
and $A^+$ is a subring of the ring $A^\circ$ of all power-bounded
elements of $A^\trg$. We shall consider exclusively affinoid
rings $A$ of {\em topologically finite type over $K$}; for such
$A$ one has always $A^+=A^\circ$. The subring $A^\circ$ is
characterized by a topological -- rather than metric -- condition.
Hence in principle the notation of this section may conflict with
\eqref{eq_notation+}. However, when $A^\trg$ is reduced, one
knows that the {\em supremum seminorm\/} $|\cdot|_{\sup}$ on
$A^\trg$ is a power-multiplicative norm (\cite[\S6.2.4, Th.1]{Bo-Gun}),
and furthermore in this case $A^\circ=\{a\in A^\trg~|~|a|_{\sup}\leq 1\}$
(\cite[\S6.2.3, Prop.1]{Bo-Gun}). For any such $A$ we
shall also write $A^\sim:=A^\circ/\fm A^\circ$.

\sset\subsubsection{}\label{subsec_confusion}
Another possible source of confusion is the following
situation. Let $A$ be a normal domain of topologically
finite type over $K$, and endow $A$ with its supremum norm
$|\cdot|_{\sup}$; let also $A\to B$ be an injective finite
ring homomorphism. According to proposition \ref{prop_spectral-norm},
$B$ is endowed with its spectral seminorm $|\cdot|_\sp$;
on the other hand, $B$ can also be endowed with its supremum
seminorm and the problem arises whether these two seminorms
coincide. According to \cite[\S3.8.1, Prop.7]{Bo-Gun}, this
turns out to be the case, provided that $B$ is torsion-free
as an $A$-module.

\begin{lemma}\label{lem_subdomain}
Let $A$ be a reduced affinoid $K$-algebra of topologically finite
type. Let $U\subset X:=\Spa\,A$ be an affinoid subdomain.
Then $\cO_X(U)$ is reduced.
\end{lemma}
\begin{proof} It suffices to show that, for every maximal ideal
$\fq\subset B:=\cO_X(U)$, the localization $B_\fq$ is reduced.
However, by a theorem of Kiehl, $B_\fq$ is excellent
(see \cite[Th.1.1.3]{Conr} for a proof), hence it suffices to
show that the $\fq$-adic completion $B_\fq^\wedge$ of $B_\fq$
is reduced. Let $\fp:=\fq\cap A$; since $U$ is a subdomain in $X$,
the natural map $A\to B$ induces an isomorphism of complete local
rings $A_\fp^\wedge\simeq B_\fq^\wedge$, so we are reduced to
showing that $A_\fp^\wedge$ is reduced, which holds because
$A_\fq$ is reduced and excellent (again by Kiehl's theorem).
\end{proof}

\sset\subsubsection{}\label{subsec_setup}
For every $a,b\in\Gamma_K$ with $a\leq b$, we denote by $\D(a)$ the
disc of radius $a$, and by $\D(a,b)$ the annulus of radii $a$ and $b$.
Say that $a=|\alpha|$ and $b=|\beta|$ for $\alpha,\beta\in K^\times$;
then 
$$
\D(a,b):=\Spa\,A(a,b)\quad\text{and}\quad\D(a):=\Spa\,A(a)
$$
where $A(a,b)$ (resp. $A(a)$) is the affinoid $K$-algebra of topologically
finite type such that
$$
A(a,b)^\trg:=K\langle\alpha/\xi,\xi/\beta\rangle\qquad
\text{(resp. $A(a)^\trg:=K\langle\xi/\alpha\rangle$)}.
$$
Hence $A(a,b)^+=A(a,b)^\circ=K^+\langle\alpha/\xi,\xi/\beta\rangle$
and $A(a)^+=A(a)^\circ=K^+\langle\xi/\alpha\rangle$.

\sset\subsubsection{}\label{subsec_classes}
Let $X$ be any adic space locally of finite type over $\Spa\,K$ and
with $\dim X=1$. The points of $X$ fall into three distinct
classes, according to whether: (I) they admit  neither a proper
generalization nor a proper specialization, or (II) they admit a
proper specialization, or else (III) they have a proper generalization.
For every point $x\in X$ of class (III), we shall denote by
$x^\flat$ the unique generization of $x$ in $X$, so $x^\flat$ is
a point of class (II). The value group $\Gamma_x$ of $|\cdot|_x$
admits a natural decomposition (see \cite[\S1.1 and Cor.5.4]{Hu3})
$$
\Gamma_x\simeq\Gamma_x^\mathrm{div}\oplus\langle\gamma_0\rangle
$$
where $\Gamma_x^\mathrm{div}$ is the maximal divisible subgroup
contained in $\Gamma_x$, and $\langle\gamma_0\rangle\simeq\Z$ is
the subgroup generated by the element $\gamma_0$ uniquely characterized
as the largest element of the subset $\Gamma_x^+\setminus\{1\}$
(notation of \eqref{subsec_gamma-plus}).

\sset\subsubsection{}
For instance, take $X:=(\A^1_K)^\ad$, the analytification of the affine
line. The topological space underlying $(\A^1_K)^\ad$ consists of the
equivalence classes of continuous valuations $v:K[\xi]\to\Gamma_v$
extending the valuation of $K$. These valuations are described in
\cite[\S5]{Hu3}: to the class (I) belong {\em e.g.} the height one
valuations of the form
$$
f(\xi)\mapsto |f(x)|\qquad\text{for all $f(\xi)\in K[\xi]$}
$$
where $x\in K=\A^1_K(K)$ is any element (these are the $K$-rational
points of $(\A^1_K)^\ad$). The valuations of classes (II) and (III)
are all of height respectively one and two. The elements of these
classes admit a uniform description, that we wish to explain.
To this aim we consider an imbedding of ordered fields
$$
(\R,<)\hookrightarrow(\R(\eps),<)
$$
where $\R(\eps)$ is a purely transcendental extension of $\R$,
generated by an element $\eps$ such that $0<\eps<r$ for every
real number $r>0$. One can view $\R(\eps)$ as a subfield of the
ordered field of hyperreal numbers ${}^*\R$ (see \cite{Gold}).
For every $x\in K$, every real number $r>0$ and
every $\omega\in\{1,1-\eps,1/(1-\eps)\}\subset\R(\eps)$
consider the valuation
$$|\cdot|_{r\cdot\omega}:K[\xi]\to\R(\eps)
\quad : \quad
\sum_{i=0}^n a_i(\xi-x)^n\mapsto
\max(|a_i|\cdot r^i\cdot\omega^i~|~i=0,\dots,n).
$$
If $\omega=1$, this is the usual Gauss (sup) norm attached to the
disc of radius $r$ centered at the point $x$; this is a valuation
of height one. For $\omega\neq 1$ we get a valuation which should
be thought of as the sup norm on a disc of radius $r\cdot\omega$,
again centered at $x$; this new kind of valuation is a specialization
of $|\cdot|_r$, and indeed all the specializations of the latter
occur in this manner. If $r\notin\Gamma_K$, then $|\cdot|_r$
belongs to the class (I); in this case the valuations
$|\cdot|_{r\cdot\omega}$ are all equivalent, regardless of $\omega$,
and therefore they induce the same point of $(\A^1_K)^\ad$. If 
$r\in\Gamma_K$ then $|\cdot|_r$ is in the class (II); in
this case the $|\cdot|_{r\cdot\omega}$ for $\omega\neq 1$ are two
inequivalent valuations of height two, hence of class (III), and
indeed all valuations of class (III) arise in this way.

\sset\subsubsection{}\label{subsec_localize}
Let $a,b\in\Gamma_K$ with $a\leq b$. For $r\in(a,b]\cap\Gamma_K$,
the valuation $|\cdot|_{r(1-\eps)}$ extends to $A(a,b)$ by continuity;
if moreover $r>a$, then the point of $(\A^1_K)^\ad$ corresponding
to this valuation lies in the open subdomain $\D(a,b)$. This point
shall be denoted henceforth by $\eta(r)$, and to lighten notation
we shall write $\kappa(r)$ (resp. $\kappa(r^\flat)$) for the residue
field of $\eta(r)$ (resp. of $\eta(r)^\flat$). Notice that $\kappa(r)$
is also the same as the stalk $\cO_{\D(a,b),\eta(r)}$.

Likewise, if $r\in[a,b)\cap\Gamma_K$, the valuation
$|\cdot|_{r/(1-\eps)}$ determines a point $\eta'(r)$ of
$(\A^1_K)^\ad$ that lies in the open subdomain $\D(a,b)$;
the residue field of $\eta'(r)$ shall be denoted $\kappa'(r)$.
Notice that $\eta(r)^\flat=\eta'(r)^\flat$.

\sset\subsubsection{}\label{subsec_direct-image}
Let $f:X\to\D(a,b)$ be a finite and flat morphism of affinoid adic
spaces of degree $d$, and suppose that $X$ is reduced ({\em i.e.}
$X=\Spa\,B$ where $B$ is a reduced flat affinoid algebra of rank $d$
as an $A(a,b)$-module). For every $r\in(a,b]\cap\Gamma_K$,
we set :
$$
\cB(r):=(f_*\cO_X)_{\eta(r)}
$$
which is a reduced finite $\kappa(r)$-algebra, in view of lemma
\ref{lem_subdomain}. We endow $\cB(r)$ with the spectral norm
$|\cdot|_{\sp,\eta(r)}$ relative to the valuation $|\cdot|_{\eta(r)}$;
it follows that
$$
\cB(r)^+=(f_*\cO^+_X)_{\eta(r)}.
$$

\begin{lemma}\label{lem_much-ado} In the situation of
\eqref{subsec_direct-image}, let $y\in\D(a,b)$ any point.
Then $f^{-1}(y)$ is the set of all the valuations on
$(f_*\cO_X)_y$ that extend the valuation $|\cdot|_y$
corresponding to $y$.
\end{lemma}
\begin{proof} Let $X=\Spa(B^\trg,B^\circ)$; the topology of
$B^\trg$ is the $A(a,b)$-module topology, {\em i.e.} the
unique one such that the family of $A(a,b)^\circ$-submodules
$(\pi B^\circ~|~\pi\in\fm\setminus\{0\})$ is a fundamental
system of neighborhoods of $0$. Let
$B_y:=B^\trg\otimes_{A(a,b)^\trg}\cO_y=(f_*\cO_X)_y$; similarly
$B_y$ has a well-defined $\cO_y$-module topology and the fibre
$f^{-1}(y)$ consists of the continuous valuations $|\cdot|'$
on $B_y$ extending the
valuation $|\cdot|_y$, and such that
\set\begin{equation}\label{eq_the-same-goes}
|s|'\leq 1\qquad\text{for all $s\in B^\circ$.}
\end{equation}
Let $|\cdot|'$ be any valuation on $B_y$ extending $|\cdot|_y$,
and let $\fp\subset B_y$ be the support of $|\cdot|'$;
the quotient topology on $E:=B_y/\fp$ is the $\kappa(y)$-module
topology, where $\kappa(y)$ is the residue field of $y$. However,
let $\Gamma_y$ and $\Gamma_E$ be the value groups of $|\cdot|_y$
and respectively $|\cdot|_E$; since $[\Gamma_E:\Gamma_y]$ is
finite, it is easy to see that $|\cdot|'$ is continuous. Hence,
continuity holds for all $|\cdot|'$ extending $|\cdot|_y$, and
since anyway $B^\circ$ is the integral closure of $A(a,b)^\circ$
in $B^\trg$, the same goes for condition \eqref{eq_the-same-goes}.
\end{proof}

\sset\subsubsection{}\label{subsec_dir-img-etale}
In the situation of \eqref{subsec_direct-image}, let
$x\in X$ be a point of class (III). The valuation $|\cdot|_x$
is an extension of the valuation $|\cdot|_{f(x)}:\kappa(f(x))\to\R(\eps)$,
hence its value group can be realized inside the multiplicative
group of the field $\R((1-\eps)^{1/d!})$, which is an algebraic
extension of $\R(\eps)$ of degree $d!$, admitting a unique
ordering extending the ordering of $\R(\eps)$
(again, one can think of all this as taking place inside the
hyperreal numbers; of course, there is no real need to introduce
this auxiliary field: it is nothing more than a suggestive notational device).
In terms of the decomposition of \eqref{subsec_classes}, we
have $\Gamma_x^\mathrm{div}\subset\R^\times_{>0}$ and
$\langle\gamma_0\rangle\subset\{(1-\eps)^i~|~i\in\frac{1}{d!}\Z\}$.
We shall consider the two projections:
$$
\Gamma_x\to\Gamma_x^\mathrm{div}\quad:\quad\gamma\mapsto\gamma^\flat
\qquad\text{and}\qquad
\Gamma_x\to\frac{1}{d!}\Z\quad:\quad\gamma\mapsto\gamma^\natural
$$
where $\gamma^\natural$ is characterized by the identity
$$
(1-\eps)^{\gamma^\natural}\cdot\gamma^\flat=\gamma
\qquad\text{for every $\gamma\in\Gamma_x$}.
$$
Sometimes it is more natural to use an additive (rather than
multiplicative) notation; in order to switch from one to the other,
of course one takes logarithms. Hence we define :
\set\begin{equation}\label{eq_switch}
\log\gamma:=\log\gamma^\flat-\gamma^\natural\cdot\eps\in\R+\eps\R
\qquad\text{for every $\gamma\in\Gamma_x$}.
\end{equation}
The composition
$$
B\to\Gamma_x^\mathrm{div}\cup\{0\}\quad:\quad s\mapsto|s|_x^\flat
$$
is a continuous rank one valuation of $B$ and determines the unique
generization $x^\flat$ of $x$ in $X$.
If we view $\R((1-\eps)^{1/d!})$ as a subfield of the
hyperreal numbers, then the projection $|s|_x^\flat$ corresponds
to the shadow of the bounded hyperreal $|s|_x$.

\sset\subsubsection{}\label{subsec_not-hensel}
The ring $\cB(r)$ is a product of finite field extensions
$F_1\times\cdots\times F_k$ of $\cO_{\eta(r)}$, and the factors
$F_j$ are in natural bijective correspondence with the
elements of the fibre $f^{-1}(\eta(r)^\flat)$ (see also
\cite[Prop.1.5.4]{Hu2}). Moreover,
$\cB(r)^+$ decomposes as the product $F_1^+\times\cdots\times F_k^+$,
where $F_j^+$ is the integral closure of $\kappa(r)^+$
in $F_j$. The valuation ring $\kappa(r)^+$ is not henselian,
hence it may occur that $F_j^+$ is not a valuation ring, but
only a normed $\kappa(r)^+$-algebra; that happens precisely
when there are distinct points $x,y\in f^{-1}(\eta(r))$ such that
$x^\flat=y^\flat$ (see example \ref{ex_summing-up}).

\begin{lemma}\label{lem_local-sit}
In the situation of \eqref{subsec_direct-image}, suppose
furthermore that the morphism $f$ is generically \'etale.
Then, for every $r\in(a,b]\cap\Gamma_K$ we have:
\begin{enumerate}
\item
The normed ring $(\cB(r)^+,|\cdot|_{\sp,\eta(r)})$
is a free cartesian $\kappa(r)^+$-module of rank $d$.
\item
$|s|_{\sp,\eta(r)}=\max(|s|_x~|~x\in f^{-1}(\eta(r)))$
\quad for all $s\in\cB(r)$.
\item
Let us view $\cB(r)^+$ as a submodule of the normed cartesian
module $(\cB(r),|\cdot|_{\sp,\eta(r)})$, so that the value
$|\cB(r)^+|_{\sp,\eta(r)}$ is defined (notation of
\eqref{subsec_A^+-submodules}). Then:
$$
2|\cB(r)^+|^\natural_{\sp,\eta(r)}=\deg(f)-\sharp f^{-1}(\eta(r))
$$
where $\sharp f^{-1}(\eta(r))$ denotes the cardinality of the fibre
$f^{-1}(\eta(r))$.
\item
If $R\subset\cB(r)^+$ is a finite $\kappa(r)^+$-algebra such that
$R\otimes_{K^+}K=\cB(r)$, then $R=\cB(r)^+$ if and only if $R/\fm R$
is a reduced $K^\sim$-algebra.
\end{enumerate}
\end{lemma}
\begin{proof} (i): under the current assumptions, the field
extension $\kappa(r)^\wedge\subset\kappa(x)^\wedge$
is finite (see remark \ref{rem_comp-hensel}(i)) and separable.
Let us set :
$$
\cB(r)^{\wedge+}:=\cB(r)^+\otimes_{\kappa(r)^+}\kappa(r)^{\wedge+}
\qquad
\cB(r)^\wedge:=\cB(r)\otimes_{\kappa(r)}\kappa(r)^\wedge.
$$
Since $\cB(r)$ is an \'etale $\kappa(r)$-algebra,
$\cB(r)^\wedge$ is an \'etale $\kappa(r)^\wedge$-algebra,
especially it is reduced; by flatness, $\cB(r)^{\wedge+}$ is
a subalgebra of $\cB(r)^\wedge$, hence it is reduced as well.
Hence lemma \ref{lem_immediate}(ii) applies, and reduces to
showing that $(\cB(r)^{\wedge+},|\cdot|^\wedge_\sp)$ is a free
cartesian $\kappa(r)^{\wedge+}$-module of rank $d$.
Furthermore, proposition \ref{prop_normal}(iii) implies that
$\cB(r)^{\wedge+}$ is normal, hence it is the direct product
\set\begin{equation}\label{eq_decomp-cB}
\cB(r)^{\wedge+}=L^+_1\times\cdots\times L^+_k
\end{equation}
of finitely many normal domains, and each $L^+_i$ is the integral
closure of $\kappa(r)^{\wedge+}$ in a finite algebraic
extension $L_i$ of $\kappa(r)^\wedge$. Notice as well, that
$\kappa(r^\flat)^{\wedge +}$ is henselian, since it is
complete and of rank one; hence
$L^+_i\otimes_{\kappa(r)^{\wedge+}}\kappa(r^\flat)^{\wedge+}$ is a
valuation ring whose valuation extends $|\cdot|^\flat_{\eta(r)}$.
On the other hand, by lemma \ref{lem_much-ado} the points
of $x\in f^{-1}(\eta(r))$ correspond to the valuations
$|\cdot|_x$ on $\cB(r)$ that extend $|\cdot|_{\eta(r)}$;
these are also the valuations on $\cB(r)^\wedge$ that
extend $|\cdot|^\wedge_{\eta(r)}$. Hence, the decomposition
\eqref{eq_decomp-cB} induces a partition
$$
f^{-1}(\eta(r))=\Sigma_1\cup\cdots\cup\Sigma_k
$$
where, for each $i\leq k$, $\Sigma_i$ is the set of valuations
of $L_i$ that extend $|\cdot|^\wedge_{\eta(r)}$.

By lemma \ref{lem_nobigdeal}(ii), we are reduced to showing :

\begin{claim} For every $i\leq k$, let $|\cdot|_{\sp,i}$
be the spectral norm of the $\kappa(r)^\wedge$-algebra
$L_i$; then $(L^+_i,|\cdot|_{\sp,i})$ is a
$\kappa(r)^{\wedge+}$-cartesian module.
\end{claim}
\begin{pfclaim} It suffices to show that the finite field
extension $\kappa(r)^\wedge\subset L_i$ fulfills
conditions (a)--(d) of proposition \ref{prop_mess}.
However, condition (a) is none else than
\cite[Lemma 5.3(ii)]{Hu3}.
Next, say that $\Sigma_i=\{x_1,\dots,x_l\}$, and let
$\Gamma_1,\dots,\Gamma_l$ be the value groups of the residue
fields $\kappa(x_i)$. By inspecting the construction, one sees
easily that $x_i^\flat=x_j^\flat$ for every $i,j\leq l$,
which means that the subgroup $\Delta_{ij}\subset\Gamma_i$
defined as in \eqref{subsec_artificial}, is the unique
convex subgroup corresponding to $x_i^\flat$, so (d) holds,
and also (b) is clear. Finally, (c) follows easily from
\cite[Cor.5.4 and Prop.1.2(iii)]{Ga-Ra}.
\end{pfclaim}

(ii): In view of lemma \ref{lem_nobigdeal}(i.a), it suffices
to show the analogous identity for
$$
(\cB(r)^{\wedge h},|\cdot|^{\wedge h}_{\sp,\eta(r)}):=
(\cB(r),|\cdot|_{\sp,\eta(r)})\otimes_{\kappa(r)}\kappa(r)^{\wedge h}.
$$
By the proof of (i) we know that
$\cB(r)^{\wedge h}=\prod_{x\in f^{-1}(\eta(r))}\kappa(x)^{\wedge h}$,
so the assertion follows from lemma \ref{lem_nobigdeal}(ii),(iii).

(iii): In light of lemma \ref{lem_immediate}(i) it suffices to show
the same identity for $|\cB(r)^{\wedge h+}|^{\wedge h\natural}_\sp$;
however, from \cite[Prop.1.2(iii) and Cor.5.4]{Hu3} we deduce that
$$
|\kappa(x)^{\wedge h+}|^{\wedge h\natural}_x=
\frac{[\kappa(x)^{\wedge h}:\kappa(r)^{\wedge h}]}{2}
$$
for every $x\in f^{-1}(\eta(r))$; then the assertion follows easily.

(iv): Suppose first that $R/\fm R$ is reduced.
According to (i), $R$ is a finitely generated submodule of a free
$\kappa(r)^+$-module of finite rank; hence it is free as a
$\kappa(r)^+$-module. Hence every $x\in R$
can be written in the form $x=ay$ for some $a\in K^+$ and an
element $y\in R$ whose image in $R/\fm R$ does not vanish.
It follows that every $x\in\cB(r)$ can be written in the form
$x=ay$ for some $a\in K$ and some $y\in R\setminus\fm R$.
Let now $x\in\cB(r)^+$, and suppose that $x=ay$ for some
$y\in R\setminus\fm R$ and $a\in K$ with $|a|>1$; clearly
$x$ is integral over $R$, so we can write
$$
x^n+b_1x^{n-1}+\cdots+b_n=0
$$
for some $b_1,\dots,b_n\in R$. Hence
$$
y^n+b_1a^{-1}y^{n-1}+\cdots+b_na^{-n}=0.
$$
In other words, $y^n\in\fm R$, whence $y\in\fm R$, since $R/\fm R$
is reduced; the contradiction shows that $R=\cB(r)^+$.
Conversely, suppose that $R=\cB(r)^+$ and let $x\in R$ whose
image in $R/\fm R$ is nilpotent; then $x^n\in\fm R$ for $n\in\N$
large enough, say $x^n=ay$ for some $a\in\fm$ and $y\in R$.
We can write $a=b^nc$ for $b,c\in\fm$, so $(x/b)^n=cy\in R$,
so $x/b\in R$, since the latter is integrally closed in
$R\otimes_{K^+}K$. Finally, $x\in\fm R$, as claimed.
\end{proof}

\begin{lemma}\label{lem_estimate}
In the situation of \eqref{subsec_direct-image}, let
$r\in(a,b]\cap\Gamma_K$, $U\subset\D(a,b)$ an open
neighborhood of $\eta(r)$, $s\in\Gamma(U,f_*\cO_X)$
and set $k:=|s|_{\sp,\eta(r)}^\natural\in\frac{1}{d!}\Z$.
Then there exists $r'\in(a,r)$ such that:
\begin{enumerate}
\item\ \ 
$\eta(t)\in U$ for every $t\in(r',r]\cap\Gamma_K$.
\item\ \ 
$|s|_{\sp,\eta(t)}=|s|_{\sp,\eta(r)}\cdot(t/r)^k$\ \  whenever
$t\in(r',r]\cap\Gamma_K$.
\end{enumerate}
\end{lemma}
\begin{proof} (i) is obvious. In order to prove (ii), consider
an integral equation
$$
s^n+g_1s^{n-1}+\cdots+g_n=0
$$
where $g_i\in\Gamma(U,\cO_{\D(a,b)})$ for every $i\leq n$.
It is easy to see that the assertion for $s$ will follow once
the same assertion is known for the sections $g_1,\dots,g_n$.
Hence, we are reduced to the case $X=\D(a,b)$. We can also
assume that there exist $\alpha_1,\dots,\alpha_m\in K$ such
that $|\alpha_i|=r$ for every $i\leq m$ and
$$
U=\D(r',r)\setminus\bigcup_{i=1}^m\E_i
$$
where $\E_i=\{p\in\D(r',r)~|~|\xi-\alpha_i|_p<r\}$ (so each $\E_i$
is a closed subset of $\D(r',r)$). Then, in view of
\cite[Prop.2.8]{Fr-VdP} we have a Mittag-Leffler decomposition
$$
s=s_0+s_1+\sum_{i=1}^mh_i
$$
where
$$
s_0=\sum_{n=0}^\infty a_n\xi^n \qquad
s_1=\sum_{n=1}^\infty b_n\xi^{-n} \qquad
h_i=\sum_{n=1}^\infty c_{i,n}(\xi-\alpha_i)^{-n}\quad (i=1,\dots,m)
$$
for coefficients $a_n,b_n,c_{i,n}\in K$ subject to the conditions:
$$
\lim_{n\to\infty}r^n|a_n|=\lim_{n\to\infty}r^{\prime -n}|b_n|=
\lim_{n\to\infty}r^{-n}|c_{i,n}|=0.
$$
Recalling the standard power series identity:
$$
(\xi-\alpha_i)^{-n}=
\sum_{k=0}^\infty\binom{-n}{k} (-\alpha_i)^{-n-k}\xi^k
$$
we deduce that
$$
|s_1|_{\eta(t)}=\sup_{n>0}|b_n|\cdot t^{-n}(1-\eps)^{-n}
\qquad\text{for every $t\in(r'r]\cap\Gamma_K$}
$$
and
$$
\left|s_0+\sum_{i=1}^mh_i\right|_{\eta(t)}=
\sup_{k\in\N}|d_k|\cdot t^k(1-\eps)^k
\quad\text{for every $t\in(r'r]\cap\Gamma_K$}.
$$
where
$$
d_k=a_k+\sum_{i=1}^m\sum^\infty_{n=1}\binom{-n}{k} (-\alpha_i)^{-n-k}c_{i,n}
\quad\text{for every $k\in\N$}.
$$
Notice that $|s_1|_{\eta(t)}\neq|s_0+\sum_{i=1}^mh_i|_{\eta(t)}$
for every $t\in(r',r]\cap\Gamma_K$, therefore 
$$
|s|_{\eta(t)}=
\max\left(|s_1|_{\eta(t)},\left|s_0+\sum_{i=1}^mh_i\right|_{\eta(t)}\right)
$$
and hence it suffices to prove assertion (ii) separately for
$s_1$ and $s_0+\sum_{i=1}^mh_i$.

The coefficients $d_k$ enjoy the following property. There exists
a smallest $k_0\in\N$ such that
$$
|d_{k_0}|=\sup_{k\in\N}|d_k|\cdot r^k.
$$
Namely, $k_0$ is the unique integer such that
$$
\left|s_0+\sum_{i=1}^mh_i\right|_{\eta(r)}=
|d_{k_0}|\cdot r^{k_0}(1-\eps)^{k_0}.
$$
Thus, when evaluating $|s_0+\sum_{i=1}^mh_i|_{\eta(t)}$ for $t<r$,
we can disregard all the terms $d_k\xi^k$ for $k>k_0$, since their
norms decrease faster than that of the leading term $d_{k_0}\xi^{k_0}$.
There remain to consider the finitely many monomials
$d_0, d_1\xi,\dots,d_{k_0-1}\xi^{k_0-1}$; however, it is clear
from the definition of $k_0$ that the norms of these terms
are still strictly smaller than the norm of the leading term
$d_{k_0}\xi^{k_0}$, as long as $t$ is sufficiently close to $r$.
Thus, we can replace the section $s_0+\sum_{i=1}^mh_i$ by its
leading monomial, for which the sought assertion trivially holds;
an analogous, though easier argument also works for $s_1$ : we
leave the details to the reader.
\end{proof}

\subsection{Convexity and piecewise linearity of the discriminant
function}\label{sec_convexity}
The assumptions and notations are as in \eqref{sec_annuli}.
The statements proven so far make use of only a few relatively
simple local properties of the sheaf $f_*\cO^+_X$; nevertheless,
they would already suffice to prove most of the forthcoming
proposition \ref{prop_precise}.
However, in order to show theorem \ref{th_big-deal}, it will be
necessary to cast a closer look at the ring of global sections
of $f_*\cO^+_X$; the following lemmata \ref{lem_finitely-pres},
\ref{lem_reduction} and proposition \ref{prop_free} will provide
us with everything we need.

\begin{lemma}\label{lem_finitely-pres} Let $A\to B$ be a finite
morphism of $K$-algebras of topologically finite type. Then $A^\circ$
is of topologically finite presentation over $K^+$ and $B^\circ$ is
a finitely presented $A^\circ$-module.
\end{lemma}
\begin{proof} By \cite[\S6.4.1, Cor.5]{Bo-Gun} we know that $B^\circ$
is a finite $A^\circ$-module; moreover, applying {\em loc.cit.} to an
epimorphism $K\langle T_1,\dots,T_n\rangle\to A$ we deduce that
$A^\circ$ is finite over $K^+\langle T_1,\dots,T_n\rangle$.
By \cite[Prop.7.1.1(i)]{Ga-Ra} we deduce that $A^\circ$ is finitely
presented as a $K^+\langle T_1,\dots,T_n\rangle$-module, and also that
$B^\circ$ is finitely presented over $A^\circ$.
\end{proof}

\begin{lemma}\label{lem_reduction}
Let $B$ be a flat, reduced $K^+$-algebra of topologically
finite type, and set $A:=B\otimes_{K^+}K$. Then $B=A^\circ$ if
and only if $B/\fm B$ is reduced.
\end{lemma}
\begin{proof} (Notice that the assertion is a global version of
lemma \ref{lem_local-sit}(iv), and indeed its proof is analogous
to that of the former.)

\begin{claim}\label{cl_Bfree}
$B/aB$ is a free $K^+/aK^+$-module for every $a\in\fm$.
\end{claim}
\begin{pfclaim} By \cite[Prop.7.1.1(i)]{Ga-Ra} we have
$B\simeq K^+\langle T_1,\dots,T_r\rangle/I$
for some $r\geq 0$ and a finitely generated ideal $I$. It follows
that $B/aB\simeq K^+/aK^+[T_1,\dots,T_r]/J$, where $J$ is the image of $I$.
We can write $K^+/aK^+=\bigcup_{\lambda\in\Lambda} R_\lambda$, the
filtered union of its noetherian local subalgebras $R_\lambda$.
By \cite[Ch.IV, Prop.8.5.5]{EGAIV-3} we can find $\lambda\in\Lambda$
and a finitely presented flat $R_\lambda$-algebra $B_\lambda$ such
that $B/aB\simeq B_\lambda\otimes_{R_\lambda}K^+/aK^+$. It suffices
to show that $B_\lambda$ is a free $R_\lambda$-module; however,
since $R_\lambda$ is artinian, this follows from \cite[Th.7.9]{Mat}.
\end{pfclaim}

\begin{claim} $A^\circ$ is the integral closure of $B$ in $A$.
\end{claim}
\begin{pfclaim} Choose a continuous surjection
$\phi:C:=K^+\langle T_1,\dots,T_r\rangle\to B$. 
It suffices then to notice that $C=(C\otimes_{K^+}K)^\circ$ and
apply \cite[\S6.3.4, Prop.1]{Bo-Gun}.
\end{pfclaim}

By claim \ref{cl_Bfree}, every $x\in B\setminus\{0\}$ can be written
in the form $x=ay$ for some $a\in K^+$ and an element $y\in B$ whose
image in $B/\fm B$ does not vanish. It follows that every
$x\in A\setminus\{0\}$ can be written in the form $x=ay$ for some
$a\in K$ and some $y\in B\setminus\fm B$. After these remarks, one
can proceed as in the proof of lemma \ref{lem_local-sit}(iv): the
details shall be entrusted to the reader.
\end{proof}

\begin{proposition}\label{prop_free} Let $(F,|\cdot|_F)$ be a
complete algebraically closed valued field extension of $K$,
such that $|\cdot|_F$ is a valuation of rank one. Let $A$ be
a normal domain of topologically finite type over $K$, such
that $A^\sim$ is a principal ideal domain, $B$ a finite,
reduced and flat $A$-algebra, and $g\in A$ such that $|g|_{\sup}=1$.
Set $A_F:=A\widehat\otimes_KF$. Then :
\begin{enumerate}
\item
$B^\circ$ is a free $A^\circ$-module of finite rank.
\item
$B\langle g^{-1}\rangle^\circ=
B^\circ\otimes_{A^\circ}A\langle g^{-1}\rangle^\circ$.
\item
$(B\otimes_AA_F)^\circ=B^\circ\otimes_{A^\circ}A_F^\circ$.
\end{enumerate}
\end{proposition}
\begin{proof} To start out, let us endow $A$ and $B$ with their
supremum norms; then by \eqref{subsec_conflict} and 
\cite[\S3.8.1, Prop.7]{Bo-Gun}, we have $A^\circ=A^+$ and
$B^\circ=B^+$. By lemma \ref{lem_finitely-pres} we deduce
that $A^+$ is of topologically finite type over $K^+$ and
$B^+$ is finitely presented over $A^+$.

\begin{claim}\label{cl_free}
$B^\sim$ is free of finite rank over $A^\sim$.
\end{claim}
\begin{pfclaim} By the foregoing we know already that $B^\sim$
is finite over $A^\sim$, hence it suffices to show that
$B^\sim$ is torsion-free as an $A^\sim$-module. However, under
the current assumptions the norm $|\cdot|_{\sup}$ on $A$ is
a valuation (\cite[\S6.2.3, Prop.5]{Bo-Gun}). It follows that
\set\begin{equation}\label{eq_max}
|b|_\sp=\max_{v}v(b)\qquad\text{for all $b\in B$}
\end{equation}
where $v$ ranges over the finitely many extensions of
the supremum valuation of $A$ to $B$ (\cite[\S3.3.1, Prop.1]{Bo-Gun}).
For each such $v$, let $\mathrm{supp}(v):=v^{-1}(0)$, which is
a prime ideal of $B$, and denote by
$B_v\subset\mathrm{Frac}(B/\mathrm{supp}(v))$ the valuation
ring of the valuation induced by $v$ on $B/\mathrm{supp}(v)$.
Since $\Gamma_K$ is divisible, it is easy to see that
$\fm B_v$ is the maximal ideal of $B_v$. From \eqref{eq_max}
it is clear that $B^\sim\subset\prod_vB_v/\fm B_v$.
Finally, for every $v$ the field $B_v/\fm B_v$ is a finite
extension of $\mathrm{Frac}(A^\sim)$, especially it is torsion-free
over $A^\sim$, and the same holds then for $B^\sim$.
\end{pfclaim}

From claim \ref{cl_free} and \cite[Ch.IV, Prop.8.5.5]{EGAIV-3} it
follows that there exists $\pi\in\fm$ such that $B^+/\pi B^+$ is
flat over $A^+/\pi A^+$. In view of \cite[Lemma 5.2.1]{Ga-Ra}
we conclude that $B^+$ is flat, hence projective over $A^+$.
Finally, a standard application of Nakayama's lemma shows that
any lifting of a basis of $B^\sim$ is a basis of the $A^+$-module
$B^+$, which proves (i).
\begin{claim}\label{cl_candidate}
The ring $C:=B^\circ\otimes_{A^\circ}A\langle g^{-1}\rangle^\circ$
is reduced.
\end{claim}
\begin{pfclaim}
From (i) we deduce that the natural map
$$
B^\circ\otimes_{A^\circ}A\langle g^{-1}\rangle^\circ\to
B\langle g^{-1}\rangle=B^\circ\otimes_{A^\circ}A\langle g^{-1}\rangle
$$
is injective. Hence it suffices to show that $B\langle g^{-1}\rangle$
is reduced, which holds by lemma \ref{lem_subdomain}.
\end{pfclaim}

In view of claim \ref{cl_candidate} and lemma \ref{lem_reduction},
assertion (ii) will follow once we know that
$C/\fm_kC$ is reduced. However, the latter is isomorphic
to $B^\sim\otimes_{A^\sim}A^\sim[\bar g{}^{-1}]$, where
$\bar g\in A^\sim$ is the image of $g$ (\cite[\S7.2.6, Prop.3]{Bo-Gun}).
Again lemma \ref{lem_reduction} ensures that $B^\sim$ is reduced.

(iii): According to \cite[Lemma 3.3.1.(1)]{Conr},
$B_F:=B\otimes_AA_F=B^\circ\otimes_{A^\circ}A_F$ is reduced.
From (i) we deduce that $D:=B^\circ\otimes_{A^\circ}A_F^\circ$
is a subalgebra of $B_F$, so it is reduced as well. Hence in order to
prove (iii) it suffices, by lemma \ref{lem_reduction}, to show that
$D/\fm_F D$ is reduced (where $\fm_F$ is the maximal ideal of $F^+$).
However, $D/\fm_F D\simeq B^\sim\otimes_{K^\sim}F^\sim$ and the
extension $K^\sim\to F^\sim$ is separable, so everything is clear.
\end{proof}

\sset\subsubsection{}
Let $a,b\in\Gamma_K$ with $a<b$ and $f:X\to\D(a,b)$ a finite,
flat and generally \'etale morphism, say of degree $d$. For every
$r\in[a,b]\cap\Gamma_K$ we define:
$$
\cB(r)^\flat:=(f_*\cO^+_X)_{\eta(r)^\flat}.
$$
It follows easily from proposition \ref{prop_free}(ii) that
\set\begin{equation}\label{eq_asymmetric}
\cB(r)^\flat=\cB(r)^+\otimes_{\kappa(r)^+}\kappa(r^\flat)^+
\qquad\text{for every $r\in(a,b]\cap\Gamma_K$}
\end{equation}
(notation of \eqref{subsec_localize} and \eqref{subsec_direct-image}).

\sset\subsubsection{}\label{subsec_asymmetry}
The apparent asymmetry between the values $a$ and $b$ can be
easily resolved. Indeed, let us consider the isomorphism
$$
q:\D(1/b,1/a)\to\D(a,b)
\quad : \quad \xi\mapsto\xi^{-1}
$$
and let $g:=q^{-1}\circ f$. For every $r\in(1/b,1/a]$, the
image $q(\eta(r))$ is the point $\eta'(1/r)$ 
(notation of \eqref{subsec_localize}).
Hence, $g^*:\cO_{\D(1/b,1/a)}\to g_*\cO_X$
endows $(f_*\cO^+_X)_{\eta'(a)}$ with a structure of
$\kappa(1/a)^+$-algebra, and since
$\eta'(a)^\flat=\eta(a)^\flat$, we deduce that an identity
analogous to \eqref{eq_asymmetric} holds also for $r=a$, provided
we replace $\eta(a)$ by $\eta'(a)$. Especially, since -- according 
to lemma \ref{lem_local-sit}(i) -- the stalk $\cB(r)^+$ is a free
$\kappa(r)^+$-module of rank $d$, we see that $\cB(r)^\flat$
is a free $\kappa(r^\flat)^+$-module of rank $d$ for every
$r\in[a,b]\cap\Gamma_K$.

\sset\subsubsection{}\label{subsec_discr-at-r}
Now, as $f$ is generically \'etale, the trace forms
$\Tr_{\cB(r)^\flat/\kappa(r^\flat)^+}$ and
$\Tr_{\cB(r)^+/\kappa(r)^+}$ induce the same perfect
pairing after tensoring with $\kappa(r^\flat)$. We set
$$
\fd_f^+(r):=\fd_{\cB(r)^+/\kappa(r)^+}\in
\kappa(r)^+
\qquad\text{for every $r\in(a,b]\cap\Gamma_K$.}
$$
(respectively:
$$
\fd_f^\flat(r):=\fd_{\cB(r)^\flat/\kappa(r^\flat)^+}\in
\kappa(r^\flat)^+
\qquad\text{for every $r\in[a,b]\cap\Gamma_K$.}
$$
Notation of \eqref{sec_discrim}.) Since $\fd_f^\flat(r)$
is well defined up to the square of an invertible element
of $\kappa(r^\flat)^+$, the real-valued function:
$$
\delta_f:[\log 1/b,\log 1/a]\cap\log\Gamma_K\to\R_{\geq 0}\qquad
-\log r\mapsto -\log|\fd_f^\flat(r)|_{\eta(r)^\flat}
$$
is well defined independently of all choices. Unless we have to
deal with more than one morphism, we shall usually drop the subscript,
and write $\delta$, $\fd^\flat$ instead of $\delta_f$, $\fd^\flat_f$.
We call $\delta$ the {\em discriminant function\/} of the morphism
$f$.

\begin{example}\label{ex_summing-up} Let $f:X\to\D(a,a^{-1})$
be a finite, flat and generically \'etale morphism, where $a:=|\pi|$
for some $\pi\in\fm$. Using the Mittag-Leffler decomposition
\cite[Prop.2.8]{Fr-VdP} one verifies easily that
$A(a,a^{-1})^\circ=K^+\langle\pi/\xi,\xi/\pi\rangle\simeq
K^+\langle S,T\rangle/(ST-\pi^2)$ (alternatively, one sees this
via lemma \ref{lem_reduction}).
We set
$h:=\Spa(\psi_K)\circ f:X\to\D(1)$, where $\psi_K:=\psi\otimes_{K^+}K$,
with $\psi$ defined as in example \ref{ex_discriminant}. A direct
computation shows that
$$
h^{-1}(\D(r,1))=f^{-1}(\D(a,a/r))\cup
f^{-1}(\D(r/a,a^{-1}))
\qquad\text{for every $r\in(a,1]\cap\Gamma_K$}.
$$
Consequently:
$$
\fd^\flat_h(r)=\fd^\flat_f(r/a)\cdot\fd^\flat_f(a/r)
\qquad\text{whenever $r\in(a,1]\cap\Gamma_K$}
$$
and therefore
\set\begin{equation}\label{eq_summing-up}
\delta_h(-\rho)=\delta_f(\rho-\log a)+\delta_f(\log a-\rho)
\qquad\text{for $\rho\in(\log a,0]\cap\log\Gamma_K$}.
\end{equation}
Incidentally, let $\eta'(a)\in\D(1)$ be defined as in
\eqref{subsec_localize}; it is easy to check that the
preimage of $\eta'(a)$ in $\D(a,a^{-1})$ under $\Spa\,\psi_K$
is the subset $\{\eta(1),\eta'(1)\}$.
\end{example}

\sset\subsubsection{}
Let $f:[r,s]\to\R$ be a piecewise linear function; for every
$\rho\in[r,s)$ we denote by $df/dt(\rho^+)$ the {\em right slope\/}
of $f$ at the point $r$, {\em i.e.} the unique real number
$\alpha$ such that $f(\rho+x)=f(\rho)+\alpha x$ for every sufficiently
small $x\geq 0$. Similarly we can define the {\em left slope\/}
$df/dt(\rho^-)$ for every $\rho\in(r,s]$. More generally, the
definition makes sense whenever $f$ is defined on a dense subset
of $[r,s]$.

\begin{example}\label{ex_asymm}
Let $f$ and $g$ be as in \eqref{subsec_asymmetry}. Then
$$
\delta_f(\rho)=\delta_g(-\rho)
\qquad\text{and}\qquad
\frac{d\delta_f}{dt}(\rho^-)=-\frac{d\delta_g}{dt}(-\rho^+)
$$
for every $\rho\in(\log 1/b,\log 1/a]\cap\log\Gamma_K$.
\end{example}

\begin{proposition}\label{prop_precise}
With the notation of \eqref{subsec_discr-at-r}, the function
$\delta$ is piecewise linear;  moreover:
$$
\frac{d\delta}{dt}(-\log r^+)=|\fd^+(r)|^\natural_{\eta(r)}-
2|\cB(r)^+|^\natural_{\sp,\eta(r)}
\qquad\text{for every $r\in(a,b]\cap\Gamma_K$.}
$$
(Notation of \eqref{subsec_A^+-submodules}.)
\end{proposition}
\begin{proof} By lemma \ref{lem_local-sit}(i) we can find
an orthogonal basis $c_1,\dots,c_d$ of $(\cB(r)^+,|\cdot|_{\sp,\eta(r)})$
over $(\kappa(r)^+,|\cdot|_{\eta(r)})$. Obviously we have
$|c_i|_{\sp,\eta(r)}^\flat=1$ for every $i=1,\dots,d$; we set
\set\begin{equation}\label{eq_for.c_i}
\gamma_i:=|c_i|^\natural_{\sp,\eta(r)}\in\frac{1}{d!}\Z
\qquad \text{for $i=1,\dots,d$}.
\end{equation} 
For every $i,j\leq d$ we can find uniquely determined
$m_{ij1},\dots,m_{ijd}\in\kappa(r)^+$ such that
$$
c_i\cdot c_j=\sum_{k=1}^d m_{ijk}c_k.
$$
For every $i,j,k\leq d$ we set
\set\begin{equation}\label{eq_interpret}
x_{ijk}:=|m_{ijk}|^\flat_{\eta(r)}
\qquad\text{and}\qquad
\mu_{ijk}:=|m_{ijk}|^\natural_{\eta(r)}
\end{equation}
so that $x_{ijk}\in[0,1]$ and $\mu_{ijk}\in\N$.
Since the $c_i$ are orthogonal, it follows easily that
\set\begin{equation}\label{eq_bound-for-mu}
\mu_{ijk}\geq\gamma_i+\gamma_j-\gamma_k\qquad\text{whenever $x_{ijk}=1$}.
\end{equation}
Lemma \ref{lem_estimate}(ii) implies that
\set\begin{equation}\label{eq_choose-again}
\frac{m_{ijk}}{a^{\mu_{ijk}}},\frac{c_i}{a^{\gamma_i}}\in\cB(sr)^+
\qquad\text{for every $i,j,k=1,\dots,d$}
\end{equation}
whenever $a\in K^+$ is an element with $s:=|a|$ sufficiently
close to $1$. Clearly we have
\set\begin{equation}\label{eq_new.product}
\frac{c_i}{a^{\gamma_i}}\cdot\frac{c_j}{a^{\gamma_j}}=
\sum_k\frac{m_{ijk}}{a^{\gamma_i+\gamma_j-\gamma_k}}\cdot
\frac{c_k}{a^{\gamma_k}}
\end{equation}
which, in view of \eqref{eq_bound-for-mu} and \eqref{eq_choose-again},
means that
\set\begin{equation}\label{eq_finite}
\cC(sr):=\sum_{i=1}^d\kappa(sr)^+\cdot\frac{c_i}{a^{\gamma_i}}
\end{equation}
is a finite (unitary) subalgebra of $\cB(sr)^+$ for every
$s:=|a|$ sufficiently close to $1$.

Let $\cB(r)^\sim:=\cB(r)^+\otimes_{K^+}K^\sim$. The map
$b\mapsto|b|^\natural_{\sp,\eta(r)}$ induces a power-multiplicative
norm on $\cB(r)^\sim$, whence a filtration $\Fil^\bullet\cB(r)^\sim$
defined by setting:
$$
\Fil^h\cB(r)^\sim:=
\{b\in\cB(r)^\sim~|~|b|^\natural_{\sp,\eta(r)}\geq h\}
\qquad\text{for every $h\in\frac{1}{d!}\Z$}.
$$
The filtration $\Fil^\bullet\cB(r)^\sim$ restricts to a
filtration $\Fil^\bullet\cO_{\eta(r)}^\sim$ on
$\cO_{\eta(r)}^\sim:=\kappa(r)^+\otimes_{K^+}K^\sim$.
The associated graded ring $\gr^\bullet\cB(r)^\sim$ can be
computed explicitly. Indeed, \eqref{eq_for.c_i} shows that
$\gr^\bullet\cB(r)^\sim=
\gr^\bullet\cO_{\eta(r)}^\sim[\bar c_1,\dots,\bar c_d]$, where
$\bar c_i$ is the class of $c_i$ in $\gr^{\gamma_i}\cB(r)^\sim$
for $i=1,\dots,d$. Furthermore, for every $i,j\leq d$ we have
the rule:
$$
\bar c_i\cdot\bar c_j=\sum_{k=1}^d\beta_{ijk}\bar c_k
$$
where $\beta_{ijk}$ is determined as follows. If either $x_{ijk}<1$
or $\mu_{ijk}>\gamma_i+\gamma_j-\gamma_k$, then $\beta_{ijk}=0$
and otherwise $\beta_{ijk}$ is the class of $m_{ijk}$ in
$\gr^{\mu_{ijk}}\cO^\sim_{\eta(r)}$.

\begin{claim}\label{cl_C(r).mod.m}
$\cC(sr)/\fm\cC(sr)\simeq\gr^\bullet\cB(r)^\sim$ for every
$s\in\Gamma^+_K\setminus\{1\}$ sufficiently close to $1$.
\end{claim}
\begin{pfclaim} With the notation of \eqref{eq_finite} we have
$$
\cC(sr)/\fm\cC(sr)=
\cO_{\eta(sr)}^\sim
\left[\frac{c_1}{a^{\gamma_1}},\dots,\frac{c_d}{a^{\gamma_d}}\right].
$$
We define an isomorphism 
\set\begin{equation}\label{eq_was-f}
\gr^\bullet\cO^\sim_{\eta(r)}\simeq\cO_{\eta(sr)}^\sim
\end{equation}
by the rule:
$$
% \xi\mapsto\xi/a.
g\mapsto g/a^n\pmod{\fm\cO^+_{\eta(sr)}}
\qquad\text{for every $g\in\gr^n\cO^\sim_{\eta(r)}$}.
$$
Via \eqref{eq_was-f}, $\cC(sr)/\fm\cC(sr)$ becomes a
free $\gr^\bullet\cO^\sim_{\eta(r)}$-module, whose rank
$d$ is the same as the rank of $\gr^\bullet\cB(r)^\sim$.
Obviously, we would like to extend the isomorphism
\eqref{eq_was-f} by setting 
\set\begin{equation}\label{eq_souhait}
\bar c_i\mapsto c_i/a^{\gamma_i}
\qquad\text{for every $i=1,\dots,d$}.
\end{equation}
 In view of \eqref{eq_new.product}, in order
to prove that \eqref{eq_souhait} yields a well-defined ring
homomorphism, it suffices to show that the class of
$m_{ijk}/a^{\gamma_i+\gamma_j-\gamma_k}$ in $\cO^\sim_{\eta(sr)}$
agrees with the image of $\beta_{ijk}$ under \eqref{eq_was-f},
whenever $s$ is sufficiently close to $1$. This can be checked
easily by inspecting the definitions; we leave the details
to the reader.
\end{pfclaim}

\begin{claim}\label{cl_C(r)} $\cC(sr)=\cB(sr)^+$ whenever
$s\in\Gamma^+_K$ is sufficiently close to $1$.
\end{claim}
\begin{pfclaim} In view of claim \ref{cl_C(r).mod.m} and
lemma \ref{lem_local-sit}(iv), it suffices to show that
$\gr^\bullet\cB(r)^\sim$ is reduced. But this is clear,
since the norm $|\cdot|^\natural_{\sp,\eta(r)}$ on
$\cB(r)^\sim$ is power-multiplicative.
\end{pfclaim}

In view of claim \ref{cl_C(r)} we have
\set\begin{equation}\label{eq_estimate.D}
\fd^+(sr)=\det\left(\Tr_{\cB(rs)^+/\kappa(rs)^+}
\left(\frac{c_i}{a^{\gamma_i}}\otimes\frac{c_j}{a^{\gamma_j}}\right)
~|~1\leq i,j\leq d\right)=
a^{-2\cdot\sum_i\gamma_i}\cdot\fd^+(r)_{\eta(rs)}
\end{equation}
whenever $s:=|a|$ is sufficiently close to $1$ (here
$\fd^+(r)_{\eta(rs)}$ denotes the image of $\fd^+(r)$ in
$\kappa(rs)^+$; this is well-defined whenever $s$
is sufficiently close to $1$). However, one deduces easily
from \eqref{eq_asymmetric} that
$|\fd^+(t)|^\flat_{\eta(t)}=|\fd^\flat(t)|_{\eta(t)^\flat}$
for every $t\in(a,b]\cap\Gamma_K$. Since we have as well:
$$
|\cB(r)^+|^\natural_{\sp,\eta(r)}=\sum_{i=1}^d\gamma_i
$$
the contention follows from lemma \ref{lem_estimate}(ii).
\end{proof}

\sset\subsubsection{}
Proposition \ref{prop_precise} expresses the slope of
the discriminant function at a given radius $r$ as a
local invariant depending only on the behaviour of the
morphism $f$ over the point $\eta(r)$. We wish now to
consider a special situation, where the slope can also
be obtained as a global invariant of the ring
$\Gamma(X,\cO^+_X)$. Namely, suppose that $g:X\to\D(1)$
is a finite, flat and generically \'etale morphism;
proceeding as in the foregoing we attach to $g$ a
discriminant function $\delta_g$, which clearly shall
be defined over the set
$[0,+\infty)\cap\log\Gamma_K=-\log\Gamma_K^+$.
However, our present aim is to compute the right slope of
$\delta_g(\rho)$ in a small neighborhood of $\rho=0$.
To this purpose, let $B^\circ:=\Gamma(X,\cO^+_X)$ and
$B:=B^\circ\otimes_{K^+}K$; according to proposition
\ref{prop_free}(i), $B^\circ$ is a free module,
necessarily of rank $d:=\deg(g)$, over the ring
$A(1)^\circ=\Gamma(\D(1),\cO^+_{\D(1)})$. Clearly
$B\otimes_{A(1)}\kappa(1)=\cB(1)$, hence the
natural map
\set\begin{equation}\label{eq_embed-glob}
B^\circ_\eta:=B^\circ\otimes_{A(1)^\circ}\kappa(1)^+
\to\cB(1)^+
\end{equation}
is injective. Let $\fd^\circ_g:=\fd_{B^\circ/A(1)^\circ}$;
combining lemmata \ref{lem_Fitt-discr} and \ref{lem_Fitting}
we deduce:
\set\begin{equation}\label{eq_local-to-glob}
|\fd_g^+(1)|^\natural_{\eta(1)}-
2\cdot|\cB(1)^+|^\natural_{\sp,\eta(1)}=
|\fd^\circ_g|^\natural_{\eta(1)}-
2\cdot|B^\circ_\eta|^\natural_{\sp,\eta(1)}.
\end{equation}
Notice that the left-hand side of this identity calculates
the right slope of $\delta_g$ at the point $\rho=0$, hence the
right-hand side is the sought global formula for this
slope. 

\sset\subsubsection{}\label{subsec_analyze}
The contribution $|B^\circ_\eta|^\natural_{\sp,\eta(1)}$
can be further analyzed. Indeed, let us set
$$
B^{\circ h}_\eta:=B^\circ\otimes_{A(1)^\circ}\kappa(1)^{\wedge h+}.
$$
The ring $B^{\circ h}_\eta$ is henselian along the ideal
$\fp_\eta B^{\circ h}_\eta$, where $\fp_\eta$ is the maximal ideal of
$\kappa(1)^+$. Let $\fq_1,\dots,\fq_k\subset B^\circ_\eta$
be the finitely many prime ideals lying over $\fp_\eta$; the
ring $B^{\circ h}_\eta$ decomposes as a direct product of henselian
local rings:
$$
B^{\circ h}_\eta=B^{\circ h}_{\fq_1}\times\cdots\times B^{\circ h}_{\fq_k}.
$$
For every $i=1,\dots,k$ set
$$
\fF(\fq_i):=\{x\in g^{-1}(\eta(1))~|~
\text{$\kappa(x)^+$ dominates $B^\circ_{\fq_i}$}\}.
$$
After completion, henselization and localization at $\fq_i$, the map
\eqref{eq_embed-glob} yields injective ring homomorphisms
(see the proof of lemma \ref{lem_local-sit}(i)) :
\set\begin{equation}\label{eq_integral-form}
B^{\circ h}_{\fq_i}\to\cB(1)^{\wedge h+}_{\fq_i}\simeq
\prod_{x\in\fF(\fq_i)}\kappa(x)^{\wedge h+}
\end{equation}
More precisely, let $\bar\kappa(\fq_i)$ (resp. $\bar\kappa(x)$) be the
residue field of $B^{\circ h}_{\fq_i}$ (resp. of $\kappa(x)^{\wedge h+}$);
the maps \eqref{eq_integral-form} induce isomorphisms
$\bar\kappa(\fq_i)\stackrel{\sim}{\to}\bar\kappa(x)$ for every
$x\in\fF(\fq_i)$, hence the image of \eqref{eq_integral-form}
lands in the {\em seminormalization\/} of $B^{\circ h}_{\fq_i}$,
{\em i.e.} the subring
$$
B^{h\nu}_{\fq_i}:=\kappa(x_1)^{\wedge h+}\times_{\bar\kappa(\fq_i)}\cdots
\times_{\bar\kappa(\fq_i)}\kappa(x_r)^{\wedge h+}
$$
(the fibre product over $\bar\kappa(\fq_i)$ of the rings
$\kappa(x_i)^{h+}$, where $\{x_1,\dots,x_r\}=\fF(\fq_i)$). Let us
set:
$$
\alpha(\fq_i):=
|F_0(B^{h\nu}_{\fq_i}/B^{\circ h}_{\fq_i})|^{\wedge h\natural}_{\eta(1)}
\qquad\text{for every $i=1,\dots,k$.}
$$

\begin{lemma}\label{lem_seminormalize}
With the notation of \eqref{subsec_analyze}, we have:
$$
2\cdot|B^\circ_\eta|^\natural_{\sp,\eta(1)}=
\deg(g)+\sum_{i=1}^k(2\alpha(\fq_i)+\sharp\fF(\fq_i)-2).
$$
\end{lemma}
\begin{proof} (Here $\sharp\fF(\fq_i)$ denotes the cardinality
of the finite set $\fF(\fq_i)$.) First of all, we remark that
$\bar\kappa(\fq_i)\simeq\bar\kappa(\eta(1))\simeq
\kappa(1)^+/\xi\kappa(1)^+$, where $\xi\in A(1)$
is an element such that $|\xi|_{\eta(1)}=1-\eps$, hence
$|F_0(\bar\kappa(\fq_i))|^\natural_{\eta(1)}=1$; it follows easily
that
$$
|F_0(\cB(1)^{\wedge h+}_{\fq_i}/B^{h\nu}_{\fq_i})|^{\wedge h\natural}_{\eta(1)}
=\sharp\fF(\fq_i)-1\qquad\text{for every $i=1,\dots,k$}.
$$
Hence
$$
\begin{array}{r@{\:=\:}l}
2\cdot|B^\circ_\eta|^\natural_{\sp,\eta(1)} &
2\cdot|B^{\circ h}_\eta|^{\wedge h\natural}_{\sp,\eta(1)}=
2\cdot(|\cB(1)^{\wedge h+}|^{\wedge h\natural}_{\sp,\eta(1)}+
|F_0(\cB(1)^{\wedge h+}/B^{\circ h}_\eta)|^{\wedge h\natural}_{\eta(1)}) \\
& 2\cdot|\cB(1)^+|^\natural_{\sp,\eta(1)}+
2\sum_{i=1}^k(
|F_0(\cB(1)^{\wedge h+}_{\fq_i}/B^{h\nu}_{\fq_i})|^{\wedge h\natural}_{\eta(1)}
+|F_0(B^{h\nu}_{\fq_i}/B^{\circ h}_{\fq_i})|^{\wedge h\natural}_{\eta(1)}) \\
& 2\cdot|\cB(1)^+|^\natural_{\sp,\eta(1)}+
2\sum_{i=1}^k(\sharp\fF(\fq_i)-1+\alpha(\fq_i)) \\
& \deg(g)-\sharp g^{-1}(\eta(1))+
2\sum_{i=1}^k(\sharp\fF(\fq_i)-1+\alpha(\fq_i))
\end{array}
$$
where the last equality holds by lemma \ref{lem_local-sit}(iii).
Since clearly
$$
\sum_{i=1}^k\sharp\fF(\fq_i)=\sharp g^{-1}(\eta(1))
$$
the assertion follows.
\end{proof}

\begin{theorem}\label{th_big-deal} With the notation of
\eqref{subsec_discr-at-r} :
\begin{enumerate}
\item
$\delta_f$ extends to a continuous, piecewise linear function
$\delta:[\log 1/b,\log 1/a]\to\R_{\geq 0}$ with integer slopes.
\item
If moreover $f$ is \'etale, then $\delta$ is convex.
\end{enumerate}
\end{theorem}
\begin{proof} Let $(F,|\cdot|_F)$ be an algebraically closed valued
field extension of $K$ with $|F|_F=\R_{\geq 0}$, and denote by
$f_F:X\times_{\Spa\,K}\Spa\,F\to\D(a,b)\times_{\Spa\,K}\Spa\,F$ the
morphism deduced by base change of $f$; using proposition
\ref{prop_free}(iii) one sees that $\delta_{f_F}$ agrees with $\delta_f$
wherever the latter is defined. Hence we can and do assume from start
that $|K|=\R_{\geq 0}$. Now, for the proof of (i) it suffices
to show that $\delta_f$ is piecewise linear in the neghborhood of
every real number $\rho:=\log 1/r\in[\log 1/b,\log 1/a]$. Using a morphism
$g$ as in example \ref{ex_asymm}, one reduces to consider the case
where $r>a$, and study the function $\delta_f$ in a small interval
$[\rho,\rho+x]$. In such situation, the more precise proposition
\ref{prop_precise} shows that the assertion holds.

Suppose next that $f$ is \'etale. In order to show (ii), we need
to study the function $\delta$ in any small neighborhood of the form
$[\log 1/r-x,\log 1/r+x]\subset[\log 1/b,\log 1/a]$. Assertion (ii)
then means that the function
$\rho\mapsto\delta(\log 1/r-\rho)+\delta(\log 1/r+\rho)$
has positive derivative in a neighborhood of $0$. We can assume that
$r=1$ and $x=-\log|\pi|$ for some $\pi\in\fm$, so we reduce to consider
an \'etale morphism $f:X\to\D(a,a^{-1})$ (for $a:=|\pi|$). In view
of \eqref{eq_summing-up} we can further reduce to studying the
morphism $h:=\Spa(\psi_K)\circ f:X\to\D(1)$, defined as in example
\ref{ex_summing-up}, and then we have to show that the left slope
of $\delta_h$ is negative in a small neighborhood $(x,0]$.
Say that $X=\Spa\,B$; in the notation of \eqref{eq_local-to-glob}
we have
\set\begin{equation}\label{eq_eval-discr}
|\fd^\circ_h|_{\eta(1)}^\natural=|(\fd_\psi)^d|_{\eta(1)}^\natural=2d
=\deg(g)
\end{equation}
by example \ref{ex_discriminant} and proposition
\ref{prop_tower-discr}. Finally, in view of \eqref{eq_local-to-glob},
\eqref{eq_eval-discr}, proposition \ref{prop_precise} and lemma
\ref{lem_seminormalize}, the sought assertion is implied by the
following:
\begin{claim}\label{cl-at-least-two} Resume the notation of
\eqref{subsec_analyze}. Then $\sharp\fF(\fq_i)\geq 2$ for every
$i=1,\dots,k$.
\end{claim}
\begin{pfclaim}[] Recall that $\fq_1,\dots,\fq_k$ are by definition
the prime ideals of $B^\circ_\eta$ lying over the maximal ideal
of $\kappa(1)^{h+}$, or what is the same, the prime
ideals of $B^\circ$ lying over the maximal ideal
$\fp:=\fm A(1)^\circ+\xi A(1)^\circ$ of $A(1)^\circ$. Now, we
have already observed (example \ref{ex_summing-up}) that
$A(a,a^{-1})^\circ\simeq K^+\langle S,T\rangle/(ST-\pi^2)$,
and $\psi$ is the map
$K^+\langle\xi\rangle\to K^+\langle S,T\rangle/(ST-\pi^2)$ such
that $\xi\mapsto S+T$. Hence
$A(a,a^{-1})^\circ/\fp A(a,a^{-1})^\circ\simeq
K^\sim[S,T]/(ST,S+T)\simeq K^\sim[S]/(S^2)$. Thus, there is exactly
one prime ideal $\fP\subset A(a,a^{-1})^\circ$ lying over
$\fp$ and necessarily $\fq_i\cap A(a,a^{-1})^\circ=\fP$
for every $i=1,\dots,k$. On the other hand, the fibre
$\psi^{-1}(\eta(1))\subset\D(a,a^{-1})$ consists of the two
valuations $\eta'(a),\eta(1/a)$, and clearly both of them dominate
the local ring $A(a,a^{-1})^\circ_\fP$. It is now a standard
fact that, for each prime ideal $\fq_i$, there are valuations
$\eta_1$, $\eta_2$ on $B$ which extend respectively $\eta'(a)$
and $\eta(1/a)$, and which dominate the local ring $B^\circ_{\fq_i}$.
By lemma \ref{lem_much-ado} we have $\eta_1,\eta_2\in\Spa\,B$,
whence the claim.
\end{pfclaim}
\end{proof}

\begin{remark} The continuity and piecewise linearity of
the function $\delta_f$ are also proved in the preprint
\cite{Schm} of T.Schmechta. He also obtains some interesting
results in the case where the base field has positive characteristic.
His methods are refinements of those of L\"utkebohmert \cite{Lut}.
(However, as far as I understand, he does not prove the convexity
of the function $\delta_f$.)
\end{remark}

\subsection{The $p$-adic Riemann existence theorem}\label{sec_Riemann}
In this section we show how to use theorem \ref{th_big-deal}
to solve the so-called $p$-adic Riemann existence problem
in case $K$ is a field of characteristic zero. We choose
an argument that maximizes the use of valuation theory;
see remark \ref{rem_alternative} for some indications of
an alternative, slightly different proof.

\sset\subsubsection{}
Recall that a finite \'etale covering $f:X\to\D(a,b)$ is
said to be {\em of Kummer type}, if there exists an integer
$n>0$ and an isomorphism $g:\D(a^{1/n},b^{1/n})\isom X$ such
that $f\circ g=\Spa\,\phi$, where $\phi:A(a,b)\to A(a^{1/n},b^{1/n})$
is the map of $K$-affinoid algebras given by the rule $\xi\mapsto\xi^n$
(notation of \eqref{subsec_setup}).

\begin{lemma}\label{lem_Kummer}
Let $f:X\to\D(a,b)$ and $g:Y\to X$ be two finite \'etale
coverings. Then $f$ and $g$ are of Kummer type, if and only
if the same holds for $f\circ g$.
\end{lemma}
\begin{proof} Left to the reader.
\end{proof}

\begin{theorem}\label{th_Riemann} Suppose that $K$ is
algebraically closed of characteristic zero, and let
$f:X\to\D(a,b)$ be a finite \'etale morphism of degree $d$.
There is a constant $c:=c(d)\in(0,1]$ such that the restriction
$f^{-1}(\D(c^{-1}a,cb))\to\D(c^{-1}a,cb)$ splits as the
disjoint union of finitely many finite coverings of Kummer type.
\end{theorem}
\begin{proof} First of all, let $f':Y\to\D(a,b)$ be the
smallest Galois \'etale covering that dominates $f$ ({\em i.e.}
such that $f'$ factors through $f$); it is well-known that
the degree of $f'$ is bounded by $d!$. Suppose now that
the theorem is known for $f'$; then we may find $c\in(0,1]$
such that the restriction of $f'$ to the preimage of
$\D(c^{-1}a,cb)$ is of Kummer type. Using lemma \ref{lem_Kummer}
we deduce that the same holds for the restriction of $f$ to
$f^{-1}(\D(c^{-1}a,cb))$. Hence, we may replace $f$ by
$f'$ and assume from start that $f$ is a Galois covering,
say of finite group $G$.

Next, we consider the function
$\delta:[\log 1/b,\log 1/a]\cap\log\Gamma_K\to\R_{\geq 0}$
corresponding to the covering $f$. To start out, lemma
\ref{lem_bound-discr} implies that $\delta$ admits an upper
bound that depends only on $d$; since $\delta$ is convex, piecewise
linear and non-negative and since its slopes are integers (theorem
\ref{th_big-deal}), it follows easily that we may find a
constant $c\in(0,1]$, depending only on the degre $d$, such
that $\delta$ is linear (indeed constant) on the interval
$[\log 1/(bc),\log c/a]\cap\log\Gamma_K$. We may therefore
assume from start that $\delta$ is linear. 
Also, we may assume that $a<1$ and $b=a^{-1}$, in which case
we let $g:=\Spa\,\psi_K:\D(a,a^{-1})\to\D(1)$, where $\psi$
is defined as in example \ref{ex_summing-up}, and
$h:=g\circ f:X\to\D(1)$.
Let $\fp:=\fm A(1)^\circ+\xi A(1)^\circ\subset A(1)^\circ$;
in the proof of claim \ref{cl-at-least-two} we have established
that there exists a unique prime ideal $\fP\subset A(a,a^{-1})^\circ$
lying over $\fp$, and both rings $\kappa(\eta'(a))^+$ and
$\kappa(b)^+$ dominate the localization $A(a,b)^\circ_\fP$;
denote also $\fq_1,\dots,\fq_k\subset B^\circ$ the finitely
many prime ideals lying over $\fp$.

The natural map $A(1)^\circ\to A(1)^{\circ\wedge}\simeq K^+[[\xi]]$
from $A(1)^\circ$ to its $\xi$-adic completion, factors through
the henselization $A(1)^{\circ h}_\fp$ of $A(1)^\circ$ along
its ideal $\fp$; hence
$B^{\circ\wedge}:=B^\circ\otimes_{A(1)^\circ}A(1)^{\circ\wedge}$
decomposes as a direct product of algebras :
$B^{\circ\wedge}\simeq C_1\times\cdots\times C_k$.
Moreover, for every $t\in\fm\setminus\{0\}$, the map
$A(1)^\circ\to A(|t|)^\circ$ induced by the open imbedding
$\D(|t|)\to\D(1)$ factors through $A(1)^{\circ\wedge}$, and
induces an isomorphism $A(1)^\circ/\fp\isom A(|t|)^\circ/\fp_t$,
where $\fp_t:=\fm A(|t|)^\circ+(\xi/t) A(|t|)^\circ$.
It follows that the prime ideals of
$$
B^\circ\otimes_{A(1)^\circ}A(|t|)^\circ\simeq
(C_1\otimes_{A(1)^{\circ\wedge}}A(|t|)^\circ)\times\cdots\times
(C_k\otimes_{A(1)^{\circ\wedge}}A(|t|)^\circ)
$$
lying over $\fp_t$ are in natural bijection with the prime
ideals $\fq_1,\dots,\fq_k$, and therefore every factor
$C_i(t):=C_i\otimes_{A(1)^\circ}A(|t|)^\circ$ contains
exactly one of these prime ideals. Notice that
$g^{-1}(\D(|t|))=\D(a/|t|,|t|/a)$, whence
natural isomorphisms :
$$
f^{-1}(\D(a/|t|,|t|/a))\simeq
\Spa(C_1(t)\otimes_{K^+}K)\amalg\cdots\amalg
\Spa(C_k(t)\otimes_{K^+}K)
$$
so that the restriction of $\delta$ to
$[\log a/|t|,\log |t|/a]\cap\log\Gamma_K$ decomposes
as a sum $\delta=\delta_1+\cdots+\delta_k$, where $\delta_i$
is the discriminant function of the restriction
$$
f_i:\Spa(C_i(t)\otimes_{K^+}K)\to\D(a/|t|,|t|/a)
$$
for every $i\leq k$. Since every such $\delta_i$ is still
convex, and their sum is linear, it follows that $\delta_i$ is
linear for every $i\leq k$. We remark as well, that each
$f_i$ is still a Galois covering, whose Galois group is
the subgroup of $G$ that stabilizes $\fq_i$, for the natural
action of $G$ on the set $\{\fq_1,\dots,\fq_k\}$.
Clearly it suffices to prove the theorem separately for every
\'etale covering $f_i$, hence we may replace from start $f$
by $f_i$, and assume additionally that $k=1$, in which case
we shall write $\fq$ instead of $\fq_1$. Define $\alpha(\fq)$,
$\fF(\fq)$ as in \eqref{subsec_analyze}. By inspection of the
proof of theorem \ref{th_big-deal} we deduce that $\delta$ is
linear precisely when :
\set\begin{equation}\label{eq_crossing}
\alpha(\fq)=0 \quad\text{and}\quad \sharp\fF(\fq)=2.
\end{equation}
Especially, the preimages $f^{-1}(\eta'(a))$ and $f^{-1}(\eta(b))$
both consist of precisely one point; let $x\in X$ be the
only point lying over $\eta(b)$. We deduce a Galois extension
of valued fields :
$$
\kappa(b)\to\kappa(x)
$$
whose Galois group is isomorphic to $G$. Since the residue
field of $K$ is algebraically closed, the residue field extension
$\bar\kappa(\eta(b))\to\bar\kappa(x)$ is trivial, and therefore
$G$ is a solvable group. Thus, we may factor $f$ as a composition
of finitely many \'etale coverings :
$$
X_n:=X\stackrel{g_n}{\longrightarrow}X_{n-1}
\stackrel{g_{n-1}}{\longrightarrow}\cdots
\stackrel{g_1}{\longrightarrow} X_0:=\D(a,b)
$$
such that the degree of $g_i$ is a prime number for every
$i\leq n$. Using lemma \ref{lem_Kummer} and an easy induction,
we may then further reduce to the case where $G$ is a cyclic
group of prime order $d$.

Such coverings are classified by the \'etale cohomology group
$H:=H^1(\D(a,b)_\et,\Z/d\Z)$ (where $\D(a,b)_\et$ denotes the
\'etale site of $\D(a,b)$, as defined in \cite{Hu2}). The
latter can be computed by the Kummer exact sequence (on the
\'etale site of $\D(a,b)$) :
$$
0\to\mu_d\to\cO^\times\stackrel{(-)^d}{\longrightarrow}\cO^\times\to 0
$$
(recall that $K$ has characteristic zero) and since the Picard
group of $\D(a,b)$ is trivial (\cite[Th.2.2.9(3)]{Fr-VdP}), one
obtains a natural isomorphism :
$$
H\isom A(a,b)^\times/(A(a,b)^\times)^d
$$
where $A(a,b)^\times$ denotes the invertible elements of
$A(a,b)$. Under this isomorphism, the Kummer coverings
of degree $d$ correspond to the equivalence classes of the
sections $\xi^j$, for $j=0,\dots,d-1$ (notation of
\eqref{subsec_setup}). Thus, we come down to verifying the following :

\begin{claim} There exists a constant $c:=c(d)\in(0,1]$
such that, for every $u\in A(a,b)^\times$, the restriction
$u':=u_{|\D(a/c,bc)}$ can be written in the form :
$u'=v^d\cdot\xi^j$ for some $v\in A(a/c,cb)^\times$
and $0\leq j\leq d-1$.
\end{claim}
\begin{pfclaim}[] Let $\alpha,\beta\in K^\times$ such that
$|\alpha|=a$, $|\beta|=b$; it is well known that every invertible
element $u$ of $A(a,b)$ can be written in the form
$u=\gamma\cdot\xi^n\cdot(1+h)$ where $\gamma\in K^\times$,
$n\in\Z$ and $h\in A(a,b)^\circ$ of the form
\set\begin{equation}\label{eq_expl-descr}
h(\xi)=\sum_{k\in\Z\setminus\{0\}}h_k\xi^k\in
K^+\langle\xi/\alpha,\beta/\xi\rangle
\end{equation}
with $|h|_{\sup}<1$. Hence we are reduced to showing that $1+h$
admits a $d$-th root, after restriction to a smaller annulus
$\D(a/c,b/c)$. This is clear in case $d\neq p$, in which case
we may even choose $c=1$. Finally, suppose that $d=p$; it is
well known that $1+h$ admits a $p$-th root as soon as
\set\begin{equation}\label{eq_estimate}
|h|_{\sup}<|p|^{1/(p-1)}.
\end{equation}
Using the explicit description \eqref{eq_expl-descr} we
may easily determine $c\in(0,1]$ such that the estimate
\eqref{eq_estimate} holds for the restriction $h_{|\D(a/c,cb)}$. 
\end{pfclaim}
\end{proof}

\begin{remark}\label{rem_alternative} Keep the notation of
the proof of theorem \ref{th_Riemann}. Alternatively, one may
deduce from \eqref{eq_crossing} that $\fq$ is an ordinary
double point of the analytic reduction of $X$, in which case
\cite[Prop.2.3]{Bo-Lut} shows that the corresponding formal
fibre is an open annulus, and then theorem \ref{th_Riemann}
follows without too much trouble.
\end{remark}

\section{Study of the conductors}\label{chap_conductors}

\subsection{Algebraization}\label{sec_algebrize}
Let $V$ be a henselian local ring, $s$ the closed point of
$\Spec\,V$ and $\kappa(s)$ its residue field. For every affine
$V$-scheme $X$ we let $X_s:=X\times_{\Spec\,V}\Spec\,\kappa(s)$.
More generally, let $X:=\Spec\,R$ be any affine scheme, and $Z\subset X$
a closed subscheme, say $Z=V(I)$ for an ideal $I\subset R$; we denote
by $R^h_I$ the henselization of $R$ along the ideal $I$ and by
$R^\wedge_I$ the $I$-adic completion of $R$.
The {\em henselization\/} of $X$ along $Z$ is the affine scheme
$X^h_{/Z}:=\Spec\,R^h_I$.

\begin{lemma}\label{lem_various-compl}
Let $A$ be a noetherian ring, $I,J\subset A$ two ideals.
\begin{enumerate}
\item The natural commutative diagram
$$
\xymatrix{ A^\wedge_{I\cap J} \ar[r] \ar[d] & A^\wedge_I \ar[d] \\
           A^\wedge_J \ar[r] & A^\wedge_{I+J}
}
$$
is cartesian.
\item
Moreover, for every $n\in\N$, set $A_n:=A/J^n$. Then there is a natural
isomorphism of $A$-algebras:
$$
A^\wedge_{I+J}\stackrel{\sim}{\to}\liminv{n\in\N}\,A^\wedge_{n,I}.
$$
\end{enumerate}
\end{lemma}
\begin{proof} (i): Using \cite[Ch.0, Lemme 19.3.10.2]{EGAIV} we see
that the natural commutative diagram
$$
\xymatrix{ A/(I^n\cap J^n) \ar[r] \ar[d] & A/I^n \ar[d] \\
           A/J^n \ar[r] & A/(I^n+J^n)
}
$$
is cartesian for every $n>0$. Set $A':=\liminv{n\in\N}\,A/(I^n+J^n)$
and $A'':=\liminv{n\in\N}\,A/(I^n\cap J^n)$. We deduce easily a cartesian
commutative diagram:
$$
\xymatrix{ A'' \ar[r] \ar[d] & A^\wedge_I \ar[d] \\
           A^\wedge_J \ar[r] & A'
}
$$
and it remains only to show that the natural maps
$A'\to A^\wedge_{I+J}$ and $A^\wedge_{I\cap J}\to A''$ are
isomorphisms. For the former, it suffices to remark that
\set\begin{equation}\label{eq_trivial-topol}
(I+J)^{2n-1}\subset I^n+J^n\subset(I+J)^n
\end{equation}
for every $n>0$. For the latter, one uses the Artin-Rees lemma
\cite[Th.8.5]{Mat} to show that for every $n\in\N$ there exists
$m\in\N$ such that
$$
I^m\cap J^n\subset I^nJ^n\subset(I\cap J)^n
$$
from which (i) follows easily. (ii) is an easy consequence of
\eqref{eq_trivial-topol}; we leave the details to the reader.
\end{proof}

\begin{theorem}\label{th_algebrize} In the situation of \eqref{sec_algebrize},
suppose that $\kappa(s)$ is a perfect field. Let $X$ be an affine
finitely presented $V$-scheme of pure relative dimension one,
$\Sigma\subset X_s$ a finite subset such that $X_s\setminus\Sigma$
is smooth over\/ $\Spec\,\kappa(s)$. Then there exists a projective
$V$-scheme $Y$ of pure relative dimension one and an open affine
subset $U\subset Y$ with an isomorphism of $V$-schemes
\set\begin{equation}\label{eq_hens-iso}
X^h_{/X_s}\simeq U^h_{/U_s}.
\end{equation}
Moreover, $U_s$ is dense in $Y_s$ and $Y$ is smooth over\/
$\Spec\,V$ at all the points of\/ $Y_s\setminus U_s$. 
\end{theorem}
\begin{proof} We begin with the following:
\begin{claim}\label{cl_compactify} There exists a projective purely
one-dimensional $\kappa(s)$-scheme $Y_0$ and a dense open imbedding
of $\kappa(s)$-schemes
\set\begin{equation}\label{eq_imbedding}
X_s\subset Y_0
\end{equation}
such that $Y_0$ is smooth over $\Spec\,\kappa(s)$ at all the points
of $Y_0\setminus X_s$.
\end{claim}
\begin{pfclaim} This is standard : one picks a projective $\kappa(s)$-scheme
$\bar X_s$ containing $X_s$ as a dense open subscheme, and let $X'_s$
be the normalization of $\bar X_s\setminus\Sigma$.  By
\cite[Ch.V, \S3.2, Th.2]{BouAC} we know that $X'_s$ is of finite
type over $\kappa(s)$, and since $X'_s$ has dimension one, we know as
well that all its local rings are regular, hence they are formally
smooth over $\kappa(s)$, in view of \cite[\S28, Lemma 1]{Mat}
and \cite[Ch.IV, Prop.17.5.3]{EGA4}. One can then glue $X_s$ and
$X'_s$ along their common open subscheme $X_s\setminus\Sigma$;
the resulting scheme $Y_0$ will do.
\end{pfclaim}

We can write as usual $V$ as the colimit of a filtered family
$(V_\lambda~|~\lambda\in\Lambda)$ of noetherian local subrings
of $V$, essentially of finite type over an excellent discrete
valuation ring, and such that the inclusion maps
$j_\lambda:V_\lambda\to V$ are local ring homomorphisms.
Then $j_\lambda$ extends to a map $V_\lambda^h\to V$ from
the henselization of $V_\lambda$, and $V$ is still the colimit
of the filtered family $(V^h_\lambda~|~\lambda\in\Lambda)$.
For some $\lambda\in\Lambda$ we can find an affine
finitely presented $V^h_\lambda$-scheme $X_\lambda$ and an isomorphism
of $V$-schemes:
$$
\Spec\,V\times_{\Spec\,V^h_\lambda}X_\lambda\stackrel{\sim}{\to} X.
$$
For every $\lambda\in\Lambda$, let $k_\lambda$ be the residue field
of $V^h_\lambda$; we can even choose $\lambda$ in such a way that
the scheme $Y_0$  provided by claim \ref{cl_compactify} descends
to a projective $k_\lambda$-scheme $Y_{0,\lambda}$, so that
$Y_0\simeq\Spec\,\kappa(s)\times_{\Spec\,k_\lambda}Y_{0,\lambda}$.
Let $\Sigma_\lambda\subset Y_{0,\lambda}$ be the image of $\Sigma$;
we can furthermore assume that $Y_{0,\lambda}$ is smooth over
$\Spec\,k_\lambda$ outside $\Sigma_\lambda$,
and that there exists an open imbedding of $k_\lambda$-schemes
\set\begin{equation}\label{eq_another}
\Spec\,k_\lambda\times_{\Spec\,V^h_\lambda}X_\lambda\subset Y_{0,\lambda}
\end{equation}
inducing \eqref{eq_imbedding}, after base change to $\Spec\,\kappa(s)$
(see \cite[Ch.IV, Prop.17.7.8]{EGA4}). At the cost of trading
the residue field $\kappa(s)$ with a non-perfect field, we can
then replace the given ring $V$ by one such $V^h_\lambda$, the scheme
$X$ by $X_\lambda$ and $\Sigma$ by $\Sigma_\lambda$; hence we can
assume that $V$ is the henselization of a ring of essentially
finite type over an excellent discrete valuation ring, and additionally,
that there exists a projective $\kappa(s)$-scheme $Y_0$ as in claim
\ref{cl_compactify}.

Let $\fn\subset V$ be the maximal ideal; denote by $V\Alg$
the category of $V$-algebras, by $\Set$ the category of sets.
We define a functor $\cF:V\Alg\to\Set$ as follows. For a $V$-algebra
$A$, $\cF(A)$ is the set of equivalence classes of data of the
form $(Z_A,Y_A,f_A,g_A,h_A)$ where:
\begin{itemize}
\item
$Z_A$ and $Y_A$ are finitely presented $A$-schemes, and $Y_A$
is projective over $\Spec\,A$.
\item
$f_A:Z_A\to X_A:=\Spec\,A\times_{\Spec\,V}X$ and $g_A:Z_A\to Y_A$
are \'etale morphisms of $A$-schemes.
\item
$h_A:Y_{A,s}:=\Spec\,A/\fn A\times_{\Spec\,A}Y_A\to
\Spec\,A/\fn A\times_{\Spec\,\kappa(s)}Y_0$ is an isomorphism.
\item
The restriction
$f_{A,s}:Z_{A,s}:=\Spec\,A/\fn A\times_{\Spec\,A}Z_A\to
X_{A,s}:=\Spec\,A/\fn A\times_{\Spec\,V}X$ is an isomorphism.
\item
The restriction $g_{A,s}:Z_{A,s}\to Y_{A,s}$ is an open imbedding
and the morphism
$$
\Spec\,A/\fn A\times_{\Spec\,\kappa(s)}\eqref{eq_imbedding}:
X_{A,s}\to Y_{A,s}
$$
agrees with $h_A\circ g_{A,s}\circ f_{A,s}^{-1}$.
\item
$Y_A$ is smooth over $\Spec\,A$ at all the points of
$Y_{A,s}\setminus g_{A,s}(Z_{A,s})$.
\end{itemize}
Two data $(Z_A,Y_A,f_A,g_A,h_A)$ and $(Z'_A,Y'_A,f'_A,g'_A,h'_A)$
are said to be equivalent if there are isomorphisms of $A$-schemes
$Z_A\stackrel{\sim}{\to}Z'_A$, $Y_A\stackrel{\sim}{\to}Y'_A$
such that the obvious diagrams commute. A map $A\to A'$ of
$V$-algebras induces an obvious base change map $\cF(A)\to\cF(A')$.

\begin{claim}\label{cl_small-fctr}
$\cF(V)\neq\emptyset$ if and only if there exists a pair $(U,Y)$
as in the theorem, with an isomorphism of $V$-schemes
$Y_s\stackrel{\sim}{\to}Y_0$.
\end{claim}
\begin{pfclaim}
Indeed, suppose we have found $(Z,Y,f,g,h)\in\cF(V)$. Then we
can set $U:=g_V(Z)\subset Y$; the morphism $f$ and $g$ induce
isomorphisms $f^h:Z^h_{/Z_s}\stackrel{\sim}{\to}X^h_{/X_s}$
and $g^h:Z^h_{/Z_s}\stackrel{\sim}{\to}U^h_{/X_s}$, and
therefore the pair $(U,Y)$ fulfills the conditions of the
theorem; furthermore $h$ yields an isomorphism
$Y_s\stackrel{\sim}{\to} Y_0$. Conversely, suppose that
a pair $(U,Y)$ has been found that fulfills the conditions
of the theorem; especially, \eqref{eq_hens-iso} induces a
morphism of $V$-schemes $g:X^h_{/X_s}\to U$. By
\cite[Ch.XI, \S2, Th.2]{Ray} $X^h_{/X_s}$ is the projective
limit of a cofiltered family of morphisms of $V$-schemes
$(f_\alpha:X_\alpha\to X~|~\alpha\in I)$ such that every $X_\alpha$
is affine and finitely presented over $\Spec\,V$, the induced
morphisms $X_{\alpha,s}\to X_s$ are isomorphisms, and for
every $x\in X_\alpha$ the induced map
$\cO_{X,f_\alpha(x)}\to\cO_{X_\alpha,x}$ is local ind-\'etale.
By \cite[Ch.IV, Th.11.1.1]{EGAIV-3} we deduce that there is an
open (quasi-compact) subset $Z_\alpha$ of $X_\alpha$ containing
$X_{\alpha,s}$ such that the restriction of $f_\alpha$ to
$Z_\alpha$ is flat; then, up to shrinking $Z_\alpha$, we can
achieve that the morphism $Z_\alpha\to X$ is an \'etale neighborhood
of $X_s$, and $X$ is still the projective limit of the $V$-schemes
$Z_\alpha$.
We can then find $\alpha\in I$ such that $g$ factors through a
morphism $g_\alpha:Z_\alpha\to U$; after further shrinking
$Z_\alpha$, the morphism $g_\alpha$ becomes \'etale. Hence,
the datum $(Z_\alpha,Y,f_\alpha,g_\alpha)$ represents an
element of $\cF(V)$.
\end{pfclaim}

Next, in view of \cite[Ch.IV, Th.8.10.5]{EGAIV-3} and
\cite[Ch.IV, Prop.17.7.8(ii)]{EGA4} one sees easily that
$\cF$ is a functor of finite presentation. Let $V^\wedge$
be the $\fn$-adic completion of $V$; according to Artin's
approximation theorem \cite[Th.I.12]{Art}, every element of
$\cF(V^\wedge)$ can be approximated arbitrarily closely
in the $\fn$-adic topology by elements of $\cF(V)$,
especially:
\set\begin{equation}\label{eq_Artin-app}
\cF(V)\neq\emptyset\Leftrightarrow\cF(V^\wedge)\neq\emptyset.
\end{equation}
So finally, in view of \eqref{eq_Artin-app} and claim
\ref{cl_small-fctr}, we can replace $V$ by $V^\wedge$ and
suppose from start that $V$ is a complete noetherian local ring.
Let $V_n:=V/\fn^{n+1}$ and $X_n:=\Spec\,V_n\times_{\Spec\,V}X$
for every $n\in\N$. We endow $V$ with its $\fn$-adic topology;
the family $(X_n~|~n\in\N)$ defines a unique affine formal
$\Spf\,V$-scheme $\fX$.

\begin{claim}
There exists a proper $\Spf\,V$-scheme $\fY$ with an open
imbedding $\fX\subset\fY$, such that the reduced fibre of $\fY$
is $V_0$-isomorphic to $Y_0$, and such that $\fY$ is formally
smooth over $\Spf\,V$ at all points of $Y_0\setminus X_0$.
\end{claim}
\begin{pfclaim} In view of \cite[Ch.III, \S3.4.1]{EGA-III} and
\cite[Ch.I, Prop.10.13.1]{EGAI}, it suffices to lift the
$V_0$-scheme $Y_0$ to a compatible family $(Y_n~|~n\in\N)$
of schemes such that
\begin{enumerate}
\item
for every $n\in\N$ there are isomorphisms of $V_n$-schemes:
$Y_n\stackrel{\sim}{\to}\Spec\,V_n\times_{\Spec\,V_{n+1}}Y_{n+1}$;
\item
moreover, \eqref{eq_imbedding} lifts to a system of open imbeddings
of $V_n$-schemes $X_n\subset Y_n$ compatible with the isomorphisms
(ii);
\item
$Y_n$ is smooth over $\Spec\,V_n$ at all points of $Y_n\setminus X_n$.
\end{enumerate}
To this aim, we may assume that $\Sigma\neq\emptyset$, in
which case $Y'_0:=Y_0\setminus\Sigma$ is smooth and affine
over $\Spec\,V_0$. One can then lift $Y'_0$ to a
compatible system of schemes $(Y'_n~|~n\in\N)$ satisfying a condition
as the foregoing (i), and such that furthermore, $Y'_n$ is smooth
over $\Spec\,V_n$ for every $n\in\N$. Again some basic deformation
theory shows that the imbedding $X_0\setminus\Sigma\subset Y'_0$
lifts to a compatible system of open imbeddings
$X_n\setminus\Sigma\subset Y'_n$ for every $n\in\N$. Therefore
one can glue $X_n$ and $Y'_n$ along their common open subscheme
$X_n\setminus\Sigma$; the resulting schemes $Y_n$ will do.
\end{pfclaim}

Next we are going to construct an invertible $\cO_\fY$-module on
$\fY$. To this aim we proceed as follows.
Let $\{y_1,\dots,y_n\}:=Y_0\setminus X_0$. By construction $\fY$ is
formally smooth over $\Spf\,V$ at the point $y_i$, for every
$i=1,\dots,n$. For every $i\leq n$, the maximal ideal of $\cO_{Y_0,y_i}$
is principal, say generated by the regular element $\bar t_i\in\cO_{Y_0,y_i}$.
The natural ring homomorphism $\cO_{\fY,y_i}\to\cO_{Y_0,y_i}$ is
surjective, hence we can lift $\bar t_i$ to an element $t_i\in\cO_{\fY,y_i}$.
\begin{claim}\label{cl_regular} $t_i$ is a regular element in $\cO_{\fY,y_i}$
for every $i=1,\dots,n$.
\end{claim}
\begin{pfclaim} By \cite[Ch.0, Th.19.7.1]{EGAIV} the ring $\cO_{\fY,y_i}$
is flat over $V$; then the claim follows from \cite[Ch.0, Prop.15.1.16]{EGAIV}.
\end{pfclaim}

In view of claim \ref{cl_regular} we can find an open affine
subscheme $\fU_i\subset\fY$ such that $y_i\in\fU_i$ and $t_i$ 
extends to a regular element of $\Gamma(\fU_i,\cO_\fY)$. Finally,
set $\fV_i:=\fY\setminus\{y_i\}$ and $\fW_i:=\fU_i\cap\fV_i$; we
define the invertible $\cO_\fY$-module $\cL_i^\wedge$ by gluing the sheaves
$\cO_{\fU_i}$ (defined on $\fU_i$) and $\cO_{\fV_i}$ (defined on $\fV_i$);
the gluing map is the isomorphism 
$$
\cO_{\fV_i|\fW_i}\stackrel{\sim}{\to}\cO_{\fU_i|\fW_i}
\qquad f\mapsto t_i\cdot f.
$$
So, the global sections of $\cL_i^\wedge$ are identified naturally
with the pairs $(f,g)$ where $f\in\Gamma(\fU_i,\cO_\fY)$,
$g\in\Gamma(\fV_i,\cO_\fY)$ and $f_{|\fW_i}=t_i\cdot g_{|\fW_i}$.
Clearly, $\cL_i^\wedge/\fn\cL_i^\wedge$ is the invertible sheaf
$\cO_{Y_0}(y_i)$ on $Y_0$; set
$\cL^\wedge:=\cL_1^\wedge\otimes_{\cO_\fY}\cdots\otimes_{\cO_\fY}\cL_n^\wedge$.
notice that, since $X$ is affine, every irreducible component
of $X_0$ meets $\{y_1,\dots,y_n\}$; it then follows from
\cite[\S7.5, Prop.5]{Liu} that $\cL^\wedge/\fn\cL^\wedge$
is ample on $Y_0$, and then \cite[Ch.III, Th.5.4.5]{EGA-III} shows that
$\fY$ is algebraizable to a projective scheme $Y$ over $\Spec\,V$,
and $\cL^\wedge$ is the formal completion of an ample invertible
$\cO_Y$-module that we shall denote by $\cL$. Especially, 
$\cL/\fn\cL$ is the ample sheaf $\cO_{Y_0}(y_1+\cdots+y_n)$ on $Y_0$.

By inspecting the construction, we see that $\cL\simeq\cO_Y(D)$,
where $D$ is an ample divisor on $Y$ whose support
$\mathrm{Supp}(D)\subset Y$ intersects $Y_0$ in precisely the
closed subset $Y_0\setminus X_0$. Therefore the open subset
$U:=Y\setminus\mathrm{Supp}(D)$ is affine and clearly
$U\cap Y_0=X_0$. Let $\fU\subset\fY$ be the formal completion
of $U$ along its closed subscheme $X_0$; it follows that the
imbedding $\fX\subset\fY$ induces an isomorphism of affine formal
$\Spf\,V$-schemes:
\set\begin{equation}\label{eq_iso-of-ofrm.affine}
\fU\stackrel{\sim}{\to}\fX.
\end{equation}
Let $Y^\sm\subset Y$ be the set of all points $x\in Y$ such
that $Y$ is smooth over $\Spec\,V$ at $x$; according to
\cite[Ch.IV, Cor.6.8.7]{EGAIV-2} and \cite[Ch.IV, Cor.17.5.2]{EGA4},
$Y^\sm$ is an open subset of $Y$.

\begin{claim} $Y_0\setminus X_0\subset Y^\sm$.
\end{claim}
\begin{pfclaim}
Indeed, for every $y\in Y_0$, the completions of the local
rings $\cO_{\fY,y}$ and $\cO_{Y,y}$ are isomorphic topological
$V$-algebras, therefore the claim follows from
\cite[Ch.IV, Prop.19.3.6]{EGAIV} and \cite[Ch.IV, Prop.17.5.3]{EGA4}.
\end{pfclaim}

\begin{claim}\label{cl_Y-sing.is.finite} The (reduced) closed
subscheme $Y^\sing:=Y\setminus Y^\sm$ is finite over
$\Spec\,V$, and moreover $Y^\sing\subset U$.
\end{claim}
\begin{pfclaim} First of all, since the morphism $Y\to\Spec\,V$
is proper (especially, universally closed), the closure of any
point of $Y$ meets the closed fibre $Y_0$. It follows easily
that the dimension of $Y^\sing$ is the maximum of the dimension
of the local rings $\cO_{Y^\sing,y}$, where $y$ ranges over all
the points of $Y^\sing\cap Y_0$. However, $Y^\sing\cap Y_0$
is precisely the set of points $y\in Y_0$ with the property that
$Y_0$ is not smooth over $\Spec\,V_0$ at $y$; in other words,
$Y^\sing\cap Y_0\subset\Sigma$, which consists of finitely
many closed points, whence:
$$
\dim\cO_{Y^\sing,y}\otimes_VV_0=0\qquad
\text{for every $y\in Y^\sing\cap Y_0$.}
$$
On the other hand, quite generally we have the inequality:
$$
\dim\cO_{Y^\sing,y}\leq\dim V+\dim\cO_{Y^\sing,y}\otimes_VV_0
$$
for every $y\in Y^\sing\cap Y_0$ (see \cite[Th.15.1]{Mat}),
therefore we conclude that the relative dimension of $Y^\sing$
over $\Spec\,V$ equals zero, hence $Y^\sing$ is finite over
$\Spec\,V$, as claimed. Since $U$ contains $Y^\sing\cap Y_0$
by construction, and since every point of $Y^\sing$ admits a
specialization to a point of $Y_0$, it is clear that $Y^\sing$
is contained in $U$.
\end{pfclaim}

For any $V$-algebra $A$ denote by $A^\wedge_\fn$ the $\fn$-adic
completion of $A$.

\begin{claim}\label{cl_sledgeham} For any $V$-algebra $R$ of finite
type, the natural morphism $\Spec\,R^\wedge_\fn\to\Spec\,R$ is regular.
\end{claim}
\begin{pfclaim} According to \cite[Ch.IV, Prop.7.4.6]{EGAIV-2},
it suffices to show that all the formal fibres of $R$ are
geometrically regular. However, this follows from
\cite[Ch.IV, Th.7.4.4(ii)]{EGAIV-2}.
\end{pfclaim}

Say that $X=\Spec\,R$ and $U=\Spec\,S$, for some $V$-algebras
$R$ and $S$ of finite type; \eqref{eq_iso-of-ofrm.affine} induces
an isomorphism of $V$-schemes
$$
\phi^\wedge:U^\wedge:=\Spec\,S^\wedge_\fn\stackrel{\sim}{\to}
X^\wedge:=\Spec\,R^\wedge_\fn.
$$

\begin{claim}\label{cl_Elkik-goes} Let $x\in U^\wedge$ be
any point whose image $y\in U$ lies outside $Y^\sing$.
Let $x':=\phi^\wedge(x)$ and denote by $y'$ the image
of $x'$ in $X$. Then $X$ is smooth over $\Spec\,V$ at
the point $y'$.
\end{claim}
\begin{pfclaim} It suffices to show that the induced morphism
$\Spec\,\cO_{X,y}\to\Spec\,V$ is regular
(\cite[Ch.IV, Cor.17.5.2]{EGA4}). However, by assumption
$\Spec\,\cO_{Y,y}\to\Spec\,V$ is regular, hence so is the
morphism $\Spec\,\cO_{Y^\wedge,x}\to\Spec\,V$, in view of claim
\ref{cl_sledgeham}. Consequently the morphism
$\Spec\,\cO_{X^\wedge,x'}\to\Spec\,V$ is regular, and since the
map $\cO_{X,y'}\to\cO_{X^\wedge,x'}$ is faithfully flat, the
claim follows from \cite[Th.32.1]{Mat}.
\end{pfclaim}

Next, let $J\subset S$ be an ideal with $V(J)=Y^\sing$.
In view of claim \ref{cl_Y-sing.is.finite}, the quotient
$S/J^n$ is finite over $V$ for every $n\in\N$; especially,
it is complete for the $\fn$-adic topology. We deduce from
lemma \ref{lem_various-compl}(ii) that the natural map
$S^\wedge_J\to S^\wedge_{\fn S+J}$ is an isomorphism (notation
of \eqref{sec_algebrize}), and then lemma \ref{lem_various-compl}(i)
implies that the natural map $S^\wedge_{\fn J}\to S^\wedge_\fn$ is an
isomorphism as well, whence, by composing $\phi^\wedge$ with
the natural map $X^\wedge\to X$, a morphism of $V$-schemes
$\psi:\Spec\,S^\wedge_{\fn J}\to X$. The latter can be seen as
a section
$$
\sigma^\wedge:\Spec\,S^\wedge_{\fn J}\to U\times_{\Spec\,V}X
$$
of the $U$-scheme $U\times_{\Spec\,V}X$ (namely, $\sigma^\wedge$
is the unique section such that $\pi\circ\sigma^\wedge=\psi$, where
$\pi:U\times_{\Spec\,V} X\to X$ is the natural projection).
It follows from claim \ref{cl_Elkik-goes} that the section
$\sigma^\wedge$ fulfills the assumptions of
\cite[Ch.II, Th.2 bis]{Elk}, so that there exists a section
$$
\sigma:\Spec\,S^h_{\fn J}\to U\times_{\Spec\,V}X.
$$
such that $\sigma$ and $\sigma^\wedge$ agree on the closed subscheme
$\Spec\,S/\fn J$. Let $S^h_\fn$ be the henselization of $S$
along the ideal $\fn S$; finally we define a morphism of $V$-schemes
$$
\beta:\Spec\,S^h_\fn\stackrel{\omega}{\to}\Spec\,S^h_{\fn J}
\stackrel{\sigma}{\to}U\times_{\Spec\,V}X\stackrel{\pi}{\to}X
$$
where $\omega$ is the natural morphism. The theorem is
now a straightforward consequence of the following :

\begin{claim} $\beta$ induces an isomorphism
$\Spec\,S^h_\fn\stackrel{\omega}{\to} X^h_{/X_0}$.
\end{claim}
\begin{pfclaim}[] By construction, $\beta$ and $\phi^\wedge$ agree
on the closed subscheme $\Spec\,S/\fn S$.
Let $\beta^\wedge:U^\wedge\to X^\wedge$ be the morphism induced
by $\beta$ and denote by $\gr^\bullet_\fn R$ and $\gr^\bullet_\fn S$
the graded rings associated to the $\fn$-preadic filtrations of $R$
and respectively $S$; we deduce that the morphisms $\beta$ and
$\phi^\wedge$ induce the same homomorphism
$\gr^\bullet_\fn R\to\gr^\bullet_\fn S$ of graded rings.
Since $\phi^\wedge$ is an isomorphism, it then follows that
$\beta^\wedge$ is an isomorphism as well
(\cite[Ch.III, \S2, n.8, Cor.3]{BouAC}). Next, according to
\cite[Ch.XI, \S2, Th.2]{Ray}, the ring $S^h_\fn$ is the filtered
colimit of a family of \'etale $S$-algebras
$(S_\lambda~|~\lambda\in\Lambda)$ such that $S_\lambda/\fn S_\lambda$
is $S$-isomorphic to $S/\fn S$ for every $\lambda\in\Lambda$.
Especially, the $\fn$-adic completion $S^\wedge_{\lambda,\fn}$
is $S$-isomorphic to $S^\wedge_\fn$, and we can find
$\lambda\in\Lambda$ such that $\beta$ descends to a morphism
$\beta_\lambda:\Spec\,S_\lambda\to X$.
Consequently, the map $R\to S^\wedge_{\lambda,\fn}$ is formally
\'etale for the $\fn$-adic topology, therefore the same holds
for the map $R\to S_\lambda$ induced by $\beta_\lambda$.
Finally, let $\fp\in\Spec\,S_\lambda/\fn S_\lambda$, and set
$\fq:=\fp\cap R\in X_0$; {\em a fortiori\/} we see that $S_{\lambda,\fq}$
is formally \'etale over $R_\fq$ for the $\fq$-preadic
and $\fp$-preadic topologies; by \cite[Ch.IV, Prop.17.5.3]{EGA4}
we deduce that $S_\lambda$ is \'etale over $R$ at all the points
of $\Spec\,S_\lambda/\fn S_\lambda$. Since
$S^h_\fn\simeq S^h_{\lambda,\fn}$, the claim follows.
\end{pfclaim}
\end{proof}

\sset\subsubsection{}\label{subsec_alter-data}
Suppose that $V$ is a valuation ring whose field of fractions $K$
is algebraically closed, and let $S:=\Spec\,V$. We conclude this
section with a result stating the existence of semi-stable
models for curves over the generic point of $S$. The proof
consists in reducing to the case where $V$ is noetherian and
excellent, to which one can apply de Jong's method of
alterations \cite{deJ}. Recall that a {\em semi-stable $S$-curve\/}
is a flat and proper morphism $g:Y\to S$ such that all the geometric
fibres of $g$ are connected curves having at most ordinary
double points as singularities. Denote by $\eta$ the generic
point of $S$. We consider a projective finitely presented
morphism $f:X\to S$ such that $f^{-1}(\eta)$ is irreducible
of dimension one. We also assume that $X$ is an integral scheme,
and $G$ is a given finite group of $S$-automorphisms of $X$.

\begin{proposition}\label{prop_alterate}
In the situation of \eqref{subsec_alter-data}, there exists a
projective and birational morphism $\phi:X'\to X$ of $S$-schemes
such that :
\begin{enumerate}
\alphaenu
\item
The structure morphism $f':X'\to S$ is a semi-stable $S$-curve
whose generic fibre $f^{\prime -1}(\eta)$ is irreducible and
smooth over $\Spec\,K$.
\item
$G$ acts on $X'$ as a group of $S$-automorphisms, and $\phi$
is $G$-equivariant.
\romanenu
\end{enumerate}
\end{proposition}
\begin{proof} Let us write $V$ as the colimit of the filtered
family $(V_i~|~i\in I)$ of its excellent noetherian subrings.
By \cite[Ch.IV, Th.8.8.2(i),(ii)]{EGAIV-3}, we may find $i\in I$
such that both $f$ and the action of $G$ descend to, respectively,
a finitely presented morphism $f_i:X_i\to S_i:=\Spec\,V_i$ and
a finite group of $S_i$-automorphisms of $X_i$. By
\cite[Ch.IV, Th.8.10.5]{EGAIV-3} we may even suppose that
$f_i$ is projective, and -- up to replacing $I$ by a cofinal
subset -- we may suppose that the latter property holds for
all $i\in I$. For every $i\in I$, let $\eta_i$ be the
generic point of $S_i$; we apply 
\cite[Ch.IV, Prop.8.7.2  and Cor.8.7.3]{EGAIV-3} to the
projective system of schemes $(f_i^{-1}(\eta_i)~|~i\in I)$
to deduce that there exists $i\in I$ such that $f_i^{-1}(\eta_i)$
is geometrically irreducible over $\Spec\,\kappa(\eta_i)$.
Let $Z_i\subset X_i$ be the Zariski closure (with reduced
structure) of $f_i^{-1}(\eta_i)$; then $Z_i\times_{S_i}S$
is a closed subscheme of $X$ containing $f^{-1}(\eta)$, so
it coincides with $X$, since the latter is an integral scheme.
Moreover, since $f_i$ is a closed morphism and $\eta_i\in f_i(Z_i)$,
we have $f_i(Z_i)=S_i$. Furthermore, the $G$-action on $X_i$
restricts to a finite group of $S_i$-automorphisms of $Z_i$.
The restriction $Z_i\to S_i$ of $f_i$ fulfills the conditions
of \cite[Th.2.4]{deJ}, hence we may find a commutative diagram :
$$
\xymatrix@R-4pt{
Z'_i \ar[r]^-{\phi_i} \ar[d]_{f'_i} & Z_i \ar[d]^{f_i} \\
S'_i \ar[r]^-{\psi_i} & S_i
}$$
such that $Z'_i$ and $S'_i$ are integral and excellent,
$\phi_i$ and $\psi_i$ are projective, dominant and generically finite,
$f'_i$ is a semi-stable projective $S'_i$-curve and moreover a finite
group $G'$ acts by $S'_i$-automorphisms on $X'_i$ in such a way that
$f'_i$ is $G'$-equivariant; also there is a surjection $G'\to G$ so
that $\phi_i$ is equivariant for the induced $G'$-action (condition
(v) of \cite[Th.2.4]{deJ}). Furthermore, if $\eta'_i$ denotes the
generic point of $S'_i$, the induced morphism
\set\begin{equation}\label{eq_gen-fibres}
f^{\prime -1}_i(\eta'_i)\to
f_i^{-1}(\eta_i)\times_{\Spec\,\kappa(\eta_i)}\Spec\,\kappa(\eta'_i)
\end{equation}
is birational (condition  of \cite[Th.2.4(vii)(b) and Rem.2.3(v)]{deJ}).
After taking the base change $S\to S_i$ we arrive at the
commutative diagram :
$$
\xymatrix{
Z':=Z'_i\times_{S_i}S \ar[r]^-{\phi_{i,S}} \ar[d]_{f'_{i,S}} &
X \ar[d]^f \\
S':=S'_i\times_{S_i}S \ar[r]^-{\psi_{i,S}} & S
}$$
where again $\psi_{i,S}$ is generically finite and projective.
Since $\kappa(\eta)$ is algebraically closed, it follows that
the induced map $\kappa(\eta)\to\kappa(\eta')$ is bijective
for every $\eta'\in\psi_{i,S}^{-1}(\eta)$; by the valuative
criterion for properness (\cite[Ch.II, Th.7.3.8(b)]{EGAII})
we deduce that $\psi_{i,S}$ admits a section
$\sigma:S\to S'$. We set $X':=Z'\times_{S'}S$ (the base change
along the morphism $\sigma$) and denote by $\phi:X'\to X$
the restriction of $\phi_{i,S}$. By construction, $f':X'\to S$
is a semi-stable $S$-curve. From \eqref{eq_gen-fibres}
we deduce easily that the induced morphism
$f^{\prime -1}(\eta)\to f^{-1}(\eta)$ is birational,
and then, since $f'$ is flat, it follows that $X'$ is integral
and $\phi$ is birational.
Moreover, the action of $G'$ is completely determined by
its restriction to the generic fibre $f^{\prime -1}(\eta)$,
and since $\phi$ is equivariant, it follows that this action
factors through $G$.
\end{proof}

\subsection{Vanishing cycles}
We resume that notation of \eqref{sec_annuli}; we let $S:=\Spec\,K^+$,
and denote by $s$ the closed point of $S$. According to
\cite[\S4.2]{Hu2}, to every $S$-scheme $\sX$ and every abelian
sheaf $F$ on $\sX_\et$ one attaches a complex of abelian sheaves
$$
R\Psi_{\sX/S}(F)\in\sD(\sX_{s,\et},\Z)
$$
where $\sX_s:=\sX\times_S\Spec\,K^\sim$ and $\sD(\sX_{s,\et},\Z)$
denotes the derived category of the category of abelian sheaves on
$\sX_{s,\et}$. Its stalk at a point $x\in\sX_s$ can be computed as
follows. Denote by $\sX^h_{/\{x\}}$ the spectrum of the henselization
of the local ring $\cO_{\sX,x}$; then there is a natural isomorphism
in the derived category of complexes of abelian groups:
\set\begin{equation}\label{eq_stalks-nearby}
R\Psi_{\sX/S}(F)_x\simeq R\Gamma((\sX^h_{/\{x\}}\times_S\Spec\,K)_\et,F).
\end{equation}
We let also $\sX_K:=\sX\times_S\Spec\,K$. There is a natural map
$$
R\Gamma(\sX_K,F)\to R\Gamma(\sX_s,R\Psi_{\sX/S}(F))
$$
which is an isomorphism when $\sX$ is proper over $S$. On the
other hand, one has a natural morphism in $\sD(\sX_{s,\et},\Z)$
\set\begin{equation}\label{eq_van-cycle}
F_{|\sX_s}[0]\to R\Psi_{\sX/S}(F)
\end{equation}
and the cone of \eqref{eq_van-cycle} is a complex on $\sX_s$
called the {\em complex of vanishing cycles\/} of the sheaf
$F$, and denoted by $R\Phi_{\sX/S}(F)$.

Any $S$-morphism $\phi:\sX\to\sY$ induces a natural map of complexes
\set\begin{equation}\label{eq_vanish-map}
\phi^*R\Phi_{\sY/S}(F)\to R\Phi_{\sX/S}(\phi^*F)
\end{equation}
for every abelian sheaf $F$ on $\sY_\et$.

\sset\subsubsection{}\label{subsec_keep}
Let $\pi\in\fm$ be any non-zero element, and set $\sA:=K^+[S,T]/(ST-\pi^2)$;
the $\pi$-adic completion of $\sA$ is the subring $A^\circ$ of the
affinoid ring $A:=A(a,a^{-1})$, where $a:=|\pi|$ (cp. example
\ref{ex_discriminant}).
Let $\sA^h$ be the henselization of $\sA$ along its ideal $\pi\sA$;
by proposition \ref{prop_normal}(i), the base change functor :
$$
\sA^h\Alg_{\fpet/K}\to A^\circ\Alg_{\fpet/K}\quad :\quad
\sB\mapsto A^\circ\otimes_{\sA^h}\sB.
$$
is an equivalence.

\sset\subsubsection{}\label{subsec_sans-serif}
Let $f:X:=\Spa\,B\to\D(a,a^{-1})$ be a finite \'etale morphism,
so that $B$ is a finite \'etale $A$-algebra, and suppose moreover
that a finite group $G$ acts freely on $X$ in such a way that $f$
becomes a $G$-equivariant morphism, provided we endow $\D(a,a^{-1})$
with the trivial $G$-action. This situation includes the basic case
where $f$ is a Galois (\'etale) morphism with Galois group $G$,
but we also allow the case where $G$ is the trivial
group. Under these assumptions, $B$ is normal, therefore the same
holds for $B^\circ$. Moreover, by lemma \ref{lem_finitely-pres},
$B^\circ$ is a finitely presented $A^\circ$-module;
we denote by $\sB$ the unique (up to unique isomorphism)
$\sA^h$-algebra corresponding to $B^\circ$ under the
equivalence of \eqref{subsec_keep}. By proposition \ref{prop_normal}(iii),
$\sB$ is normal, and clearly the $G$-action on $B$ translates
into a $G$-action on $\sB$, fixing $\sA^h$. Since
$\sA^h/\fm\sA^h\simeq A^\circ/\fm A^\circ$, we can view the prime
ideal $\fP\subset A^\circ$ defined in the proof of claim
\ref{cl-at-least-two}, as an element of $\Spec\,\sA^h/\fm\sA^h$.
Similarly, the finitely many prime ideals
$\fq_1,\dots,\fq_n\subset B^\circ$ lying over $\fP$ can be viewed
as elements of $\Spec\,\sB/\fm\sB$. We also obtain a natural
action of $G$ on the set $\{\fq_1,\dots,\fq_n\}$. For every
$i=1,\dots,n$, we let:
\begin{itemize}
\item
$St(\fq_i)\subset G$ be the stabilizer of $\fq_i$,
\item
$\sA^h_\fP$ (resp. $\sB_i^h$) the henselization of the local
ring $\sA_\fP$ (resp. of $\sB_{\fq_i}$),
\item
$\sA_{\fP,K}^h:=\sA_\fP^h\otimes_{K^+}K$,
$\sB_{i,K}^h:=\sB_i^h\otimes_{K^+}K$,
$\sT^h_K:=\Spec\,\sA_{\fP,K}^h$, $\sX^h_{i,K}:=\Spec\,\sB_{i,K}^h$
\item
$\Lambda$ a finite local ring such that $\ell^n\Lambda=0$, for
a prime number $\ell\neq\chara\,K^\sim$ and some integer $n>0$.
\end{itemize}
With this notation we define:
$$
\Delta(X,\fq_i,F):=R\Gamma((\sX_{i,K}^h)_\et,F)
\quad\text{and}\quad
\Delta(X,F):=\bigoplus^n_{i=1}\Delta(X,\fq_i,F)
$$
for every sheaf of $\Lambda$-modules $F$ on the \'etale site of
$\sX^h_{i,K}$. Moreover, if $C^\bullet$ is a bounded complex of
$\Lambda$-modules and $H^\bullet C^\bullet$ is a finite $\Lambda$-module,
we shall denote by $\chi(C^\bullet)$ the
{\em Euler-Poincar\'e characteristic\/} of $C^\bullet$, which is
defined by the rule :
$$
\chi(C^\bullet):=
\frac{\sum_{i\in\Z}(-1)^i\cdot\mathrm{length}_\Lambda(H^iC^\bullet)}
{\mathrm{length}_\Lambda(\Lambda)}.
$$
In case $C^\bullet$ is a complex of free $\Lambda$-modules
(especially, $C^\bullet$ is perfect), we have also the identity:
$$
\chi(C^\bullet)=
\sum_{i\in\Z}(-1)^i\cdot\rk_\Lambda(C^i).
$$
As usual (\cite[Exp.X, \S1]{SGA5}) one can view the constant
sheaf $\Lambda_{\sX_{i,K}^h}$ on $(\sX_{i,K}^h)_\et$ as a sheaf of
$G$-modules with trivial $G$-action, and then by functoriality,
$\Delta(X,\fq_i,\Lambda)$ is a complex of $\Lambda[St(\fq_i)]$-modules
in a natural way. Furthermore, $\Delta(X,\Lambda)$ is a complex
of $\Lambda[G]$-modules, whose structure can be analyzed as follows.
Let $O_1,\cup\cdots\cup O_k=\{\fq_1,\dots,\fq_n\}$
be the decomposition into orbits under the $G$-action, and
for every $i\leq k$ let us pick a representative $\fq_i\in O_i$;
then :
\set\begin{equation}\label{eq_induce-from-St}
\Delta(X,\Lambda)\simeq
\sum_{i=1}^k\Ind^G_{St(\fq_i)}\Delta(X,\fq_i,\Lambda)
\qquad\text{in $\sD(\Lambda[G]\Mod)$}.
\end{equation}
The proof is left to the reader.

\begin{lemma}\label{lem_chi}
With the notation of \eqref{subsec_sans-serif}, the following holds:
\begin{enumerate}
\item
The complex of $\Lambda$-modules $\Delta(X,\fq_i,F)$ is perfect
of amplitude $[0,1]$, for every constructible sheaf $F$
of free\/ $\Lambda$-modules on $(\sX_{i,K}^h)_\et$.
\item
For every $n\in\N$, there is a natural isomorphism in\/
$\sD^+(\Z/\ell^n\Z\Mod)$ :
$$
\Z/\ell^n\Z\derotimes_{\Z/\ell^{n+1}\Z}\Delta(X,\fq_i,\Z/\ell^{n+1}\Z)
\isom\Delta(X,\fq_i,\Z/\ell^n\Z).
$$
\item
For every locally constant sheaf $F$ of\/ $\Lambda$-modules on
$(\sT_K^h)_\et$ we have :
$$
\chi(\Delta(\D(a,a^{-1}),\fP,F))\leq 0.
$$
\end{enumerate}
\end{lemma}
\begin{proof} To start out, notice that the scheme
$\Spec\,\sB_{\fq,K}^h$ is a cofiltered limit of one-dimensional
affine $K$-schemes of finite type; then (ii) is an easy
consequence of \cite[Th. finitude, Cor.1.11]{SGA4half}.
We also deduce that the cohomological dimension of
$\Spec\,\sB_{i,K}^h$ is $\leq 1$; furthermore, it follows from
\eqref{eq_stalks-nearby} and \cite[Prop.4.2.5]{Hu2} that
$H^n\Delta(X,\fq_i,F)$ has finite length for every $n\in\N$
and every constructible sheaf $F$ of $\Lambda$-modules,
hence assertion (i) follows from :

\begin{claim}\label{cl_perf-crit}
Let $R$ be a (not necessarily commutative) right noetherian
ring with center $R_0\subset R$, $\phi:\sZ\to\sY$ a morphism
of schemes, and $C^\bullet$ a complex in $\sD^-(\sZ_\et,R)$.
Suppose that the functor
$$
R\phi_*:\sD^+(\sZ_\et,R_0)\to\sD^+(\sY_\et,R_0)
$$
has finite cohomological dimension (thus, $R\phi_*$ extends
to the whole of $\sD(\sZ_\et,R_0)$). Then :
\begin{enumerate}
\item
$C^\bullet$ is pseudo-coherent if and only if $H^n C^\bullet$
is coherent for every $n\in\Z$.
\item
$C^\bullet$ is perfect if and only if it is pseudo-coherent
and has locally finite $\Tor$-dimension.
\item
If the $\Tor$-dimension of $C^\bullet$ is $\leq d$
(for some $d\in\Z$), then the $\Tor$-dimension of
$R\phi_*C^\bullet\in\mathrm{Ob}(\sD^-(\sY_\et,R))$ is $\leq d$.
\end{enumerate}
\end{claim}
\begin{pfclaim} (i) (resp. (ii), resp. (iii)) is a special case of
\cite[Exp.I, Cor.3.5]{Bert} (resp. \cite[Exp.I, Cor.5.8.1]{Bert},
resp. \cite[Exp.XVII, Th.5.2.11]{SGA4-3}).
\end{pfclaim}

(iii): It follows from (i) that the Euler-Poincar\'e
characteristic of $\Delta(\D(a,a^{-1}),\fP,F)$ is well defined
when $\Lambda$ is a finite field of characteristic $\ell$.
For the general case, let $\fm_\Lambda\subset\Lambda$ be the
maximal ideal; we consider the descending filtration
$F\supset\fm_\Lambda F\supset\fm_\Lambda^2 F\supset
\cdots\supset\fm_\Lambda^rF=0$, whose graded subquotients
are sheaves of modules over the residue field $\kappa(\Lambda)$
of $\Lambda$; since the expression that we have to evaluate is
obviously additive in $F$, we are then reduced to the case where
$F$ is an irreducible locally constant sheaf of
$\kappa(\Lambda)$-modules. In such case,
if $F$ is not constant the assertion is clear, so we can further
suppose that $F$ is the constant sheaf with stalks isomorphic
to $\kappa(\Lambda)$. However, $\sA_\fP^h$ is a normal domain,
therefore $\sT_K^h$ is connected, so
$H^0\Delta(\D(a,a^{-1}),\fP,\kappa(\Lambda))=\kappa(\Lambda)$.
Finally, from \cite[Exp.XV, \S2.2.5]{SGA7} we derive
$H^1\Delta(\D(a,a^{-1}),\fP,\kappa(\Lambda))=\kappa(\Lambda)$.
({\em loc.cit.} considers the vanishing cycle functor relative to
a family defined over a strictly henselian discrete valuation ring,
but by inspecting the proofs it is easy to see that the same
argument works {\em verbatim\/} in our setting as well.)
\end{proof}

\sset\subsubsection{}\label{subsec_not-necess}
Let $R$ be a (not necessarily commutative) ring. We denote by
$K^0(R)$ (resp. by $K_0(R)$) the Grothendieck group of finitely
generated projective (resp. finitely presented) left $R$-modules.
Any perfect complex $C^\bullet$ of $R$-modules determines a class
$[C^\bullet]\in K^0(R)$. Namely, one chooses a quasi-isomorphism
$P^\bullet\stackrel{\sim}{\to}C^\bullet$ from a bounded complex
of finitely generated projective left $R$-modules, and sets
$[C^\bullet]:=\sum_{i\in\Z}(-1)^i\cdot[P^i]$; a standard
verification shows that the definition does not depend on the
chosen projective resolution.

\begin{proposition}\label{prop_perf-with-H}
In the situation of \eqref{subsec_sans-serif} we have:
\begin{enumerate}
\item
$\Delta(X,\fq_i,\Lambda)$ is a perfect complex of
$\Lambda[St(\fq_i)]$-modules of amplitude $[0,1]$.
\item
$[\Delta(X,\fq_i,\F_\ell)[1]]\in K^0(\F_\ell[St(\fq_i)])$ is the
class of a projective $\F_\ell[St(\fq_i)]$-module.
\item
$[\Delta(X,\F_\ell)[1]]\in K^0(\F_\ell[G])$ is the
class of a projective $\F_\ell[G]$-module.
\end{enumerate}
\end{proposition}
\begin{proof} (i): This is well known: let $\sff_K:\sX^h_{i,K}\to\sT^h_K$
be the natural morphism; one shows as in the proof of
\cite[Exp.X, Prop.2.2]{SGA5} that $\sff_{K*}\Lambda$ is a flat sheaf
of $\Lambda[St(\fq_i)]$-modules, and then claim \ref{cl_perf-crit}(iii)
implies that the complex of $\Lambda[St(\fq_i)]$-modules
$\Delta(X,\fq_i,\Lambda)\simeq\Delta(\D(a,a^{-1}),\fP,\sff_{K*}\Lambda)$ is of
finite Tor-dimension. It then follows from claim \ref{cl_perf-crit}(i),(ii)
and lemma \ref{lem_chi}(i) that $\Delta(X,\fq_i,\Lambda)$ is a perfect complex
of $\Lambda[St(\fq_i)]$-modules of amplitude $[0,1]$.

(ii): Let us choose a complex $\Delta^\bullet:=(\Delta^0\to\Delta^1)$
of finitely generated projective $\F_\ell[St(\fq_i)]$-modules with a
quasi-isomorphism
$\Delta^\bullet\stackrel{\sim}{\to}\Delta(X,\fq_i,\F_\ell)$; if $M$ is any
$\F_\ell[St(\fq_i)]$-module of finite length, we deduce a quasi-isomorphism
(cp. the proof of lemma \ref{lem_chi}(ii))
$$
M\otimes_{\F_\ell[St(\fq_i)]}\Delta^\bullet\stackrel{\sim}{\to}
R\Gamma(\sT^h_K,M\otimes_{\F_\ell[St(\fq_i)]}\sff_{K*}\F_\ell).
$$
Whence $\chi(M\otimes_{\F_\ell[St(\fq_i)]}\Delta^\bullet)\leq 0$, in view
of lemma \ref{lem_chi}(iii). In other words :
$$
\rk_{\F_\ell} M\otimes_{\F_\ell[St(\fq_i)]}\Delta^0\leq
\rk_{\F_\ell} M\otimes_{\F_\ell[St(\fq_i)]}\Delta^1
\qquad\text{for every $M$ of finite length}.
$$
On the other hand, $K_0(\F_\ell[St(\fq_i)])$ is endowed with an
involution (\cite[Exp.X, \S3.7]{SGA5})
$$
K_0(\F_\ell[St(\fq_i)])\to K_0(\F_\ell[St(\fq_i)])\quad:\quad
[M]\mapsto [M]^*:=[M^*]:=[\Hom_{\F_\ell}(M,\F_\ell)].
$$
We have a natural isomorphism
$$
\Hom_{\F_\ell[St(\fq_i)]}(N,M^*)\simeq(N\otimes_{\F_\ell[St(\fq_i)]}M)^*
$$
for every $\F_\ell[St(\fq_i)]$-modules of finite length $M$ and $N$
(\cite[Exp.X, Prop.3.8]{SGA5}).
Since clearly $\rk_{\F_\ell}M=\rk_{\F_\ell}M^*$ for every such $M$,
we conclude that
\set\begin{equation}\label{eq_envelopes}
\rk_{\F_\ell}\Hom_{\F_\ell[St(\fq_i)]}(\Delta^0,M)\leq
\rk_{\F_\ell}\Hom_{\F_\ell[St(\fq_i)]}(\Delta^1,M)
\quad\text{for every $M$ of finite length}.
\end{equation}
By \cite[\S14.3, Cor.1,2]{Se2} the projective modules $\Delta^0$
and $\Delta^1$ are direct sums of projective envelopes of simple
$\F_\ell[St(\fq_i)]$-modules; however, \eqref{eq_envelopes}
implies that the multiplicity in $\Delta^0$ of the projective
envelope $P_N$ of any simple module $N$ is $\leq$ the
multiplicity of $P_N$ in $\Delta^1$, whence the assertion.

(iii) is an easy consequence of (ii) and \eqref{eq_induce-from-St}.
\end{proof}

\sset\subsubsection{}\label{subsec_go-to-adic}
Keep the assumptions of proposition \ref{prop_perf-with-H}.
For every $n\in\N$ we set $\Lambda_n:=\Z/\ell^n\Z$; in view
of lemma \ref{lem_chi}(ii) we derive natural isomorphisms
in $\sD^+(\Lambda_n[St(\fq_i)]\Mod)$ :
\set\begin{equation}\label{eq_-adic-system}
\Lambda_n[St(\fq_i)]\derotimes_{\Lambda_{n+1}[St(\fq_i)]}
\Delta(X,\fq_i,\Lambda_{n+1})
\isom\Lambda_n\derotimes_{\Lambda_{n+1}}\Delta(X,\fq_i,\Lambda_{n+1})
\isom\Delta(X,\fq_i,\Lambda_n).
\end{equation}
Then, according to \cite[Exp.XIV, \S3, n.3, Lemme 1]{SGA5}
we may find :
\begin{itemize}
\item
An inverse system $(\Delta^\bullet_n(X,\fq_i)~|~n\in\N)$, such that
$\Delta^\bullet_n(X,\fq_i)$ is a complex of projective
$\Lambda_n[St(\fq_i)]$-modules of finite rank, concentrated in
degrees $0$ and $1$, and the transition maps are
isomorphisms of complexes of $\Lambda_n[St(\fq_i)]$-modules :
\set\begin{equation}\label{eq_strict-adic-syst}
\Lambda_n\otimes_{\Lambda_{n+1}}\Delta^\bullet_{n+1}(X,\fq_i)
\isom\Delta^\bullet_n(X,\fq_i)
\qquad\text{for every $n\in\N$.}
\end{equation}
\item
A system of isomorphisms :
$\Delta^\bullet_n(X,\fq_i)\isom\Delta(X,\fq_i,\Lambda_n)$
in $\sD^+(\Lambda_n[St(\fq_i)]\Mod)$, compatible with the isomorphisms
\eqref{eq_-adic-system} and \eqref{eq_strict-adic-syst}.
\end{itemize}
(Actually, {\em loc.cit.} includes the assumption that the
coefficient rings are commutative. This assumption is not verified
by our system of rings $\Lambda_n[St(\fq_i)]$; however, by inspecting
the proof, one sees easily that the commutativity is not needed.)

We let $\Delta^\bullet_\infty(X,\fq_i)$ be the inverse limit of
the system $(\Delta^\bullet_n(X,\fq)~|~n\in\N)$; this is a complex
of projective $\Z_\ell[St(\fq_i)]$-modules of finite rank, concentrated
in degrees $0$ and $1$, and we have isomorphisms of complexes of
$\Lambda_n[St(\fq_i)]$-modules :
$$
\Lambda_n[St(\fq_i)]\otimes_{\Z_\ell[St(\fq_i)]}\Delta^\bullet_\infty(X,\fq_i)
\simeq\Lambda_n\otimes_{\Z_\ell}\Delta^\bullet_\infty(X,\fq_i)
\simeq\Delta^\bullet_n(X,\fq_i)\qquad\text{for every $n\in\N$}.
$$
Likewise, we set
$$
\Delta^\bullet_n(X):=\bigoplus^n_{i=1}\Delta^\bullet_n(X,\fq_i)
$$
and the analogue of \eqref{eq_induce-from-St} holds for
$\Delta^\bullet_n(X)$, especially the latter is a complex
of finitely generated projective $\Lambda_n[G]$-modules
and the inverse limit $\Delta^\bullet_\infty(X)$ of the
system $(\Delta^\bullet_n(X)~|~n\in\N)$ is a complex of
finitely generated projective $\Z_\ell[G]$-modules.

\begin{lemma}\label{lem_go-to-adic} In the situation
of \eqref{subsec_go-to-adic} :
\begin{enumerate}
\item
The element
$[\Delta^\bullet_\infty(X,\fq_i)[1]]\in K^0(\Z_\ell[St(\fq_i)]\Mod)$
is the class of a finitely generated projective
$\Z_\ell[St(\fq_i)]$-module.
\item
The element
$[\Delta^\bullet_\infty(X)[1]]\in K^0(\Z_\ell[G]\Mod)$
is the class of a finitely generated projective
$\Z_\ell[G]$-module.
\end{enumerate}
\end{lemma}
\begin{proof} (i) follows from proposition
\ref{prop_perf-with-H}(ii) and \cite[\S14.4, Cor.3]{Se2}, and
(ii) follows from (i).
\end{proof}

\sset\subsubsection{}\label{subsec_explain-l-adic}
In view of lemma \ref{lem_go-to-adic}, the element
$[\Delta^\bullet_\infty(X,\fq_i)\otimes_\Z\Q]\in K^0(\Q_\ell[St(\fq_i)])$
is the class of a finite-dimensional $\ell$-adic
representation of $St(\fq_i)$. For such representations $\rho$,
it makes sense to ask whether the associated character
takes only rational values, {\em i.e.} whether the class
$[\rho]$ lies in the subgroup
$\bar R_\Q(St(\fq_i))\subset K^0(\Q_\ell[St(\fq_i)])$
(notation of \cite[\S12.1]{Se2}), and a complete charaterization
of $\Q\otimes_\Z\bar R_\Q(St(\fq_i))$ is provided by the criterion of
\cite[\S13.1]{Se2}. The availability of that criterion
is the main reason why we are interested in $\ell$-adic
representations (rather than just $\ell$-torsion ones).

\begin{theorem}\label{th_rationality}
With the notation of \eqref{subsec_explain-l-adic} :
\begin{enumerate}
\item
The class $[\Delta^\bullet_\infty(X,\fq_i)\otimes_\Z\Q]$ lies in
$\bar R_\Q(St(\fq_i))$.
\item
The class $[\Delta^\bullet_\infty(X)\otimes_\Z\Q]$ lies in
$\bar R_\Q(G)$.
\end{enumerate}
\end{theorem}
\begin{proof} Of course it suffices to show (i).
We begin with the following :

\begin{claim}\label{cl_apply-semistable} There exist :
\begin{itemize}
\item
a projective birational morphism $\phi:\sY\to\sX$ of integral
projective $S$-schemes, where $\sY\to S$ is a semistable $S$-curve
with smooth connected generic fibre $\sY_K\to\Spec\,K$;
\item
group homomorphisms $St(\fq_i)\to\Aut_S\sY$ and $St(\fq_i)\to\Aut_S\sX$
such that $\phi$ is $St(\fq_i)$-equivariant;
\item
a point $x\in\sX_s$, fixed by $St(\fq_i)$, with an $St(\fq_i)$-equivariant
isomorphism : $\cO^h_{\sX,x}\isom\sB^h_i$;
\item
an open neighborhood $\sU\subset\sX$ of $x$, which is a connected
normal scheme.
\end{itemize}
\end{claim}
\begin{pfclaim} To start out, $\sB^h_i$ is the colimit of a
filtered family $(\sB_\mu~|~\mu\in J)$ of \'etale $\sB$-algebras,
and since since $\sB^h_i$ is a normal domain, we may assume
that the same holds for every $\sB_\mu$.
We may find $\mu\in J$ such that the action of $St(\fq_i)$ on
$\sB^h_i$ descends to an action by $K^+$-automorphisms on
$\sB_\mu$ (\cite[Ch.IV, Th.8.8.2(i)]{EGAIV-3}).
If $x\in\sU:=\Spec\,\sB_\mu$ denotes the contraction of the ideal
$\fq_i\subset\sB$, we have an $St(\fq_i)$-equivariant isomorphism
$\cO^h_{\sU,x}\isom\sB^h_i$. Let $V\subset\sB_\mu$ be a
finitely generated $K^+$-submodule, say of rank $n+1$, that generates
$\cO_\sU$; $V$ determines a morphism $\psi:\sU\to\P^n_S$, and by choosing
$V$ large enough, we may achieve that $\psi$ is a closed immersion;
moreover we may suppose that $V$ is stable under the natural
$St(\fq_i)$-action, in which case $\psi$ is $St(\fq_i)$-equivariant.
Denote by $\sX$ the Zariski closure of $\psi(\sU)$ (with reduced structure).
Then the action of $St(\fq_i)$ extends to $\sX$ and $\sff:\sX\to S$ is a
projective finitely presented morphism; moreover, since $\sU$
is a normal scheme, the generic fibre $\sff^{-1}(\eta)$ is irreducible.
To conclude, it suffices to invoke proposition \ref{prop_alterate}.
\end{pfclaim}

Hence, let $\phi:\sY\to\sX$ and $x\in\sU\subset\sX_s$ be as in
claim \ref{cl_apply-semistable}, and $\phi_s:\sY_s\to\sX_s$ (resp.
$\phi_\eta:\sY_K\to\sX_K$) the restriction of $\phi$;
applying the proper base change theorem
(\cite[Exp.XII, Th.5.1]{SGA4-3}), we derive a natural
$St(\fq_i)$-equivariant isomorphism :
\set\begin{equation}\label{eq_proper-b-change}
R\phi_{s*}R\Psi_{\sY/S}\Lambda\isom
R\Psi_{\sX/S}R\phi_{\eta*}\Lambda
\end{equation}
in $\sD(\sX_{s,\et},\Lambda)$, for every ring $\Lambda$
as in \eqref{subsec_sans-serif}. Set $Z:=\phi_s^{-1}(x)$;
taking
the stalk over $x$ of the map \eqref{eq_proper-b-change} 
yields an $St(\fq_i)$-equivariant isomorphism :
\set\begin{equation}\label{eq_ordin-sing}
R\Gamma(Z,R\Psi_{\sY/S}\Lambda)\isom
(R\Psi_{\sX/S}R\phi_{\eta*}\Lambda)_x.
\end{equation}
Moreover, $\sU_K$ is smooth over $\Spec\,K$, hence $\phi_\eta$
restricts to an isomorphism on $\phi^{-1}\sU_K$, hence:
\set\begin{equation}\label{eq_ordin-sing2}
(R\Psi_{\sX/S}R\phi_{\eta*}\Lambda)_x\simeq
(R\Psi_{\sX/S}\Lambda)_x\simeq\Delta(X,\fq_i,\Lambda)
\end{equation}
where again these isomorphisms are $St(\fq_i)$-equivariant.
On the other hand, \cite[Exp.XV, \S2.2]{SGA7} yields natural
isomorphisms :
\set\begin{equation}\label{eq_annalyze}
R^0\Psi_{\sY/S}\Lambda\simeq\Lambda_{|\sY_s}\qquad
R^1\Psi_{\sY/S}\Lambda\simeq i_*\Lambda(-1)_{|\sY^\sing_s}\qquad
R^j\Psi_{\sY/S}\Lambda=0\quad\text{for $j>1$}
\end{equation}
where $i:\sY_s^\sing\to\sY_s$ is the closed immersion
of the singular locus of $\sY_s$ (which consists of finitely
many $\Spec\,K^\sim$-rational points) and $(-1)$ denotes
the Tate twist. (Actually, {\em loc.cit.} considers the case
where $S$ is a henselian discrete valuation ring, but
by inspecting the proof one sees easily that the same
argument works in our situation as well.) Since $Z$ is proper
over $\Spec\,K^\sim$, one may apply
\cite[Exp.XVII, Th.5.4.3]{SGA4-3} to deduce that
\eqref{eq_ordin-sing} and \eqref{eq_ordin-sing2} still hold
after we replace $\Lambda$ by $\Q_\ell$ and $\Delta(X,\fq_i,\Lambda)$
by $\Delta_\infty^\bullet(X,\fq_i)\otimes_\Z\Q$.
Hence, $[\Delta_\infty^\bullet(X,\fq_i)\otimes_\Z\Q]$
is the difference of the classes :
$$
R_1:=[R\Gamma(Z,\Q_{\ell})] \qquad R_2:=[\Gamma(Z^\sing,\Q_\ell(-1))].
$$
where the $St(\fq_i)$-actions on $R_1$ and $R_2$ are deduced by functoriality
from the actions of $St(\fq_i)$ on the sheaves $\Q_\ell$ and $\Q_\ell(-1)$,
and the latter are defined via \eqref{eq_annalyze}. So the theorem
follows from the following :

\begin{claim}\label{cl_firdt-cl} $R_1,R_2\in\bar R_\Q(St(\fq_i))$.
\end{claim}
\begin{pfclaim}[] Concerining $R_1$ : first of all, notice that the
action of $St(\fq_i)$ on $\Lambda_{|\sY_s}$ (resp. on $\Q_{\ell|\sY_s}$)
induced by the isomorphism \eqref{eq_annalyze}, is the trivial one
(this isomorphism is the map \eqref{eq_van-cycle}). 
Let $\rho:Z'\to Z$ be the normalization morphism,
and say that $W_1,\dots,W_k$ are the irreducible components of $Z'$;
a standard {\em d\'evissage\/} shows that 
$$
R_1=[R\Gamma_c(Z',\Q_\ell)]-[\Gamma(\rho^{-1}Z^\sing,\Q_\ell)]+
[\Gamma(Z^\sing,\Q_\ell)]
$$
where the $St(\fq_i)$-actions on the terms appearing on the right-hand
side are deduced, by functoriality, from the trivial $St(\fq_i)$-actions
on the constant $\ell$-adic sheaves $\Q_\ell$ on the scheme $Z'$.
Let $g\in St(\fq_i)$ be any element, $g':Z'\to Z'$ the
unique $K^\sim$-automorphism that lifts the action of $g$ on $Z$,
and $g'':g^{\prime *}\Q_\ell\isom\Q_\ell$ the isomorphism that
defines the trivial $St(\fq_i)$-action on $\Q_{\ell|Z'}$.
We have a natural decomposition :
$R\Gamma_c(Z',\Q_\ell)\simeq\bigoplus_{i=1}^kR\Gamma_c(W_i,\Q_\ell)$,
and $R\Gamma_c(Z',g'')$ restricts to isomorphisms :
$$
\omega_i:R\Gamma_c(g'(W_i),\Q_\ell)\isom R\Gamma_c(W_i,\Q_\ell)
\qquad\text{for every $i=1,\dots,k$}.
$$
It follows that the trace of $R\Gamma_c(Z',g'')$ is the sum
of the traces of the maps $\omega_i$ such that $g'(W_i)=W_i$.
Each $W_i$ is either a point or a smooth projective $K^\sim$-curve.
In case $W_i$ is a point, $\omega_i$ is the identity map; to
determine the trace of $\omega_i$ in case $W_i$ is a smooth
curve, we may apply the Lefschetz fixed point formula
\cite[Rapport, Th.5.3]{SGA4half}, and it follows easily that
$[R\Gamma_c(Z',\Q_\ell)]\in\bar R_\Q(St(\fq_i))$.
Next, consider the term $[\Gamma(Z^\sing,\Q_\ell)]$;
a similar argument shows that, in order to compute the trace of
the automorphism induced by $g$ on $\Gamma(Z^\sing,\Q_\ell)$, we
may neglect all the points of $Z^\sing$ that are not fixed by the
action of $g$; if $z\in Z^\sing$ is fixed by $g$, then clearly the
trace of $\Gamma(\{z\},g)$ equals $1$, so we get as well
$[\Gamma(Z^\sing,\Q_\ell)]\in\bar R_\Q(St(\fq_i))$.
Finally, we consider $[\Gamma(\rho^{-1}Z^\sing,\Q_\ell)]$ : let
$(Z^\sing)^g$ be the set of points of $Z^\sing$ that are fixed
by $g$; an argument as in the foregoing shows that the trace of
$\Gamma(\rho^{-1}Z^\sing,g'')$ is the same as the trace of
$\Gamma(\rho^{-1}(Z^\sing)^g,g'')$. For every $z\in(Z^\sing)^g$,
the fibre $\rho^{-1}(z)$ consists of two points $z'_1$ and $z'_2$,
and clearly $g$ either exchanges them, in which case the corresponding
contribution to the trace is $0$, or else $g$ fixes them, in which
case the contribution is $2$. Hence
$[\Gamma(\rho^{-1}Z^\sing,\Q_\ell)]\in\bar R_\Q(St(\fq_i))$,
so the claim holds for $R_1$.

Concerning $R_2$ : let again $g\in St(\fq_i)$ be any element;
the action of $g$ on $R_2$ is induced by an action of $g$ on
$(R^1\Psi_\eta\Q_\ell)_{|\sY^\sing}$, {\em i.e.} by an isomorphism
$g':g^*\Q_\ell(-1)_{|\sY^\sing}\isom\Q_\ell(-1)_{|\sY^\sing}$.
Arguing as in the foregoing case, we see that the trace of
$\Gamma(Z^\sing,g')$ is the same as the trace of
$\Gamma((Z^\sing)^g,g')$. Hence, for our purposes, it suffices
to determine the automorphism $g'_z$ of the stalk over any
$z\in Z^\sing$ which is fixed by $g$. By \cite[Exp.XV, \S2.2]{SGA7},
Poincar\'e duality yields a perfect pairing :
$$
(R^1\Psi_\eta\Q_\ell)_z\times H^1_{\{z\}}(\sY_s,R\Psi_\eta\Q_\ell(1))
\to\Q_\ell.
$$
Hence it suffices to show that $[H^1_{\{z\}}(\sY_s,R\Psi_\eta\Q_\ell(1))]$
lies in $\bar R_\Q(St(\fq_i))$. However, according to
\cite[Exp.XV, lemme 2.2.7]{SGA7}, we have a natural
short exact sequence :
$$
0\to H^0(\{z\},\Q_\ell)(1)\to H^0(\rho^{-1}\{z\},\Q_\ell)(1)\to
H^1_{\{z\}}(\sY_s,R\Psi_\eta\Q_\ell(1))\to 0
$$
which becomes $St(\fq_i)$-equivariant, provided we endow the
$\ell$-adic sheaves $\Q_{\ell|\{z\}}$ and $\Q_{\ell|\rho^{-1}\{z\}}$
with their trivial actions. It follows that the trace of 
$g'_z$ equals $1$ if $g$ fixes the points of $\rho^{-1}\{z\}$,
and equals $-1$ if $g$ exchanges these two points.
\end{pfclaim}
\end{proof}

\sset\subsubsection{}
Keep the notation of \eqref{subsec_sans-serif} and let
$H\subset G$ be any subgroup; since $A$ is a Japanese ring
(\cite[\S6.1.2, Prop.4]{Bo-Gun}) we see easily that the subring
$B^H$ of elements fixed by $H$ is an affinoid algebra; we can
then consider the morphism
$f_H:X/H:=\Spa\,B^H\to\D(a,a^{-1})$; clearly $f_H$ is
again \'etale (indeed, this can be checked after an \'etale
base change, especially, after base change to $X$, in which
case the assertion is obvious). Moreover, obviously $(B^H)^\circ=(B^\circ)^H$; 
under the equivalence of \eqref{subsec_keep}, the finitely
presented $A^\circ$-algebra $(B^H)^\circ$ corresponds to a
unique (up to unique isomorphism) finitely presented
$\sA^h$-algebra $\sC$ such that $\sC\otimes_{K^+}K$ is \'etale over
$\sA^h_K:=\sA^h\otimes_{K^+}K$. By lemma \ref{lem_future} the
natural map $\sA^h\to A^\circ$ is faithfully flat; by considering
the left exact sequence of $\sA^h$-modules
$$
\xymatrix{
0 \ar[r] & \sB^H \ar[r] &\sB \ar[rr]^-{\oplus_{h\in H}(1-h)} & &
\bigoplus_{h\in H}\sB
}$$
one deduces easily that
\set\begin{equation}\label{eq_compat}
\sC=\sB^H.
\end{equation}

\sset\subsubsection{}\label{subsec_invariance}
Suppose next, that the subgroup $H\subset G$ is contained in
$St(\fq_i)$. We denote by $\fq_i^H$ the image of $\fq_i$ in
$\Spec\,(B^H)^\circ$. In view of \eqref{eq_compat} the induced map
$$
\sg:\sX^h_{i,K}\to\sY^h_{i,K}:=\Spec\,\sC_{\fq_i^H,K}^h
$$
is a Galois \'etale covering with Galois group $H$.

\begin{lemma}\label{lem_H-action}
In the situation of \eqref{subsec_invariance}, we have a natural
isomorphism in $\sD(\Z_\ell\Mod)$ :
$$
\Delta^\bullet_\infty(X/H,\fq^H_i)\stackrel{\sim}{\to}
\Delta^\bullet_\infty(X,\fq_i)^H.
$$
\end{lemma}
\begin{proof} (Notice that, since in general $H$ is not a
normal subgroup of $St(\fq_i)$, the only group
surely acting on $\Delta(X/H,\fq_i^H,\Lambda_n)$ is the trivial
one, so $\Delta^\bullet_\infty(X/H,\fq_i^H)$ is to be meant as a
complex of free $\Z_\ell$-modules.)
To start with, let $\Lambda$ be any ring as in
\eqref{subsec_sans-serif}; the functor
$F\mapsto\underline\Gamma^H(F):=F^H$ on sheaves of
$\Lambda[H]$-modules on $(\sY^h_{i,K})_\et$ induces a derived functor
$$
R\underline\Gamma^H:\sD^+((\sY^h_{i,K})_\et,\Lambda[H])
\to\sD^+((\sY^h_{i,K})_\et,\Lambda).
$$
Likewise, we have a derived functor :
$$
R\Gamma^H:\sD^+(\Lambda[H]\Mod)\to\sD^+(\Lambda\Mod).
$$
Especially, consider the inverse system
$(\Delta^\bullet_n(X,\fq_i)~|~n\in\N)$ of \eqref{subsec_go-to-adic};
since each $\Delta^\bullet_n(X,\fq_i)$ is a complex of projective
$\Lambda_n[St(\fq_i)]$-modules, the natural map :
$$
\Gamma^H\Delta^\bullet_n(X,\fq_i)\to R\Gamma^H\Delta^\bullet_n(X,\fq_i)
$$
is an isomorphism in $\sD(\Lambda_n\Mod)$, for every $n\in\N$.
Similarly, since $\sg_*\Lambda_{n,\sX^h_{i,K}}$ is a sheaf of projective
$\Lambda_n[H]$-modules, we have natural isomorphisms of sheaves
on $(\sY_{i,K}^h)_\et$ :
\set\begin{equation}\label{eq_invariance}
\Lambda_{n,\sY_{i,K}^h}\stackrel{\sim}{\to}
R\underline\Gamma^H\sg_*\Lambda_{n,\sX_{i,K}^h}.
\end{equation}
Now, by applying to \eqref{eq_invariance} the triangulated functor
$$
R\Gamma:\sD((\sY_{i,K}^h)_\et,\Lambda_n)\to\sD(\Lambda_n\Mod)
$$
and using the obvious isomorphism of triangulated functors
$$
R\Gamma\circ R\underline\Gamma^H\simeq
R\Gamma^H\circ R\Gamma:
\sD^+((\sY_{i,K}^h)_\et,\Lambda_n[H])\to\sD^+(\Lambda_n\Mod)
$$
we deduce natural isomorphisms :
$$
\Delta^\bullet_n(X/H,\fq_i^H)\isom R\Gamma((\sY^h_{i,K})_\et,\Lambda_n)
\isom\Gamma^H\Delta^\bullet_n(X,\fq_i)\qquad\text{for every $n\in\N$}.
$$
The assertion then follows after taking inverse limits.
\end{proof}

\sset\subsubsection{}
In the situation of \eqref{subsec_sans-serif},
let $\delta:[\log a,-\log a]\cap\log\Gamma_K\to\R_{\geq 0}$ be
the discriminant function of the morphism $f$. We saw in the course of
the proof of theorem \ref{th_big-deal} how one can calculate the variation
of the slope of $\delta$ at the point $\rho=0$ -- that is, the slope around
$\rho=0$ of the function $\rho\mapsto\delta(-\rho)+\delta(\rho)$.
The expression is a sum of contributions indexed by the prime ideals
$\fq_1,\dots,\fq_n$.
Proposition \ref{prop_vanish} explains how these localized
contributions can be read off from the complexes
$\Delta(X,\fq_i,\F_\ell)$ (where $\F_\ell$ is the finite
field with $\ell$ elements).

\begin{proposition}\label{prop_vanish}
Resume the notation of \eqref{subsec_analyze} and
\eqref{subsec_sans-serif}. Then :
\set\begin{equation}\label{eq_was-enum}
2\alpha(\fq_i)+\fF(\fq_i)-2=\chi(\Delta(X,\fq_i,\F_\ell)[1]).
\qquad\text{for every $i=1,\dots,n$}.
\end{equation}
\end{proposition}
\begin{proof} We shall give a global argument: first, according to
\cite[Ch.IV, Prop.17.7.8]{EGA4} and \cite[Ch.XI, \S2, Th.2]{Ray},
we can find :
\begin{itemize}
\item
an \'etale map $\sA\to\sA'$ of $K^+$-algebras of finite presentation
such that the induced map
$\sA\otimes_{K^+}K^\sim\to\sA'\otimes_{K^+}K^\sim$  is an isomorphism;
\item
a finitely presented $\sA'$-algebra $\sB$ with an isomorphism
$\sA^h\otimes_{\sA'}\sB'\stackrel{\sim}{\to}\sB$. Especially, the induced
morphism $\sA'\otimes_{K^+}K\to\sB'\otimes_{K^+}K$ is still \'etale.
\end{itemize}
Set $\sT:=\Spec\,\sA'$ and $\sX:=\Spec\,\sB'$. Using
\eqref{eq_stalks-nearby} one deduces natural isomorphisms
in the derived category of complexes of $\F_\ell$-vector spaces:
\set\begin{equation}\label{eq_natural-f}
\Delta(X,\fq_i,\F_\ell)\simeq R\Psi_{\sX/S}(\F_{\ell,\sX})_{\fq_i}
\end{equation}
for every $i=1,\dots,n$. The special fibre $\sX_s:=\Spec\,\sB'/\fm\sB'$
of the $S$-scheme $\sX$ is of pure dimension one, since it
is finite over $\Spec\,\sA/\fm\sA$, and it is reduced, in view of
lemma \ref{lem_reduction}. Hence $\sX_s$ is generically smooth over
$\Spec\,K^\sim$; denote by $\sX^\nu_s$ and $\sX^\mathrm{n}_s$
respectively the seminormalization and normalization of $\sX_s$
(cp. \cite[\S7.5, Def.13]{Liu}). There are natural finite morphisms
$$
\sX^\mathrm{n}_s\stackrel{\pi_1}{\to}\sX^\nu_s
\stackrel{\pi_2}{\to}\sX_s
$$
and the quotient $\cO_{\sX_s}$-modules
$$
Q_1:=(\pi_2\circ\pi_1)_*\cO_{\sX^\mathrm{n}_s}/\pi_{2*}\cO_{\sX^\nu_s}
\quad\text{and}\quad
Q_2:=\pi_{2*}\cO_{\sX^\nu_s}/\cO_{\sX_s}$$
are torsion sheaves concentrated on the singular locus
$\sX_s^\sing\subset\sX_s$. By inspecting the definition,
one verifies easily that
\set\begin{equation}\label{eq_Q_1}
\alpha(\fq_i)=\dim_{K^\sim}Q_{2,\fq_i}\qquad\text{for every $i=1,\dots,n$}.
\end{equation}
Similarly, using \cite[Th.10.1]{Mat} and lemma \ref{lem_much-ado}
one finds a natural bijection between the points of $\fF(\fq_i)$
and the points of $\sX^\mathrm{n}_s$ lying over the point
$\fq_i\in\sX_s$ (cp. the proof of \cite[Th.6.3]{Hu3}).
This leads to the identity:
\set\begin{equation}\label{eq_Q_2}
\sharp\fF(\fq_i)=1+\dim_{K^\sim}Q_{1,\fq_i}
\qquad\text{for every $i=1,\dots,n$}.
\end{equation}
Now, let us fix one point $\fq:=\fq_i$ and choose an affine open
neighborhood $\sV\subset\sX$ of $\fq$ such that $\sX\setminus\sV$
contains $\sX_s^\sing\setminus\{\fq\}$; by further restricting
$\sV$ we can even achieve that the special fibre $\sV_s$ is connected.
One can then apply theorem \ref{th_algebrize} to produce a projective
$S$-scheme $\sY$ of pure relative dimension one containing an open
subscheme $\sU$ such that $\sY_s$ is connected, $\sY$ is smooth over $S$
at the points of $\sY_s\setminus\sU_s$, and furthermore
$\sU^h_{/\sU_s}\simeq\sV^h_{/\sV_s}$.
By construction, the generic fibre $\sY_K$ is a smooth projective
curve over $\Spec\,K$, and $\sY^\sing_s\subset\{\fq\}$. It is also
clear that the morphism $\sY\to S$ is flat, whence an equality
of Euler-Poincar\'e characteristics:
$$
\chi(\sY_K,\cO_{\sY_K})=\chi(\sY_s,\cO_{\sY_s}).
$$
On the other hand, let $c$ be the number of irreducible components
of $\sY_s$; \cite[\S7.5, Cor.33]{Liu} yields the identity
$$
\begin{array}{r@{\:=\:}l}
\dim_{\F_\ell}H^1(\sY_{s,\et},\F_\ell) &
\dim_{\F_\ell}H^1(\sY^\mathrm{n}_{s,\et},\F_\ell)+
\dim_{K^\sim}(\pi_*\cO_{\sY^\mathrm{n}_s}/\cO_{\sY^\nu_s})-c+1 \\
& \dim_{\F_\ell}H^1(\sY^\mathrm{n}_{s,\et},\F_\ell)+\fF(\fq)-c
\end{array}
$$
where $\sY^\mathrm{n}_s\stackrel{\pi}{\to}\sY^\nu_s$ is the natural
morphism from the normalization to the seminormalization of $\sY_s$.
Since $\dim_{\F_\ell}H^0(\sY^\mathrm{n}_{s,\et},\F_\ell)=c$, we deduce
$$
\chi_\et(\sY_s,\F_\ell)=\chi_\et(\sY^\mathrm{n}_s,\F_\ell)-\fF(\fq)+1.
$$
Furthermore, in light of \eqref{eq_Q_1} and \eqref{eq_Q_2} we can
write
$$
\chi(\sY_s,\cO_{\sY_s})=
\chi(\sY^\mathrm{n}_s,\cO_{\sY_s})-\alpha(\fq)-\sharp\fF(\fq)+1.
$$
By Riemann-Roch we have
$$
\chi_\et(\sY^\mathrm{n}_s,\F_\ell)=2\chi(\sY^\mathrm{n}_s,\cO_{\sY_s})
\qquad\text{and}\qquad
\chi_\et(\sY_K,\F_\ell)=2\chi(\sY_K,\cO_{Y_s}).
$$
Putting everything together we end up with the identity:
$$
\chi_\et(\sY_K,\F_\ell)=
\chi_\et(\sY_s,\F_\ell)-\sharp\fF(\fq)+1-2\alpha(\fq).
$$
Now, the complex $R\Phi_{\sY/S}(\F_\ell)$ is concentrated at 
$\sY^\sing_s\subset\{\fq\}$; it follows that
$R\Phi_{\sY/S}(\F_\ell)\simeq j_!R\Phi_{\sU/S}(\F_\ell)$, where
$j:\sU_s\to\sY_s$ is the open imbedding. Since furthermore
the henselizations of $\sU$ and $\sV$ are isomorphic, we have
a natural identification
$R\Phi_{\sU/S}(\F_\ell)_\fq\simeq R\Phi_{\sV/S}(\F_\ell)_\fq\simeq
R\Phi_{\sX/S}(\F_\ell)_\fq$. Consequently
\set\begin{equation}\label{eq_for-X}
\chi(R\Phi_{\sX/S}(\F_\ell)_\fq)=
\chi_\et(\sY_K,\F_\ell)-\chi_\et(\sY_s,\F_\ell)=
1-\sharp\fF(\fq)-2\alpha(\fq).
\end{equation}
Clearly
$\chi(R\Psi_{\sX/S}(\F_\ell)_\fq)=1+\chi(R\Phi_{\sX/S}(\F_\ell)_\fq)$,
whence the contention.
\end{proof}

\subsection{Conductors}\label{sec_conductors}
We consider now a finite Galois \'etale covering $f:X\to\D(a,b)$
of Galois group $G$.
For given $r\in(a,b]\cap\Gamma_K$, pick any $x\in f^{-1}(\eta(r))$;
$G$ acts transitively on the set $f^{-1}(\eta(r))$ and the stabilizer
subgroup $St_x\subset G$ of $x$ is naturally isomorphic to the Galois
group of the extension of henselian valued fields
$\kappa(r)^{\wedge h}\subset\kappa(x)^{\wedge h}$ (see \cite[\S5.5]{Hu3}).
By \cite[Prop.1.2(iii) and Cor.5.4]{Hu3}, $\kappa(x)^{\wedge h+}$
is a free $\kappa(r)^{\wedge h+}$-module of rank equal to the
order $o(St_x)$ of $St_x$; hence the different $\cD_{x/\eta(r)}^+$ of
the ring extension $\kappa(r)^{\wedge h+}\subset\kappa(x)^{\wedge h+}$
is well-defined and, by arguing as in \eqref{subsec_different}
one sees that it is principal. Moreover, it is shown in
\cite[Lemma 2.1(iii)]{Hu3} that for every $\sigma\in St_x$
there is a value $i_x(\sigma)\in\Gamma^+_x\cup\{0\}$ such that
$$
|t-\sigma(t)|^{\wedge h}_x=i_x(\sigma)
$$
for all $t\in\kappa(x)^{\wedge h}$ such that $|t|^{\wedge h}_x$
is the largest element of $\Gamma_x^+\setminus\{0\}$.
The same argument as in the case of discrete valuations
shows the identity
\set\begin{equation}\label{eq_diff-char}
|\cD^+_{x/\eta(r)}|^{\wedge h}_x=
\prod_{\sigma\in St_x\setminus\{1\}}i_x(\sigma).
\end{equation}
One defines the {\em higher ramification subgroups\/} of $St_x$ by
setting :
$$
P_\gamma:=\{\sigma\in St_x~|~i_x(\sigma)<\gamma\}
\qquad\text{for every $\gamma\in\Gamma^+_x$}
$$
and one says that $\gamma\in\Gamma^+_x$ is a {\em jump\/} in the
family $(P_\gamma~|~\gamma\in\Gamma^+_x)$ if
$P_{\gamma'}\neq P_\gamma$ for every $\gamma'<\gamma$.
When $\gamma<1$, the subgroup $P_\gamma$ is contained in the
unique $p$-Sylow subgroup $St_x^{(p)}$ of $St_x$.

Furthermore, one has Artin and Swan characters; to explain this,
let us introduce the {\em total Artin conductor}:
$$
\sa_x:St_x\to\Gamma_x
\qquad
\sa_x(\sigma):=\left\{\begin{array}{ll}
                   i_x(\sigma)^{-1}
                   & \qquad\text{if $\sigma\neq 1$} \\
                   & \\
                   \displaystyle\prod_{\tau\in St_x\setminus\{1\}}
                   i_x(\tau)
                   & \qquad\text{if $\sigma=1$.}
               \end{array}\right.
$$
It is convenient to decompose this total conductor into two
(normalized) factors:
$$
\sa_x^\natural(\sigma):=o(St_x)\cdot\sa_x(\sigma)^\natural
\quad\text{and}\quad
\sa_x^\flat(\sigma):=-o(St_x)\cdot\log\sa_x(\sigma)^\flat
\qquad\text{for all $\sigma\in St_x$}.
$$
and as usual the Swan character is
$\ssw^\natural_x:=\sa^\natural_x-\su_{St_x}$,
where $\su_{St_x}:=\sreg_{St_x}-1_{St_x}$ is the augmentation character,
{\em i.e.} the regular character $\sreg_{St_x}$ minus the constant
function $1_{St_x}(\sigma):=1$ for every $\sigma\in St_x$. Huber shows
that $\sa^\natural_x$ (resp. $\ssw^\natural_x$) is the character
of an element of $K_0(\Q_\ell[St_x])$ (resp. of $K_0(\Z_\ell[St_x])$ :
see \cite[Th.4.1]{Hu3});
these are respectively the Artin and Swan representations.
In general, these elements are however only virtual representations
(whereas it is well known that in the case of discrete valuation
rings one obtains actual representations). 

\sset\subsubsection{}\label{subsec_in-light-of}
The identities obtained in \cite[\S19.1]{Se2} also generalize as
follows. Let
$\gamma_0:=1\in\Gamma_x$ and $\gamma_1>\gamma_2>\cdots>\gamma_n$
be the jumps in the family $(P_\gamma~|~\gamma\in\Gamma^+_x)$
which are $<1$. Then directly from the definitions we deduce
the identities :
\set\begin{equation}\label{eq_double-wham}
\sa^\flat_x=\sum^n_{i=1}
o(P_{\gamma_i})\cdot\log\frac{\gamma^\flat_{i-1}}{\gamma_i^\flat}
\cdot\Ind^{St_x}_{P_{\gamma_i}}\su_{P_{\gamma_i}} \qquad
\ssw^\natural_x=\sum^n_{i=1}
o(P_{\gamma_i})\cdot(\gamma_i^\natural-\gamma_{i-1}^\natural)\cdot
\Ind^{St_x}_{P_{\gamma_i}}\su_{P_{\gamma_i}}
\end{equation}
where $\su_{P_{\gamma_i}}$ is the augmentation character of
the group $P_{\gamma_i}$. Especially, notice that there exist
two $\R$-valued and, respectively, $\Q$-valued class functions
$\sa^\flat_{St_x^{(p)}}$ and $\ssw^\natural_{St_x^{(p)}}$ of the
$p$-Sylow subgroup $St_x^{(p)}$,
such that:
$$
\sa^\flat_x=\Ind^{St_x}_{St_x^{(p)}}\sa^\flat_{St_x^{(p)}}
\qquad\text{and}\qquad
\ssw^\natural_x=\Ind^{St_x}_{St_x^{(p)}}\ssw^\natural_{St_x^{(p)}}
$$
(the induced class functions from the subgroup $St_x^{(p)}$ to $St_x$;
here we view $\sa^\flat_x$ as an $\R$-valued class function).
Next, we define
$$
\sa^\natural_G(r^+):=\mathrm{Ind}^G_{St_x}\,\sa^\natural_x
\qquad
\sa^\flat_G(r^+):=\mathrm{Ind}^G_{St_x^{(p)}}\,\sa^\flat_{St^{(p)}_x}
\qquad
\ssw^\natural_G(r^+):=\mathrm{Ind}^G_{St_x^{(p)}}\,\ssw^\natural_{St_x^{(p)}}.
$$
Notice that $\sa^\natural_G(r^+)$ and $\ssw^\natural_G(r^+)$
do not depend on the choice of the point $x\in f^{-1}(\eta(r))$.

Moreover, for any element $\chi\in K_0(\C[G])$ we let :
$$
\sa^\flat_G(\chi,r^+):=\langle\sa^\flat_G(r^+),\chi\rangle_G
\qquad
\ssw^\natural_G(\chi,r^+):=\langle\ssw^\natural_G(r^+),\chi\rangle_G
$$
where $\langle\cdot,\cdot\rangle_G$ is the natural scalar product
of $\R\otimes_\Z K_0(\C[G])$ (\cite[\S7.2]{Se2}). For future
reference we point out :

\begin{lemma}\label{lem_positivity}
{\em(i)}\ \ $\sa^\flat_G(\chi,r^+)\in\R$ for every
$\chi\in K_0(\C[G])$.
\begin{enumerate}
\addenu
\item
Moreover, $\sa^\flat_G(\chi,r^+)\geq 0$ whenever
$\mathrm{Res}^G_{St^{(p)}_x}\chi$ is a positive element
of $K_0(\C[St^{(p)}_x])$.
\end{enumerate}
\end{lemma}
\begin{proof}
Both assertions follow easily from \eqref{eq_double-wham}.
\end{proof}

\sset\subsubsection{}
Now, suppose that $f':X'\to\D(a,b)$ is another finite Galois \'etale
covering which dominates $X$, {\em i.e.} such that $f'$ factors
through $f$ and an \'etale morphism $g:X'\to X$. Then $g(X')$ is
a union of connected components of $X$. Let $G'$ be the
Galois group of $f'$; we assume as well that $g$ is equivariant
for the $G'$-action on $X'$ and the $G$-action on $X$, {\em i.e.}
there is a group homomorphism
$$
\phi:G'=\Aut(X'/\D(a,b))\to G=\Aut(X/\D(a,b))
$$
such that $\phi(\sigma)\circ g=g\circ\sigma$ for every $\sigma\in G'$.
Pick $x'\in X'$ lying over $x$; there follows a commutative diagram
of group homomorphisms :
$$
\xymatrix{
St_{x'} \ar[r] \ar[d] & St_x \ar[d] \\
G' \ar[r]^-\phi & G
}$$
whose vertical arrows are injections and whose top horizontal
arrow is a surjection. One shows as in \cite[Ch.VI, \S2, Prop.3]{Se}
that :
$$
\sa^\flat_x=\mathrm{Ind}^{St_x}_{St_{x'}}\sa^\flat_{x'}
\qquad
\sa^\natural_x=\mathrm{Ind}^{St_x}_{St_{x'}}\sa^\natural_{x'}
$$
whence the identities :
\set\begin{equation}\label{eq_bigger-group}
\sa^\flat_G(r^+)=\mathrm{Ind}^G_{G'}\sa^\flat_{G'}(r^+)
\qquad
\sa^\natural_G(r^+)=\mathrm{Ind}^G_{G'}\sa^\natural_{G'}(r^+).
\end{equation}

\sset\subsubsection{}\label{subsec_base-change}
Later we shall also need to know that the conductors are
invariant under changes of base field. Namely, let
$(K,|\cdot|)\to(F,|\cdot|_F)$ be a map of algebraically closed
valued fields of rank one; say that $X=\Spa\,B$, and set :
$$
X_F:=\Spa\,B\hat\otimes_KF \qquad
\D(a,b)_F:=\Spa\,A(a,b)\hat\otimes_KF.
$$
The natural map of adic space $X_F\to X$ is surjective and
$G$-equivariant, hence we may choose $x'\in X_F$ lying over $x$.
Let $f':X_F\to\D(a,b)_F$ be the morphism deduced by base change
from $f$. Then we may define Artin and Swan conductors for the
extension $\kappa(f'(x'))^\wedge\subset\kappa(x')^\wedge$.

\begin{lemma}\label{lem_base-change-cond}
In the situation of \eqref{subsec_base-change}, the following holds :
\begin{enumerate}
\item
The natural map $\kappa(x)\to\kappa(x')$ induces isomorphisms
on value groups and Galois groups :
\set\begin{equation}\label{eq_ident-b-c}
\Gamma_x\isom\Gamma_{x'} \qquad St_{x'}\isom St_x.
\end{equation}
\item
Under the identification \eqref{eq_ident-b-c}, we have :
$$
i_x(\sigma)=i_{x'}(\sigma)\qquad\text{for every $\sigma\in St_x$}.
$$
\end{enumerate}
\end{lemma}
\begin{proof} Set $\eta(r)_F:=f'(x')$ and $\kappa(r)_F:=\kappa(\eta(r)_F)$;
one checks easily that the natural commutative diagram :
$$
\xymatrix{
\kappa(r) \ar[r] \ar[d] & \cB(r):=(f_*\cO_X)_{\eta(r)} \ar[d] \\
\kappa(r)_F \ar[r] & \cB(r)_F:=(f'_*\cO_{X_F})_{\eta(r)_F}
}$$
is cocartesian. By lemma \ref{lem_much-ado}, the set
$f^{-1}(\eta(r))$ (resp. $f^{\prime-1}(\eta(r)_F)$) is in natural
bijection with the set of valuations on $\cB(r)$ (resp. $\cB(r)_F$)
that extend $|\cdot|_{\eta(r)}$ (resp. $|\cdot|_{\eta(r)_F}$).
It then follows by standard valuation theory
(see \cite[Ch.VI, \S2, Exerc.2]{BouAC}) that the map
$f^{\prime-1}(\eta(r)_F)\to f^{-1}(\eta(r))$ is a surjection;
let $N$ (resp. $N'$) be the cardinality of $f^{-1}(\eta(r))$
(resp. $f^{\prime-1}(\eta(r)_F)$).
Moreover, $\cB(r)_F$ is reduced, by \cite[Lemma 3.3.1.(1)]{Conr};
let $r$ be its rank over $\kappa(r)_F$, which is also the rank
of $\cB(r)$ over $\kappa(r)$.
Since $\kappa(\eta(r))$ and $\kappa(\eta(r)_F)$ are defectless
in every finite separable extension (\cite[Lemma 5.3(iii)]{Hu3}),
we deduce that :
$$
r=N\cdot(\Gamma_x:\Gamma_{\eta(r)})=
N'\cdot(\Gamma_{x'}:\Gamma_{\eta(r)_F}).
$$
We know already that $N'\geq N$, and since
$\Gamma_{\eta(r)}=\Gamma_{\eta(r)_F}$, it is also clear that
$(\Gamma_{x'}:\Gamma_{\eta(r)_F})\geq(\Gamma_x:\Gamma_{\eta(r)})$,
whence (i). Next, in light of (i), for every $\sigma\in St_{x'}$
we may compute $i_{x'}(\sigma)$ as $|\sigma(t)-t|_{x'}$, where
$t\in\kappa(x)$ is any element such that $|t|_x$ is the largest
element in $\Gamma^+_x\setminus\{1\}$.
Assertion (ii) is then an immediate consequence.
\end{proof}

\begin{lemma}\label{lem_Fuhrer} For every subgroup $H\subset G$,
denote by $f_H:X/H\to\D(a,a^{-1})$ the morphism deduced from $f$.
The following identities hold:
$$
\delta_{f_H}(-\log r)=\sa^\flat_G(\C[G/H],r^+)
\qquad\text{for every $r\in(a,b]\cap\Gamma_K$}
$$
and
$$
\frac{d\delta_{f_H}}{dt}(-\log r^+)=\ssw^\natural_G(\C[G/H],r^+)
\qquad\text{for every $r\in(a,b]\cap\Gamma_K$}.
$$
\end{lemma}
\begin{proof} (Here $\C[G/H]=\mathrm{Ind}^G_H 1_H$, where
$1_H$ is the trivial character of $H$.) First of all, one
applies \eqref{eq_diff-char} to derive, as in
\cite[Ch.VI, \S3, Cor.1]{Se} that
$$
|\fd_{f_H}^+(r)|^\natural_{\eta(r)}=
\langle\sa^\natural_G(r^+),\mathrm{Ind}^G_H 1_H\rangle_G
\qquad\text{and}\qquad
-\log|\fd_{f_H}^+(r)|^\flat_{\eta(r)}=
\langle\sa^\flat_G(r^+),\mathrm{Ind}^G_H 1_H\rangle_G
$$
(notation of \eqref{subsec_discr-at-r}), which already implies
the first of the sought identities.
Moreover we deduce:
$$
\begin{array}{r@{\:=\:}l}
\langle\ssw^\natural_G(r^+),\mathrm{Ind}^G_H 1_H\rangle_G &
|\fd_{f_H}^+(r)|^\natural_{\eta(r)}-
\langle\mathrm{Ind}^G_{St_x}\su_{St_x},\mathrm{Ind}^G_H 1_H\rangle_G \\
& |\fd_{f_H}^+(r)|^\natural_{\eta(r)}+\sharp(H\backslash G/St_x)-(G:H).
\end{array}
$$
which is equivalent to the second stated identity, in view
of lemma \ref{lem_local-sit}(iii) and proposition \ref{prop_precise}.
\end{proof}

\sset\subsubsection{}\label{subsec_was-a-rem}
Of course, one can also repeat the same discussion with the
point $\eta'(r)$ instead of $\eta(r)$ (notation of \eqref{subsec_localize});
then one obtains characters $\sa^\natural_G(r^-)$,
$\sa^\flat_G(r^-)$, $\ssw^\natural_G(r^-)$ and -- in view of
example \eqref{ex_asymm} -- the identity:
\set\begin{equation}\label{eq_was-a-rem}
-\frac{d\delta_f}{dt}(-\log r^-)=
\langle\ssw^\natural_G(r^-),\sreg_G\rangle_G
\qquad\text{for every $r\in[a,b)\cap\Gamma_K$}.
\end{equation}
Moreover we have :

\begin{lemma}\label{lem_missing} $\sa^\flat_G(r^-)=\sa^\flat_G(r^+)$.
\end{lemma}
\begin{proof} In view of \eqref{eq_double-wham} and its analogue
for $\sa^\flat_G(r^-)$, we see that both sides of the sought
identity are elements of $\R\otimes_\Z R_\Q(G)$ (notation of
\cite[\S12.1]{Se2}). By \cite[\S13.1, Th.30]{Se2}, it then
suffices to check that
$\sa^\flat_G(\C[G/H],r^-)=\sa^\flat_G(\C[G/H],r^+)$ for every
(cyclic) subgroup $H\subset G$. The latter is clear, in light of
lemma \ref{lem_Fuhrer}.
\end{proof}

\begin{proposition}\label{prop_Fuehrer}
For every $\chi\in K_0(\C[G])$ the function
$$
\delta_{f,\chi}:(\log 1/b,\log 1/a]\cap\log\Gamma_K\to\R
\quad :\quad
-\log r\mapsto\sa^\flat_G(\chi,r^+)
$$
is the restriction of a piecewise linear continuous function
defined on $(\log 1/b,\log 1/a]$, and the following identity holds :
$$
\frac{d\delta_{f,\chi}}{dt}(-\log r^+)=
\ssw_G^\natural(\chi,r^+)
\qquad\text{for every $r\in(a,b]\cap\Gamma_K$}.
$$
\end{proposition}
\begin{proof} Due to Artin's theorem \cite[\S9.2, Cor.]{Se2}
we may assume that $\chi$ is induced from a character
$\rho:H\to\C^\times$ of a cyclic subgroup $H\subset G$. 
Denote by $\sa^\flat_{|H}$ (resp. $\ssw_{|H}^\natural$) the
restriction to $H$ of $\sa^\flat_G$ (resp. of $\ssw_G^\natural$);
by Frobenius reciprocity, we are reduced to showing that
the map :
$$
-\log r\mapsto\langle\sa^\flat_{|H}(r^+),\rho\rangle_H
$$
extends to a piecewise linear and continuous function, with
right slope $\langle\sa^\natural_{|H}(r^+),\rho\rangle_H$.
Let $k$ be the order of $\rho$, {\em i.e.} the smallest integer
such that $\rho^k=1_H$; we shall argue by induction on $k$.
For $k=1$, $\rho$ is the trivial character, and then the
assertion follows from lemma \ref{lem_Fuhrer}.
Hence, suppose that $k>1$ and that the assertion is known
for all the characters of $H$ whose order is strictly smaller
than $k$. There exists a unique subgroup $L\subset H$ with
$(H:L)=k$, and $\C[H/L]\subset \C[H]$ is the direct sum
of all characters of $H$ whose orders divide $k$. By lemma \ref{lem_Fuhrer}
the sought assertion is known for this direct sum of
characters, and then our inductive assumption implies that
the assertion is also known for the sum
$\rho':=\rho_1\oplus\cdots\oplus\rho_n$ of all characters
of $H$ whose order equals $k$. The latter are permuted under the
natural action of $(\Z/k\Z)^\times$ on $K_0(\C[H])$
(cp. \cite[\S9.1, Exerc.3]{Se2}). Moreover, for any
$j\in(\Z/k\Z)^\times$ let $\spsi_j:K_0(\C[H])\to K_0(\C[H])$
be the corresponding operator.
\begin{claim}\label{cl_spsi} $\sa^\flat_{|H}=\spsi_j(\sa^\flat_{|H})$
for every $j\in(\Z/k\Z)^\times$.
\end{claim}
\begin{pfclaim} From \cite[Lemma 2.6]{Hu3}
it follows that $\sa_x(\sigma)=\sa_x(\tau)$ whenever $\sigma$
and $\tau$ generate the same subgroup of $St_x$.
The claim is a direct consequence.
\end{pfclaim}

Using claim \ref{cl_spsi} we compute :
$$
\langle\sa^\flat_{|H}(r^+),\rho\rangle=
\langle\spsi_j(\sa^\flat_{|H}(r^+)),\spsi_j(\rho)\rangle=
\langle\sa^\flat_{|H}(r^+),\spsi_j(\rho)\rangle
$$
for every $r\in(a,b]\cap\Gamma_K$ and every $j\in(\Z/k\Z)^\times$. Thus :
$$
\sa_G^\flat(\chi,r^+)=\frac{\sa^\flat_G(\rho,r^+)}{o((\Z/k\Z)^\times)}.
$$
A similar argument yields a corresponding identity for
$\ssw_G^\natural(\chi,r^+)$, and concludes the proof of
the proposition.
\end{proof}

The following is the main result of this chapter.

\begin{theorem}\label{th_Fuehrer} Suppose that $b=a^{-1}$,
so that we can define the complex of $\Z_\ell[G]$-modules
$\Delta^\bullet_\infty(X)$ as in \eqref{subsec_go-to-adic}.
Then we have the identity :
$$
[\Q\otimes_\Z\Delta^\bullet_\infty(X)[1]]=
\ssw^\natural_G(1^+)+\ssw^\natural_G(1^-).
$$
\end{theorem}
\begin{proof} We begin with the following:
\begin{claim}\label{cl_H-invariants} For every abelian
subgroup $H\subset G$ we have a natural identification:
$$
\Delta^\bullet_\infty(X)^H=\Delta^\bullet_\infty(X/H).
$$
\end{claim}
\begin{pfclaim} Indeed, let $\fp_1$,\dots,$\fp_k$ be the prime
ideals of $(B^\circ)^H$ lying over $\fP$ (notation of
\eqref{subsec_sans-serif}); for every $j:=1,\dots,k$ let
$S_j$ be the set of the prime ideals of $B^\circ$ lying
over $\fp_j$. Clearly $\bigcup_{j=1}^kS_j=\{\fq_1,\dots,\fq_n\}$,
and for every $j\leq k$ the subgroup $H$ stabilizes the
direct sum
$$
\Delta^\bullet_\infty(X,S_j):=
\bigoplus_{\fq\in S_j}\Delta^\bullet_\infty(X,\fq).
$$
Let $L_j\subset H$ be the stabilizer of any (hence each)
element of $S_j$; by lemma \ref{lem_H-action} we have
a natural identification:
$$
\Delta^\bullet_\infty(X,S_j)^{L_j}=
\bigoplus_{\fq\in S_j}\Delta^\bullet_\infty(X/L_j,\fq^{L_j})=
\bigoplus_{\fq\in S_j}\Delta^\bullet_\infty(X/H,\fp_j).
$$
However, the quotient $H/L_j$ permutes transitively the
summands of $\Delta^\bullet_\infty(X,S_j)^{L_j}$,
whence the claim.
\end{pfclaim}

\begin{claim}\label{cl_for-cyclic-gps}
For every cyclic subgroup $H\subset G$ we have:
\set\begin{equation}\label{eq_equal-representation}
\langle[\Q\otimes_\Z\Delta^\bullet_\infty(X)[1]],\C[G/H]\rangle_G=
\langle \ssw^\natural_G(1^+)+\ssw^\natural_G(1^-),\C[G/H]\rangle_G.
\end{equation}
\end{claim}
\begin{pfclaim} We use the morphism $h:=\Spa\,\psi\circ f$
defined in example \ref{ex_summing-up}. For $H\subset G$ any
cyclic subgroup, set $h_H:=\Spa\,\psi\circ f_H$. In view of
proposition \ref{prop_vanish} and claim \ref{cl_H-invariants},
and arguing as in the proof of theorem \ref{th_big-deal}, we
see that the left-hand side of \eqref{eq_equal-representation}
computes $(-1)$ times the left slope at the point $0$
of the discriminant function of $h_H$. But combining lemma
\ref{lem_Fuhrer}, \eqref{eq_was-a-rem} and example \ref{ex_asymm}
we conclude that also the right-hand side admits the same
interpretation.
\end{pfclaim}

Since $\ssw^\sharp_G(1^+)+\ssw^\sharp_G(1^-)$ is a rational
valued function on $G$, the theorem follows from claim
\ref{cl_for-cyclic-gps}, theorem \ref{th_rationality}(ii) and
\cite[\S13.1, Th.30]{Se2}.
\end{proof}

\begin{corollary}\label{cor_Fuehrer} Keep the notation of proposition
{\em\ref{prop_Fuehrer}}. For every complex representation
$\rho$ of $G$, the function $\delta_{f,\rho}$ is the restriction
of a non-negative, piecewise linear, continuous, and convex
real-valued function defined on $(\log 1/b,\log 1/a]$, with integer
slopes.
\end{corollary}
\begin{proof} Continuity and piecewise linearity are already
known from proposition \ref{prop_Fuehrer}, and
$\delta_{f,\rho}$ takes values in $\R_{\geq 0}$, by lemma
\ref{lem_positivity}(i),(ii).
In view of \cite[Th.4.1]{Hu3}, proposition \ref{prop_Fuehrer}
also implies that the slopes of $\delta_{f,\rho}$ are integer.
Finally, let $r\in (a,b]\cap\Gamma_K$ be a radius such
that the left and right slope of $\delta_{f,\rho}$ are
different at the value $-\log r$; up to restricting the covering
$f$, and rescaling the coordinate, we may assume that $r=1$
and $b=a^{-1}$, in which case theorem \ref{th_Fuehrer}, lemma
\ref{lem_missing} and \eqref{eq_was-a-rem} yield the identity :
\set\begin{equation}\label{eq_diff-slopes}
\langle[\Q\otimes_\Z\Delta^\bullet_\infty(X)][1],\rho\rangle_G=
\frac{d\delta_{f,\rho}}{dt}(0^+)-\frac{d\delta_{f,\rho}}{dt}(0^-).
\end{equation}
However, lemma \ref{lem_go-to-adic} implies that the left-hand
side of \eqref{eq_diff-slopes} is a non-negative integer for
any representation $\rho$, whence the contention.
\end{proof}

\sset\subsubsection{}\label{subsec_modular}
Next we consider modular representations of $G$. Namely, let
$\Lambda$ be a complete discrete valuation ring with residue
field $\bar\Lambda$ of positive characteristic $\ell\neq p$,
and field of fractions $\Lambda_\Q:=\Lambda[1/\ell]$ of characteristic
zero. We assume that we are also given a fixed imbedding
of $\Lambda_\Q$ into the field of complex numbers :
\set\begin{equation}\label{eq_fixed-imbed}
\Lambda_\Q\hookrightarrow\C.
\end{equation}
As usual (cp. \eqref{subsec_not-necess}), for any group $H$
we denote by $K_0(\bar\Lambda[H])$ (resp. by $K^0(\bar\Lambda[H])$)
the Grothendieck group of the category of $\bar\Lambda[H]$-modules
of finite rank over $\bar\Lambda$ (resp. of projective
$\bar\Lambda[H]$-modules of finite rank).
We shall also consider $K_0(\Lambda_\Q[H])=K^0(\Lambda_Q[H])$.
The tensor product (over $\Lambda$) induces a ring structure
on these groups, and according to \cite[\S15.5, Prop.43]{Se2}
there is a commutative diagram of ring homomorphisms :
$$
\xymatrix{
K^0(\bar\Lambda[H]) \ar[rr]^-{c_H} \ar[dr]_{e_H} & & K_0(\bar\Lambda[H]) \\
& K_0(\Lambda_\Q[H]) \ar[ru]_{d_H}
}$$
such that $d_H$ and $e_H$ are adjoint maps for the natural
bilinear pairing
$$
\langle\cdot,\cdot\rangle_H:
K^0(\bar\Lambda[H])\times K_0(\bar\Lambda[H])\to\Z
$$
{\em i.e.} we have the identity
(\cite[\S15.4]{Se2}) :
\set\begin{equation}\label{eq_adjoints}
\langle\rho,d_H(\chi)\rangle_H=\langle e_H(\rho),\chi\rangle_H
\qquad\text{for every $\chi\in K_0(\Lambda_\Q[H])$ and
$\rho\in K^0(\bar\Lambda[H])$.}
\end{equation}

\sset\subsubsection{}\label{subsec_produce}
For instance, for every $r\in(a,b]\cap\Gamma_K$ there exists a unique
element $\bar\ssw{}^\natural_G(r^+)\in K^0(\bar\Lambda[G])$
such that $\ssw^\natural_G(r^+)=e_G(\bar\ssw{}^\natural_G(r^+))$.
Likewise, by inspecting \eqref{eq_double-wham} we see that
there exist elements
$$
\bar\sa{}^\flat_G(r^+)\in\R\otimes_\Z K^0(\bar\Lambda[G])
\qquad
\bar\sa{}^\flat_{St^{(p)}_x}(r^+)\in\R\otimes_\Z K^0(\bar\Lambda[St^{(p)}_x])
$$
such that :
$$
\sa^\flat_{St^{(p)}_x}(r^+)=e_{St^{(p)}_x}(\bar\sa{}^\flat_{St^{(p)}_x}(r^+))
\qquad
\sa^\flat_G(r^+)=e_G(\bar\sa{}^\flat_G(r^+))
\qquad
\bar\sa{}^\flat_G(r^+)=\Ind^G_{St^{(p)}_x}\bar\sa{}^\flat_{St^{(p)}_x}(r^+).
$$
Now, let $\bar\chi\in K_0(\bar\Lambda[G])$ be any element;
we define the function :
$$
\delta_{f,\bar\chi}:(\log 1/b,\log 1/a]\cap\log\Gamma_K\to\C
\qquad -\log r\mapsto\langle\bar\sa{}^\flat_G(r^+),\bar\chi\rangle_G.
$$

\begin{proposition}\label{prop_modular}
If $\bar\chi$ is the class of a $\bar\Lambda$-linear
representation ({\em i.e.} a positive element of
$K_0(\bar\Lambda[G])$), then $\delta_{f,\bar\chi}$ is
the restriction of a non-negative, piecewise linear,
convex and continuous real-valued function defined on
$(\log 1/b,\log 1/a]$, and moreover :
\set\begin{equation}\label{eq_and-moreover}
\frac{d\delta_{f,\bar\chi}}{dt}(-\log r^+)=
\langle\bar\ssw{}^\natural_G(r^+),\bar\chi\rangle_G
\end{equation}
\end{proposition}
\begin{proof} According to \cite[\S16.1, Th.33]{Se2}, we
may find $\chi\in K_0(\Lambda_\Q[G])\subset K_0(\C[G])$
such that $d_G(\chi)=\bar\chi$. For every $r\in(a,b]\cap\Gamma_K$,
we compute :
$$
\delta_{f,\bar\chi}(-\log r)=
\langle\bar\sa{}^\flat_G(r^+),d_G(\chi)\rangle_G=
\langle\sa^\flat_G(r^+),\chi\rangle_G
$$
and then piecewise linearity, continuity, as well as
\eqref{eq_and-moreover} follow from proposition \ref{prop_Fuehrer}.
Since $\mathrm{Res}^G_{St^{(p)}_x}(\chi)=
d_{St^{(p)}_x}^{-1}\circ\mathrm{Res}^G_{St^{(p)}_x}(\bar\chi)$
is a positive element of $K_0(\C[St^{(p)}_x])$, lemma
\ref{lem_positivity}(i),(ii) implies that
$\delta_{f,\bar\chi}(-\log r)\in[0,+\infty)$.

As for the convexity, notice that we cannot apply directly corollary
\ref{cor_Fuehrer}, since the element $\chi\in K_0(\C[G])$
may fail to be positive.
Nevertheless, the argument proceeds along the same lines :
let $r\in(a,b]\cap\Gamma_K$ be a radius such that the
left and right slopes of $\delta_{f,\bar\chi}$ are different
at the value $-\log r$; we may assume that $b=a^{-1}$ and $r=1$,
in which case the class
$[\Delta^\bullet_\infty(X)[1]]\in K^0(\Z_\ell[G])$ is
well defined and positive. Clearly we have :
$$
[\Q\otimes_\Z\Delta^\bullet_\infty(X)[1]]=
e_G([\Delta^\bullet(X,\bar\Delta)[1]])
$$
(notation of \eqref{subsec_sans-serif}). Then we compute
using theorem \ref{th_Fuehrer}, lemma \ref{lem_missing} and
\eqref{eq_was-a-rem} :
\set\begin{equation}\label{eq_positive-again}
\frac{d\delta_{f,\bar\chi}}{dt}(0^+)-\frac{d\delta_{f,\bar\chi}}{dt}(0^-)=
\langle[\Q\otimes_\Z\Delta^\bullet_\infty(X)][1],\chi\rangle_G=
\langle[\Delta^\bullet(X,\bar\Delta)[1]],\bar\chi\rangle_G
\end{equation}
where again the last identity follows from \eqref{eq_adjoints}.
To conclude we apply proposition \ref{prop_perf-with-H}(iii)
to deduce the sought positivity of the right-most term in
\eqref{eq_positive-again}.
\end{proof}

Next, we wish to investigate the continuity properties of the
higher ramification filtration. These are gathered in the following :

\begin{theorem}\label{th_continuities}
In the situation of \eqref{sec_conductors}, there exists $r'\in(a,r)$
and, for every $s\in(r',r]\cap\Gamma_K$, a point
$x(s)\in f^{-1}(\eta(s))$ such that :
\begin{enumerate}
\item
The stabilizer $St_{x(s)}\subset G$ of\/ $x(s)$ under the natural
$G$-action on $f^{-1}(\eta(s))$, is a subgroup independent of $s$.
\item
The length of the higher ramification filtration
of $St_{x(s)}=\Gal(\kappa(x(s))^{\wedge h}/\kappa(s)^{\wedge h})$ :
$$
P_{\gamma_n(s)}\subset\cdots\subset P_{\gamma_1(s)}\subset St_{x(s)}^{(p)}
$$
is independent of $s\in(r',r]\cap\Gamma_K$.
\item
Set $\gamma_k^\natural:=\gamma_k(r)^\natural$\ \  for every
$k\leq n$. Then :
$$
\gamma_k(s)=(s/r)^{\gamma_k^\natural}\cdot\gamma_k(r)
\qquad\text{for every $s\in(r',r]\cap\Gamma_K$}.
$$
\end{enumerate}
\end{theorem}
\begin{proof} Choose $r'\in(a,r)\cap\Gamma_K$ such that the
discriminant function $\delta_f$ has constant left slope on
$[r',r]\cap\Gamma_K$. Up to rescaling the coordinates, we
may assume that $r=a$ and $r'=a^{-1}$ with $a:=|\pi|$ for some
$\pi\in\fm$. Then $\D(r',r)=\D(a,a^{-1})=\Spa\,A(a,a^{-1})$, and
$A^\circ:=A(a,a^{-1})^\circ=K^+\langle S,T\rangle/(S\cdot T-\pi^2)$.
Also, there is no harm in replacing $X$ by its restriction
to $f^{-1}(\D(r',r))$. Let $B$ be an affinoid $K$-algebra
such that $X=\Spa\,B$, and denote by $\fP\subset A^\circ$
the unique maximal ideal such that $S,T\in\fP$. Let
$\fq_1,\dots,\fq_k\subset B^\circ$ be the maximal ideals
lying over $\fP$. 
\begin{claim}\label{cl_one-pt}
For every $s\in(a,a^{-1}]\cap\Gamma_K$ and every $i\leq k$
there is exactly one point $x_i(s)\in f^{-1}(\eta(s))$ such
that $\kappa(x_i(s))^+$ dominates $B^\circ_{\fq_i}$.
\end{claim}
\begin{pfclaim} One argues as in the proof of theorem
\ref{th_Riemann}, using \eqref{eq_crossing} : the
details shall be left to the reader.
\end{pfclaim}

By standard arguments we may find an \'etale ring homomorphism
$\phi:A^\circ\to C$ such that :

\begin{itemize}
\item
$\Spec\,C$ is connected, $\fP':=\fP C$ is a maximal ideal of $C$,
and the induced map $A^\circ/\fP\to C/\fP'$ is an isomorphism.
\item
$\Spec\,C\otimes_{A^\circ}B^\circ$ is the union of $k$ connected
components, which are therefore in natural bijection with the set
$\{\fq_1,\dots,\fq_k\}$.
\end{itemize}
Let $C^\wedge$ be the $\pi$-adic completion of $C$; then
$C^\wedge_K:=C^\wedge\otimes_{K^+}K$ is an affinoid $K$-algebra,
and the induced morphism of affinoid spaces
$g:Y:=\Spa\,C^\wedge_K\to\D(a,a^{-1})$ is \'etale
(\cite[Cor.1.7.3(iii)]{Hu2}).

\begin{claim}\label{cl_two-pt}
(i)\ \ $(C^\wedge_K)^\circ=C^\wedge$.
\begin{enumerate}
\addenu
\item
$g$ restricts to an isomorphism
$Y_t:=g^{-1}(\D(a/|t|,|t|/a))\isom\D(a/|t|,|t|/a)$,
for every $t\in\fm\setminus\{0\}$.
\item
For every $s\in(a,a^{-1}]\cap\Gamma_K$, the preimage
$g^{-1}(\eta(s))$ consists of one point.
\end{enumerate}
\end{claim}
\begin{pfclaim} (i): Since $C$ is $K^+$-flat, the same holds
for $C^\wedge$ (cp. the proof of lemma \ref{lem_future});
moreover $C/\fm C$ is reduced, since it is \'etale over
the reduced ring $A^\sim$ (notation of \eqref{subsec_conflict}).
Then the assertion follows from lemma \ref{lem_reduction}.

Next, denote by $A^{\circ h}_\fP$ (resp. $C_{\fP'}^h$)
the henselization of $A^\circ$ (resp. of $C$) along $\fP$
(resp. along $\fP'$); for every $t\in\fm\setminus\{0\}$, the
map $A^\circ\to A_t^\circ:=A(a/|t|,|t|/a)^\circ$ (dual to the
open immersion $\D(a/|t|,|t|/a)\subset\D(a,a^{-1})$) factors
through a natural morphism $A^{\circ h}_\fP\to A^\circ_t$.
Since $\phi$ induces an isomorphism $A^{\circ h}_\fP\isom C_{\fP'}^h$,
assertion (ii) follows easily. This implies already assertion
(iii), for every $s\in(a,a^{-1})\cap\Gamma_K$. It remains to
show that there exists exactly one point in $y\in Y$ lying over
$\eta(a^{-1})$; in view of (i), this is the same as showing
that $C\otimes_{K^+}K$ admits exactly one continuous valuation
$|\cdot|_y:C\otimes_{K^+}K\to\Gamma_y\cup\{0\}$ such that
$$
|t|_y=|t|_{\eta(a^{-1})}\quad\text{for every $t\in A^\circ$}
\qquad\text{and}\qquad|t|_y\leq 1\quad\text{for every $t\in C$.}
$$
These valuations are in natural bijection with the set
of valuations $|\cdot|_{\bar y}$ on $C/\fm C$ such that
$|\bar t|_{\bar y}=|t|^\natural_{\eta(a^{-1})}$ for the
class $\bar t\in A^\sim$ of every element
$t\in A^\circ\setminus\fm A^\circ$.
But since $C/\fm C$ is \'etale over $A^\sim$ and $\fP'$
is a prime ideal, there is exactly one such valuation
$|\cdot|_{\bar y}$.
\end{pfclaim}

Let $Z:=X\times_{\D(a,a^{-1})}Y$; then
$Z=\Spa\,C^\wedge_K\otimes_{A^\circ}B^\circ$. Set
$D:=C^\wedge\otimes_{A^\circ}B^\circ$.

\begin{claim}\label{cl_identify} $\cO^+_Z(Z)=D$ and $D$ is henselian
along its ideal $\fm D$.
\end{claim}
\begin{pfclaim} Obviously the second assertion follows from
the first. Since $B^\circ$ is a finitely presented $A^\circ$-module
(lemma \ref{lem_finitely-pres}), $D$ is the $\pi$-adic completion of
$C\otimes_{A^\circ}B^\circ$, and since the latter is flat
over $K^+$, the same holds for $D$ (cp. the proof of lemma
\ref{lem_future}). Furthermore,
$D/\fm D=(C/\fm C)\otimes_{A^\sim}B^\sim$. By lemma
\ref{lem_reduction}, the ring $B^\sim$ is reduced, hence
the same holds for $D/\fm D$, and then the claim follows
by a second application of lemma \ref{lem_reduction}.
\end{pfclaim}

It follows from claim \ref{cl_identify} and the construction
of $C$, that $Z$ splits as the disjoint union of $k$ open and
closed subsets $Z_{\fq_i}$ ($i=1,\dots,k$), in natural bijection
with the set $\Sigma:=\{\fq_1,\dots,\fq_k\}$, and such that
precisely one prime ideal of $\cO^+_Z(Z_{\fq_i})$ lies over
$\fP'$, for every $i\leq k$. Choose arbitrarily $\fq\in\Sigma$,
and let $Z_\fq\subset Z$ be the corresponding open and closed subset.
The restriction $h:Z_\fq\to Y$ of $f\times_{\D(a,a^{-1})}Y$
is a finite Galois \'etale morphism, whose Galois group is the
stabilizer $St(\fq)\subset G$ of $\fq$ for the natural $G$-action
on $\Sigma$. Furthermore, combining claims \ref{cl_one-pt},
\ref{cl_two-pt} we deduce that, for every $s\in(a,a^{-1}]\cap\Gamma_K$
there exists precisely one point $z(s)\in Z_\fq$ such that
$\{h(z(s))\}=g^{-1}(\eta(s))$; set $y(s):=h(z(s))$ for every
such $s$. Let $g':Z_\fq\to X$ be the restriction of
$X\times_{\D(a,b)}g$, and set $x(s):=g'(z(s))$ for every
$s\in(a,a^{-1}]\cap\Gamma_K$. Since $g'$ is equivariant
for the action of $St(\fq)$, we see that this family
$(x(s)~|~s\in(a,a^{-1}]\cap\Gamma_K)$ fulfills condition (i);
indeed we have $St_{x(s)}=St(\fq)$ for each point $x(s)$.

\begin{claim}\label{cl_equiv-ident}
For every $s\in(a,a^{-1}]\cap\Gamma_K$, the following holds :
\begin{enumerate}
\item
The natural maps $\kappa(s)\to\kappa(y(s))$ and
$\kappa(x(s))\to\kappa(z(s))$ induce $St(\fq)$-equivariant
identifications :
$$
\kappa(s)^{\wedge h+}\isom\kappa(y(s))^{\wedge h+}
\qquad\text{and}\qquad
\kappa(x(s))^{\wedge h+}\isom\kappa(z(s))^{\wedge h+}
$$
\item
The map $\kappa(y(s))\to\kappa(z(s))$ is a finite
field extension, and $\kappa(z(s))=(h_*\cO_{Z_\fq})_{y(s)}$.
\end{enumerate}
\end{claim}
\begin{pfclaim} (i): For $s<a^{-1}$ this is already clear from
claim \ref{cl_two-pt}(ii); indeed in this case one does not
need to complete, nor to henselize, to obtain isomorphisms.
For the case where $s=a^{-1}$, set
$D:=\kappa(a^{-1})^{\wedge h+}\otimes_{A^\circ}C$, and let
$D^\nu\subset D$ be the normalization of
$\kappa(a^{-1})^{\wedge h+}$ in $D$. Then $\Spec\,D$ and
$\Spec\,D^\nu$ have the same number of irreducible components;
moreover $D^\nu$ is a product of finitely many valuation rings,
each of which is a finite extension of $\kappa(a^{-1})^{\wedge h+}$.
Since $C$ admits only one valuation extending $|\cdot|_{\eta(a^{-1})}$,
we deduce that precisely one of the irreducible components of
$\Spec\,D$ is finite over $\kappa(a^{-1})^{\wedge h+}$.
Let $V$ be the direct factor of $D$ corresponding to this
irreducible component; then $V\cap\mathrm{Frac}(C)$ is the valuation
ring of the field of fractions $\mathrm{Frac}(C)$, corresponding to
$y(a^{-1})\in Y$. On the one hand, the extension
$\kappa(a^{-1})^{\wedge h+}\to V$ is \'etale, hence an isomorphism;
on the other hand, the natural map $V\to\kappa(y(a^{-1}))^{\wedge h+}$
induces a bijection on value groups, so the map
$\kappa(a^{-1})^{\wedge h+}\to\kappa(y(a^{-1}))^{\wedge h+}$
is an isomorphism, as stated. The remaining assertion concerning
$\kappa(x(a^{-1}))^{\wedge h+}$ is an easy consequence; the
details shall be left to the reader.

(ii): It suffices to show that $z(s)$ admits a cofinal system
of open neighborhoods $U\subset Z_\fq$ such that the restriction
$h_{|U}:U\to h(U)$ is finite. However, let $V\subset Z_\fq$ be
any open neighborhood of $z(s)$; since $\{z(s)\}=h^{-1}(h(z(s)))$,
the subset $U:=h^{-1}(Y\setminus h(Z_\fq\setminus V))\subset V$ is
also an open neighborhood of $z(s)$, and clearly $h_{|U}$
is finite.
\end{pfclaim}

\begin{claim}\label{cl_once-again}
Let $U\subset Z_\fq$ be an open neighborhood of $z(a^{-1})$,
$t\in\cO_Z(U)$ any section, and set $k:=|t|^\natural_{z(a^{-1})}$.
Then there exists $r\in(a,a^{-1})$ such that :
\begin{enumerate}
\item
$z(s)\in U$ for every $s\in[r,a^{-1}]\cap\Gamma_K$.
\item
$|t|_{z(s)}=|t|_{z(a^{-1})}\cdot(s\cdot a)^k$\ \ 
for every $s\in[r,a^{-1}]\cap\Gamma_K$.
\end{enumerate}
\end{claim}
\begin{pfclaim} Notice that claim \ref{cl_equiv-ident}(ii), and lemma
\ref{lem_nobigdeal}(iii) imply that the spectral norm induced
on the $\kappa(y(s))$-algebra $\kappa(z(s))$ agrees with the
valuation $|\cdot|_{z(s)}$, for every $s\in(a,a^{-1}]\cap\Gamma_K$.
Then we may argue as in the proof of lemma \ref{lem_estimate}, to
reduce the claim to the special case where $Z_\fq=Y$ and $h$ is
the identity map, therefore $U\subset Y$ is an open neighborhood
of $y(s)$ and $t\in\cO_Y(U)$. In such case, in view of claim
\ref{cl_equiv-ident}(i), we may find $t'\in\kappa(a^{-1})$
such that $|t/t'|_{y(a^{-1})}=1$.
Then $t/t'$ extends to a section of $\cO^+_Y$ over an open subset
$U'\subset U$ containing $y(a^{-1})$, and up to restricting
$U'$ we may assume that $t/t'$ is invertible in $\cO^+_Y(U')$.
Applying again claim \ref{cl_equiv-ident}(i) we see that
$|t|_{y(s)}=|t'|_{\eta(s)}$ for every $s\in[s,a^{-1}]\cap\Gamma_K$
such that $y(s)\in U'$, and then both assertions follow from
lemma \ref{lem_estimate}.
\end{pfclaim}

Let $d:=o(St_\fq)$, and pick $t\in\kappa(z(a^{-1}))$ such that
$$
|t|_{z(a^{-1})}=(1-\eps)^{1/d}.
$$
By claim \ref{cl_once-again}, there exists an open neighborhood
$U\subset Z_\fq$ of $z(a^{-1})$ and $r\in(a,a^{-1})$ such that
$t\in\cO_Z(U)$, $z(s)\in U$ for every $s\in[r,a^{-1}]\cap\Gamma_K$
and :
\set\begin{equation}\label{eq_rate}
|t|_{z(s)}=|t|_{z(a^{-1})}\cdot(s\cdot a)^{1/d}=
(1-\eps)^{1/d}\cdot(s\cdot a)^{1/d}
\qquad\text{for every $s\in[r,a^{-1}]\cap\Gamma_K$}.
\end{equation}
Let $\sigma\in St(\fq)$ be any element; by definition, we have :
$$
i_{z(a^{-1})}(\sigma)=|t-\sigma(t)|_{z(a^{-1})}.
$$
Now, let $s\in[r,a^{-1}]\cap\Gamma_K$, and choose $c\in K^\times$
such that $|c|=(s\cdot a)^{1/d}$; in view of \eqref{eq_rate}, we also
obtain the identity :
\set\begin{equation}\label{eq_rate-i}
i_{z(s)}(\sigma)=|c^{-1}\cdot t-\sigma(c^{-1}\cdot t)|_{z(s)}=
(s\cdot a)^{-1/d}\cdot|t-\sigma(t)|_{z(s)}.
\end{equation}
Set $k:=i_{z(a^{-1})}(\sigma)^\natural$; if we apply claim
\ref{cl_once-again} to $t-\sigma(t)$, we can rewrite \eqref{eq_rate-i}
in the form :
$$
i_{z(s)}(\sigma)=(s\cdot a)^{k-1/d}\cdot|t-\sigma(t)|_{z(a^{-1})}=
(s\cdot a)^{k-1/d}\cdot i_{z(a^{-1})}(\sigma)
$$
which easily yields assertions (ii) and (iii).
\end{proof}

\sset\subsubsection{}\label{subsec_continuities}
Clearly, the analogue of theorem \ref{th_continuities} holds
also for the fibres over the points $\eta'(s)$ (see
\eqref{subsec_was-a-rem}). Namely, suppose that $r\in[a,b)$;
then there exists $r'\in(r,a^{-1})$ and for every
$s\in[r,r')\cap\Gamma_K$ a point $x'(s)\in f^{-1}(\eta'(s))$
such that :
\begin{itemize}
\item
The stabilizer subgroup $St_{x'(s)}\subset G$ of\/ $x'(s)$
under the natural $G$-action on $f^{-1}(\eta'(s))$, is independent
of $s$.
\item
The length of the higher ramification filtration :
$$
P_{\beta_m(s)}\subset\cdots\subset P_{\beta_1(s)}\subset St_{x'(s)}^{(p)}
$$
of $St_{x'(s)}=\Gal(\kappa(x'(s))^{\wedge h}/\kappa'(s)^{\wedge h})$
is independent of $s\in[r,r')\cap\Gamma_K$.
\item
Set $\beta_k^\natural:=\beta_k(r)^\natural$\ \  for every
$k\leq m$. Then :
$$
\beta_k(s)=(s/r)^{-\beta_k^\natural}\cdot\beta_k(r)
\qquad\text{for every $s\in[r,r')\cap\Gamma_K$}.
$$
\end{itemize}
In other words, at the left (resp. at the right), of every
$r\in(a,b)\cap\Gamma_K$, the higher ramification filtrations
of the points lying over $\eta(s)$ (resp. $\eta'(s)$) for $s$
sufficiently close to $r$, change in a continuous -- indeed
linear -- fashion. To get the complete picture, we must also
analyze what happens when we switch from the left to the right
of a given radius $r$, {\em i.e.} we need to understand how
the filtrations $(P_{\gamma_i(r)}~|~i=1,\dots,n)$ and
$(P_{\beta_i(r)}~|~i=1,\dots,m)$ are related. The key is to
compare both ramification filtrations to a third one, attached
to the finite Galois extension $\kappa(r^\flat)\subset
\kappa(x(r)^\flat)=\kappa(x'(r)^\flat)$ (see \cite[Prop.1.5.4]{Hu2}).
To this aim, we make the following :

\begin{definition}\label{def_central}
Let $f:X\to\D(a,b)$ be a Galois finite \'etale
covering, with Galois group $G$, $x\in X$ a point of type (III)
and $x^\flat$ its unique proper generization (notation of
\eqref{subsec_classes}).
Let $St^\flat_x\subset G$ be the stabilizer subgroup of
the point $x^\flat$, under the natural action of $G$ on
$f^{-1}f(x^\flat)$. Then $St^\flat_x$ is naturally identified
with the Galois group
$\Gal(\kappa(x^\flat)^\wedge/\kappa(f(x^\flat))^\wedge)$.
For any given $c\in K^+\setminus\{0\}$ we set :
$$
P^\flat_\gamma:=
\Ker(St^\flat_x\to\Aut(\kappa(x^\flat)^{\wedge+}\otimes_{K^+}K^+/cK^+))
\qquad
\text{where $\gamma:=|c|$}.
$$
Clearly $P^\flat_\gamma$ is a normal subgroup of $St^\flat_x$
for every $\gamma\in\Gamma^+_K$. and the sequence
$(P^\flat_\gamma~|~\gamma\in\Gamma^+_K)$ is called the
{\em higher ramification filtration\/} of $St^\flat_x$.
\end{definition}

\begin{proposition}\label{prop_superjump}
Let $x\in X$ be as in \eqref{sec_conductors}, and let
$(P_\gamma~|~\gamma\in\Gamma^+_x)$ (resp.
$(P^\flat_\gamma~|~\gamma\in\Gamma^+_K)$) be the higher ramification
filtration of $St_x$ (resp. of $St^\flat_x$). Then
$St_x\subset St^\flat_x$ and moreover :
$$
P^\flat_\gamma=\bigcup_{n\in\Z}P_{\gamma_0^n\cdot\gamma}
\qquad\text{for every $\gamma\in\Gamma^+_K\setminus\{1\}$}.
$$
(Here $\gamma_0$ is the largest element of\/ $\Gamma^+_x\setminus\{1\}$.)
\end{proposition}
\begin{proof} The first assertion is obvious. Next, let
$\{x_1,\dots,x_k\}:=f^{-1}f(x)$ be the orbit of $x=x_1$ under the
$G$-action; in light of lemma \ref{lem_completing-kappa}, we see that
$\kappa(x_i)^\wedge=\kappa(x^\flat)^\wedge$ for every $i\leq k$,
and the rings $\kappa(x_i)^{\wedge+}$ are the valuation rings
of the field $\kappa(x^\flat)^\wedge$ that dominate
$\kappa(f(x))^{\wedge+}$. Let $C$ be the integral closure
of $\kappa(f(x))^{\wedge+}$ in $\kappa(x^\flat)^\wedge$;
$C$ is semilocal ring, whose localizations at the maximal
ideals are the valuation rings $\kappa(x_i)^{\wedge+}$;
applying \cite[Th.1.4]{Mat} we may then find $t\in C$ such that :
\set\begin{equation}\label{eq_ribo}
|t|^\wedge_x=1\qquad\text{and}\qquad
|t|^\wedge_{x_i}<1\quad\text{for every $i=2,\dots,k$}.
\end{equation}
Suppose now that $\sigma\in St^\flat_x\setminus St_x$; \eqref{eq_ribo}
implies easily that $|\sigma(t)-t|^\wedge_{x^\flat}=
(|\sigma(t)-t|^\wedge_x)^\flat=1$, hence
$\sigma\notin P^\flat_\gamma$ whenever $\gamma<1$, in other words :
\set\begin{equation}\label{eq_dueto}
P^\flat_\gamma\subset St_x\qquad\text{for every $\gamma<1$}.
\end{equation}
Thus, suppose $\sigma\in P^\flat_\gamma$ for some $\gamma<1$,
and choose $t\in\kappa(x)^\wedge$ such that $|t|_x=\gamma_0$;
due to \eqref{eq_dueto}, we have :
$i_x(\sigma)^\flat=|\sigma(t)-t|_x^\flat\leq\gamma$, whence
$P^\flat_\gamma\subset P'_\gamma:=\bigcup_{n\in\Z}P_{\gamma_0^n\cdot\gamma}$.

Conversely, suppose $\sigma\in P_\gamma$ for some $\gamma\in\Gamma^+_x$
such that $\gamma^\flat<1$; by definition, this means that
$|\sigma(t)-t|_x<\gamma$ for every $t\in\kappa(x)^{\wedge+}$
(\cite[Lemma 2.1(ii)]{Hu3}).
Let $s\in\kappa(x^\flat)^{\wedge+}$; then $c\cdot s\in\kappa(x)^{\wedge+}$
for every $c\in\fm$, hence :
$$
|c|\cdot|\sigma(s)-s|=|\sigma(c\cdot s)-c\cdot s|_x<\gamma
$$
therefore $|\sigma(c\cdot s)-c\cdot s|_x<\gamma\cdot|c|^{-1}$
for every $c\in\fm\setminus\{0\}$ and consequently
$|\sigma(s)-s|_x^\flat\leq\gamma^\flat$, {\em i.e.}
$\sigma\in P^\flat_{\gamma^\flat}$, which shows that
$P'_\gamma\subset P^\flat_\gamma$ for every
$\gamma\in\Gamma^+_K\setminus\{1\}$, as stated.
\end{proof}

\begin{remark}\label{rem_superjump}
If we apply proposition \ref{prop_superjump} to a pair of
points $x,x'\in X$ lying over $\eta(r)$ and respectively
$\eta'(r)$, and such that $x^\flat=x^{\prime\flat}$, we
see that both ramification filtrations ``to the left''
and ``to the right'' of the radius $r$ can be compared
to the same ``central'' ramification filtration for
$St^\flat_x=St^\flat_{x'}$. This expresses the sought
continuity property for the jumps of the ramification filtrations.
\end{remark}

\section{Local systems on the punctured disc}\label{chap_local-sys}

In this chapter we fix a complete and algebraically closed
valued field $(K,|\cdot|)$ of rank one and of zero characteristic.
We shall use the standard notation of \eqref{sec_annuli},
and we suppose that the characteristic of the residue field
$K^\sim$ is $p>0$.

\subsection{Break decomposition}\label{sec_break}
Let $\Lambda$ be an artinian local $\Z[1/p]$-algebra whose
residue field $\bar\Lambda$ has positive characteristic $\ell\neq p$;
we assume that the group $\Lambda^\times$ of invertible elements
of $\Lambda$ contains a subgroup isomorphic to 
$\mu_{p^\infty}:=\bigcup_{n>0}\mu_{p^n}$ (where $\mu_{p^n}$
denotes the group of $p^n$-th roots of one contained in $K^\times$),
and we fix such an imbedding :
\set\begin{equation}\label{eq_imbed-mu}
\mu_{p^\infty}\subset\Lambda^\times.
\end{equation}
Moreover, we shall also suppose that $\Lambda$ is the filtered
union of its finite subrings. This latter condition is motivated
by the following :

\begin{lemma}\label{lem_easy-down}
Let $X$ be a quasi-compact analytic adic space over $K$, $F$
a locally constant constructible $\Lambda$-module on the \'etale
site $X_\et$ of $X$. Then there exists a finite subring
$\Lambda'\subset\Lambda$ and a locally constant constructible
$\Lambda'$-module $F'$ on $X_\et$ such that
$F\simeq F'\otimes_{\Lambda}\Lambda'$.
\end{lemma}
\begin{proof} This lemma -- stated in the language of Berkovich's
non-archimedean analytic varieties -- appears in \cite[Lemma 4.1.8]{Ra}.
We sketch here the argument in the case of adic spaces. First,
using \cite[Lemma 1.4.7, Cor.1.7.4]{Hu2} we find finitely many
affinoid open subsets $U_1,\dots,U_n\subset X$ covering $X$,
and for each $i\leq n$ a finite \'etale morphism $f_i:V_i\to U_i$
such that $f^*_iF_{|U_i}$ is a constant $\Lambda$-module, whose
stalk (at some chosen geometric point of $V_i$) we denote $M_i$.
Then $V_i$ is an affinoid adic space for every $i\leq n$
(\cite[\S1.4.4]{Hu2}), hence the same holds for $W_i:=V_i\times_{U_i}V_i$,
especially the set $\pi_0(W_i)$ of connected components of $W_i$
is finite. Then the descent datum for $F_{|U_i}$ relative to the
morphism $f_i$ amounts to a finite set of $\Lambda$-automorphisms
$(\phi_{ij}~|~j\in\pi_0(W_i))$ of $M_i$, fulfilling a certain
cocycle condition. Since $M_i$ is of finite type and $\Lambda$
is noetherian, we may find a finite subring $\Lambda_i\subset\Lambda$,
a finite $\Lambda_i$-module $M'_i$ and a set of automorphisms
$(\phi'_{ij}~|~j\in\pi_0(W_i))$ such that
$M_i\simeq M'_i\otimes_{\Lambda_i}\Lambda$ and
$\phi_{ij}=\phi'_j\otimes_{\Lambda_i}\one_\Lambda$ for every
$i\leq n$ and $j\in\pi_0(W_i)$. Furthermore, after replacing
$\Lambda_i$ by some larger finite subring, we can achieve that
the cocycle conditions are still fulfilled by the system
$(\phi'_{ij}~|~j\in\pi_0(W_i))$; hence the latter furnishes
a descent datum for $M'_i$ relative to $f_i$, whence a
$\Lambda_i$-module $F'_i$ on $U_{i,\et}$ such that
$F'_i\otimes_{\Lambda_i}\Lambda\simeq F_{|U_i}$. Next, let
$U_{ij}:=U_i\cap U_j$ for every $i,j\leq n$, so that $F$
is defined by a cocycle system of isomorphisms
$(F'_i\otimes_{\Lambda_i}\Lambda)_{|U_{ij}}\isom
(F'_j\otimes_{\Lambda_j}\Lambda)_{|U_{ij}}$. Again, these
isomorphisms are already defined over some larger subring
$\Lambda_{ij}\supset\Lambda_i+\Lambda_j$, and the claim
follows easily.
\end{proof}

\begin{remark}\label{rem_easy-down}
(i)\  Keep the situation of lemma \ref{lem_easy-down};
an easy corollary is the following fact. There exists a finite
\'etale covering $f:Y\to X$ such that $f^*F$ is a constant
$\Lambda$-module on $Y_\et$. 

(ii)\ This is in general false, if $\Lambda$ is not a filtered
union of finite subrings : for a counterexample, consider an
elliptic curve $E$ over $K$, with bad reduction over $K^+$; it
is well known that the associated analytic space $E^\an$ can
be uniformized by the analytic torus $\G_{m,K}^\an$, and the
corresponding \'etale covering $\G_{m,K}^\an\to E^\an$ is Galois
with group $G\simeq\Z$. If we take $\Lambda:=\F_\ell(T)$, we may
define an action $\chi:G\to\End_\Lambda(L)$ on a one-dimensional
$\Lambda$-vector space $L$, by letting a generator $\sigma\in G$
act as multiplication by $T$. The character $\chi$ defines a
locally constant constructible $\Lambda$-module on $E^\an_\et$
that trivializes on $\G_{m,K}^\an$, but does not trivialize
on any finite covering of $E^\an$.
\end{remark}

\sset\subsubsection{}\label{subsec_Lambda-mod}
Let $F$ a constructible, locally constant and locally free sheaf of
$\Lambda$-modules on the \'etale site of the {\em punctured disc\/}
of radius one, {\em i.e.} the analytic adic space
$$
\D(1)^*:=\D(1)\setminus\{0\}
$$
where $0\in\D(1)$ is the closed point corresponding to the
valuation given by the rule : $f(\xi)\mapsto|f(0)|$ for every
$f\in A(1)$ (notation of \eqref{subsec_setup}). Following
\cite[\S8]{Hu3}, for every $r\in\Gamma_K^+$
one defines the Swan conductor of the stalk $F_{\eta(r)}$,
viewed as a representation of the Galois group of the algebraic
closure of the henselian field $\kappa(r)^{\wedge h}$ : 
\set\begin{equation}\label{eq_pione-r}
\pi_1(r):=\Gal((\kappa(r)^{\wedge h})^a:\kappa(r)^{\wedge h}).
\end{equation}
This conductor is a (possibly negative) integer which we
shall denote by :
$$
\ssw^\natural(F,r^+).
$$
(Huber denotes this quantity by $\alpha_x(F)$ in {\em loc.cit.},
with $x:=\eta(r)$.)

\sset\subsubsection{}
In this chapter we wish to investigate the properties of
the function
\set\begin{equation}\label{eq_log-swan}
-\log\Gamma^+_K\to\Z \quad:\quad
-\log r\mapsto\ssw^\natural(F,r^+).
\end{equation}
To start out, the identity (\cite[Lemma 8.1(iii)]{Hu3}) :
$$
\ssw^\natural(F,r^+)=\mathrm{length}_\Lambda\Lambda\cdot
\ssw^\natural(F\otimes_\Lambda\bar\Lambda,r^+)
$$
reduces the study of $\ssw^\natural(F,r^+)$ to that of
$\ssw^\natural(F\otimes_\Lambda\bar\Lambda,r)$. Now, for every
$r\in\Gamma^+_K$ let us fix a finite Galois \'etale
covering $f_r:X_r\to\D(r,1)$ such that $f_r^*F_{|\D(r,1)}$ is
a constant sheaf (since $\D(r,1)$ is a quasi-compact open subspace
of $\D(1)^*$, the existence of $f_r$ is ensured by remark
\ref{rem_easy-down}(i)).
Let $G_r$ be the Galois group of the covering $f_r$; the sheaf
$(F\otimes_\Lambda\bar\Lambda)_{|\D(r,1)}$ can be regarded in
the usual way as a finite dimensional $\bar\Lambda$-linear
representation of $G$, in other words, as a positive element
$\bar\chi\in K_0(\bar\Lambda[G])$.
By inspecting the definitions we find the identity :
$$
\ssw^\natural(F,s^+)=\mathrm{length}_\Lambda\Lambda\cdot
\langle\bar\ssw{}^\natural_{G_r}(s^+),\bar\chi\rangle_{G_r}
\qquad\text{for every $s\in(r,1]\cap\Gamma_K$}.
$$
This leads us to set :
$$
\delta_F(-\log s):=\mathrm{length}_\Lambda\Lambda\cdot
\delta_{f_r,\bar\chi}(-\log s)\qquad
\text{for every $s\in(r,1]\cap\Gamma_K$}
$$
(notation of \eqref{subsec_produce}). Suppose now that
$f'_r:X'_r\to\D(r,1)$ is another Galois covering that dominates
$f'_r$ ({\em i.e.} such that $f'_r$ factors through $f_r$).
Let $G'$ be the Galois group of $f'_r$, and set
$\bar\chi':=\mathrm{Res}^G_{G'}\bar\chi$; it follows easily from
\eqref{eq_bigger-group} that
$\delta_{f_r,\bar\chi}=\delta_{f'_r,\bar\chi{}'}$. Since any two
Galois \'etale coverings are dominated by a common one, we deduce
that the function $\delta_F$ thus defined is independent of
the choice of $f$. Especially, let $r'<r$ be another positive
real number in $\Gamma_K$, choose a Galois \'etale covering
$f_{r'}:X'\to\D(r',1)$ trivializing $F$, let $\bar\chi{}'$ be
the $\bar\Lambda[G_{r'}]$-module corresponding to
$(F\otimes_\Lambda\bar\Lambda)_{|\D(r',1)}$ and define the function
$\delta'_F:=\mathrm{length}_\Lambda\Lambda\cdot
\delta_{f_{r'},\bar\chi{}'}:(0,\log 1/r']\cap\log\Gamma_K\to\R$;
it follows easily that $\delta'_F$ agrees with $\delta_F$
wherever the latter is defined. Hence, by patching these
locally defined function, we obtain a well defined mapping :
$$
\delta_F:-\log\Gamma_K^+\to\R.
$$

\begin{proposition}\label{prop_invoke}
In the situation of \eqref{subsec_Lambda-mod}, the mapping
$\delta_F$ is the restriction of a non-negative, piecewise
linear, continuous and convex real-valued function. Moreover :
$$
\frac{d\delta_F}{dt}(-\log r^+)=\ssw^\natural(F,r^+)
\qquad
\text{for every $r\in\Gamma^+_K$}.
$$
\end{proposition}
\begin{proof} All the hard work has already been done, and
it remains only to invoke proposition \ref{prop_modular}.
\end{proof}

\begin{corollary}\label{cor_monotone}
In the situation of \eqref{subsec_Lambda-mod}, the function :
$$
-\log\Gamma_K^+\to\Z\qquad
-\log r\mapsto\ssw^\natural(F,r^+)
$$
is monotonically non-decreasing, and moreover :
$$
\ssw^\natural(F,0^+):=
\lim_{r\to 0^+}\ssw^\natural(F,r^+)\in\N\cup\{+\infty\}.
$$
\end{corollary}
\begin{proof} The monotonicity follows from the convexity
of $\delta_F$. The limit value $\ssw^\natural(F,0^+)$ cannot
be negative, since $\delta_F$ is non-negative.
\end{proof}

\sset\subsubsection{}\label{subsec_breaks}
To advance further, we use the break decomposition of \cite[\S8]{Hu3}.
We choose the elegant presentation of N.Katz in \cite[Ch.1]{Kat},
which makes it transparent that this is really a general
representation-theoretic device. Indeed, suppose that $H$ is a finite
group with a unique (hence normal) $p$-Sylow subgroup $P$, and assume
that $P$ admits a finite descending filtration :
$$
P_n:=\{1\}\subset P_{n-1}\subset\cdots P_1\subset P_0:=P
$$
consisting of subgroups $P_i$ normal in $H$ for every $i\leq n$.
Let $R$ be any $\Z[1/p]$-algebra, and for every $i\leq n$ let us
define the element :
$$
e_i:=\frac{1}{o(P_i)}\sum_{g\in P_i}g\in R[H].
$$
Since $P_i$ is normal in $H$, every $e_i$ is a central idempotent
element in $R[H]$. One verifies easily that the central idempotents :
$$
e_0,\quad e_1\cdot(1-e_0),\quad e_2\cdot(1-e_1),\quad\cdots\quad,
e_n\cdot(1-e_{n-1})=1-e_{n-1}
$$
are orthogonal and sum to $1$, hence define a natural decomposition
of $R[H]$ in $n+1$ direct factors. Correspondingly, every $R[H]$-module
$M$ admits a {\em break decomposition} :
$$
M=M_{-1}\oplus M_0\oplus\cdots\oplus M_{n-1}
$$
into $R[H]$-submodules such that :
$$
M_{-1}=M^P\qquad M_i^{P_i}=0 \quad\text{for every $i\geq 0$}
\qquad\text{and}\qquad M_i^{P_j}=M_i\quad\text{whenever $j>i$.} 
$$
Furthermore, for every $i\leq n$ the rule $M\mapsto M_i$ defines
an exact functor : $R[H]\Mod\to R[H]\Mod$, and for every pair of
$R[H]$-modules $M$, $N$ we have :
$$
\Hom_{R[H]}(M_i,N_j)=0\qquad\text{whenever $i\neq j$}.
$$
One deduces easily that :
\set\begin{equation}\label{eq_estimate-break}{
\begin{array}{r@{\:\subset\:}ll}
M_i\otimes_R N_j & (M\otimes_R N)_{\max(i,j)} & \text{if $i\neq j$} \\
M_i\otimes_R N_i & \sum_{j\leq i}(M\otimes_R N)_j &
\text{for every $i=-1,\dots,n$} \\
\Hom_R(M_i,N_j) & \Hom_R(M,N)_{\max(i,j)} & \text{if $i\neq j$} \\
\Hom_R(M_i,N_i) & \sum_{j\leq i}\Hom_R(M,N)_j &
\text{for every $i=-1,\dots,n$}.
\end{array}}
\end{equation}
See \cite[Lemma 1.3]{Kat} for details. Moreover, the break
decomposition is invariant under arbitrary base-change $R\to R'$,
{\em i.e.} we have :
\set\begin{equation}\label{eq_base-change-brk}
(M\otimes_RR')_i=M_i\otimes_RR'
\qquad\text{for every $i\leq n$}.
\end{equation}

\sset\subsubsection{}\label{subsec_upper-numb}
The generalities of \eqref{subsec_breaks} shall be applied to
the group $H:=St_x$ of \eqref{sec_conductors} and its higher
ramification filtration, and with $R:=\Lambda$. However, for
book-keeping purposes, it is convenient to replace the
lower-numbering indexing of this filtration, by a upper-numbering
one. With our current notation, this is defined as follows. First,
one considers the order-preserving bijection :
$$
\phi:\Q\otimes_\Z\Gamma_x\to\Q\otimes_\Z\Gamma_x \qquad
\gamma\mapsto\prod_{g\in St_x}\max(\gamma,i(g)/\gamma_0)
$$
where $\gamma_0\in\Gamma_x^+$ is defined as in
\eqref{subsec_dir-img-etale}. Notice that $\phi$ maps
$(\Q\otimes_\Z\Gamma_x)^+$ bijectively onto itself. Next we let :
$$
P^\gamma:=P_{\phi^{-1}(\gamma)} \qquad\text{for every
$\gamma\in\Gamma_x^+$}.
$$
If $\gamma_1>\cdots>\gamma_{n-1}>\gamma_n$
are the jumps in the family $(P_\gamma~|~\gamma\in\Gamma_x^+)$ that
are less than $1$, we obtain therefore a finite filtration of $St_x$ :
$$
\{1\}\subset P^{\phi(\gamma_n)}\subset P^{\phi(\gamma_{n-1})}
\subset\cdots\subset P^{\phi(\gamma_1)}\subset P
$$
where $P$ is the $p$-Sylow subgroup of $St_x$. If now $M$ is any
$\Lambda$-module, we derive a break decomposition as in \eqref{subsec_breaks} :
$$
M=M(1)\oplus M(\phi(\gamma_1))\oplus\cdots\oplus M(\phi(\gamma_n))
$$
such that $M(1)=M^P$ and :
\set\begin{equation}\label{eq_characterize-brks}
M(\phi(\gamma_i))^{P^{\phi(\gamma_i)}}=0 \quad\text{for every $i\leq n$}
\quad\text{and}\quad
M(\phi(\gamma_i))^{P^{\phi(\gamma_j)}}=M(\phi(\gamma_i))
\quad\text{whenever $j>i$}.
\end{equation}
The values $\phi(\gamma_i)$ such that $M(\phi(\gamma_i))\neq 0$ are
called {\em the breaks\/} of $M$.

\sset\subsubsection{}\label{subsec_break-stalk}
Especially, let $F$ be as in \eqref{subsec_Lambda-mod}, and
$r\in\Gamma_K^+$. Choose $a\in(0,r)\cap\Gamma_K$
and a finite Galois \'etale covering $f_a:X_a\to\D(a,1)$
trivializing $F_{|\D(a,1)}$. Denote by $G_a$ the Galois group of $f_a$,
pick any point $x\in f^{-1}_a(\eta(r))$, and let $St_x\subset G_a$
be the stabilizer of $x$. The stalk $F_r:=F_{\eta(r)}$ is a
$\Lambda[G_a]$-module of finite type, hence a $\Lambda[St_x]$-module,
by restriction. Then the upper numbering filtration of $St_x$ yields
a break decomposition :
\set\begin{equation}\label{eq_break-decomp}
F_r=F_r(\beta_0(r))\oplus F_r(\beta_1(r))\oplus\cdots\oplus F_r(\beta_n(r))
\end{equation}
with $1=\beta_0(r)>\beta_1(r)>\cdots>\beta_n(r)$.
Of course, this is also a $\pi_1(r)$-equivariant decomposition
(notation of \eqref{eq_pione-r}).

\begin{lemma}\label{lem_sum-of-brks}
With the notation of \eqref{subsec_break-stalk}, we have the identity :
$$
\delta_F(-\log r)=-\sum_{i=1}^n\log\beta_i(r)^\flat\cdot
\mathrm{length}_\Lambda(F_r(\beta_i(r))).
$$
\end{lemma}
\begin{proof} First, \eqref{eq_base-change-brk} allows to reduce
to the case where $\Lambda=\bar\Lambda$. Then the  sought identity
is derived from \eqref{eq_double-wham} by a standard calculation
(cp. the proof of \cite[Prop.8.2 and Cor.8.4]{Hu3}).
\end{proof}

\begin{lemma}\label{lem_continu-brk}
In the situation of \eqref{subsec_break-stalk}, choose $r'\in(a,r)$
such that the conditions {\em (i)} and {\em (ii)} of theorem
{\em\ref{th_continuities}} are fulfilled, with $f:=f_a$.
Then, for every $s\in(r',r]\cap\Gamma_K$ :
\begin{enumerate}
\item
The length of the break decomposition of $F_s$ :
$$
F_s=F_s(\beta_0(s))\oplus F_s(\beta_1(s))\oplus\cdots\oplus F_s(\beta_n(s))
$$
is independent of $s$.
\item
For every other $s'\in(r',r]\cap\Gamma_K$, the isomorphism
$\omega_{s,s'}$ of theorem {\em\ref{th_continuities}(i)} induces
equivariant isomorphisms $F_s(\beta_k(s))\simeq F_{s'}(\beta_k(s'))$,
for every $k=0,\dots,n$.
\item
Moreover, set $\beta_k^\natural:=\beta_k(r)^\natural$\ \  for every
$k\leq n$. Then :
$$
\beta_k(s)=(s/r)^{\beta_k^\natural}\cdot\beta_k(r)
\qquad\text{for every $s\in(r',r]\cap\Gamma_K$}.
$$
\end{enumerate}
\end{lemma}
\begin{proof} This is an exercise in translating from lower
to upper numbering. Indeed, assertions (i) and (ii) are clear
from theorem \ref{th_continuities}, and it remains only to
show (iii). However, let $\{x(s)~|~s\in(r',r]\cap\Gamma_K\}$
be a family of points as in theorem \ref{th_continuities}
(with $f:=f_a$), and 
$P_{\gamma_n(s)}\subset\cdots\subset P_{\gamma_1(s)}\subset St^{(p)}_{x(s)}$
the lower numbering ramification filtration of $St_{x(s)}$.
Let also $\gamma_0$ be the largest element in
$\Gamma^+_{x(s)}\setminus\{1\}$. Then, for every $k\leq n$
we may compute :
$$
\begin{array}{r@{\:=\:}l}
\beta_k(s):=\phi(\gamma_k(s)) & \gamma_k(s)^{o(P_{\gamma_k(s)})}\cdot
\displaystyle\prod_{g\in St_{x(s)}\setminus P_{\gamma_k(s)}} i(g)/\gamma_0 \\
& \gamma_k(s)^{o(P_{\gamma_k(s)})}\cdot
\displaystyle\prod_{1\leq t<k}
\gamma_t(s)^{o(P_{\gamma_t(s)})-o(P_{\gamma_{t+1}(s)})}
\end{array}
$$
so the sought identities follow from theorem \ref{th_continuities}(iii).
\end{proof}

\subsection{Local systems with bounded ramification}\label{sec_bound-ram}
We keep the notation of \eqref{sec_break}.
Corollary \ref{cor_monotone} suggests the following :

\begin{definition} Let $F$ be a locally constant and locally
free $\Lambda$-module of finite rank on the \'etale site of $\D(1)^*$.
We say that $F$ has {\em bounded ramification\/} if
$\ssw^\natural(F,0^+)\in\N$.
\end{definition}

The class of sheaves with bounded ramification includes
that of meromorphically ramified $\Lambda$-modules from
\cite{Ra}. The first result concerning these sheaves is :

\begin{theorem}\label{th_tannaka}
{\em(i)}\ \  If $F$ and $F'$ are two $\Lambda$-modules with
bounded ramification on $\D(1)^*_\et$, then $F\otimes_\Lambda F'$
and $\mathscr{H}om_\Lambda(F,F')$ have also bounded ramification.
\begin{enumerate}
\addenu
\item
Especially, if $\Lambda$ is a field, the full subcategory of
the category of $\Lambda$-modules on $\D(1)^*_\et$, consisting
of all $\Lambda$-modules with bounded ramification, is tannakian.
\end{enumerate}
\end{theorem}
\begin{proof} Of course (ii) follows from (i). To show assertion (i)
for $F\otimes_\Lambda F'$, since we know {\em a priori\/} that
$\delta_{F\otimes_\Lambda F'}$ is piecewise linear, continuous
and convex (proposition \ref{prop_invoke}), it suffices to provide
a rough estimate on the rate of growth of the latter mapping, in terms of
$\delta_F$ and $\delta_{F'}$. However, for given $r\in\Gamma_K^+$
let
$$
F_r=F_r(1)\oplus F_r(\gamma_1)\oplus\cdots\oplus F_r(\gamma_n)
\qquad\text{and}\qquad
F'_r=F'_r(1)\oplus F'_r(\gamma'_1)\oplus\cdots\oplus F'_r(\gamma'_m)
$$
be the break decompositions of the stalks of $F$ and $F'$ over
the point $\eta(r)$. We set :
$$
\lambda_i:=\mathrm{length}_\Lambda F_r(\gamma_i) \qquad
x_i:=-\log\gamma_i^\flat\qquad\text{for every $i\leq n$}
$$
and
$$
\lambda'_j:=\mathrm{length}_\Lambda F'_r(\gamma'_j) \qquad
y_j:=-\log\gamma_j^{\prime\flat}\qquad\text{for every $j\leq m$}.
$$
Clearly $(F\otimes_\Lambda F')_{\eta(r)}=
\oplus_{i,j}F_r(\gamma_i)\otimes_\Lambda F'_r(\gamma'_j)$,
and using \eqref{eq_estimate-break} and lemma \ref{lem_sum-of-brks}
we arrive at the inequality :
$$
\delta_{F\otimes_\Lambda F'}(-\log r)\leq
\sum_{ij} x_iy_j\cdot\max(\lambda_i,\lambda'_j)\leq
(\rk_\Lambda F)\cdot(\rk_\Lambda F')\cdot
\max(\delta_F(-\log r),\delta_{F'}(-\log r))
$$
as required. A similar argument takes care of $\mathscr{H}om_\Lambda(F,F')$
and concludes the proof of the theorem.
\end{proof}

\begin{example}\label{ex_Kumm-chi} Choose a coordinate $T$ on
$\A^1_K$, (so that $\A^1_K=\Spec\,K[T]$), and for
every $m\in\N$, denote by $f_m:\A^1_K\to\A^1_K$ the morphism
such that $f_m^*(T)=T^m$. The restriction of $f_m$ to
$\G_{m,K}:=\Spec\,K[T,T^{-1}]$ is a torsor in the \'etale
topology of $\G_{m,K}$ for the group $\bmu_m\subset K^\times$,
hence the analytification $f^\ad_m$ is a $\bmu_m$-torsor in
the \'etale topology of $\G_{m,K}^\ad$.
For every character $\chi:\bmu_m\to\Lambda^\times$, we let
$\cK(\chi)$ be the locally free rank one $\Lambda$-module on
$\D(1)^*_\et$ which is induced, via $\chi$, by the restriction
of the torsor $f^\ad_m$.
Let $n\leq m$ be the order of $\chi$ ({\em i.e.} the smallest
$k\in\N$ such that $\chi^k$ is the trivial character), and
denote by $\mu$ a chosen generator of $\bmu_n$. It follows
from \cite[Ex.8.8]{Hu3} that :
$$
\ssw^\natural(\cK(\chi),r^+)=0
\qquad\text{for every $r\in\Gamma_K^+$}.
$$
Especially, $\cK(\chi)$ is a $\Lambda$-module with bounded
ramification. Moreover, if $\chi$ is not trivial, the (unique)
break of $\cK(\chi)$ equals $|1-\mu|$ for every
$r\in\Gamma^+_K$.
\end{example}

\begin{example}\label{ex_bricks}
Keep the notation of example \ref{ex_Kumm-chi}, and let :
$$
\D(1^-):=\bigcup_{r\in\Gamma_K^+\setminus\{1\}}\D(r).
$$
The morphism :
$$
\log:\D(1^-)\to(\A^1_K)^\ad\quad\text{such that}\quad
\log^*(T)=\sum^\infty_{n=1}\frac{(-1)^{n-1}}{n}\cdot T^n
$$
is a torsor in the \'etale topology of $(\A^1_K)^\ad$ for the
constant group $\mu_{p^\infty}$ (see \cite[Lemma 9.4]{Hu3} or
\cite[Lemma 6.1.1]{Ra}). We denote by $\cL$ the locally free
$\Lambda$-module on $(\A^1_K)^\ad_\et$ which is induced by this
torsor via the inclusion \eqref{eq_imbed-mu}. The sheaf $\cL$
has been studied at length in \cite{Ra}. For instance, one
can show that for any two morphisms
$\phi,\psi:(\A^1_K)^\ad\to(\A^1_K)^\ad$, there exists
a natural isomorphism :
\set\begin{equation}\label{eq_atlength}
(\phi^*\cL)\otimes_\Lambda(\psi^*\cL)\simeq(\phi+\psi)^*\cL
\end{equation}
where $\phi+\psi$ is the addition of $\phi$ and $\psi$, regarded
as sections of the structure sheaf of $(\A^1_K)^\ad$.

Let $g:\G_{m,K}\to\A^1_K$ be the morphism such that
$g^*(T)=T^{-1}$, and for every $q\in\Q_{\geq 0}$,
let $(m,n)\in\N^2$ be the unique pair of relatively prime integers
such that $q=n/m$; we set :
$$
\cL(q):=f_{m*}\circ f^*_n\circ g^*\cL
$$
where $f_m$ and $f_n$ are defined as in example \ref{ex_Kumm-chi}.
The following lemma \ref{lem_bricks} shows that the sheaves
$\cL(q)_{|\D(1)^*}$ have bounded ramification.
\end{example}

\begin{lemma}\label{lem_bricks}
Let $q\in\Q_{>0}$, and write $q=n/m$, with $n,m\in\N$ and $(n,m)=1$;
moreover, write $n=p^a N$, $m=p^bM$, with $a,b\geq 0$ and $(N,p)=(M,p)=1$
(of course, either $a=0$ or $b=0$).
Set $\lambda:=|p|^{1/(p-1)}$ and $l:=\mathrm{length}_\Lambda\Lambda$.
The following holds (notice that $\delta_{\cL(q)}(\rho)$ is defined for
every $\rho\in\log\Gamma_K$) :
\begin{enumerate}
\item
$\delta_{\cL(q)}$ is the unique continuous piecewise linear
function such that :
\begin{enumerate}
\item
$\delta_{\cL(q)}(\rho)=0$\ \ whenever $\rho\leq q^{-1}\log\lambda$.
\item
The right slope of $\delta_{\cL(q)}$ equals :
\begin{itemize}
\item
$lN$ on the interval $[q^{-1}\log\lambda,q^{-1}\log|1/p||)\cap\log\Gamma_K$;
\item
$lp^jN$ on the interval
$[jq^{-1}\log|1/p|,(j+1)q^{-1}\log|1/p|)\cap\log\Gamma_K$,
for every $j=1,\dots,a-1$.
\item
$ln$ on the half-line $[aq^{-1}\log|1/p|,+\infty)\cap\log\Gamma_K$.
\end{itemize}
\end{enumerate}
\item
$\ssw^\natural(\cL(q),0^+)=ln$.
\item
For every $q\in\Q_{\geq 0}$, the sheaf $\cL(q)_{|\D(1)^*}$ is
indecomposable in the category of locally free $\Lambda$-modules
on $\D(1)^*_\et$.
\item
More precisely, for every $r\in\Gamma_K$, let $\pi_1(r)^{(p)}$
be the unique $p$-Sylow subgroup of $\pi_1(r)$. Set :
$$
r_0:=\left\{\begin{array}{ll}
               |p|^{(b-1)/q}     & \text{if\ \  $b\neq 0$} \\
               \lambda^{-1/q}    & \text{if\ \  $b=0$}.
             \end{array}\right.
$$
Then, for every $r\leq r_0$, the stalk $\cL(q)_r$ is an indecomposable
$\Lambda[\pi_1(r)^{(p)}]$-module (notation of \eqref{subsec_break-stalk});
especially, $\cL(q)_r$ has a unique break $\beta(q,r)$.
\item
Suppose that $a=b=0$. Then :
$$
\beta(q,r)=\left\{\begin{array}{ll}
                    r^q\cdot(1-\eps)^q\cdot\lambda\cdot|p|^{-a} &
                            \text{for\ \ $r\leq\lambda^{-1/q}$} \\
                    1 & \text{otherwise}.
                  \end{array}\right.
$$
\end{enumerate}
\end{lemma}
\begin{proof} Without loss of generality, we may assume that
$\Lambda$ is a field, hence $l=1$. Notice that $\cL$ is trivial
on the disc $\D(\lambda^-):=\bigcup_{r<\lambda}\D(r)$, hence :
$$
\delta_{\cL(1)}(\rho)=0 \qquad\text{whenever $\rho<\log\lambda$}.
$$
In view of proposition \ref{prop_invoke}, $\delta_{\cL(1)}$ is
then completely determined, once we know its right derivative
$\ssw^\natural(\cL(1),\cdot)$. However, \cite[Lemma 9.4]{Hu3}
shows that :
$$
\ssw^\natural(\cL(1),r^+)=1\qquad\text{whenever $r\leq\lambda^{-1}$}
$$
which gives (i), for $q=1$. Suppose now that $n=p^aN>1$ is an
integer, and set $P:=p^a$; according to \cite[Ex.8.8(i)]{Hu3},
we have :
$$
\ssw^\natural(\cL(n),r^+)=\ssw^\natural(f_{k*}\cL(n),r^{n+})
\qquad\text{for every $r\in\Gamma_K$ and every $k\in\N$}
$$
which translates as the identity :
\set\begin{equation}\label{eq_push-forw}
\frac{d\delta_{\cL(n)}}{dt}(\rho^+)=
\frac{d\delta_{f_{k*}\cL(n)}}{dt}(k\rho^+)
\qquad\text{for every $\rho\in\log\Gamma_K$ and every $k\in\N$}.
\end{equation}
We shall apply \eqref{eq_push-forw} with $k=n$, so we need to
calculate the conductor of $f_{n*}\cL(n)$. The projection formula
yields :
$$
f_{n*}\cL(n)\simeq\cL(1)\otimes_\Lambda f_{n*}\Lambda.
$$
\begin{claim} There is a natural isomorphism :
$$
f_{n*}\Lambda\simeq f_{N*}\Lambda\otimes_\Lambda f_{P*}\Lambda.
$$
\end{claim}
\begin{pfclaim} The morphism $f_n$ induces an inclusion
of fundamental groups
$\phi:H:=\pi_1(\G_{m,K},x)\to G:=\pi_1(\G_{m,K},f_n(x))$
(for any choice of a geometric point $x$) whose image is a normal
subgroup with cokernel isomorphic to $\Z/n\Z$. The constant
sheaf $\Lambda$ on $\G_{m,K}$ corresponds to the trivial
representation of $H$, and $f_{n*}\Lambda$ is the induction
of this representation along the inclusion $\phi$. Likewise
we may describe $f_{N*}\Lambda$ and $f_{P*}\Lambda$.
However :
$$
\Lambda\otimes_{\Lambda[H]}\Lambda[G]\simeq\Lambda[G/H]=\Lambda[\Z/n\Z]
$$
so $f_{n*}\Lambda$ is also the induction of the trivial
representation of the trivial group $\{0\}$, along the inclusion
$\{0\}\subset\Z/n\Z$, and likewise for $f_{N*}\Lambda$ and
$f_{P*}\Lambda$. So finally, the sought isomorphism boils
down to the $\Lambda$-algebra isomorphism :
$\Lambda[\Z/n\Z]\simeq\Lambda[\Z/P\Z]\otimes_\Lambda\Lambda[\Z/N\Z]$.
\end{pfclaim}

Since $\mu_{p^\infty}\subset\Lambda^\times$, we have :
$$
f_{P*}\Lambda\simeq\bigoplus_{\chi:\bmu_P\to\mu_{p^\infty}}\cK(\chi)
$$
where the sum runs over the $P$ different characters of
$\bmu_P\subset K^\times$ (notation of example \ref{ex_Kumm-chi}).
Thus, $f_{n*}\cL(n)$ decomposes as the direct sum of $P$ terms
of the form
$$
\cM(\chi):=\cL(1)\otimes_\Lambda\cK(\chi)\otimes_\Lambda f_{M*}\Lambda.
$$
Let $p^j>1$ be the order of the character $\chi$; according to
\cite[Ex.8.8(ii)]{Hu3}, the unique break of $\cK(\chi)_r$ is
independent of $r$, and equals $|p|^j\cdot\lambda$. On the other
hand, the unique break $\beta(r)$ of $\cL(1)_r$ can be computed from
$\delta_{\cL(1)}$ using proposition \ref{prop_invoke} : we
get $\beta(r)=1$ for $r>\lambda^{-1}$ and $\beta(r)=r\lambda\cdot(1-\eps)$
for $r\leq\lambda^{-1}$. From \eqref{eq_estimate-break} we deduce
that the unique break of $(\cL(1)\otimes_\Lambda\cK(\chi))_r$ equals 
$1$ for $r>|p|^j$ and $r\lambda\cdot(1-\eps)$ for $r\leq|p|^j$.
Next, since $(N,p)=1$, the only possible break of the stalk
$(f_{N*}\Lambda)_r$ equals $1$, hence the stalks of
$\cL(1)\otimes_\Lambda\cK(\chi)$ and $\cM(\chi)$ have the same breaks.
Consequently :
\set\begin{equation}\label{eq_contrib}
\frac{d\delta_{\cM(\chi)}(\rho^+)}{dt}=
\left\{\begin{array}{ll}
         0 & \text{for $\rho<j\log|1/p|$} \\
         N & \text{otherwise}
       \end{array}\right.
\end{equation}
Since $\delta_{f_{n*}\cL(n)}=\sum_\chi\delta_{\cM(\chi)}$, assertion
(i) for $q=n$ follows from \eqref{eq_push-forw} (with $k:=n$) and
\eqref{eq_contrib}.
Finally, let $q:=n/m$, with $n,m$ two relatively prime integers;
in order to determine the right slope of $\delta_{\cL(q)}$, it
suffices to apply \eqref{eq_push-forw} with $k:=m$. This completes
the proof of (i).

Assertion (ii) is an immediate consequence of (i); also (iii)
follows directly from (iv), and (v) follows from (i) and (iv).
Hence it remains only to show (iv) when $\Lambda$ is a field,
which we may assume to be algebraically closed.
The assertion is obvious if $q$ is an integer, since in that case
$\cL(q)_r$ has rank one. For the general case $q=n/m$, notice
that the action of $\pi_1(s)$ (resp. $\pi_1(s^m)$) on $\cL(n)_s$
(resp. $\cL(q)_{s^m}$) factors through a finite quotient
$H_s$ (resp. $G_s$), and :
$$
\cL(q)_{s^m}\simeq\Ind^{G_s}_{H_s}\cL(n)_s
\qquad\text{for every $s\in\Gamma_K$}.
$$
The morphism $f_m$ is a torsor for the group $\bmu_m$, and
we have a natural identification $G_s/H_s\simeq\bmu_m$.
We shall apply Mackey's irreducibility criterion (this is
shown in \cite[\S7.4, Cor.]{Se2} in case the base field has
characteristic zero, but the result holds whenever the
characteristic of the algebraically closed field $\Lambda$
does not divide the order of $G_s$; this latter condition
is clearly fulfilled here). To this aim, we have to show
that, for every $\mu\in\bmu_m\setminus\{1\}$, the conjugate
representation $\cL(n)_s^\mu$ is not isomorphic to $\cL(n)_s$.
However, we have a natural identification :
$$
\cL(n)_s^\mu\simeq\mu^*\cL(n)_s
$$
where $\mu:\G_{m,K}:=\Spec K[T,T^{-1}]\to\G_{m,K}$
is the morphism such that $\mu^*(T)=\mu\cdot T$.
According to \eqref{eq_atlength}, we have a natural isomorphism :
$$
\mu^*\cL(n)\otimes_\Lambda\cL(n)^{-1}\simeq g^*\cL
$$
where $g:\G_{m,K}\to\G_{m,K}$ is the morphism such that
$g^*(T)=(1-\mu^{-n})T^{-n}$. Since $\cL_r$ is not trivial
whenever $r\geq\lambda$, it follows that $(g^*\cL)_s$ is not
trivial whenever $|1-\mu^{-n}|\cdot s^{-n}\geq\lambda$.
However, $\mu$ is a primitive $p^b$-root of unity for some
$j=1,\dots,b$, we have :
$$
|1-\mu|=\lambda\cdot|p|^{j-1}
$$
so in this case, $(g^*\cL)_s$ is not trivial for $s\leq|p|^{(j-1)/n}$.
If $b=0$, then $|1-\mu|=1$, and then $(g^*\cL)_s$ is not
trivial for $s\leq\lambda^{-1/n}$. Letting $r:=s^m$, we obtain the
contention, in either case. As an immediate
consequence, we  see that $\cL(q)_r$ admits a single break
$\beta(q,r)$ for $r\leq r_0$; this break can be determined by
evaluating $\delta_{\cL(q)}(-\log r)$, since the latter
must equal $-m\log\beta(q,r)^\flat$ (lemma \ref{lem_sum-of-brks});
we leave to the reader the elementary calculation.
\end{proof}

\sset\subsubsection{}\label{subsec_ideally}
Let $F$ be a locally free $\Lambda$-module on $\D(1)^*$ with
bounded ramification. We wish to define the {\em breaks of $F$
around the origin}. Ideally, one would like to define a stalk
$F_{\eta(0)}$ that captures the behaviour of $F$ in arbitrarily
small punctured open discs $\D(\eps)^*$ centered at the origin;
then the sought breaks should be numerical invariants associated
to this stalk.
To make sense of this, one would like to complete somehow
the sequence of points $(\eta(r)~|~r\in\Gamma^+_K)$ with
a limit point $\eta(0)$; however, such a limit point seems
to elude the grasp of the formalism of adic spaces, hence we
have to proceed in a rather more indirect fashion. But the
ideal picture should be kept in mind, as it motivates much
of what we are trying to do in the remainder of this work.

To begin with, for given $r,\alpha\in\Gamma_K^+$ we set :
\set\begin{equation}\label{eq_flat-breaks}
F_r^\flat(\alpha):=\bigoplus_{\beta_i(r)^\flat=\alpha}F_r(\beta_i(r))
\end{equation}
where $1=\beta_0(r)>\cdots>\beta_n(r)$ are the breaks of $F_r$,
so that we have the break decomposition \eqref{eq_break-decomp}.
Say that $\rk_\Lambda F=d$ and $-\log r=\rho$; we consider the
unique sequence of real numbers:
\set\begin{equation}\label{eq_redundant}
0\leq f_1(\rho)\leq f_2(\rho)\leq\cdots\leq f_d(\rho)
\end{equation}
in which, for every $\beta\in\Gamma^+_K$, the value $-\log\beta$
appears with multiplicity equal to $\rk_\Lambda F_r^\flat(\beta)$.

\begin{lemma}\label{lem_surprise} The functions $f_1,\dots,f_d$
extend to piecewise linear continuous maps $\R_{\geq 0}\to\R_{\geq 0}$.
\end{lemma}
\begin{proof} Using lemma \ref{lem_continu-brk}, we deduce
already that $f_1,\dots,f_d$ extend to continuous, linear functions
on every small segment of the form $[-\log r,-\log r')$.

It remains to show that for every $\rho>0$ there
is some small segment $(\rho',\rho]$ on which the functions
$f_i$ are continuous. To this aim, we remark that all the
considerations of \eqref{subsec_ideally} and
\eqref{subsec_upper-numb} can be repeated for the
family of stalks $F_{\eta'(r)}$ (instead of $F_r:=F_{\eta(r)}$).
We obtain in this way a break decomposition
$F_{\eta'(r)}(\delta_0)\oplus\cdots\oplus F_{\eta'(r)}(\delta_l)$
for $F_{\eta'(r)}$, and we may define the submodules
$F_{\eta'(r)}^\flat(\alpha)$ for every $\alpha\in\Gamma^+_K$,
just as in \eqref{eq_flat-breaks}.
Using the ranks of the modules $F_{\eta'(r)}^\flat(\alpha)$,
we may finally construct a non-decreasing sequence
$0\leq f'_1(\rho)\leq f'_2(\rho)\leq\cdots\leq f'_d(\rho)$
analogous to \eqref{eq_redundant}. Making use of
\eqref{subsec_continuities} (rather than theorem \ref{th_continuities}),
one can then show the analogue of lemma \ref{lem_continu-brk}
which expresses the continuity of the breaks $\delta_i$; from
the latter, we see that the functions $f'_i$ are continuous on
segments of the form $(\rho',\rho]$. To conclude it suffices to
show that $f_i=f'_i$ for every $i\leq d$. This boils down to
the following :

\begin{claim} $F_{\eta'(r)}^\flat(\alpha)\simeq F_r^\flat(\alpha)$
for every $\alpha\leq 1$.
\end{claim}
\begin{pfclaim}[] We do not only assert the existence of an
isomorphism in the category $\Lambda\Mod$, but more precisely,
that the two modules are equivariantly isomorphic, in the
following sense.
Say that $r\in(a,b)$ for some $a,b\in\Gamma^+_K$,
and pick a Galois \'etale covering $f:X\to\D(a,b)$ that trivializes
$F_{|\D(a,b)}$; choose also points $x,x'$ lying over repectively
$\eta(r)$ and $\eta'(r)$, such that $x^\flat=x^{\prime\flat}$.
We shall use the ramification filtration
$(P^\flat_\gamma~|~\gamma\in\Gamma^+_K)$ of $St^\flat_x$,
given by definition \ref{def_central}. According to proposition
\ref{prop_superjump}, the $p$-Sylow subgroup $St^{\flat(p)}_x$ is
naturally a subgroup of $St^{(p)}_x$ and $St^{(p)}_{x'}$, and the
claim amounts to a $St^{\flat(p)}_x$-equivariant identification of
$F_{\eta'(r)}^\flat(\alpha)$ and $F_r^\flat(\alpha)$.

Proceeding as
in \eqref{subsec_upper-numb}, we replace the lower-numbering
indexing by the upper-numbering
$(P^{\flat,\gamma}~|~\gamma\in\Gamma^+_K)$. Say that
$(P^\gamma~|~\gamma\in\Gamma^+_x)$ (resp.
$(Q^\gamma~|~\gamma\in\Gamma^+_{x'})$) is the upper-numbering
ramification filtration for $St_x$ (resp. for $St_{x'}$).
Then proposition \ref{prop_superjump} yields the identities :
\set\begin{equation}\label{eq_superjumps}
\bigcup_{n\in\Z}Q^{\gamma_0^n\cdot\gamma}=P^{\flat,\gamma}=
\bigcup_{n\in\Z}P^{\gamma_0^n\cdot\gamma}
\qquad\text{for every $\gamma\in\Gamma^+_K\setminus\{1\}$}.
\end{equation}
Recall that, by construction, for each
$\beta\in\{\beta_1(r),\dots,\beta_n(r)\}$, the direct summand
$F_r(\beta)$ is of the form $e_\beta\cdot F_r$, where
$e_\beta$ is a certain central idempotent in $\Lambda[St^{(p)}_x]$.
Likewise, for every $\delta\in\{\delta_1,\dots,\delta_l\}$,
we have $F_{\eta'(r)}(\delta)=e'_\delta\cdot F_{\eta'(r)}$ for a certain
central idempotent $e'_\delta\in\Lambda[St^{(p)}_{x'}]$.
Then, by inspecting the definitions, we see that $e_\beta$
(resp. $e'_\delta$) actually lies in the subring $\Lambda[St^{\flat(p)}_x]$,
whenever $\beta^\flat<1$ (resp. $\delta^\flat<1$). 
(Here $St^{\flat(p)}_x$ is the unique $p$-Sylow subgroup of
$t^\flat_x$.) Moreover, \eqref{eq_superjumps} leads to the
identities :
$$
\sum_{\beta^\flat=\alpha}e_\beta=\sum_{\delta^\flat=\alpha}e'_\delta
\qquad\text{for every $\beta\in\Gamma^+_K\setminus\{1\}$}.
$$
This already shows the claim for $\alpha<1$. The case for
$\alpha=1$ can similarly be dealt with, using \eqref{eq_superjumps}
and the characterization \eqref{eq_characterize-brks} of the
breaks : the details shall be left to the reader.
\end{pfclaim}
\end{proof}

\sset\subsubsection{}
We shall denote by $\Delta(F)\subset\R_{\geq 0}\times\R_{\geq 0}$
the union of the graphs of the functions $f_1,\dots,f_d$ defined
in \eqref{eq_redundant}.
Let now $(K',|\cdot|')$ be an algebraically closed valued
field extension of $(K,|\cdot|)$ whose value group is
$\R_{>0}$ ({\em e.g.} we can take a maximally complete
field containing $K$). The given $\Lambda$-module $F$ pulls
back to a locally constant $\Lambda$-module $F'$ on the adic
space $\D(1)^*\times_{\Spa\,K}\Spa\,K'$. In view of lemma
\ref{lem_base-change-cond}, we see that for every
$r\in\Gamma_K^+$ the breaks of $F'_r$ are the same
as that of $F_r$, therefore the subset $\Delta(F')$ is none else
than the topological closure of $\Delta(F)$. Hence for the
considerations that follow we may replace $K$ by $K'$ and
$F$ by $F'$, and assume that $\Gamma_K=\R_{>0}$. Simple operations
on $F$ can be translated into corresponding changes
for the subset $\Delta(F)$. For instance, for any $s\in K^+\setminus\{0\}$,
let $\mu_s:\D(1)^*\to\D(1)^*$ be the ``shrinking'' morphism
such $\mu_s^*(\xi)=s\cdot\xi$. We have $(\mu^*_sF)_x\simeq F_{\mu_s(x)}$
for every $x\in\D(1)^*$, so that
$$
\Delta(\mu^*_sF)=(\log|s|,0)+\Delta(F):=\R^2_{\geq 0}\cap
\{(x+\log|s|,y)~|~(x,y)\in\Delta(F)\}.
$$
We are now ready to make the following :

\begin{definition}\label{def_break-funct}
Assume that $\Gamma_K=\R_{>0}$ and let $F$ be a locally constant
and locally free $\Lambda$-module on $\D(1)^*_\et$ of finite rank.
Then the {\em break function\/} of $F$ is the mapping :
$$
\beta(F,\cdot):\Q_{\geq 0}\to\R_{\geq 0}\cup\{\infty\}
$$
defined by the rule :
$$
\beta(F,q):=\frac{1}{\rk_\Lambda\cL(q)}\cdot
\sup\,\{\ssw^\natural(F\otimes_\Lambda\mu_s^*\cL(q),0^+)~|~
s\in K^+\setminus\{0\}\}
\qquad\text{for every $q\in\Q_{\geq 0}$}.
$$
\end{definition}

\begin{remark}\label{rem_station}
Suppose that $F$ is the restriction of a sheaf $F'$ of
$\Lambda$-modules on $(\A^1_K)^\ad_\et$, and consider the
case $q=1$ : the cohomology complex $R\Gamma_c(F'\mu_s^*\cL(1))$
is none else than the stalk over the point $s^{-1}\in K=\A^1_K(K)$
of the Fourier transform $\cF(F')$ of $F$, as defined in \cite{Ra}.
In view of the (analogue of the) Grothendieck-Ogg-Shafarevich formula
(see \cite[Th.10.2]{Hu3}), we see that the function $\beta(F,1)$
essentially calculates the Euler-Poincar\'e characteristic of
$\cF(F)$ in a neighborhood of the point $\infty\in(\P^1_K)^\ad$.
This is the sort of quantities that appear in the method of
stationary phase (see the introduction, \S\ref{sec_dulcis}),
and indeed this sort of sheaf-theoretic harmonic analysis has
motivated the definition of the function $\beta$. 
\end{remark}

For any rational number $q$, we define the {\em denominator\/}
of $q$ as the smallest positive integer $n$ such that $nq\in\Z$.

\begin{theorem}\label{th_Newton}
Let $F$ be as in definition {\em\ref{def_break-funct}},
and suppose moreover that $F$ has bounded ramification. Then :
\begin{enumerate}
\item
For every $q\in\Q_{\geq 0}$, $\beta(F,q)$ is a positive rational
number, whose denominator divides the denominator of $q$.
\item
$\beta(F,q)\geq q\cdot\rk_\Lambda F$ for every $q\in\Q_{\geq 0}$, and
the inequality is an equality for every sufficiently large $q$
(so $\beta(F,\cdot)$ is eventually linear).
\item
The break function $\beta(F,\cdot)$ is the restriction of a function
$\R_{\geq 0}\to\R_{\geq 0}$ which is convex, continuous,
non-decreasing and piecewise linear whose slopes are integers.
\end{enumerate}
\end{theorem}
\begin{proof} Without loss of generality, we may assume that
$\Lambda$ is a field. We begin by introducing some notation :
we let $\SS\subset\Q$ be the subset of the numbers of the form
$n/m$ where $n$, $m$ are relatively prime positive integers,
such that $(p,n)=(p,m)=1$. Also, for every $s\in K^\times$ and
$q\in\R$, set
$$
c(q,s):=(p-1)^{-1}\log|p|+q\log|s|.
$$
By lemma \ref{lem_bricks}(v), for any $q\in\SS$, and every
$s\in K^\times$, the subset $\Delta(\mu^*_s\cL(q))$ is the graph
of the function :
$$
\rho\mapsto b(q,\rho,s):=\max\{0,q\rho-c(q,s)\}
\qquad\text{for every $\rho\in\R_{\geq0}$}.
$$
In the study of the ``stalk over $\eta(0)$'', we are allowed to
disregard the behaviour of our sheaf $F$ outside any given punctured
disc $\D(\eps)^*$, {\em i.e.} we may disregard the part of
$\Delta(F)$ that lies in a vertical band of the form $[0,c]\times\R$;
hence, let us define $\Sigma(q)$ as the subset of all $c(q,s)\in\R$
such that
$$
\Delta(\mu_s^*\cL(q))\cap\Delta(F)\cap([q^{-1}c(q,s),+\infty)\times\R)
$$
is a set whose cardinality is at most countable. Notice that
$\R\setminus\Sigma(q)$ has at most countable cardinality;
especially, $\Sigma(q)$ is dense in $\R$, for every $q\in\SS$.

Let $s\in K^+\setminus\{0\}$, $q\in\SS$, set $\cL':=\mu_s^*\cL(q)$,
and suppose that $c(q,s)\in\Sigma(q)$; this means that for every
$\rho\geq \max(q^{-1}c(q,s),0)$ :
\begin{itemize}
\item
either $(\rho,b(q,\rho,s))\notin\Delta(F)$,
\item
or else the right and left slope of $b(q,\cdot,s)$
at the point $\rho$ are different from the slopes of each of
the functions $f_i$ as in \eqref{eq_redundant}, such that
$f_i(\rho)=b(q,\rho,s)$.
\end{itemize}
However, say that $\rho=-\log r$, let $\gamma$ be the unique
break of $\cL'_r$, and $\beta_1,\dots,\beta_k$
the finitely many breaks of $F_r$; then by definition,
$\Delta(F)\cap(\{r\}\times\R)$ consists of the values
$-\log\beta_j^\flat$ (for $j=1,\dots,k$), and
$-\log\gamma^\flat=b(q,\rho,s)$. Furthermore,
by theorem \ref{th_continuities} and \eqref{subsec_continuities}, the
(right and left) slopes of the functions $f_i$ at the point $\rho$
are none else than the values $\beta_j^\natural$
(and likewise for the slope of $b(q,\cdot,s)$).
We conclude that $\gamma\notin\{\beta_1,\dots\beta_k\}$, and then
\eqref{eq_estimate-break} implies that the breaks of
$(F\otimes_\Lambda\cL')_r$ are the values
\set\begin{equation}\label{eq_good-breaks}
\beta'_j:=\min(\gamma,\beta_j)\qquad\text{for $j=1,\dots,k$}.
\end{equation}
Moreover, let
$$
M(\beta_1)\oplus\cdots\oplus M(\beta_k)
\qquad
M'(\beta'_1)\oplus\cdots\oplus M'(\beta'_k)
$$
be the break decompositions of $F_r$ and respectively
$(F\otimes_\Lambda\cL')_r$; then :
\set\begin{equation}\label{eq_ranks}
\rk_\Lambda M'(\beta'_j)=\rk_\Lambda\cL(q)\cdot\rk_\Lambda M(\beta_j)
\qquad\text{for every $j\leq k$}.
\end{equation}
Set $d:=\rk_\Lambda F$; combining lemma \ref{lem_bricks}(iii),
\eqref{eq_good-breaks} and \eqref{eq_ranks} we arrive at the identity :
$$
\delta_{F\otimes_\Lambda\cL'}(\rho)=\rk_\Lambda\cL(q)\cdot
\sum_{i=1}^d\max(f_i(\rho),b(q,\rho,s)).
$$
Notice that, since $f_i\geq 0$ for every $i\leq d$, the foregoing
identity persists also for $\rho<q^{-1}c(q,s)$. Recall also that
these functions $f_i$ are continuous and piecewise linear
(lemma \ref{lem_surprise}). This motivates the following :

\begin{claim}\label{cl_eventually}
For every $q\in\R_{\geq 0}$ and $c\in\R$, consider the function
$$
f_{q,c}:\R_{\geq 0}\to\R_{\geq 0}\qquad
\rho\mapsto\sum_{i=1}^d\max(f_i(\rho),q\rho-c).
$$
Then :
\begin{enumerate}
\item
$f_{q,c}$ is continuous, convex and piecewise linear.
\item
The (right and left) slopes of $f_{q,c}$ are of the form $qa+b$, where
$a\in\{0,1,\dots,d\}$, $b\in\Z$.
\item
Moreover, $f_{q,c}$ is eventually linear ({\em i.e.} of the form
$\rho\mapsto\rho x+y$ for every sufficiently large $\rho$).
\item
More precisely, if $q\in\Q_{\geq 0}$, then for every sufficiently
large $\rho\in\R_{\geq 0}$ the left and right slope of $f_{q,c}$
coincide, and their common value is a rational number whose denominator
divides the denominator of $q$.
\item
For every $\rho\in\R_{\geq 0}$, the function
$\R_{\geq 0}\to\R_{\geq 0}$ : $q\mapsto f_{q,c}(\rho)$
is non-decreasing, convex, continuous and piecewise linear.
\end{enumerate}
\end{claim}
\begin{pfclaim} By construction, the function
$\rk_\Lambda\cL(q)\cdot f_{q,c}$ is the mapping
$\delta_{F\otimes_\Lambda\cL'}$, whenever $q\in\SS$,
$\cL'=\mu_s^*\cL(q)$ and $c=c(q,s)\in\Sigma(q)$.
Hence, for $c\in\Sigma(q)$, convexity and continuity (and piecewise
linearity) of $f_{q,c}$ follow from proposition \ref{prop_invoke}.
Now, if $c$ and $c'$ are any two positive real numbers,
it is clear that :
$$
|f_{q,c}(\rho)-f_{q,c'}(\rho)|\leq q\cdot|c-c'|
\qquad\text{for every $\rho\in\R_{\geq 0}$}.
$$
Since $\Sigma(q)$ is dense in $\R$, it follows easily
that $f_{q,c}$ is convex and continuous for every $c\in\R$.
Similarly, for $q,q'\in\R_{\geq 0}$, we may bound the difference
$|f_{q,c}-f_{q',c}|$ in terms of $|q-q'|$, on every bounded subset
of $\R_{\geq 0}$; since $\SS$ is dense in $\R_{\geq 0}$, we
deduce continuity and convexity of $f_{q,c}$ for every $q\in\R_{\geq 0}$
and $c\in\R$.
Next, for given $\rho_0:=-\log r_0\geq 0$, let $\beta_1,\dots,\beta_k$
be the breaks of $F_{r_0}$, so that we have the break decomposition
$F_{r_0}=F_{r_0}(\beta_1)\oplus\cdots\oplus F_{r_0}(\beta_k)$. Set
$m_j:=\rk_\Lambda F_{r_0}(\beta_j)$ for every $j\leq k$.
We may find a segment $[\rho_0,\rho_1]$, and for every
$i\leq d$, an integer $j_i\leq k$ such that :
$$
f_i(\rho)=f_i(\rho_0)+(\rho-\rho_0)\cdot\beta^\natural_{j_i}
\qquad\text{for every $\rho\in[\rho_0,\rho_1]$}.
$$
It follows easily that there exists $\rho_2\in(\rho_0,\rho_1]$
such that :
\set\begin{equation}\label{eq_explicit-f}
f_{q,c}(\rho)=(qa+b)\cdot\rho+c'
\qquad\text{for every $\rho\in[\rho_0,\rho_2]$}
\end{equation}
where $a:=\max\{i\leq d~|~-\log\beta_{j_i}\leq(q\rho_0-c)+q\eps\}$
(notation of \eqref{eq_switch}), and :
$$
b:=\sum_{j>j_a}^k m_j\beta_j^\natural
\qquad
c':=ca+\sum^k_{j>j_a}m_j(\beta^\flat_j-\rho_0\beta^\natural_j).
$$
This shows the piecewise linearity of $f_{q,c}$. We deduce as
well that $b\in\Z$, since each term $m_j\beta_j^\natural$ is
the Swan conductor of the Galois module $F_{r_0}(\beta_j)$ (denoted
$\alpha(F_{r_0}(\beta_j))$ in \cite[\S8]{Hu3}).

This shows that (i) and (ii) hold. Moreover, \eqref{eq_explicit-f}
also easily implies (v). Assertion (iii) is already
known for every pair $(q,c)$ with $q\in\SS$ and $c\in\Sigma(q)$
(theorem \ref{th_tannaka}(i)). Next, if $c'\in\R$ is arbitrary,
since the distance between $f_{q,c}$ and $f_{q,c'}$ is bounded, and
$f_{q,c'}$ is convex and piecewise linear, it is easy to deduce that
$f_{q,c'}$ is also eventually linear. Finally, if $q'\leq q$ is any
positive real number, it is clear that $f_{q',c'}\leq f_{q,c'}$;
since $f_{q,c'}$ is eventually linear and $f_{q',c'}$ is convex,
it follows that right derivative $\rho\mapsto df_{q',c'}/dt(\rho^+)$
is non-decreasing and bounded; but from (ii) we see that the set of
possible slopes for $f_{q'c'}$ does not admit accumulation points,
hence the right derivative of $f_{q',c'}$ must be eventually constant.
This concludes the proof of (iii).

Assertion (iv) is clear from (ii).
\end{pfclaim}

Claim \ref{cl_eventually}(iii) says that, for every $q\in\R_{\geq 0}$
and $c\in\R$, the limit :
\set\begin{equation}\label{eq_asymptotic}
s(q):=\lim_{\rho\to+\infty}f_{q,c}(\rho)/\rho
\end{equation}
exists and is a rational number independent of $c$, whose denominator
divides the denominator of $q$. Now, suppose $q\in\SS$, $c\in\Sigma(q)$
and $c'>c$ is some real number; one sees easily that
$$
\ssw^\natural(F\otimes_\Lambda\mu^*_s\cL(q),r^+)\geq
\ssw^\natural(F\otimes_\Lambda\mu^*_{s'}\cL(q),r^+)
$$
for every $r\in(0,1]$ and every $s,s'\in K^+$ with $\log|s|=c$
and $\log|s'|=c'$. It follows that
$$
\beta(F,q)=
\sup\{\ssw^\natural(F\otimes_\Lambda\mu^*_s\cL(q),0^+)~|~
\text{$s\in K^+$ and $\log|s|\in\Sigma(q)$}\}=s(q)
$$
for every $q\in\SS$.

\begin{claim} The function $\Q_{\geq 0}\to\R$ : $q\mapsto\beta(F,q)$
is non-decreasing.
\end{claim}
\begin{pfclaim} It suffices to show that, if $q'<q<q''$ with
$q',q''\in\SS$ and $q\in\Q$, then $\beta(F,q')\leq\beta(F,q)\leq\beta(F,q'')$.
However, choose $s'\in K^+\setminus\{0\}$ such that $c(q',s')\in\Sigma(q')$;
from lemma \ref{lem_bricks} we see that there exists $\rho_0\in\R$
such that $\Delta(\cL(q))\cap([\rho_0,+\infty)\times\R_{\geq0})$
is the graph of a linear map. We may then find $s\in K^\times$
such that $|s|<|s'|$ and such that
$\Delta(\mu^*_s\cL(q))\cap([\rho_0,+\infty)\times\R_{\geq0})\cap\Delta(F)$
is a countable subset. Since $\Delta(\mu^*_s\cL(q))$ lies above
$\Delta(\mu^*_{s'}\cL(q'))$ in the region
$[\rho_0,+\infty)\times\R_{\geq0}$, an argument as in the foregoing
shows that $\delta_{F\otimes\mu_s^*\cL(q)}(\rho)>
\delta_{F\otimes\mu_{s'}^*\cL(q')}(\rho)$ for $\rho\geq\rho_0$.
But since $q'\in\SS$, we have seen that the slope of
$\delta_{F\otimes\mu_{s'}^*\cL(q')}$ equals $\beta(F,q')$, hence
$\beta(F,q)\geq\beta(F,q')$, as required. The proof of the other
inequality is similar, and shall be left to the reader.
\end{pfclaim}

From claim \ref{cl_eventually}(v) we deduce that $s$ is a non-decreasing
function. Since $s$ and $\beta(F,\cdot)$ agree on the dense subset $\SS$,
they must coincide for all $q\in\Q_{\geq0}$.
Combining with claim \ref{cl_eventually}(iv), this proves assertion (i).

(ii): From the definition of $f_{q,c}$, it is obvious that
$f_{q,0}(\rho)\geq q\rho\cdot\rk_\Lambda F$ for every
$\rho\in\R_{\geq 0}$, hence $s(q)\geq q\cdot\rk_\Lambda F$.
Furthermore, since $\sum_{i=1}^df_i(\rho)=\delta_F(\rho)$
is a convex function which is eventually linear of slope
$q_0:=\ssw^\natural(F,0^+)$, one sees easily that there exists
$c\in\R$ such that :
$$
f_i(\rho)\leq q_0\rho+c\qquad\text{for every $\rho\in\R_{\geq 0}$
and every $i\leq d$}.
$$
Hence $f_{q,0}(\rho)=q\rho\cdot\rk_\Lambda F$ for every $q>q_0$,
provided $\rho$ is large enough. Consequently
$s(q)=q\cdot\rk_\Lambda F$ for every $q>q_0$, so (ii) holds.

(iii): Claim \ref{cl_eventually}(v) implies that the
function $q\mapsto s(q)$ is convex and non-decreasing.
Next, if $q.q'$ are any two positive real numbers, it is clear
that $|s(q)-s(q')|\leq d\cdot|q-q'|$, so the mapping $s$ is also
continuous. Moreover, the convexity of $s$ implies that
the right derivative $ds/dt(\rho^+)$ exists for every
$\rho\in\R_{>0}$ and is non-decreasing, and $s$ is a
primitive of its right derivative.

\begin{claim}\label{cl_integral-slope}
$ds/dt(\rho^+)\in\Z$ for every $\rho\in\Q$.
\end{claim}
\begin{pfclaim} Write $\rho=a/b$ with relatively prime
positive integers $a,b$. By definition, we have :
\set\begin{equation}\label{eq_ration}
\frac{ds}{dt}(\rho^+)=\lim_{n\to+\infty}
bn\cdot\left\{s\left(\rho+\frac{1}{bn}\right)-s(\rho)\right\}.
\end{equation}
Now, if we let $n$ run over the positive integers, the right-hand
side of \eqref{eq_ration} is the limit of a sequence of integers,
since the denominators of both $(\rho+1/(bn))$ and $s(\rho)$
divide $bn$. 
\end{pfclaim}

We deduce from claim \ref{cl_integral-slope} that the right
derivative of $s$ is a non-decreasing step function (constant
on segments of the form $[a,b)$). Hence $s$ is piecewise linear
with integral slopes, which concludes the proof of (iii) and
of the theorem.
\end{proof}

\sset\subsubsection{}\label{subsec_Newton}
Let $F$ be as in theorem \ref{th_Newton}.
The idea is that the graph of $\beta(F,\cdot)$ should be
the Newton polygon associated to the sought break decomposition
of the stalk $F_0$ of $F$ over the missing point $\eta(0)$ (see
\eqref{subsec_ideally}). According to this picture, the
breaks of $F_0$ are the values $q_i\in\R_{>0}$
such that $ds/dt(q^-)\neq ds/dt(q^+)$ (where $s$ is defined as in
\eqref{eq_asymptotic});
naturally we call these the {\em break points\/} of $\beta(F,\cdot)$.
The first observation is that there are only finitely many
break points, and all of them are rational; indeed, this is
a straightforward consequence of theorem \ref{th_Newton}.
Let $0<q_1<q_2<\cdots<q_n$ be these break points, and set
$q_0:=0$. Since $\beta(F,\cdot)$ is piecewise linear and
non-negative, we may find unique $\mu_0,\mu_1,\dots,\mu_n\geq 0$
such that :
\set\begin{equation}\label{eq_on-the-mu}
\beta(F,q)=\sum_{j=0}^n\mu_j\cdot\max(q_j,q)
\qquad\text{for every $q\in\Q_{\geq 0}$}.
\end{equation}
Indeed, by deriving both sides of \eqref{eq_on-the-mu}, we find :
\set\begin{equation}\label{eq_deriva}
\sum^i_{j=0}\mu_j=
\left.\frac{d\beta(F,q)}{dq}\right|_{q=q_i^+}\qquad\text{for every $j\leq n$}.
\end{equation}
And since $\beta(F,\cdot)$ has integer slopes (theorem
\ref{th_Newton}(iii)), we deduce that $\mu_j\in\N$ for
every $j\leq n$. For every $j\leq n$, the integer
$\mu_j$ should be nothing else than the rank of the direct
factor of $F_0$ which is pure of break $q_j$. This is borne
out by the identity :
$$
\rk_\Lambda F=\sum_{j=0}^n\mu_j
$$
which holds, since $\beta(F,\cdot)$ is eventually linear of
slope $\rk_\Lambda F$ (theorem \ref{th_Newton}(ii)).
For this reason, we shall say that $\mu_i$ is the
{\em multiplicity\/} of the break $q_i$, for every $i=0,\dots,n$.
Now, let :
$$
0<\tau_1\leq\tau_2\leq\cdots\leq\tau_d
$$
be the unique sequence of rational numbers in which the value
$q_i$ appears with multiplicity $\mu_i$, for every $i=0,\dots,n$.
We have :

\begin{theorem}\label{th_missing-stalk}
Keep the notation of \eqref{subsec_Newton}, and let
$d:=\rk_\Lambda F$. Then there exist $\rho_0\geq 0$,
and real numbers $c_1,c_2,\dots,c_d$ such that :
$$
f_i(\rho)=\tau_i\cdot\rho+c_i\qquad\text{for every $\rho\geq\rho_0$
and every $i=1,\dots,d$}
$$
where $f_1,f_2,\dots,f_d:\R_{\geq 0}\to\R_{\geq 0}$ are defined
as in \eqref{eq_redundant}.
\end{theorem}
\begin{proof} Recall (claim \ref{cl_eventually})
that $f_{q,0}(\rho):=\sum_{i=1}^d\max(f_i(\rho),q\rho)$ for every
$\rho\in\R_{\geq 0}$. For every $\rho,y\in\R_{\geq 0}$, let us set :
$$
I_{\rho,y}:=
\{i\in\N~|~\text{$1\leq i\leq d$\ \ and\ \ $f_i(\rho)\leq y$}\}
\qquad N_{\rho,y}:=\sharp I_{\rho,y}
$$
(where, as usual, for any set $I$, we denote by $\sharp I$ the
cardinality of $I$). Let $q'>q$ be any two real numbers; we may
compute :
$$
f_{q',0}(\rho)-f_{q,0}(\rho)=q'\rho N_{\rho,q'\rho}-
q\rho N_{\rho,q\rho}-\sum_{i\in J}f_i(\rho)
\qquad\text{for every $\rho\geq0$}
$$
where :
$$
J:=\{i\in\N~|~\text{$1\leq i\leq d$\ \ and\ \ $q\rho<f_i(\rho)<q'\rho$}\}.
$$
It follows easily that :
\set\begin{equation}\label{eq_comine}
g(q,q',\rho):=\frac{f_{q',0}(\rho)-f_{q,0}(\rho)}{(q'-q)\rho}\in
[N_{\rho,q\rho},N_{\rho,q'\rho}]
\qquad\text{for every $\rho>0$}.
\end{equation}
Now, let $q>0$ be any real number which is not a break point for
$\beta(F,\cdot)$; let $k\in\{0,1,\dots,n\}$ be the largest
integer such that $q_k<q$. If $k=n$, set $q_{n+1}:=q+1$, so
that in any case $q\in(q_k,q_{k+1})$. We notice that the function
$s$ defined as in \eqref{eq_asymptotic} is linear on the interval
$[q_k,q_{k+1}]$, since the proof of theorem \ref{th_Newton} shows
that $s$ agrees with $\beta(F,\cdot)$ on $\Q_{\geq 0}$.
Especially :
$$
\lim_{\rho\to +\infty}g(q,q_{k+1},\rho)=
\lim_{\rho\to +\infty}g(q_k,q,\rho)=
\left.\frac{d\beta(F,q)}{dq}\right|_{q=q_k^+}.
$$
But recall that the slopes of $\beta(F,\cdot)$ are integers;
therefore, combining with \eqref{eq_comine} we deduce :
\set\begin{equation}\label{eq_relationshippy}
\left.\frac{d\beta(F,q)}{dq}\right|_{q=q_k^+}\in
[N_{\rho,q_k\rho},N_{\rho,q\rho}]\cap[N_{\rho,q\rho},N_{\rho,q_{k+1}\rho}]=
\{N_{\rho,q\rho}\}
\quad\text{for all large $\rho$}.
\end{equation}
The meaning of \eqref{eq_relationshippy} is that, if $q$ is not a break,
then the points $(\rho,f_i(\rho))$ tend to ``move away'' from the
line $\{(x,y)~|~qx=y\}$; indeed, \eqref{eq_relationshippy} shows
that if $q'\in(q_k,q)$ is any other real number, then
$I_{\rho,q\rho}=I_{\rho,q'\rho}$ provided $\rho$ is large enough.
Fix any $\eps>0$ such that :
$$
2\eps<\min\{q_{k+1}-q_k~|~k=0,\dots,n-1\}
$$
and set :
$$
J_k(\rho):=I_{\rho,(q_k+\eps)\rho}\setminus I_{\rho,(q_k-\eps)\rho}
\qquad\text{for every $k\leq n$ and every $\rho\geq 0$}.
$$
Notice that, since the functions $f_i$ are continuous (lemma
\ref{lem_surprise}), each set $J_k(\rho)$ will be eventually
independent of $\rho$ ({\em i.e.} for large values of $\rho$),
and we shall therefore denote it simply by $J_k$.
Summing up, so far we have exhibited a natural partition :
$$
\{1,2,\dots,d\}=J_0\amalg J_1\amalg\cdots\amalg J_n
$$
such that, for every $k\leq d$, the values
$T_k(\rho):=\{(\rho,f_i(\rho))~|~i\in J_k\}$ ``cluster'' around
a straight line of slope $q_k$. Explicitly, for every $\eps>0$ and
for every large $\rho$, the points of $T_k(\rho)$ lie
in the cone $C_\eps(k):=\{(x,y)\in\R^2_{>0}~|~|y/x-q_k|<\eps\}$.
Next we show that, for every large $\rho$, the set $T_k(\rho)$ actually
lies in a band of slope $q_k$ and fixed bounded width.
To this aim, for every $k=0,1,\dots,n$ and every $c\in\R$, set :
$$
h_k(\rho):=\sum_{i\in J_k}f_i(\rho)
\qquad
h^*_{k,c}(\rho):=\sum_{i\in J_k}\max(f_i(\rho),q_k\rho-c)\qquad
\text{for every $\rho\in\R_{\geq 0}$}.
$$
\begin{claim}\label{cl_properties-of-h}
(i)\ \ For every $k\leq n$ the following holds :
\begin{enumerate}
\item
the functions $h_k$ and $h^*_{k,c}$ are eventually linear.
\item
$\sharp J_k=\mu_k$.
\item
$\displaystyle\lim_{\rho\to+\infty}h_k(\rho)/\rho=q_k\mu_k=
\displaystyle\lim_{\rho\to+\infty}h^*_{k,c}(\rho)/\rho$.
\end{enumerate}
\end{claim}
\begin{pfclaim} Suppose $q\in(q_k,q_{k+1})$. We can then write :
$$
f_{q,0}(\rho)=q\rho\cdot\sum_{t\leq k}\sharp J_t+
\sum_{t=k+1}^n h_t(\rho)
\qquad\text{for every sufficiently large $\rho$}.
$$
Since the function $f_{q,0}(\rho)$ is eventually linear,
we deduce that, for every $k\leq n$, the sum
$\sum_{t=k+1}^n h_t$ is eventually linear, so the same holds
for each term $h_t$. Let
$C_t:=\displaystyle\lim_{\rho\to+\infty}h_t(\rho)/\rho$.
In view of \eqref{eq_asymptotic} we deduce :
$$
s(q)=q\cdot\sum_{t\leq k}\sharp J_t+\sum_{t=k+1}^n C_t
\qquad\text{for every $q\in(q_k,q_{k+1})$ and every $k\leq n$}.
$$
Now suppose that $q'\in(q,q_{k+1})$. Taking into account
\eqref{eq_deriva}, we find :
$$
(q'-q)\cdot\sum_{t\leq k}\sharp J_t=s(q')-s(q)=
(q'-q)\cdot\sum_{j\leq k}\mu_j
$$
from which (ii) follows easily, arguing by induction on $k$.
Finally, on the one hand we know that $h_k$ is eventually linear;
on the other hand, for every $\eps>0$, each of its summands
$f_i$ (for $i\in J_k$) is eventually contained in the cone
$C_\eps(k)$, so assertion (iii) for $h_k$ follows easily from (ii).
Next, we look at the identity :
$$
f_{q_k,c}(\rho)=(q_k\rho-c)\cdot\sum_{t<k}\mu_t+\sum_{t=k+1}^n h_t(\rho)+
h^*_{k,c}(\rho)
$$
which holds for every $k\leq n$ and every large enough $\rho$,
in view of (ii). Since $f_{q_k,0}$ and $h_{k+1},\dots,h_n$ are
eventually linear functions, we see that the same holds for
$h^*_{k,c}$, for every $k\leq n$. This shows (i), and also
the remaining assertion (iii) for $h^*_k$ follows easily.
\end{pfclaim}

We now fix $k\in\{0,1,\dots,n\}$, and write just $q$, $\mu$ $J$,
$h$ and $h_c^*$ instead of $q_k$, $\mu_k$, $J_k$, $h_k$, $h^*_{k,c}$.

\begin{claim}\label{cl_boundedness} For every $i\in J$, the function
$$
\rho\mapsto|f_i(\rho)-q\rho|
$$
is bounded.
\end{claim}
\begin{pfclaim} It follows easily from claim
\ref{cl_properties-of-h} that both functions :
$$
\sum_{i\in J}\{\max(f_i(\rho),q\rho)-q\rho\}
\quad\text{and}\quad
\sum_{i\in J}\{\max(f_i(\rho),q\rho)-f_i(\rho)\}
$$
are eventually constant. Since these summands are always non-negative,
we deduce that, for every $i\in J$, the terms :
$$
\max(f_i(\rho),q\rho)-q\rho
\quad\text{and}\quad
\max(f_i(\rho),q\rho)-f_i(\rho)
$$
are bounded, which is the claim.
\end{pfclaim}

\begin{claim}\label{cl_limi-vals}
For every $i\in J$ there exists $a_i\in\R$ such that :
$$
\lim_{\rho\to+\infty}f_i(\rho)-q\rho=a_i.
$$
\end{claim}
\begin{pfclaim}
Say that $J=\{i_0,\dots,i_0+\mu-1\}$. We prove, by induction on
$t$, that $a_{i_0+t}$ with the desired property exists for every
$t<\mu$. For $t<0$, there is nothing to prove. Next, suppose
that $t\geq 0$ and that the assertion is already known for every
integer $<t$; we set :
$$
g(\rho):=\sum_{i=i_0+t}^{i_0+\mu-1}f_i(\rho)
\quad
g^*_c(\rho):=\sum_{i=i_0+t}^{i_0+\mu-1}\max(f_i(\rho),q\rho-c)
\qquad\text{for every $\rho\in\R_{\geq 0}$ and $c\in\R$}.
$$
Using the inductive assumption, and claim
\ref{cl_properties-of-h}, we see that there exists
$C\in\R$ with :
\set\begin{equation}\label{eq_with}
\lim_{\rho\to+\infty} g(\rho)-\rho q(\mu-t)=C.
\end{equation}
Set :
$$
a:=\liminf_{\rho\to+\infty} f_{i_0+t}(\rho)-q\rho
\qquad
b:=\limsup_{\rho\to+\infty} f_{i_0+t}(\rho)-q\rho.
$$
Notice that :
\set\begin{equation}\label{eq_lower-bound}
a\geq a_{i_0},\dots,a_{i_0+t-1}
\end{equation}
since $f_i\leq f_{i+1}$ for every $i=1,\dots,d-1$. Suppose $a<b$,
pick $x\in(a,b)$ and set $c:=-x$; in view of \eqref{eq_lower-bound},
we have :
$$
h_c^*(\rho)=t(q\rho-c)+g^*_c(\rho)
\qquad\text{for every large enough $\rho$}.
$$
Then claim \ref{cl_properties-of-h} implies that $g^*_c$ is
eventually linear of slope $q(\mu-t)$. Combining with
\eqref{eq_with}, we deduce that there exists $C'\in\R$ such that :
\set\begin{equation}\label{eq_lim-exists}
\lim_{\rho\to+\infty} g^*_c(\rho)-g(\rho)=C'.
\end{equation}
However, due to our choice of $x$, for every $\rho\geq 0$
and every $\eps>0$ we may find $\rho',\rho''\geq\rho$ such that :
$$
g^*_c(\rho')=g(\rho')
\qquad\text{and}\qquad
g^*_c(\rho'')-g(\rho'')>x-a+\eps
$$
which contradicts \eqref{eq_lim-exists}. Hence $a=b$, and the
common value is a real number, due to claim \ref{cl_boundedness}.
This concludes the inductive step.
\end{pfclaim}

To conclude the proof of the theorem, we shall show that
the function $f_{i_0+t}$ is eventually linear, whenever
$i_0+t\in J=\{i_0,\dots,i_0+\mu-1\}$.
We shall proceed by induction on $t$.
If $t<0$, there is nothing to prove. Hence, suppose that the
assertion is known for every integer $<t$. Set $a:=a_{i_0+t}$, where
$a_{i_0},\dots,a_{i_0+\mu-1}$ are the real numbers whose existence
is ensured by claim \ref{cl_limi-vals}. Let $J(a):=\{i\in J~|~a_i=a\}$;
we shall show simultaneously that all the functions $f_i$ with $i\in J(a)$
are eventually linear, hence we may suppose that $i_0+t$ is the smallest
element of $J(a)$. In this case, using the inductive assumption and
claim \ref{cl_properties-of-h}, we see that both functions $g$ and
$g^*_{-a}$ introduced in the proof of claim \ref{cl_limi-vals} are
eventually linear. Moreover, it is also clear that :
$$
\lim_{\rho\to+\infty}g^*_{-a}(\rho)-g(\rho)=0.
$$
It follows that $g(\rho)=g^*_{-a}(\rho)$ for every large enough $\rho$,
hence 
\set\begin{equation}\label{eq_estim-below}
f_i(\rho)\geq qp+a
\qquad\text{for every $i\in J(a)$ and every large $\rho$}.
\end{equation}
On the other hand, let $i_1$ be the largest element of $J(a)$;
if $i_1<i_0+\mu-1$, choose $b\in(a,a_{i_1+1})$; otherwise,
set $b:=a_{i_1}+1$. In either case, we may write :
$$
h^*_{-b}(\rho)=i_1(q\rho+b)+\sum_{i>i_1}f_i(\rho)
\qquad\text{for every large enough $\rho$}.
$$
Hence $\sum_{i>i_1}f_i$ is eventually linear, and therefore
the same holds for $g-\sum_{i>i_1}f_i=\sum_{i\in J(a)}f_i$.
Clearly :
$$
\lim_{\rho\to+\infty}\sum_{i\in J(a)}(f_i(\rho)-q\rho-a)=0.
$$
Combining with \eqref{eq_estim-below}, we deduce the contention.
\end{proof}

\sset\subsubsection{}
Let $(\Gamma_0,\leq)$ be the abelian group $\Q\times\Gamma_K$, endowed
with the ordering such that :
$$
(q,c)\leq(q',c')\qquad
\text{if and only if either $q'<q$ or else $q=q'$ and $c\leq c'$}.
$$
(This is the lexicographic ordering, except that the ordering
on $\Q$ is the reverse of the usual one.)
For given $r\in\Gamma^+$, let $\Gamma_r$ be the value group of
the valuation $|\cdot|_{\eta(r)}$. The mapping :
$$
\Gamma_0\to\Q\otimes_\Z\Gamma_r
\quad :\quad (q,c)\mapsto c\cdot r^q\cdot(1-\eps)^q
$$
is an isomorphism of groups which does not respect the orderings
(indeed, the ordering on $\Gamma_0$ induced by this isomorphism
is also lexicographic, but the two factors $\Q$ and $\Gamma_K$ are
swapped). Nevertheless, we may interpret theorem \ref{th_missing-stalk},
by saying that the ``missing stalk $F_0$'' admits a break decomposition
which is naturally indexed by elements of $\Gamma_0^+$.
More precisely, we have :

\begin{theorem}\label{th_decompose}
Let $F$ be as in theorem {\em\ref{th_Newton}}. Then
there exist $r_0\in\Gamma^+_K$, a connected open subset
$U\subset\D(1)^*$ and a decomposition :
$$
F_{|U}=\bigoplus_{(q,c)\in\Gamma_0^+}M(q,c)
$$
where each summand $M(q,c)$ is a locally constant $\Lambda$-module
on $U_\et$, such that :
\begin{enumerate}
\item
$U\cap\D(\eps)\neq\emptyset$ for every $\eps\in\Gamma^+_K$.
\item
For every $r\leq r_0$, we have $\eta(r)\in U$ and :
$$
M(q,c)_{\eta(r)}=F_r(c\cdot r^q\cdot(1-\eps)^q).
$$
\end{enumerate}
\end{theorem}
\begin{proof} Set $E:=\cE{nd}_\Lambda(F)$, the sheaf of $\Lambda$-linear
endomorphisms of $F$. We have to exhibit $U\subset\D(1)^*$ fulfilling
(i), and for each $(q,c)\in\Gamma_0^+$, a projector
$\pi\in E(U)$ such that
$$
(\Img\,\pi)_{\eta(r)}=F_r(c\cdot r^q\cdot(1-\eps)^q).
$$
By theorem \ref{th_missing-stalk}, we may find $\rho_0\geq 0$
such that the functions $f_i$ are linear on the half-line
$[\rho_0,+\infty)$; up to replacing $\rho_0$ by a larger real
number, we may achieve that, for every $i,j\leq d$, the graphs
of the functions $f_i, f_j$ are either disjoint or equal.
Say that $\rho_0=-\log r_0$. The open subset $U$ shall be
constructed by removing from $\D(r_0)^*$ infinitely many
closed discs. Indeed, suppose $r\leq r_0$, and let
$\beta_1(r),\dots,\beta_k(r)$ be the breaks of $F_r$; we
may assume that $c\cdot r^q\cdot(1-\eps)^q=\beta_1(r)$.
From our choice of $\rho_0$, it follows that
$\beta_i(r)^\flat=\beta_j(r)^\flat$ if and only if $i=j$;
in other words, the decomposition \eqref{eq_break-decomp}
is the same as the decomposition
$$
F_{\eta(r)^\flat}=\bigoplus_{\alpha\in\Gamma^+_K}F^\flat_r(\alpha)
$$
(notation of \eqref{eq_flat-breaks}). Hence, the stalk
$E_{\eta(r)^\flat}$ contains a projector $\pi_{\eta(r)^\flat}$
that cuts out the summand $F_r(\beta_1(r))$.
However, since $F$ is locally constant on $\D(1)^*_\et$, the same
holds for $E$; especially, $E$ is overconvergent, in the sense
of \cite[Def.8.2.1]{Hu2}. Let :
$$
\D(1)^*_\et\stackrel{\mu}{\longrightarrow}\D(1)^*
\stackrel{\nu}{\longrightarrow}\D(1)^*_{\mathrm{p.p}}
$$
be the natural morphisms of sites, where $\D(1)^*_{\mathrm{p.p}}$
denotes the topological space $\D(1)^*$ endowed with its partially
proper topology (\cite[Def.8.1.3]{Hu2}).
On the one hand, using \cite[Prop.1.5.4]{Hu2}, we deduce that
the natural map :
$$
(\mu_*E)_{\eta(r)^\flat}\to E_{\eta(r)^\flat}
$$
is a bijection; on the other hand, according to
\cite[Prop.8.1.4(a)]{Hu2}, the counit of adjunction
$\nu^*\nu_*(\mu_*E)\to(\mu_*E)$ is an isomorphism, hence also the natural
map
$$
(\nu_*\mu_*E)_{\eta(r)^\flat}\to \mu_*E_{\eta(r)^\flat}
$$
is bijective. Therefore, we may find a partially proper open
neighborhood $V\subset\D(1)^*$ of $\eta(r)^\flat$, and a section
$\pi_V\in E(V)$, such that $(\pi_V)_{\eta(r)^\flat}=\pi_{\eta(r)^\flat}$.
Since $E$ is locally constant and
$\pi_{\eta(r)^\flat}^2=\pi_{\eta(r)^\flat}$, it follows that
$\pi_V$ is a projector in $E(V)$, and its stalk $(\pi_V)_{\eta(r)}$
cuts out the direct summand $F_r(\beta_1(r))$.

Next, for every $r'<r$, let $\D(r',r^-):=\bigcup_{r'<s<r}\D(r',s)$;
by inspecting the proof of theorem \ref{th_continuities}
we see that, provided $r'$ is sufficiently close to $r$,
there exists a decomposition :
$$
F_{|\D(r',r^-)}\simeq\bigoplus_{j=1}^k G_j
$$
consisting of locally constant $\Lambda$-modules $G_j$ on
$\D(r',r^-)_{\et}$, such that :
$$
(G_j)_{\eta(s)}=F_s(\beta_j(s))
\qquad\text{for every $j\leq k$ and every $s\in[r',r)$}.
$$
Thus, we may find a unique projector $\pi_{\D(r',r^-)}\in E(\D(r',r^-))$
that cuts out the direct summand $G_1$. Finally, on the intersection
$V\cap\D(r',r^-)$, we must have $\pi_V=\pi_{\D(r',r^-)}$, hence
we obtain a section $\pi_{W(r)}$ of $E$ on the open subset
$W(r):=V\cup\D(r',r^-)$. Up to removing some small closed subset,
we may assume that $W(r)$ is of the form
$\D(r',r)\setminus\bigcup_{i=1}^m\E(a_i,\rho_i)$, where $\E(a_i,\rho_i)$
denotes the closed disc with center $a_i\in\D(r,r)$ and radius
$\rho_i<r$. By a standard compactness argument, we see that
there exists a sequence $(r_n~|~n\in\N)$ of elements of $\Gamma_K$,
with $r_n\leq r_0$ for every $n\in\N$, and $\lim_{n\to+\infty}r_n=0$,
such that $U:=\bigcup_{n\in\N}W(r_n)$ contains the subset
$\{\eta(s)~|~s\in\Gamma_K;\ \ s\leq r_0\}$. Clearly, this open subset
$U$ will do.
\end{proof}

\end{document}